\documentclass[11pt,twoside]{article}

\usepackage{subfiles}
\usepackage[top=1in, bottom=1in, left=1in, right=1in]{geometry}

\usepackage{subcaption}

\usepackage{amssymb}
\usepackage{amsmath,amsthm,amsfonts,dsfont,bm}
\usepackage{url}
\usepackage{rotating}
\usepackage{multirow}
\usepackage{color}
\usepackage{xcolor}
\usepackage{enumitem}

\usepackage{hyperref}

\usepackage[showonlyrefs]{mathtools}

\hypersetup{
  colorlinks   = true, 
  urlcolor     = blue, 
  linkcolor    = blue, 
  citecolor   = red 
}

\newtheorem{theorem}{Theorem}

\newtheorem{lemma}[theorem]{Lemma}

\newtheorem{corollary}[theorem]{Corollary}

\newtheorem{assumption}{Assumption}
\newtheorem{remark}{Remark}
\newtheorem{proposition}[theorem]{Proposition}

\numberwithin{equation}{section}
\numberwithin{theorem}{section}

\renewcommand{\P}{\operatorname{\mathbb{P}}}

\newcommand{\E}{\operatorname{\mathbb{E}}}

\newcommand{\R}{\mathbb{R}}

\usepackage[mathscr]{euscript}
\DeclareSymbolFont{rsfs}{U}{rsfs}{m}{n}
\DeclareSymbolFontAlphabet{\mathscrsfs}{rsfs}

\newcommand{\one}{\mathbf{1}}
\newcommand{\eps}{\varepsilon}
\newcommand{\rmd}{\mathrm{d}}
\newcommand{\cA}{\mathcal{A}}
\newcommand{\cS}{\mathbb{S}}
\newcommand{\cC}{\mathcal{C}}
\newcommand{\cD}{\mathcal{D}}
\newcommand{\cU}{\mathcal{U}}
\newcommand{\cN}{\mathcal{N}}
\newcommand{\btheta}{\theta}
\newcommand{\AREA}{\mathrm{Area}}
\newcommand{\VOL}{\text{Vol}}
\newcommand{\Csol}{C_{\textsc{sol}}}
\newcommand{\expressionI}{\textsc{(I)}}
\newcommand{\expressionII}{\textsc{(II)}}
\newcommand{\expressionIII}{\textsc{(III)}}
\newcommand{\expressionIV}{\textsc{(IV)}}
\newcommand{\expressionV}{\textsc{(V)}}

\newcommand{\gmm}{\textsc{\tiny GMM}}
\newcommand{\logistic}{\textsc{\tiny Logistic}}
\newcommand{\cDlogistic}{\cD_{\textsc{\tiny logistic}}}
\newcommand{\cDgmm}{\cD_{\textsc{\tiny gmm}}}
\newcommand{\GAUSSIANMEASURE}{\beta}
\newcommand{\MARGIN}{\kappa}
\newcommand{\DIMENSION}{d}
\newcommand{\SAMPLES}{n}
\newcommand{\SAMPLEIDX}{i}
\newcommand{\DIMENSIONIDX}{j}
\newcommand{\REPLICAIDX}{\ell}
\newcommand{\REPLICAS}{L}
\newcommand{\SIGNALDISORDER}{S}
\newcommand{\QUENCHEDDISORDER}{Z}
\newcommand{\ANNEALEDDISORDER}{W}
\newcommand{\NORMALIZEINTERPOLATION}{Z_{t,n-1,d}}
\newcommand{\INVERSEMILLSEXPRESSION}{V}
\newcommand{\DENSITYBOUND}{\alpha_0}
\newcommand{\SIGNALSHORTDEPENDENCE}{\lambda}
\newcommand{\SIGNALMEDIUMDEPENDENCE}{\lambda}
\newcommand{\EXPONENTEXPSCALE}{\gamma}
\newcommand{\APPROXCONSTRAINT}{\eta}
\newcommand{\VARPROBLEMSPHERICAL}{\bar\delta}

\newcommand{\OPNORM}{\textsc{op}}
\newcommand{\SIGNALGAUSSIAN}{\beta_{\textsc{signal}}}
\newcommand{\MILLSBOUND}{K_{\textsc{MILLS}}}
\newcommand{\signalsubexp}{K_{\textsc{signal}}}
\newcommand{\pdfnormal}{\varphi_{0}}
\newcommand{\pdfsignal}{\varphi_{1}}
\newcommand{\err}{\text{err}}
\newcommand{\grad}{\nabla}
\newcommand{\upref}[1]{\textup{\ref{#1}}}
\newcommand{\equpref}[1]{\textup{\eqref{#1}}}

\def\RHS{\mbox{\small{\rm RHS}}}
\def\LHS{\mbox{\small\rm LHS}}

\def\RSG{\textsc{RSG}}
\def\RS{\textsc{RS}}
\def\underH{\underline{H}}
\def\overH{\overline{H}}

\def\sLP{\mbox{\rm\tiny LP}}
\def\sGD{\mbox{\rm\tiny GD}}
\def\sBayes{\mbox{\rm\tiny Bayes}}
\def\err{\mbox{\rm err}}

\begin{document}

\title{
How abundant are good interpolators?
}

\author{
August Y. Chen \thanks{Cornell University, Department of Computer Science. Ithaca, NY, USA. Email: \texttt{ayc74@cornell.edu}}
\and
Ahmed El Alaoui \thanks{Cornell University, Department of Statistics and Data Science. Ithaca, NY, USA. Email: \texttt{elalaoui@cornell.edu}. 
}
}

\date{}
\maketitle

\vspace*{-.3cm} 
\begin{abstract}
Let $S$ be the set of unit norm linear classifiers $\theta \in \R^d$ which correctly classify every point of a labeled dataset $(X_i,y_i)_{i=1}^n$, $X_i \in \R^d$, $y_i \in \{-1,+1\}$, with a possibly negative margin $\kappa$ fixed in advance. Under two natural data-generating distributions of the $(X,y)$ pairs -- a Gaussian mixture model and a logistic model with Gaussian features -- and in the proportional regime $n/d \to \alpha$ with small enough $\alpha$, we establish a large deviation principle on the event that a point $\theta$ chosen uniformly at random from $S$ achieves a given generalization error, with high probability over the choice of the data. The associated large deviation rate function is deterministic and describes the proportion, at the exponential scale in $d$, of interpolating classifiers having a given desired performance. As a consequence, we establish the following concentration phenomenon: all but an exponentially small fraction of interpolating classifiers have approximately the same generalization performance given by the unique maximizer of this rate function. 

We numerically compare this maximizer to the performance of empirical risk minimization by gradient descent and to the performance of a natural linear program, both finding a point in $S$, and deduce that in the overparametrized regime of small $\alpha$, these efficient procedures outperform the vast majority of interpolators, pointing to their nontrivial benign overfitting in this setting. 
\end{abstract}

{
  \hypersetup{linkcolor=black}
  \tableofcontents
}

\section{Introduction}
Overparametrized statistical models have far more parameters than training points. Their empirical success in the context of deep learning has elicited a flurry of theoretical activity on the reasons such models are able to generalize well. With so many parameters, these models can perfectly fit the training data with simple and efficient algorithms minimizing an empirical loss without hurting generalization performance~\cite{belkin2019reconciling}. With such an overparametrization one expects the existence of a large set of points in the parameter space with zero training error, suggesting that good generalization is a property of the specific algorithm used to minimize the loss function, rather than a consequence of an explicit control on the model's complexity, as might be done by explicit regularization~\cite{zhang2021understanding}. This phenomenon, known as \emph{implicit regularization}, has been the topic of much research showing that gradient descent and its variants tend to find solutions of some minimum norm, thereby selecting a model in a statistically favorable way; see for instance~\cite{neyshabur2017exploring,soudry2018implicit}. 

In this paper we are interested in the generalization properties of the set of \emph{interpolators}: points in parameter space with zero training error. We want to understand the \emph{full profile} of possible generalization performances, beyond the ones algorithmically selected via empirical risk minimization, and investigate whether there is a typical performance achieved by ``most" interpolators.      
       
We consider the simple setting of linear binary classification with fixed margin $\kappa$, in $d$ dimensions and with $n$ samples, under two distinct data-generating distributions of constant signal-to-noise ratio $\lambda$: a mixture of two Gaussians, and the logistic model with Gaussian features. 
In each case we provide a precise asymptotic characterization of the fraction of interpolators having any desired generalization performance $x$, in the proportional limit $n \to\infty$ and $n/d \to \alpha$, and $\alpha$ small enough. Our result takes the form of a ``quenched" large deviation principle (LDP) on the event that a point chosen uniformly at random from the set of interpolators achieves a performance approximately $x$. The rate function of this LDP is given by a variational problem whose maximum is achieved at a value $x_{\star} = x_{\star}(\alpha,\kappa,\lambda)$, the performance of a `typical' interpolator. As a byproduct, we establish that with high probability over the data, \emph{all but an exponentially small fraction} of interpolators have a performance near $x_{\star}$. 

We next remark that this typical performance is lower than the ones algorithmically achieved: we compare the value of $x_{\star}$ to the performances $x_{\sLP}$ and $x_{\sGD}$ of two interpolators respectively found by a linear programming procedure proposed in~\cite{montanari2024tractability}, and gradient descent on the logistic loss function analyzed in~\cite{deng2022model}. We find that $x_{\star}$ is smaller than $x_{\sLP}$ and $x_{\sGD}$ for a large range of parameters. 
We deduce that ``good" interpolators --- in the sense of being at least as good as the two aforementioned efficient baseline algorithms --- are exceedingly rare in this high-dimensional, constant SNR, overparametrized regime.      

For further comparison we also derive an asymptotic characterization of the Bayes optimal performance $x_{\sBayes}$ -- that of a parameter vector drawn from the posterior distribution given the data -- under both data-generating distributions. This is similarly done by establishing an LDP on the generalization performance under the corresponding posterior distributions.      

Our result confirms a conjecture of Theisen, Klusowski and Mahoney~\cite{theisen2021good} that the generalization error of a uniformly random interpolator should concentrate around a typical value. They address in their work the question of `abundance of good interpolators' where they conduct a heuristic analysis of this typical error in the setting of the Gaussian mixture model, leading them to claim that ``good classifiers are abundant in the interpolation regime", seemingly at odds with our findings. But their analysis implicitly assumes a concentration condition about the Gram matrix of the data which can only hold at high SNR: the separation between the Gaussian centers must be growing with the dimension. In this paper we reach a more nuanced conclusion, namely that good interpolators are exceedingly rare in the constant SNR regime -- making implicit regularization a genuinely nontrivial phenomenon -- and that the typical performance $x_\star$ increases with the SNR $\lambda$ (and with the sampling ratio $\alpha$) leading to a good typical performance in the large SNR limit, in accordance with the results of~\cite{theisen2021good} in that limit.

\section{Settings} 
We consider the setting of binary classification with data $(X_i,y_i)_{i=1}^n$ where $X_i \in \R^d$ is a vector of $d$ covariates and $y_i \in \{-1,+1\}$ is the class label. Let $\theta_{\star} \in \cS^{d-1}(\sqrt{d})$ be a distinguished `signal' direction. We consider two data-generating distributions:  
\begin{enumerate}
\item \textbf{The Gaussian mixture model} $\cDgmm$:  independently for $i=1,\cdots,n$ let $y_i = \pm 1$ with equal probability, and 
\begin{equation}\label{eq:gmm_dist}
X_i \,|\, y_i \sim N\left(\sqrt{\frac{\lambda}{d}}\, y_i \theta_\star\,,\,I_d\right)\,
\end{equation}
for a fixed signal-to-noise ratio (SNR) parameter $\lambda > 0$.   
\item \textbf{The logistic model} $\cDlogistic$: independently for $i=1,\cdots,n$ let $X_i \sim N(0,I_d)$ and conditionally on $X_i$,
\begin{equation}\label{eq:logistic_dist}
y_i =
\begin{cases}
 +1 &\mbox{with prob.} ~~\varphi \big(\sqrt{\lambda/d}\,\langle X_i ,\theta_{\star}\rangle\big)\,, \\
 -1 &\mbox{with prob.} ~~1 - \varphi \big(\sqrt{\lambda/d}\,\langle X_i ,\theta_{\star}\rangle\big)\,, 
\end{cases}
\end{equation}
for a fixed SNR $\lambda > 0$.
Here $\varphi: \R \to [0,1]$ is an increasing link function which we fix throughout the paper to have the sigmoid form
 \begin{equation}\label{eq:sigmoid}
 \varphi(x) =  \frac{1}{1+e^{-x}}\,. 
 \end{equation}
 We will denote the Rademacher distribution on $\{-1,+1\}$ with probability $p$ assigned to $+1$ by $\text{Rad}(p)$, hence we can write \[y_i \, |\, X_i \sim \text{Rad}\Big(\varphi \big(\sqrt{\lambda/d}\,\langle X_i ,\theta_{\star}\rangle\big)\Big)\,.\] 
\end{enumerate} 

Let $X \in \R^{n \times d}$ be the matrix with the vectors $X_i$ as its rows, and let $\kappa \in \R$ be a fixed margin parameter. We consider the set $S_{\kappa}(X,y) \subset \cS^{d-1}(\sqrt{d})$ of directions $\theta$ that correctly classifies all data points up to margin $\kappa$:
\begin{equation}\label{eq:interpolators}
S_{\kappa}(X,y) := \Big\{\theta \in \cS^{d-1}(\sqrt{d}) ~:~ y_i \langle X_i, \theta \rangle \ge \kappa \sqrt{d} ~\text{ for all }1 \le i \le n\Big\}\,. 
\end{equation} 
This is the set of linear interpolators with margin $\kappa$.

The generalization error of an estimator of $\theta_\star$, a measurable function $\hat{\theta}_n : (\R^d \times \{-1,+1\})^n \to \cS^{d-1}(\sqrt{d})$, is the probability of misclassifying a new data point: 
\begin{equation*} 
\mathcal{E}(\hat{\theta}_n) = \P\big(y\langle X, \hat{\theta}_n \rangle < 0\big), 
\end{equation*}
where $(X,y) \sim \cD \in\{\cDgmm,\cDlogistic\}$. We remark that by rotational invariance of the Gaussian distribution, the above error only depends on the inner product $\langle \hat{\theta}_n, \theta_\star \rangle/d$ of which it is a monotonically decreasing function. Therefore we can measure performance by amount of overlap instead of generalization error.    

\paragraph{Large deviations for interpolating linear classifiers} Let $\mu_{X,y}$ denote the uniform measure over $S_{\kappa}(X,y)$, Eq.~\eqref{eq:interpolators}, whenever this set is non empty (and define $\mu_{X,y}$ arbitrarily otherwise). We want to establish an asymptotic, high-probability formula for the random quantity 
\begin{equation} \label{eq:log-ratio}
\frac{1}{d} \log \mu_{X,y}\big(\{\theta : |\langle \theta, \theta_\star \rangle/d - x| \le \eps\}\big) \, ,
\end{equation}
where $n$ and $d$ are proportionally large and $\eps$ is an arbitrary small number, capturing the fraction of interpolators with overlap near $x$ with $\theta_\star$ in the exponential scale. Equivalently we would like to establish a `quenched' large deviation principle for the event that $\langle \theta, \theta_\star \rangle/d \approx x$, where $\theta$ is drawn uniformly at random from $S_{\kappa}(X,y)$. From~\eqref{eq:log-ratio} this boils down to estimating the logarithm of the area of $S_{\kappa}(X,y)$ (as a measurable subset of the sphere) and that of its slices 
\begin{equation}\label{eq:slices}
S_{x,\eps}(X,y) := S_{\kappa}(X,y) \cap \{\theta : |\langle \theta , \theta_\star\rangle/d - x| \le \eps\}\,,~~~ x \in (-1,1)\,,~~\eps>0\,, 
\end{equation}
corresponding to all points having overlap approximately $x$ with the signal direction $\theta_\star$. The quantity in~\eqref{eq:log-ratio} is then equal to $(1/d) \log \big(|S_{x,\eps}(X,y)
| / |S_{\kappa}(X,y)|\big)$, where $|A|$ denotes the $d-1$ dimensional surface area of any Lebesgue measurable subset of $\cS^{d-1}(\sqrt{d})$. This will be the content of Theorem~\ref{thm:mainresult}.

\paragraph{Large deviations for Bayes-optimal classifiers} We also consider a Bayesian setting where the signal direction $\theta_\star$ is drawn uniformly at random from the sphere $\cS^{d-1}(\sqrt{d})$. If the data $(X,y)$ conditional on $\theta_\star$ are drawn from $\cD \in \{\cD_{\gmm}, \cD_{\logistic}\}$, the density of the posterior distribution (relative to the normalized area measure on $\cS^{d-1}(\sqrt{d})$) of $\theta_\star$ given $(X,y)$ takes the form
\begin{equation}\label{eq:posterior}
p_{X,y}(\theta) = \frac{1}{Z_{X,y}} \prod_{i=1}^n \exp u\left(\frac{y_i\big\langle  X_i ,\theta\big\rangle}{\sqrt{d}}\right)\,, ~~~ \theta \in \cS^{d-1}(\sqrt{d})\,,
\end{equation}
where $u(x) = \sqrt{\lambda} \, x$ for $\cD = \cD_{\gmm}$, and $u(x) = \log \varphi(\sqrt{\lambda}\, x)$ for $\cD = \cD_{\logistic}$. The denominator $Z_{X,y}$ is the normalizing constant called the partition function of $p_{X,y}$. 

We will similarly be interested in quantifying the probability of the large deviation event $\{|\langle \theta, \theta_\star \rangle/d - x| \le \eps\}$, when $\theta$ is drawn from the posterior distribution~\eqref{eq:posterior}, and this amounts to establishing an asymptotic formula for
\begin{equation} \label{eq:log-posterior}
\frac{1}{d} \log p_{X,y}\big(\{\theta : |\langle \theta, \theta_\star \rangle/d - x| \le \eps\}\big)\,.
\end{equation}
This will be the content of Theorem~\ref{thm:mainresultposterior}. 

\section{Main results}
\label{sec:main}
Fix the parameters $\alpha,\kappa,\lambda$, a function $u: (-1,1) \to [-\infty,\infty)$, and a random variable $S$. Consider the following function defined on $(-1,1) \times [0,1) \to \R$:
\begin{equation}\label{eq:phi_master}
\begin{aligned}
\Phi(x,q) 
&=   \alpha \E \log \E_{W} \exp u\Big( x S + \sqrt{(1-x^2)q}Z + \sqrt{(1-x^2)(1-q)} W \Big) \\
& \hspace{6.5cm}+ \frac{1}{2} \frac{q}{1-q} + \frac{1}{2} \log\big((1-x^2)(1-q)\big) \,,
\end{aligned}
\end{equation}
where $S, Z, W$ are mutually independent, and $Z, W \sim N(0,1)$. $\E_W$ is the expectation with respect to $W$ while $\E$ averages over $Z$ and $S$.

\subsection{LDP for a uniform interpolator}
\label{sec:LDPinterp}
In this section we state our main results on the large deviations of the overlap under the uniform measure $\mu_{X,y}$ on the set of interpolators. We let
\begin{equation}\label{eq:complementary_cdf_normal}
\cN(x) = \int_{x}^{+\infty} e^{-z^2/2} \frac{\rmd z}{\sqrt{2\pi}}
\end{equation}
be the complementary c.d.f.\ of the standard normal distribution, and let 
\begin{equation}
u(x) = 
\begin{cases}
0 &\mbox{if} ~~x \ge \kappa\,,\\
-\infty & \mbox{otherwise}\,.
\end{cases}
\end{equation}
The expression of $\Phi$ simplifies to
\begin{align}\label{eq:phi_interp}
\Phi(x,q) &= \alpha \E \log \cN\Biggl( \frac{\MARGIN - x S  - \sqrt{(1-x^2)q}Z}{\sqrt{(1-x^2)(1-q)}} \Biggr)  + \frac{1}{2} \frac{q}{1-q} + \frac{1}{2} \log\big((1-x^2)(1-q)\big)\,.
\end{align}

The law of the random variable $S$ will depend on the data-generating distribution:  
\begin{itemize}
\item If $\big\{ (X_i, y_i) \big\}_{i=1}^n \sim \cD_{\gmm}$ we let 
\begin{equation}\label{eq:defSgmm}
S = \sqrt{\lambda} + G\,,~~~~G \sim N(0,1)\,.
\end{equation}
\item If $\big\{ (X_i, y_i) \big\}_{i=1}^n \sim \cD_{\logistic}$ we let  
\begin{equation}\label{eq:defSlogistic}
S = YG\,,~~~~G \sim N(0,1)\,,~~~~ Y|G \sim \text{Rad}\big(\varphi(\sqrt{\lambda}G)\big)\,.
\end{equation}
\end{itemize}

We are now in position to state our quenched large deviation principle for the uniform measure $\mu_{X,y}$ on $S_{\kappa}(X,y)$. 
Let
\begin{equation}\label{eq:minmax}
\phi = \sup_{x \in (-1,1)}\inf_{q \in [0,1)}\Phi(x,q)\,,
\end{equation}
and for an interval $I \subseteq (-1,1)$,
\begin{equation}\label{eq:ratefct}
\bar{\phi}(I) = \sup_{x \in I} \inf_{q \in [0,1)}\Phi(x,q) - \phi\,.
\end{equation}
Further let $\kappa_+ = \max\{\kappa, 0\}$.
\begin{theorem}\label{thm:mainresult} 
Fix $\lambda>0$, $\kappa \in \R$ and let $n = \lfloor \alpha d \rfloor$. There exists $\alpha_0 = \alpha_0(\lambda,\kappa_+)>0$ such that the following holds for all $\alpha < \alpha_0$ in both cases~\eqref{eq:defSgmm} and~\eqref{eq:defSlogistic}. We have $\phi >-\infty$, and for any interval $I \subseteq [-1+\delta, 1-\delta]$ with $\delta = \delta(\alpha,\lambda, \MARGIN_+)>0$,
there exists $\eps_0>0$ such that for $0 < \eps \le \eps_0$ there exists $d_0, K >0 $ such that for $d \ge d_0$ we have 
\begin{equation}\label{eq:areagmm} 
\left|\frac{1}{d} \log \big| S_{\kappa}(X,y)\big| - \phi - \log\sqrt{2\pi e}\right| \le \eps\,,
\end{equation}
\begin{equation}\label{eq:ldpgmm} 
\left|\frac{1}{d} \log \mu_{X,y} \Bigl(\big\{\theta: \langle \theta, \theta_\star\rangle / d \in I\big\}\Bigr) - \bar{\phi}(I)\right| \le 2\eps\,,
\end{equation}
and 
\begin{equation}\label{eq:extreme_bd_gmm} 
\frac1d \log \mu_{X,y}\Big( \big\{ \theta: \langle \theta, \theta_\star\rangle / d \ge 1-\delta \big\} \Big) < -c\,,
\end{equation}
with probability at least $1-e^{-d/K}$ over the data $(X,y)$, where $c>0$ is a universal constant. (Here $\eps_0$ depends on $\alpha, \kappa, \lambda, |I|$, and $d_0, K$ depend additionally on $\eps$.)
\end{theorem}

The quantity $\delta=\delta(\alpha, \lambda, \MARGIN_+)$ is defined in (\ref{eq:deltaalphafulldef}); for $\alpha \le \alpha_0$, we have $\delta \le 1/10$. 
We furthermore note that as a corollary of our proof, the above results hold for any prescribed $\delta > 0$ if now $\alpha \le \alpha_0(\delta, \lambda, \MARGIN_+)$.
The above theorem characterizes the area of the set of interpolators and the rate function of the LDP for the overlap with the signal direction under $\mu_{X,y}$. This shows that the normalized overlap $\langle \theta , \theta_\star\rangle/d$ is exponentially concentrated around the maximizers of the function $x \mapsto \inf_{q} \Phi(x,q)$ when $\theta \sim \mu_{X,y}$, with high probability over the data. (Note the bound (\ref{eq:extreme_bd_gmm}) shows that the set of interpolators with normalized overlap outside $[-1+\delta, 1-\delta]$ makes an exponentially small contribution.)

We remark that $\alpha_0$ depends on the margin $\kappa$ only when the latter is positive; we will show in Appendix~\ref{sec:interpolationfreeenergy} that one can take $\alpha_0(\lambda,\kappa) = c(\lambda)/(1+\kappa^2)$ in this case. For negative margin the result holds as long as $\alpha$ is small enough as a function of the SNR $\lambda$ only.    

We further show that for $\alpha$ small enough the `sup-inf' saddle point problem~\eqref{eq:minmax} is achieved at a unique pair $(x,q) \in [0,1)^2$ which satisfies a system of two nonlinear equations.
We state the results for the Gaussian mixture model and the logistic model in two separate propositions.  

Let $\cA$ be the inverse of Mills' ratio 
\begin{equation}\label{eq:mills}
\cA(x) = - \cN'(x)/\cN(x) = \frac{1}{\sqrt{2 \pi}}\frac{e^{-x^2/2}}{\cN(x)}\,.
\end{equation}
We will say that the saddle point problem $\sup\inf \Phi$ is uniquely achieved if there is a unique $x_{\star} \in [0,1)$ achieving $\sup_{x \in (-1,1)} \inf_{q \in [0,1)} \Phi(x,q)$, and if there is a unique $q_{\star} \in [0,1)$ achieving $\inf_{q \in [0,1)}\Phi(x_{\star},q)$.

\begin{proposition}[GMM]\label{lem:RS_biv_eq_gmm}
Let $\kappa \in \R$, $\lambda>0$. There exists $\alpha_0 = \alpha_0(\lambda,\kappa_+)$ such that for all $\alpha < \alpha_0$, the saddle point problem~\eqref{eq:minmax} where $\Phi$ is defined via~\eqref{eq:phi_master} with $S = \sqrt{\lambda} + G$, $G \sim N(0,1)$, is uniquely achieved at a pair $(x_{\star},q_{\star})$ solving the system of equations
\begin{align}
    \frac{x}{\sqrt{1-x^2}} &= \alpha \sqrt{\lambda(1-q)}\E \left[\cA\left(\frac{\kappa -  xS - \sqrt{(1-x^2)q}Z}{\sqrt{(1-x^2)(1-q)}}\right)\right]\,, \label{eq:RS_biv_eq_gmm_x} \\
    \frac{q}{1-q} &= \alpha \E \left[\cA\left(\frac{\kappa - xS -\sqrt{(1-x^2)q}Z}{\sqrt{(1-x^2)(1-q)}}\right)^2\right]\,, \label{eq:RS_biv_eq_gmm_q}
\end{align}
where $Z \sim N(0,1)$ independently of $S$.
\end{proposition}

\begin{proposition}[Logistic]\label{lem:RS_biv_eq_logistic}
Let $\kappa \in \R$, $\lambda>0$. There exists $\alpha_0 = \alpha_0(\lambda,\kappa_+)$ such that for all $\alpha < \alpha_0$, the saddle point problem~\eqref{eq:minmax} where $\Phi$ is defined via~\eqref{eq:phi_master} with $S = YG$, $G \sim N(0,1)$ and $Y|G \sim \text{Rad}\big(\varphi(\sqrt{\lambda}G)\big)$, is uniquely achieved at a pair $(x_{\star},q_{\star})$ solving the system of equations
\begin{align}
    \frac{x}{\sqrt{1-x^2}} &= \alpha \sqrt{\lambda(1-q)}\E \left[\varphi\big(-\sqrt{\lambda}S\big)\cA\left(\frac{\kappa - x S  - \sqrt{(1-x^2)q}Z}{\sqrt{(1-x^2)(1-q)}}\right)\right]\,, \label{eq:RS_biv_eq_logistic_x} \\
    \frac{q}{1-q} &= \alpha \E \left[\cA\left(\frac{\kappa - x S  - \sqrt{(1-x^2)q}Z}{\sqrt{(1-x^2)(1-q)}}\right)^2\right]\,, \label{eq:RS_biv_eq_logistic_q}
\end{align}
where $Z \sim N(0,1)$ independently of $(G,Y)$.
\end{proposition}
Furthermore, in both the above cases, $(x_{\star}, q_\star) \in [0, K\alpha]^2$ where $K>0$ depends on $\lambda$ and $\MARGIN_+$.

\begin{corollary}
    The uniqueness of the maximizer of the function $\inf_{q} \Phi(\cdot,q)$, achieved at $x^\star \in [-1+\delta, 1-\delta]$, combined with the large deviation bound~\eqref{eq:ldpgmm} implies that all but an exponentially small fraction of points $\theta \in S_{\kappa} (X,y)$ have nearly the same overlap $x_\star$ with the signal direction $\theta_\star$, with high probability over the data $(X,y)$.  
\end{corollary}

\paragraph{Physical interpretation of $q$:} While the parameter $x_\star$ represents the asymptotic amount of overlap a direction $\theta \sim \mu_{X,y}$ has with $\theta_\star$, the parameter $q_\star$ has an equally concrete interpretation: let $\theta$ and $\theta'$ be two independent copies from $\mu_{X,y}$ conditional on the same realization of the data $(X,y)$, and denote by $\bar{\theta}$ and $\bar{\theta}'$ their projections (of norm $\sqrt{d}$) orthogonal to $\theta_\star$. Then $q_\star$ represents the asymptotic value of the overlap $\langle \bar{\theta},\bar{\theta}' \rangle/d$. A crucial step in the proofs of the above formulas is to show that this overlap indeed concentrates around $q_\star$. 

\subsection{LDP for the posterior distribution}
   \label{sec:LDPpost}
We state a similar result for the posterior distribution~\eqref{eq:posterior} in the Bayesian setting. 
Recall the function $\Phi$ from~\eqref{eq:phi_master}. In this Bayesian case both the function $u$ and the random variable $S$  depend on the data-generating distribution. 

\begin{itemize}
\item If $\big\{ (X_i, y_i) \big\}_{i=1}^n \sim \cD_{\gmm}$ we let 
\begin{equation}\label{eq:defuSgmm}
S = \sqrt{\lambda} + G~~~ \mbox{ and }~~~ u(x) = \sqrt{\lambda}\,x\,. 
\end{equation}
\item If $\big\{ (X_i, y_i) \big\}_{i=1}^n \sim \cD_{\logistic}$ we let  
\begin{equation}\label{eq:defuSlogistic}
S = YG ~~~ \mbox{ and }~~~ u(x) = \log \varphi(\sqrt{\lambda}\,x)\,,
\end{equation}
where as before, $G \sim N(0,1)$ and $Y|G \sim \text{Rad}\big(\varphi(\sqrt{\lambda}G)\big)$.
\end{itemize}

We have the following quenched LDP for the posterior distribution $p_{X,y}$. Let 
\begin{equation}\label{eq:minmaxBayes}
\phi = \sup_{x \in (-1,1)}\inf_{q \in [0,1)}\Phi(x,q)\,,
\end{equation}
and for $I \subset (-1,1)$,
\begin{equation}\label{eq:ratefctBayes}
\bar{\phi}(I) = \sup_{x \in I} \inf_{q \in [0,1)}\Phi(x,q) - \phi\,.
\end{equation}

\begin{theorem}\label{thm:mainresultposterior}
Fix $\lambda>0$ and let $n = \lfloor \alpha d \rfloor$. There exists $\alpha_0 = \alpha_0(\lambda)>0$ such that the following holds for all $\alpha < \alpha_0$ in both cases~\eqref{eq:defuSgmm} and~\eqref{eq:defuSlogistic}. 
For any interval $I \subseteq [-1+\delta, 1-\delta]$ with $\delta = \delta(\alpha,\lambda)>0$, there exist $\eps_0 >0$ such that for all $0 < \eps \le \eps_0$ there exists $d_0, K>0$ such that for all $d \ge d_0$ we have
\begin{equation} 
\left|\frac{1}{d} \log Z_{X,y} -  \phi - \log \sqrt{2\pi e}\right| \le \eps\,,
\end{equation}
\begin{equation} \label{eq:ldp_posterior}
\left|\frac{1}{d} \log p_{X,y} \Bigl(\big\{\theta: \langle \theta, \theta_\star\rangle / d \in I\big\}\Bigr) - \bar{\phi}(I)\right| \le 2\eps\,,
\end{equation}
and 
\begin{equation}\label{eq:extreme_bd_posterior} 
\frac1d \log p_{X,y}\Big( \big\{ \theta: \langle \theta, \theta_\star\rangle / d \ge 1-\delta \big\} \Big) < -c\,,
\end{equation}
with probability at least $1-e^{-d/K}$ over the data $(X,y)$, where $c>0$ is a universal constant. (Again, $\eps_0$ depends on $\alpha, \lambda, |I|$, and $d_0, K$ depend additionally on $\eps$.)
\end{theorem}
Here, $\delta=\delta(\alpha, \lambda)$ is defined in (\ref{eq:posteriordeltadef}); for $\alpha \le \alpha_0$, we have $\delta \le 1/10$. 
Again, the above results hold for any prescribed $\delta > 0$ if $\alpha \le \alpha_0(\delta, \lambda)$.
(The bound (\ref{eq:extreme_bd_posterior}) shows that the set of interpolators with normalized overlap outside $[-1+\delta, 1-\delta]$ makes an exponentially small contribution.)

Similarly in this case, the sup-inf saddle point problem is uniquely achieved:    

\begin{proposition}[GMM]\label{lem:RS_biv_bayes_gmm} 
For $\lambda>0$, there exists $\alpha_0 = \alpha_0(\lambda)$ such that for all $\alpha < \alpha_0$, the saddle point problem~\eqref{eq:minmaxBayes} where $\Phi$ is defined via~\eqref{eq:phi_master} with $u(x) = \sqrt{\lambda} x$ and $S = \sqrt{\lambda} + G$, $G \sim N(0,1)$, is achieved at a unique pair $ (x_{\sBayes},q_{\sBayes})$, which solves the system of equations 
\begin{align}\label{eq:explicitxq0}
    x = \frac{q}{1-q}\,,~~~~\mbox{and}~~~~
    \alpha \lambda\, (1-x^2) = \frac{q}{(1-q)^2}\,.
\end{align}
These equations admit the explicit solutions 
\begin{align}\label{eq:explicitxq}
x_{\sBayes} = \frac{\alpha \lambda}{1+\alpha \lambda}\,,~~~~\mbox{and}~~~~
    q_{\sBayes} = \frac{\alpha \lambda}{1+2\alpha \lambda}\,.
\end{align}
\end{proposition}

\begin{proposition}[Logistic]\label{lem:RS_biv_bayes_logistic}
Define the random variable 
\begin{align}\label{eq:R}
R(x,q) := \frac{\E_{W} \varphi'(\sqrt{\lambda}V)}{\E_{W}\varphi(\sqrt{\lambda}V)}\,,~~
\mbox{where} ~~ V := xYG+\sqrt{(1-x^2)q}Z + \sqrt{(1-x^2)(1-q)}W\,,
\end{align}
and $G, Z ,W \sim N(0,1)$ are mutually independent and $Y|G \sim \text{Rad}\big(\varphi(\sqrt{\lambda}G)\big)$ independently of everything else.
For $\lambda>0$ there exists $\alpha_0 = \alpha_0(\lambda)$ such that for all $\alpha < \alpha_0$, the saddle point problem~\eqref{eq:minmaxBayes} where $\Phi$ is defined via~\eqref{eq:phi_master} with $u(x) = \log \varphi(\sqrt{\lambda} x)$  and $S = YG$, is uniquely achieved at a pair $(x_{\sBayes},q_{\sBayes})$, which solves the system of equations 
\begin{align}%
\begin{aligned}\label{eq:system_bayes}
    \alpha\lambda (1-q) \E \big[\varphi\big(-\sqrt{\lambda}YG\big) R(x,q)\big] &= \frac{x}{1-x^2}\,,\\
    \mbox{and}~~~~~~~\alpha \lambda (1-x^2) \E \big[R(x,q)^2\big] &= \frac{q}{(1-q)^2}\,.
    \end{aligned}
\end{align}
\end{proposition}

The uniqueness results combined with~\eqref{eq:ldp_posterior} shows exponential concentration of the normalized overlap under the posterior distribution, and characterizes the Bayes-optimal performance in estimating the signal direction $\theta_\star$. Again, in both of the above cases, $(x_{\star}, q_\star) \in [0, K\alpha]^2$ where $K>0$ depends on $\lambda$.

\paragraph{A special relation between $x$ and $q$:}
    Keeping in mind the interpretation of $q$ mentioned in the previous section, the present Bayesian setting imposes a special relation between $x$ and $q$: for two independent copies $\theta,\theta'$ from $p_{X,y}$ we have by the law of iterated expectations 
    \begin{equation}\label{eq:nishimori}
    \E \langle \theta, \theta_\star\rangle = \E \langle \theta, \theta'\rangle \,.
    \end{equation}
    Therefore by concentration of the random variables inside each of the above expectations (see Appendix~\ref{subsec:logconcaveconcentration}) the Pythagorean theorem would impose the relation    
    \begin{equation}
        x = x^2 + (1-x^2) q\,,
        ~~~\mbox{i.e.,}~~~ x = \frac{q}{1-q}\,.
    \end{equation} 
    This identity is apparent in the Gaussian mixture case, see Eq.~\eqref{eq:explicitxq0}, but is less clear in the logistic case. We show a similar result:   
    \begin{lemma}\label{lem:nishimori_logistic}
        For $\lambda>0$, there exists $\alpha_0=\alpha_0(\lambda)>0$ such that for every $\alpha\in(0,\alpha_0)$, any solution $(x,q)\in[0,1)^2$  to the system \eqref{eq:system_bayes} satisfies
    \[x=\frac{q}{1-q}\,.\]
    More precisely, the solution is of the form
    \begin{equation} 
    x=\frac{s}{1+s}\,,\qquad q=\frac{s}{1+2s}\,,
    \end{equation}
    where $s>0$ is the unique solution to the scalar fixed-point equation
    \begin{equation}\label{eq:scalar_fixed_point_s}
    s=\alpha\lambda\,
    \E\left[R\Big(\frac{s}{1+s},\frac{s}{1+2s}\Big)^2\right]\,.
    \end{equation}
    \end{lemma}
    
    This allows to search for the pair $(x_{\sBayes},q_{\sBayes})$ along a univariate curve, usually referred to as the Nishimori line, see~\cite[Chapter 4]{NishimoriBook} and~\cite{lelarge2016fundamental}. No such simplification occurs in general outside the Bayesian setting. In particular the system of equations governing the pair $(x_\star,q_\star)$ concerning the uniform measure $\mu_{X,y}$ treated in the previous section remains genuinely bivariate.  

    The proofs of Theorems \ref{thm:mainresult} and \ref{thm:mainresultposterior} can be found in Appendix~\ref{sec:proof_main_results}. The proofs of Propositions \ref{lem:RS_biv_eq_gmm}, \ref{lem:RS_biv_eq_logistic}, \ref{lem:RS_biv_bayes_gmm}, \ref{lem:RS_biv_bayes_logistic}, and Lemma \ref{lem:nishimori_logistic} can be found in Appendix~\ref{subsec:univariate_maximizer_pf}.
    
\subsection{Comparisons and numerical simulations}
\label{sec:numerical}

\begin{figure}[h!]
    \centering
    \begin{subfigure}[b]{0.48\textwidth}
        \centering
    \includegraphics[width=\textwidth]{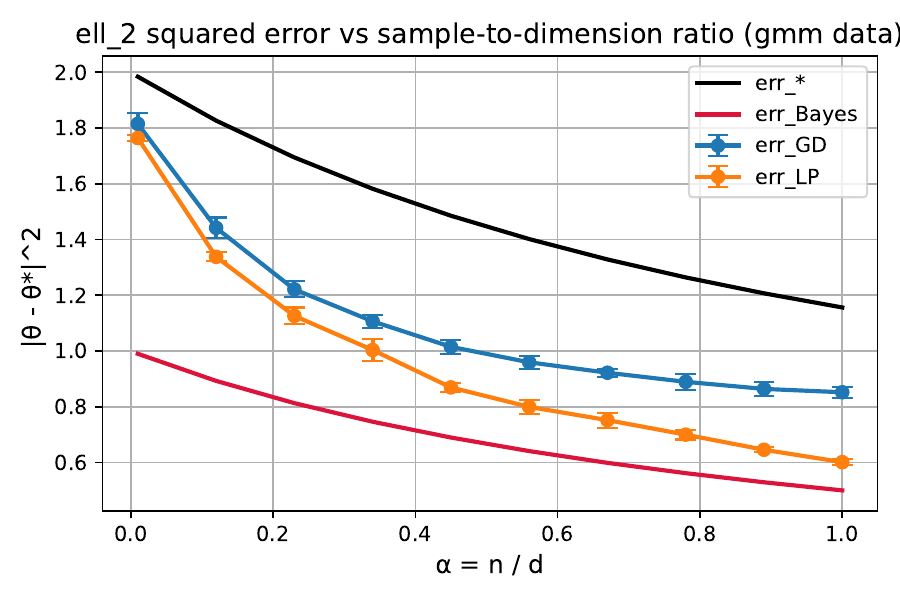}
        \caption{Data from $\cD_{\gmm}$}
        \label{fig:plot1}
    \end{subfigure}
    \hfill 
    \begin{subfigure}[b]{0.48\textwidth}
        \centering
        \includegraphics[width=\textwidth]{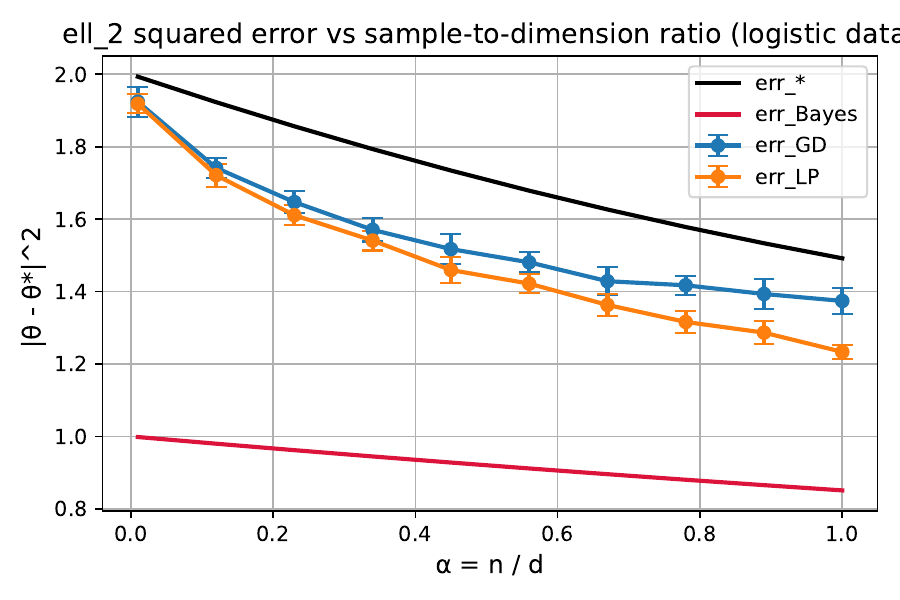}
        \caption{Data from $\cDlogistic$}
        \label{fig:plot2}
    \end{subfigure}
    \caption{The squared $\ell_2$ error $\|\hat{\theta} - \theta_\star\|_2^2$ at parameters $\kappa = 0$, $\lambda=1$. The black line is the error $\err_\star$ of an interpolator chosen uniformly at random. The red curve is the Bayes-optimal error $\err_{\sBayes}$. The blue and orange curves are the mean errors of the GD and LP estimators averaged across $T=20$ Monte Carlo trials. The error bars represent one standard deviation.} 
    \label{fig:two_plots}
\end{figure}

In this section we compare the performance of a typical interpolator relative to that of a draw from the posterior and to two efficient algorithms for finding points in $S_{\kappa}(X,y)$. (Note that the posterior distribution is not supported on the set of interpolators; a draw from $p_{X,y}$ will typically not belong to $S_{\kappa}(X,y)$ for any $\kappa$.) For convenience we rescale the parameter space by $\sqrt{d}$ in this section, so that $\theta_{\star}$ and all its estimators have unit $\ell_2$ norm.

The first algorithm is a second-order cone programming procedure  
\begin{equation}\label{eq:lp}
\mathrm{LP}~:~~~\max_{\theta \in \R^d} ~ \frac{1}{n}\sum_{i=1}^n \big\langle y_i X_i , \theta\big\rangle ~~\mbox{subj.\ to }~~ \|\theta\|_2 \le 1\,,~~ \big\langle y_i X_i , \theta \big\rangle \ge \kappa \,,~\forall i \in [n]\,. 
\end{equation}

Note that the spherical constraint $\|\theta\|_2 = 1$ is relaxed to a unit ball norm constraint $\|\theta\|_2 \le 1$ making the optimization problem convex. By abuse of terminology we refer to this procedure as a \emph{linear} program (LP), emphasizing the linear nature of the objective and the data-related constraints. This LP maximizes the correlation with the vector $v = (1/n)\sum_{i=1}^n y_iX_i$ representing a proxy of the signal direction, while staying on the correct side of the hyperplane corresponding to each data point. In~\cite[Theorem 5.3]{montanari2024tractability} it was shown that when the data is drawn from $\cD_{\logistic}$, the norm constraint is saturated, resulting in an estimator $\hat{\theta}_{\sLP}$ lying on the unit sphere, $\|\hat{\theta}_{\sLP}\|_2=1$ with high probability. Thus the algorithm produces a point in $S_{\kappa}(X,y)$ with high probability.

The second algorithm is empirical risk minimization of a smooth loss function via gradient descent: 
\begin{equation}\label{eq:erm}
\min_{\theta \in \R^d} ~ \Big\{R(\theta) := \sum_{i=1}^n \ell\big(\big\langle y_i X_i , \theta\big\rangle - \kappa \|\theta\|_2\big)\Big\}\,,
\end{equation}
where $\ell : \R \to \R_{\ge 0}$ is a decreasing function with $\lim_{x \to +\infty} \ell(x) = 0$; we only consider the logistic loss $\ell(x) = \log(1+e^{-x})$ for concreteness. The gradient descent iteration on $R$ reads 
\begin{align}\label{eq:gd}
\begin{aligned}
\mathrm{GD}~:~~~
\theta^{t+1} &= \theta^t - \eta\nabla R(\theta^t) \\
&= \theta^t -\eta \sum_{i=1}^n \big(y_i X_i - \kappa \theta^t/\|\theta^t\|\big) \ell'\big(\big\langle y_i X_i , \theta^t\big\rangle - \kappa \|\theta\|_2\big) \,,
\end{aligned}
\end{align}
where $\eta>0$ is a fixed stepsize and $t \in \{0,1,\cdots\}$. 
It was shown in~\cite[Theorem 6.1]{montanari2024tractability} that under smoothness assumptions on $\ell$ (satisfied by the logistic loss) and for small enough stepsize $\eta$, this iteration produces a margin-$\kappa$ classifier: 
\begin{equation}
\lim_{t\to \infty} \min_{1\le i \le n} \big\langle y_i X_i , \theta^t/\|\theta^t\|_2 \big\rangle \ge \kappa\,,
\end{equation} 
as long as $\|\theta^{t}\|_2 \to \infty$ as $t \to \infty$.  

For both data distributions $\cDlogistic$ and $\cD_{\gmm}$ and both the LP and GD algorithms we perform a Monte Carlo estimate of the squared $\ell_2$ distance $\|\hat{\theta}- \theta_\star\|_2^2 = 2(1-\langle \hat{\theta}, \theta_\star\rangle)$ over $T=20$ trials where the data is freshly resampled.   
We fix the dimension to $d=3000$, the margin to $\kappa = 0$, the signal-to-noise ratio to $\lambda = 1$ and vary the number of samples-to-dimension ratio $\alpha$ in the range $\{0.01,\cdots,1\}$ in increments of 0.1.

In the case of LP we use a specialized numerical solver to obtain a solution $\hat{\theta}_{\sLP}$ to~\eqref{eq:lp} and verify that the normalization condition $\|\hat{\theta}_{\sLP}\|_2=1$ is satisfied to $10^{-2}$ accuracy in all instances.

 In the case of GD we run the iteration~\eqref{eq:gd} with step size $\eta = 1.5$ until either the $\ell_2$ norm of the gradient of $R$ is smaller than $10^{-8}$ or a maximum number of iterations of $500$ is reached. We verify that the margin condition $\min_{1\le i \le n} \langle y_i X_i , \hat{\theta}_{\sGD} \rangle \ge \kappa$ is satisfied in all instances, where $\hat{\theta}_{\sGD} = \theta^t/\|\theta^t\|_2$ at the terminal time step $t$. 

 We compare these empirically obtained values to our theoretical predictions of the squared $\ell_2$ distances achieved by a typical interpolator $\theta \sim \mu_{X,y}$:
 \[\err_{\star} := \lim_{d \to \infty} \E\|\theta- \theta_\star\|_2^2  = 2(1 - x_{\star}),\]
 where $x_{\star}$ is described in Lemmas~\ref{lem:RS_biv_eq_gmm} and~\ref{lem:RS_biv_eq_logistic} for the data-generating distributions $\cD_{\gmm}$ and $\cD_{\logistic}$ respectively,  
 and the error achieved by the posterior mean $\hat{\theta}^{\sBayes}$, 
 \[\err_{\sBayes} := \lim_{d \to \infty}\E\|\hat{\theta}^{\sBayes}- \theta_\star\|_2^2  = 1 - x_{\sBayes},\]
 where $x_{\sBayes}$ is described in Lemmas~\ref{lem:RS_biv_bayes_gmm} and~\ref{lem:RS_biv_bayes_logistic} for these same distributions.
 The corresponding fixed point equations are numerically solved via bisection.
The numerical results are summarized in Figure~\ref{fig:two_plots}. We see that both algorithms significantly outperform a typical interpolator under both distributions and across all tested values of $\alpha$, with the Bayes error being smallest overall. Since this error concentrates exponentially in $d$ under the uniform measure on $S_{\kappa}(X,y)$ (Theorem~\ref{thm:mainresult}), only an exponentially small fraction of points in $S_{\kappa}(X,y)$ perform at least as well as LP or GD.

 This numerical finding can supported by theory: for data generated from the logistic model, \cite{montanari2024tractability} established that the LP (\ref{eq:lp}) achieves correlation $x_{\sLP} = \Theta(\sqrt{\alpha})$ in the limit $\MARGIN \rightarrow -\infty$ with $\lambda$ fixed; see Lemmas 22, 23, 24 therein. In contrast we have $x_{\star} = O(\alpha)$ for $\MARGIN \le 0$. This comparison again establishes that for small $\alpha$, only an exponentially small fraction of points in $S_{\MARGIN}(X, y)$ perform as well as the LP, in line with the results of our numerical simulations.

We remark that this effect is purely due to overparametrization: the performance gap between the above estimators vanishes as the sampling ratio $\alpha$ is increased. Heuristically this can be seen from the fact that the area of the set $S_{\kappa}(X,y)$ is a decreasing function of $\alpha$ since the map $\alpha \mapsto \phi$ (Eq.~\eqref{eq:minmax}) is decreasing, leading to consistent estimation for large enough $\alpha$ (or $\lambda \to \infty$). This is in agreement with the results of~\cite{montanari2024tractability} who provided bounds on the maximum and minimum estimation performance of the points in $S_{\kappa}(X,y)$. These bounds converge towards a common value as $\alpha$ increases (see Figure 4 therein).

\section{Related work and discussion}
\paragraph{Statistical literature and benign overfitting} 
There has been an impressive amount of work on overparametrization and implicit regularization over the last decade suggesting that overparametrization not only does not hurt statistical performance, but is essential for algorithmic tractability. We refer to~\cite{bartlett2021deep,belkin2021fit} for a partial synthesis of the literature. Most existing work analyzes the performance of specific algorithms and statistical procedures in the overparametrized regime and establish their success under suitable data-generating assumptions. A subset of this literature is concerned with obtaining exact high-dimensional asymptotics in (generalized) linear models under idealized data distributions, allowing a precise tracking of the generalization error and other metrics as the parameters of the model are changed. In this line of work we mention the analysis of maximum likelihood and M-estimation in logistic regression~\cite{candes2020phase,deng2022model}, Bayes-optimal estimation in GLMs~\cite{barbier2019optimal}, boosting and minimum $\ell_1$ norm interpolators~\cite{liang2021interpolating,li2021minimum,liang2022precise}, maximum margin interpolators~\cite{soudry2018implicit,montanari2025generalization}, among many others. 

Closer to our setting, \cite{montanari2024tractability} identify sharp bounds on \emph{the interpolation threshold} for data generated from the logistic model: the maximum value of $\alpha = n/d$ under which interpolators exist with positive probability, and above which the set of interpolators is empty with high probability. They  prove bounds on the generalization performance of any element of $S_{\kappa}(X,y)$ and leave open the question of whether there exists a value within those bounds which is typical. They additionally analyze the performance of the linear program~\eqref{eq:lp}, while \cite{deng2022model} analyze the ERM~\eqref{eq:erm}, which are both compared against in the previous section.      

 Our work is closely related to~\cite{theisen2021good,harel2024provable} who address the question of abundance of good interpolators and tempered overfitting of randomly selected interpolating neural networks. As mentioned earlier \cite{theisen2021good} conducts a heuristic analysis in the same setting as ours, under the Gaussian mixture model, and conclude that most interpolators have good generalization performance, in the sense of being close to Bayes-optimal. A careful reading of their argument reveals a  seemingly innocuous assumption (Eq.\ (12) therein) which can only hold at high SNR: the separation between the Gaussian centers need to grow faster than $\sqrt{d} \log n$, and is otherwise violated at constant separation. 
 \cite{harel2024provable} considers the setting of fully connected multi-layer neural networks (NN) under binary features $x$ and clean labels $h_\star(x) \in \{-1,1\}$ corrupted with bit flip noise; see also~\cite{buzaglo2024uniform}. When the clean labels are generated by a ``teacher" network $h_\star$ with $N_\star$ number of neurons, they prove that the generalization of a randomly picked interpolating ``student" neural network (of potentially much larger size) with $N$ neurons is close to Bayes-optimal when the number of samples $n$ is roughly larger than $N_\star^4 \vee N$. Importantly, the number of weights of the student network is absent from the sample complexity, which can be very large without hurting generalization, leading to tempered overfitting. But the term $N_\star^4$ in the sample complexity is roughly quadratic in the number of weights of the \emph{teacher} network, and an analogy with our setting where the teacher and the student have the same number of parameters $d$ suggests that this is still a ``large $\alpha$ regime" as $n$ needs to grow at least like $d^2$.         

\paragraph{High dimensional spin glasses and the perceptron model} 
The set of interpolators is also known as the set of solutions to the spherical perceptron model (or the half-space model when $\kappa=0$): if one draws $n$ uniformly random half-spaces of margin $\kappa$ in $d$ dimensions, when is their common intersection with the sphere typically non-empty? And what is the typical area of this intersection? Motivated by questions of learning with threshold functions, Cover~\cite{cover1965geometrical} determined the interpolation threshold in this model, also known as the \emph{storage capacity} for $\kappa=0$: the set is typically non-empty if $n \le 2(1-o(1))d$ and typically empty if otherwise $n\ge 2(1+o(1))d$, with additional finer understanding of the critical window. Gardner~\cite{gardner1988space,gardner1988optimal} proposed formulas which would later bear her name in the literature for the storage capacity and for the logarithm of the surface area, based on the non-rigorous replica method, for general $\kappa \ge0$. 

These formulas were proved by Shcherbina and Tirozzi~\cite{shcherbina2003rigorous} leveraging convexity and concentration of measure. Stojnic~\cite{stojnic2013another} found an elegant and simple Gaussian comparison argument establishing Gardner's storage capacity for $\kappa \ge 0$. Talagrand's books~\cite{talagrand2010mean,talagrand2011advanced} present a unified treatment of a larger family of models and incorporate additional important techniques such as Guerra's interpolation method~\cite{guerra2002thermodynamic,guerra2003broken} into their study.  

The behavior of the spherical perceptron model when $\kappa<0$ is much more delicate and few results are available. \cite{stojnic2013negative} showed that Gardner's formula obtained for $\kappa \ge0$ is an upper bound when $\kappa < 0$. \cite{montanari2024tractability} prove upper and lower bounds on the capacity for $\kappa \le 0$ which match as $\kappa \to -\infty$, and a detailed study at the physics level of rigor can be found in~\cite{franz2017universality}. An algorithm for finding a feasible point was proposed in~\cite{el2022algorithmic} and was shown to work for $\kappa<0$ under an unproved analytic condition regarding the structure of a certain variational problem.      

The binary version of the problem where the sphere is replaced by the binary hypercube has also been intensely studied. The storage capacity question has recently been resolved: the existence of a sharp threshold sequence was established in~\cite{nakajima2023sharp,xu2021sharp} and the capacity formula predicted by Krauth and M\'ezard~\cite{krauth1989storage} was established   in~\cite{ding2019capacity,huang2024capacity}. In addition, peculiar structural phenomena of freezing and preponderance of isolated solutions were established in a symmetric version of this model~\cite{aubin2019storage,perkins2021frozen,abbe2022proof} but remain open otherwise, and efficient algorithms for finding (a connected set of) solutions were proposed in~\cite{abbe2022binary}.

\section{Proof ideas and roadmap}
\label{sec:proofideas}
The proofs are executed in a unified way by considering the random probability density on the sphere,
\begin{equation}\label{eq:gibbs_measure0}
\mu(\theta) = \frac{1}{Z} \one\big\{|\theta_1/\sqrt{d} -x | \le \eps \big\} \prod_{i=1}^n \exp u\big(\langle g_i ,\theta \rangle / \sqrt{d}\big)\,, ~~~ \theta \in \cS^{d-1}(\sqrt{d})\,,
\end{equation}
where $x$ and $\eps>0$ are fixed, $u : \R \to \R$ is a smooth concave function with derivatives of controlled growth, and the $g_i$'s are i.i.d.\ random vectors whose first coordinates are copies of a random variable $S$ drawn from a general strongly log-concave density, and the remaining coordinates are standard normal. Without loss of generality we can set $\theta_\star = \sqrt{d} e_1$, and the density $\mu$ corresponds to each case of interest by appropriately choosing $u$ and the law of $S$; see Eqs.~\eqref{eq:defSgmm},~\eqref{eq:defSlogistic} for the (formal) correspondence with $\mu_{X,y}$ and Eqs.~\eqref{eq:defuSgmm},~\eqref{eq:defuSlogistic} for the correspondence with $p_{X,y}$. 

The goal is then to estimate the normalized logarithm of the restricted partition function $Z$, called the free energy.  We follow Talagrand's approach in proving the Gardner formula for the spherical perceptron model~\cite{talagrand2010mean,talagrand2011advanced}, 
with additional technical work required for handling the `signal part' of the data: the vectors $y_i X_i$ have a non-Gaussian component with law $S$ in the direction of $\theta_\star$ and are otherwise Gaussian on the complement subspace. This effectively turns the model into a spherical perceptron with Gaussian disorder and random inhomogeneous margins. 
While Talagrand's treatment is restricted to $\kappa \ge 0$, we are able to extend his techniques to our case. By constraining the overlap $\langle \theta, \theta_\star\rangle/d = \theta_1/\sqrt{d} \simeq x$ and leveraging the log-concavity of the law of $S$ we can control the fluctuations of the first coordinate to any desired accuracy. This allows us to conduct our analysis for $\kappa<0$ as long as $\alpha$ is suitably small.
We additionally remark that as a corollary of our proof, we prove the Gardner formula for the spherical perceptron model for all $\MARGIN \le 0$ and $\alpha \le \alpha_0$ smaller than a universal constant.

We mention that our analysis does not operate on the Nishimori line, in contrast to the recent activity in Bayesian settings; see for instance~\cite{lelarge2016fundamental,barbier2016mutual,barbier2019optimal}. As already mentioned in Section~\ref{sec:LDPpost}, Eq.~\eqref{eq:nishimori}, there is no special relation between the one-replica and two-replica overlaps $\langle \theta,\theta_\star\rangle$, $\langle \theta^1,\theta^2\rangle$ outside the Bayesian setting. Such a relation can simplify interpolation and concentration arguments, and usually lead to a simpler variational formula. The fact that we are outside this special line when studying the uniform measure $\mu_{X,y}$ on interpolators is responsible for the `sup-inf' structure of the variational formulas we obtain. Finally, we note that although our arguments are lengthy, all known proofs of the Gardner formula are technically involved \cite{talagrand2010mean, talagrand2011advanced, shcherbina2003rigorous, bolthausen2021gardner}.

In the remaining of this section we present at a high-level the main steps of the proof.

\subsection{Relaxation to the Gaussian reference measure}
In more detail, the first step is to relax the spherical constraint to full space by switching to a Gaussian reference measure, and we let     
\begin{equation}\label{eq:gibbs_measure1}
\mu_{\beta}(\theta) = \frac{1}{Z_\beta} \one\big\{|\theta_1/\sqrt{d} -x | \le \eps\big\}  \exp \left(\sum_{i=1}^n u\big(\langle g_i ,\theta \rangle / \sqrt{d}\big) -\beta \|\theta\|^2\right)\,, ~~~ \theta \in \R^d\,,
\end{equation}
now defined relative the Lebesgue measure on $\R^d$, and where $\beta >0$ is a parameter. The quadratic term acts as a confining potential.   

The measure $\mu_\beta$ is a log-concave density which has strong concentration properties. In particular the normalized squared norm $\|\theta\|^2/d$ of a draw $\theta \sim \mu_{\beta}$, and the normalized overlap $\langle \theta^1, \theta^2\rangle /d$ of two independent draws from $\mu_{\beta}$ both concentrate around their expectations. 
We then use the Gaussian interpolation method to prove a Gardner formula for $Z_\beta$ of the form 
\begin{equation}\label{eq:gardner_relaxed}
\frac{1}{d} \log Z_{\beta}  = \inf_{0 \le q \le \rho} \Big\{ \bar{\Phi}(x,q,\rho) - \beta(\rho + x^2)\Big\} + O(\eps) + o_{\P}(1)\,,
\end{equation}
where $o_{\P}(1) \to 0$ in probability in the limit $d\to \infty$ and $n/d \to \alpha$, and 
\begin{align}\label{eq:phi_master_rho}
\bar{\Phi}(x,q,\rho) 
&:=   \alpha \E \log \E_{W} \exp u\Big( x S + \sqrt{q}Z + \sqrt{\rho-q} W \Big)  + \frac{1}{2} \frac{q}{\rho-q} + \frac{1}{2} \log\big(\rho-q\big) \,,
\end{align}
This is the content of Appendix \ref{sec:interpolationfreeenergy}. Central to the successful application of this technique are 1) strong concentration of measure statements of the interpolating Gibbs measures and 2) uniqueness of the solution to the system of Replica-Symmetric (RS) equations Eq.\ (\ref{eq:niceRSeqs}). 
We demonstrate these in Appendix \ref{subsec:logconcaveconcentration} and Appendix \ref{sec:RSequationsunique} respectively. 

The unique minimizers $(q_\star,\rho_\star)$ of the formula in~\eqref{eq:gardner_relaxed} are respectively the asymptotic values of
\[\E \Big[\E_{(\theta^1,\theta^2) \sim \mu_\beta \times \mu_\beta} \big[\langle \bar{\theta}^1, \bar{\theta}^2\rangle /d\big]\Big]\,,~~~~ \mbox{and}~~~~ \E \Big[\E_{\theta \sim \mu_\beta} \big[\|\bar{\theta}\|^2/d\big]\Big]\,,\] 
where $\bar{\theta} = (\theta_2,\cdots,\theta_d)$ omits the first coordinate of $\theta$, and similarly for $\bar{\theta}^1, \bar{\theta}^2$. 

We remark that in the Bayesian setting with data generated from $\cD_{\gmm}$ we have $u(x) = \sqrt{\lambda} x$ which allows for explicit  computations of the log-partition function via a contour integral, and one can obtain the desired formula without resorting to these techniques. Nevertheless we do not pursue this and provide a unified treatment here.   

\subsection{Turning off the mollification} 
To obtain an asymptotic expression for the area of the set of interpolators with fixed overlap with the signal direction, $S_{\kappa}(X,y) \cap \{|\langle\theta,\theta_\star\rangle/d - x|\le \eps\}$, we turn off the mollification $\exp u(x)$ appearing in \eqref{eq:gibbs_measure1} and~\eqref{eq:phi_master_rho} into a hard step function $\one\{ x \ge \kappa\}$, under the Gaussian reference measure, and then convert back to the spherical setting. 
Regarding the first step, showing that the desired statement~\eqref{eq:gardner_relaxed} survives this passage to a hard step function is done by adding one constraint at a time. (Note the $i$-th datapoint imposes the constraint $\langle g_i, \theta \rangle \ge \kappa \sqrt{d}$ on the set of possible interpolators $\theta \in \cS^{d-1}(\sqrt{d})$, and henceforth we refer to the datapoints also as constraints.) Let
\begin{align}\label{eq:Ui_init_def} 
U_i = \big\{\theta \in \cS^{d-1}(\sqrt{d}): \langle g_i,\theta\rangle \ge \kappa \sqrt{d} \big\}\,,
\end{align}
and
\begin{equation}\label{eq:a_m}
a_m = \frac{1}{d}\E\log \int_{\bigcap_{i < m} U_i} \one\big\{|\theta_1/\sqrt{d} -x | \le \eps\big\}  \exp \left(\sum_{i=m}^n u\big(\langle g_i ,\theta \rangle / \sqrt{d}\big) -\beta \|\theta\|^2\right) \rmd \theta\,,~~~ m \le n\,.
\end{equation}
We want to estimate the difference $a_{m+1} - a_m$ which can be written as 
\begin{equation}\label{eq:difference}
a_{m+1} - a_m  = \frac{1}{d}\E \log G_{m}(U_{m}) - \frac{1}{d}\E \log G_{m}\big(e^{u(\langle g_m ,\theta \rangle / \sqrt{d})}\big)\,,
\end{equation}
where $G_m$ is the measure with density proportional to 
\[\one\Big\{\theta \in \bigcap_{i<m} U_i,\, |\theta_1/\sqrt{d}-x|\le \eps\Big\}\, \exp\left(\sum_{i>m}u(\langle g_i,\theta\rangle / \sqrt{d}) - \beta \|\theta\|^2\right)\,,\] 
which does not involve the $m$-th constraint. If $\exp u$ approximates the step function to accuracy $\eps'$ then we can hope to establish $|a_{m+1}-a_m | \le \eps'/d$, then sum over $m \in [n]$. The main difficulty is that  $G_m(U_m)$ may apriori be too small compared to the second term in~\eqref{eq:difference}; for instance the additional constraint $U_m$ may have an empty intersection with $\cap_{i<m}U_i$. Talagrand's approach is to truncate the logarithm in the definition of $a_m$, Eq.~\eqref{eq:a_m} to a minimum floor value $e^{-cd}$, $c>0$, and control the probability that $G_m(U_m)$ falls below this threshold. We will show that $G_m(U_m)$ is very unlikely to be exponentially small by extracting an exponential number of nearly orthogonal directions in $\cap_{i<m}U_i$, conditional on this set being large (this is provided by the truncation), and show that the fresh constraint $U_m$ has a small probability of being violated by all these directions simultaneously. This is a geometric argument executed via the tools of Gaussian processes. A minor adaptation is needed in our case to handle the non-Gaussian, but still strongly log-concave, first coordinate. This argument will also be used to establish concentration of the logarithms of the relevant partition functions, and is executed in Appendix \ref{sec:gaussianconversion}. We point out a version of this argument in \cite{nakajima2023sharp} which handles arbitrary i.i.d.\ sub-Gaussian vectors $g_i$, at the price of a suboptimal concentration rate that scales like a fractional power of $d$.           

\subsection{Conversion to the spherical reference measure}
We convert back from the Gaussian measure to the spherical measure by exploiting the concentration of $\bar{\theta}$ on a thin shell of width $o(\sqrt{d})$ around a sphere of radius approximately $\sqrt{\rho_\star d}$ when $\theta \sim \mu_{\beta}$. The parameter $\beta$ is then chosen to satisfy the unit radius condition $x^2 + \rho_\star = 1$. A formula for $(1/d)\log Z$, Eq.~\eqref{eq:gibbs_measure0} can then be obtained from~\eqref{eq:gardner_relaxed} via a volume-to-surface comparison argument. This eliminates the parameters $\beta$ and $\rho$ from the equation~\eqref{eq:phi_master_rho}, leading to a variational formula with $q$ as the only parameter. The overlap constraint on the first coordinate can be eliminated by taking a maximum over $x$, resulting in the `max-min' structure of the variational formula. 
This allows us to finish the proof of the large deviation result for the posterior distributions in the Bayesian setting (in Section~\ref{sec:LDPpost}) since the function $u$ already satisfies the required smoothness properties.
To obtain the LDP for uniform measure on the set of interpolators (in  Section~\ref{sec:LDPinterp}) we apply this argument after the mollification has been turned off. 
This is the content of Appendix \ref{sec:sphericalconversion}; see Theorems \ref{thm:sphericalmeasurersformula} and \ref{thm:posteriorrsformula}. This completes the proof of Theorems \ref{thm:mainresult} and \ref{thm:mainresultposterior}.

\paragraph{Statement on use of AI:} The authors used ChatGPT 5.1 to produce the code used for the numerical simulations presented in Section~\ref{sec:numerical}. The authors instructed the AI to write Julia code for the Monte Carlo experiment and the numerical solvers for the fixed point equations in $(x,q)$ appearing in Section~\ref{sec:main}. The code was thoroughly checked by the authors. The authors used ChatGPT 5.5 and Gemini 3.1 and 3.5 to establish Part~\ref{it:part2} in the proof of Propositions \ref{lem:RS_biv_eq_gmm}, \ref{lem:RS_biv_eq_logistic}, \ref{lem:RS_biv_bayes_gmm}, \ref{lem:RS_biv_bayes_logistic}, in writing the fixed point argument for the map $T_\alpha$, Eq.~\eqref{eq:fixed_point_map}, in the proof of Lemma~\ref{lem:nishimori_logistic}, and to prove Proposition \ref{prop:boundedapproxRSeq} and Lemma \ref{lem:posterior_spherical_continuity}. 

\paragraph{Acknowledgements:} AC would like to thank Shuangping Li for stimulating discussions. AE was supported by the National Science Foundation grant DMS-2450867.


\newcommand{\etalchar}[1]{$^{#1}$}

\newpage
\appendix
\section{Computing the expected log-partition function}\label{sec:interpolationfreeenergy}
Recall the definition of $Z_{\beta}$ from \equpref{eq:gibbs_measure1}. The goal of this section is to compute the free energy -- the expected log-partition function $(1/d) \E \log Z_{\beta}$ -- in the asymptotic regime $d \rightarrow \infty, \frac{n}{d}\rightarrow \alpha > 0$ for $\alpha \le \DENSITYBOUND$, where $\DENSITYBOUND$ is defined in \equpref{eq:alpha0def}.
Our expression depends on the solution to an explicit two-variable system of equations known as the \textit{Replica-Symmetric (RS) equations}, defined in \equpref{eq:niceRSeqs}.
Our precise results for the free energy are then given in Theorems \upref{thm:mollifiedhamiltonianfreeenergy}, \upref{thm:posterior_mollifiedhamiltonianfreeenergy}.  

\paragraph{Notation.} By rotational symmetry, we assume without loss of generality that $\theta_\star = \sqrt{d} e_1$ as per Section \ref{sec:proofideas}.
We then let $s_i=y_i X_{i1}$ for $1 \le i \le n$, the component of $y_i X_i$ in the span of $\theta_\star$.
Letting $\cD$ denote the density of the $s_i$, we impose the following technical assumption on $\cD$ which is satisfied when $\cD \in \{\cD_{\logistic}, \cD_{\gmm}\}$.
\begin{assumption}\label{ass:disorder_technical_assumption}
$S \sim \cD$ satisfies the following properties:
\begin{enumerate}
    \item $S$ has a $\SIGNALGAUSSIAN$ strongly log-concave density for some $\SIGNALGAUSSIAN>0$.
    \item $S$ is $\signalsubexp^{1/2}$ sub-Gaussian for $\signalsubexp \ge 1$, and the second and fourth moments of $S$ are upper bounded by $\signalsubexp$.
\end{enumerate}
\end{assumption}
Note $\SIGNALGAUSSIAN$, $\signalsubexp$ depend only on $\lambda$ (and the link function $\varphi$ when $\cD \equiv \cD_{\logistic}$, which we suppress for notational simplicity). We let $\pdfnormal$ denote the PDF of a standard univariate normal $N(0,1)$, and let $\pdfsignal$ denote the PDF of $S \sim \cD$.

We now formally define the quantities we are interested in understanding. 
Following the notation in Section \ref{sec:proofideas}, for any $\theta\in\R^d$, $\bar\theta = (\theta_2, \ldots, \theta_d)$ omits its first coordinate.
For all $1 \le i \le n$, we let
\begin{align}\label{eq:Si_def}
S_i = S_i(\theta) := \frac{1}{\sqrt{d}}\Big( s_i \theta_1 + \langle g_i, \bar \theta \rangle \Big)\,,
\end{align}
where $s_i \sim \cD$, $g_i \sim N(0, I_{d-1})$ are i.i.d., and where $g_{i,j}$ denotes the $j-1$-st coordinate of $g_i$. Hence we consider $2 \le j \le d$.
This notation is not to be confused with the definition of $g_i$ in Section \ref{sec:proofideas}, and for the rest of this paper, we define $g_i \sim N(0, I_{d-1})$.

Next, for an arbitrary function $u:\R \rightarrow \R$ such that the following integrals \equpref{eq:restricted_partition_function_mollified}, \equpref{eq:restricted_free_energy_mollified} exist, we define the Hamiltonian 
\begin{align}
H(\btheta) = H_{\SAMPLES, \DIMENSION}(\btheta) := \sum_{1 \le \SAMPLEIDX \le \SAMPLES} u(S_i) \,.\label{eq:pfs_hamiltonian_def}
\end{align}
Dependence on $u(\cdot)$ is suppressed here in Appendix \upref{sec:interpolationfreeenergy} and will always be clear from context. Since $\theta_{\star} = \sqrt{d} e_1$, $s_1 \theta_1 + \langle g_\SAMPLEIDX , \bar{\btheta}\rangle = y_i \langle X_i, \theta \rangle$. Thus \equpref{eq:pfs_hamiltonian_def} is identical in law to $\sum_{i=1}^n u\big(\langle y_i X_i \rangle/\sqrt{d}\big)$.
Finally, for any $\GAUSSIANMEASURE>0$ and any interval $I \subseteq \R$, we define the restricted quantities
\begin{align}
Z_{\beta}(I) &= Z_{\beta, \SAMPLES, \DIMENSION}(I) := \int \one\{\theta_1/\sqrt{\DIMENSION} \in I\} \, \exp\Big( H_{\SAMPLES, \DIMENSION}(\btheta) - \GAUSSIANMEASURE \|\btheta\|^2 \Big) \, \rmd \btheta\,, \label{eq:restricted_partition_function_mollified}\\
\phi_{\beta}(I) &= \phi_{\beta, \SAMPLES, \DIMENSION}(I) := \frac{1}{\DIMENSION}\E\log Z_{\beta, \SAMPLES, \DIMENSION}(I)\,. \label{eq:restricted_free_energy_mollified} 
\end{align}
In the above definitions, we omit dependence on $\SAMPLES, \DIMENSION$ when they are clear from context. 

We are now in position to state the main results of this section, regarding the free energy of the interpolators and posterior respectively. For the interpolators, our result applies for $u(\cdot)$ such that $\exp u(x)$ is as a sufficiently accurate approximation to $\one\{x \ge \MARGIN\}$.
\begin{theorem}[For interpolators]\label{thm:mollifiedhamiltonianfreeenergy}
Let $\alpha_0 = \alpha_0(\lambda, \MARGIN_+)$ be from \equpref{eq:alpha0def} and consider $n$ such that $1 \le n \le \DENSITYBOUND d$. 
Consider any $\eps>0$. 
Then there is a $\eps'>0$ depending on $\eps, \beta_0, \beta_1, \MARGIN, \lambda$ such that for any concave $u \le 0$ satisfying
\begin{enumerate}
    \item $u' \ge 0$ with $u' \neq 0$ on a set of positive Lebesgue measure,
    \item $u(x) = 0$ for $x \ge \MARGIN$,
    \item $\exp u(x) \le \eps'$ for all $x \le \MARGIN - \eps'$,
    \item $u$ is four times differentiable and $|u^{(l)}| \le D$ for $l = 1,2,3,4$,
\end{enumerate}
there is a $\eps_{\RS_0}>0$ depending on $\eps, u, \beta_0, \beta_1, \MARGIN, \SIGNALMEDIUMDEPENDENCE$ such that the following holds (note dependence on $\eps'$ comes in implicitly through $u$). For any $x \in [-1,1]$, $0 < \eps'' \le \eps_{\RS_0}$, and $\beta \in [\beta_0, \beta_1]$, letting $I = [x-\eps'', x +\eps'']$, for $\DIMENSION$ large enough in terms of $\eps, \eps'', \beta_0, \beta_1, \MARGIN, \SIGNALMEDIUMDEPENDENCE$, we have 
\begin{align}\label{eq:freeenergy_result}
\big| \phi_{\beta, \SAMPLES, \DIMENSION}(I) - \RSG_I(n/d, \GAUSSIANMEASURE, x) \big| \le \eps\,,
\end{align}
where $\phi_{\beta,n,d}(I)$ is from \equpref{eq:restricted_free_energy_mollified} and $\RSG_I$ is from \equpref{eq:RSGdef}.
\end{theorem}
For the posterior, $u(\cdot)$ is now fixed based on the assumed data generation process; our result applies for $u(\cdot)$ from \equpref{eq:defuSgmm} in the GMM case or \equpref{eq:defuSlogistic} in the logistic case.
\begin{theorem}[For posterior]\label{thm:posterior_mollifiedhamiltonianfreeenergy}
Let $\alpha_0 = \alpha_0(\lambda)$ be from \equpref{eq:alpha0def} and consider $n$ such that $1 \le n \le \DENSITYBOUND d$. Consider $u(x)$ from \equpref{eq:defuSgmm}, \equpref{eq:defuSlogistic} in the GMM and logistic cases respectively. Then for any $\eps>0$, there is $\eps_{\RS_0}>0$ depending on $\eps, \DENSITYBOUND, \beta_0, \beta_1, \lambda$ such that the following holds. For any $x \in [-1,1]$, $0 < \eps'' \le \eps_{\RSG}$, and $\beta \in [\beta_0, \beta_1]$, letting $I=[x-\eps'', x+\eps'']$, for $d$ large enough in terms of $\eps, \eps'', \beta_0, \beta_1, \lambda$, we have 
\begin{align}\label{eq:posteriorfreeenergy_result}
\big| \phi_{\beta, \SAMPLES, \DIMENSION}(I) - \RSG_I(n/d, \GAUSSIANMEASURE, x) \big| \le \eps\,,
\end{align}
where $\phi_{\beta, \SAMPLES, \DIMENSION}(I)$ is defined in \equpref{eq:restricted_free_energy_mollified} and $\RSG_I$ is from \equpref{eq:RSGdef}, both with this $\exp u(x)$.
\end{theorem}
In both Theorems, $\RSG_I$ is well-defined thanks to Proposition \upref{prop:nicersuniquesol}. Also, note that $d$ does not need to be large enough in terms of $\DENSITYBOUND$. This is because we apply the bound $\DENSITYBOUND \le 2$ in the bounds arising in the following proofs, and hence dependence on $\DENSITYBOUND$ in $d$ is captured by a universal constant. We now define $\alpha_0$ and $\RSG_I$ that appear in the above Theorems. \newline

\textbf{Defining $\alpha_0$:} For both the interpolators and the posterior, we let
\begin{align}
\GAUSSIANMEASURE_0 = \frac1{8}\, , \GAUSSIANMEASURE_1 = 10\, ,\Csol = \max\Big\{ \frac{1}{\GAUSSIANMEASURE_0}, 4\GAUSSIANMEASURE_1 \Big\}+1\,.\label{eq:csoldef}
\end{align}
For the interpolators, we now consider $\bar \alpha_0, c_0 \in (0, 2]$ small enough in terms of $\beta_0, \Csol, \SIGNALSHORTDEPENDENCE$, and define
\begin{align}
\DENSITYBOUND = \DENSITYBOUND(\lambda, \MARGIN_+) := \begin{cases} \bar \alpha_0 &: \MARGIN < 0\,, \\ \frac{c_0}{\MARGIN^2+1} &: \MARGIN \ge 0 \,. \end{cases}\label{eq:alpha0def}
\end{align}
For the posterior, consider $u(\cdot)$ from \equpref{eq:defuSgmm} in the GMM case or \equpref{eq:defuSlogistic} in the logistic case. In both cases, note $|u^{(l)}| \le D = D(\lambda)$ for $l=1,2,3,4$. We let $\Csol$ be large enough in terms of $D, \lambda$. Our results now apply for $\DENSITYBOUND \in (0,2]$ small enough in terms of $\beta_0, \beta_1, D, \Csol$, and therefore small enough in terms of $\beta_0, \beta_1, \lambda$. \newline

\textbf{Defining $\RSG$ and the RS equations:} We now define $\RSG$. To do so, we must define the \emph{Replica-Symmetric} (RS) equations.
Let $S \sim \cD$, $Z, W \sim N(0,1)$ be independent. Recall the definition of $\bar\Phi(x,q,\rho)$ for any $0 \le q < \rho$ and $x \in \R$ from \equpref{eq:phi_master_rho}. In what follows consider $\bar\Phi(x,q,\rho)$ with $u(x)$ therein given by $\exp u(x) = \one\{x \ge \MARGIN\}$ for the interpolators, and given by $u(\cdot)$ from \equpref{eq:defuSgmm} in the GMM case or \equpref{eq:defuSlogistic} in the logistic case for the posterior. In turn, for all such $(x, q, \rho)$ and any interval $I \subseteq \R$ we define
\begin{equation}
\label{eq:Fqrhodef}
\begin{aligned}
F_I(x,q,\rho) &:= \bar\Phi(x, q, \rho) - \beta \big(\rho + r_I(x)\big)\, \\
&= \alpha \E \log \E_{\ANNEALEDDISORDER} \exp u\big( x \SIGNALDISORDER + \sqrt{q} \QUENCHEDDISORDER + \sqrt{\rho-q} \ANNEALEDDISORDER \big)  + \frac{1}{2} \frac{q}{\rho-q} + \frac{1}{2} \log\big(\rho-q\big) \\
&\qquad \qquad - \beta \big(\rho + r_I(x)\big)\,,
\end{aligned}
\end{equation}
where 
\begin{align}\label{eq:intervalrxdef}
r_I(x) &:= \begin{cases} 0 &: 0 \in I\,, \\ x^2 &: 0 \not\in I\,.\end{cases}
\end{align}
and where dependence on $\alpha, \beta$ is implicit in the definition of $F$. Fixing $x \in [-1,1]$ and $I$, we define the RS equations by
\begin{align}
\frac{\partial F_I}{\partial q} = \frac{\partial F_I}{\partial \rho}=0\, .\label{eq:niceRSeqs}
\end{align}
We will establish in Appendix \upref{subsec:nice_rs_uniquesol_pf} the following crucial property of the system \equpref{eq:niceRSeqs}, which applies for both the interpolators and posterior.
\begin{proposition}\label{prop:nicersuniquesol}
For all $(\alpha, \GAUSSIANMEASURE, x) \in [0, \DENSITYBOUND] \times [\GAUSSIANMEASURE_0, \GAUSSIANMEASURE_1] \times [-1, 1]$, the system \equpref{eq:niceRSeqs} has a unique solution $\big(q_0(\alpha, \beta, x), \rho_0(\alpha, \beta, x) \big)=(q_0, \rho_0)$ in the domain $\big\{ (q, \rho): q < \rho, 0 \le q \le \Csol, \frac1{\Csol} \le \rho \le \Csol \big\}$. Also, we have $\rho_0(\alpha, \beta, x) - q_0(\alpha, \beta, x) \ge \frac1{2\Csol}$.

Moreover, $q_0(\alpha, \beta, x), \rho_0(\alpha, \beta, x)$ are both infinitely differentiable in $(\alpha, \beta, x)$ for $(\alpha, \beta, x) \in [0, \DENSITYBOUND] \times [\beta_0, \beta_1] \times [-1, 1]$. Thus for $x \in [-1, 1]$ and $I$ both fixed, $\RS_{0,I}(\alpha, \beta, x)$ and $\RSG_I(\alpha, \beta, x)$ are infinitely differentiable in $(\alpha, \beta)$ for $(\alpha, \beta) \in [0, \DENSITYBOUND] \times [\beta_0, \beta_1]$, where $\RS_{0, I}, \RSG_I$ are defined in \equpref{eq:RS0def}, \equpref{eq:RSGdef} below.
\end{proposition} 
\begin{remark}\label{rem:sol_nointervaldependence}
Note the system \equpref{eq:niceRSeqs} and its solution $\big(q_0(\alpha, \beta, x), \rho_0(\alpha, \beta, x) \big)$ do \textit{not} depend on $I$, since dependence on $I$ in $F_I$ comes only through $\beta r_I(x)$, which does not depend on $\rho, q$. 
\end{remark}
Define in terms of $F_I$ the following quantities, which implicitly depend on $u(\cdot)$:
\begin{align}
\RS_{0, I}(\alpha, \beta, x) &:= F_I\big(x, q_0(\alpha, \beta, x), \rho_0(\alpha, \beta, x)\big)\,,\label{eq:RS0def}\\
\RSG_I(\alpha, \beta, x) &:= \RS_{0, I}(\alpha , \beta, x) + \frac12 \log(2\pi e)
\,.\label{eq:RSGdef}
\end{align}

Theorems \upref{thm:mollifiedhamiltonianfreeenergy} and \upref{thm:posterior_mollifiedhamiltonianfreeenergy} are proven with a unified strategy that uses a Gaussian interpolation argument. 
This argument proceeds by understanding the expected contribution of adding each datapoint or constraint $i$ on the free energy $(1/d) \E\log Z_{\beta}$ for each $1 \le i \le n$. 
See Proposition \upref{prop:interpolationperconstraint}.
The free energy then can be written as a sum of such expected contributions. 

This interpolation argument is carried out in Appendix \upref{subsec:interpolation1}.
We then show in Appendix \upref{subsec:RSequationseliminatemollifier} that under the conditions on $u(x)$ in Theorem \upref{thm:mollifiedhamiltonianfreeenergy}, we can replace $u(x)$ by $\one\{x \ge \MARGIN\}$ only changes the the free energy by at most $\eps$.
This yields an expression that does not depend on $u(\cdot)$ and depends explicitly on the solution to the RS equations.
Finally, we sum these expressions across all the datapoints $i, 1 \le i \le n$ and finish the proof of Theorem \upref{thm:mollifiedhamiltonianfreeenergy} for the interpolators and Theorem \ref{thm:posterior_mollifiedhamiltonianfreeenergy} in Appendix \upref{subsec:finishinterpolation} and Appendix \upref{subsec:posterior_finishinterpolation} respectively.

Crucial to executing this proof strategy is Proposition \upref{prop:nicersuniquesol} -- that the RS equations have a unique solution.
This is done by showing that the RS equations satisfy a particular convex-concave structure.
Establishing this with the signal in the data and the case $\MARGIN < 0$ does not directly follow from the arguments of Talagrand in \cite{talagrand2010mean}; we discuss our steps to handle these complications in Appendix \upref{sec:RSequationsunique}.

\subsection{Interpolating each datapoint via cavity in $n$}\label{subsec:interpolation1}
Here in Appendix \upref{subsec:interpolation1}, we fix any $\DIMENSION \ge 1$ and $\SAMPLES$ such that $1 \le \SAMPLES \le \DENSITYBOUND \DIMENSION$, and define $\alpha = \frac{\SAMPLES}{\DIMENSION}$. 
To compute the free energy $\phi_{\beta, \SAMPLES, \DIMENSION}(I)$, we consider two interpolations. 
One interpolation is over a parameter $v, 0 \le v \le 1$, and aims to elucidate the effect of the $n$-th datapoint or constraint on the original probability distribution. 
When $v=1$, the $n$-th datapoint is fully present, and when $v=0$, it is eliminated.
Hence this interpolation is called `cavity in $n$'.
The other interpolation is over a parameter $t, 0 \le t \le 1$, and aims to elucidate the effect of the $d$-th coordinate.
When $t=1$, the $d$-th coordinate is fully present, and when $t=0$, it is eliminated.
Hence this interpolation is called `cavity in $d$'. 

Specifically, the cavity in $n$ provides us with an approximate expression for $\phi_{\beta, n, d}(I) - \phi_{\beta, n-1, d}(I)$.
This expression depends on the mollifier $u(\cdot)$ and the overlaps $q_{n,d}, \rho_{n,d}$ given in \equpref{eq:interpolationrhoqdef}; see Proposition \upref{prop:interpolationperconstraint}.
Establishing this result is the main goal of Appendix \upref{subsec:interpolation1}. 
However to define the cavity in $n$, we will first need to define the cavity in $d$. 
We will then use the cavity in $d$ more extensively next in Appendix \upref{subsec:RSequationseliminatemollifier} to relate $q_{n,d}, \rho_{n,d}$ in the expression from Proposition \upref{prop:interpolationperconstraint} to $q_0, \rho_0$ from Proposition \upref{prop:nicersuniquesol}, and eliminate $u(\cdot)$ in the case of the interpolators.
This lets us rewrite the aforementioned expression for $\phi_{\beta, n, d}(I) - \phi_{\beta, n-1, d}(I)$ in terms of $\RSG_I$.

Finally, we note that many of our following bounds do not explicitly depend on $\DENSITYBOUND$; this is because we use the upper bound $\DENSITYBOUND \le 2$.

\paragraph{Cavity in $d$:} To define the cavity in $n$ in its full generality, we must first set up the cavity in $d$.
For any $1 \le \SAMPLEIDX \le \SAMPLES$ and $t \in [0,1]$, define
\begin{align}
S_{\SAMPLEIDX, t} = S_{\SAMPLEIDX, t}(\btheta) := \frac{s_\SAMPLEIDX \theta_1}{\sqrt{\DIMENSION}}+\frac1{\sqrt{\DIMENSION}} \sum_{2 \le \DIMENSIONIDX \le \DIMENSION-1} g_{\SAMPLEIDX, \DIMENSIONIDX} \theta_\DIMENSIONIDX + \sqrt{\frac{t}{\DIMENSION}} g_{\SAMPLEIDX, \DIMENSION} \theta_\DIMENSION\,.\notag
\end{align} 
In terms of $S_{\SAMPLEIDX, t}$, define
\begin{align}
H_{t, \SAMPLES, \DIMENSION}(\btheta) := \sum_{1 \le \SAMPLEIDX \le \SAMPLES} u\big(S_{\SAMPLEIDX, t}(\btheta) \big) + \theta_\DIMENSION\sqrt{(1-t)r} Y - \frac{(1-t)(r-\bar{r})}2 \theta_\DIMENSION^2\,,\label{eq:interpolateinthamiltonian}
\end{align}
where $r,\bar{r}$ will be chosen next in \equpref{eq:interpolationrrbardef1}, \equpref{eq:interpolationrrbardef2}, and where $Y \sim N(0, 1)$.
Let $\langle \,\cdot\, \rangle_t$ denote the Gibbs average w.r.t. 
\[ \one\{\theta_1/\sqrt{\DIMENSION} \in I\} \, \exp\big( H_{t, \SAMPLES, \DIMENSION}(\btheta) - \GAUSSIANMEASURE \| \btheta \|^2 \big)\,,\]
and let
\begin{align}\label{eq:nu_t_def}
\nu_t(\cdot) := \E \langle \,\cdot\, \rangle_t\,.
\end{align}
When we write $H_t$ without explicitly stating $\SAMPLES$, we will assume this refers to $H_{t,n,d}$ as written in \equpref{eq:interpolateinthamiltonian}.
Now consider the overlaps
\begin{align}\label{eq:overlap_defs_thispaper}
R_{1,2} = \frac{1}{\DIMENSION} \langle \bar{\btheta}^1,\bar{\btheta}^2\rangle\,,\,R_{1,1} = \frac{1}{\DIMENSION}\|\bar{\btheta}\|^2\,,
\end{align}
where we view $R_{1, 1}$ as a function on $\R^{d-1}$ and $R_{1, 2}$ as a function on $(\R^{\DIMENSION-1})^2$.
Let 
\begin{align}
q_t := \nu_t(R_{1,2})\, \le\, \rho_t := \nu_t(R_{1,1}) \,.\label{eq:interpolationrhoqdefgeneralt}
\end{align}
Note $\langle \,\cdot\, \rangle_1$ corresponds to exactly the Hamiltonian $H_{\SAMPLES, \DIMENSION}(\btheta) - \GAUSSIANMEASURE \| \btheta \|^2$. Indeed, specializing to $t=1$, we define
\begin{align}
q_{\SAMPLES, \DIMENSION} :=  q_1 \, ,\, \rho_{\SAMPLES, \DIMENSION} := \rho_1\,,\label{eq:interpolationrhoqdef}
\end{align}
Consequently $q_{\SAMPLES, \DIMENSION} \le \rho_{\SAMPLES, \DIMENSION}$. Now consider independent $\QUENCHEDDISORDER, \ANNEALEDDISORDER \sim N(0, 1)$, and $\SIGNALDISORDER \sim \cD$ following the law of each $s_\SAMPLEIDX$. For arguments $0 \le q < \rho$, define
\begin{align}\label{eq:approxconstraint}
\APPROXCONSTRAINT = \APPROXCONSTRAINT(q, \rho) := x \SIGNALDISORDER + \sqrt{q}\QUENCHEDDISORDER + \sqrt{\rho - q}\ANNEALEDDISORDER\, .
\end{align}
We define the functions
\begin{align}
\psi_{\alpha}(q, \rho) &:= \alpha \E \Big( \frac{\E_{\ANNEALEDDISORDER} [ u'(\APPROXCONSTRAINT) \exp u(\APPROXCONSTRAINT) ] }{ \E_{\ANNEALEDDISORDER}[ \exp u(\APPROXCONSTRAINT) ] } \Big)^2\, ,\label{eq:psidef}\\
\bar{\psi}_{\alpha}(q, \rho) &:= \alpha \E \Big[ \frac{\E_{\ANNEALEDDISORDER} [ (u''(\APPROXCONSTRAINT) + u'(\APPROXCONSTRAINT)^2) \exp u(\APPROXCONSTRAINT) ] }{ \E_{\ANNEALEDDISORDER}[ \exp u(\APPROXCONSTRAINT) ] } \Big]\, .\label{eq:psibardef} 
\end{align}
In turn, we let
\begin{align}
r &= \psi_{\alpha}\big(q,\rho\big)=\alpha\E \Bigg[ \Bigg( \frac{\E_{\ANNEALEDDISORDER}\big[ u'(\APPROXCONSTRAINT)\exp u(\APPROXCONSTRAINT) \big]}{\E_{\ANNEALEDDISORDER}\big[ \exp u(\APPROXCONSTRAINT)\big]} \Bigg)^2 \Bigg] = \frac{\alpha}{\rho-q} \E \Bigg[ \Bigg( \frac{\E_{\ANNEALEDDISORDER}\big[ \ANNEALEDDISORDER \exp u(\APPROXCONSTRAINT) \big]}{\E_{\ANNEALEDDISORDER}\big[ \exp u(\APPROXCONSTRAINT) \big]} \Bigg)^2 \Bigg]\,,\label{eq:interpolationrrbardef1} \\
\bar{r} &= \bar{\psi}_{\alpha} \big(q, \rho \big) = \alpha \E \frac{\E_{\ANNEALEDDISORDER}\big[ (u''(\APPROXCONSTRAINT) + u'(\APPROXCONSTRAINT)^2 ) \exp u(\APPROXCONSTRAINT)\big] }{\E_{\ANNEALEDDISORDER}\big[ \exp u(\APPROXCONSTRAINT)\big]} = \frac{\alpha}{\rho - q} \E \frac{\E_{\ANNEALEDDISORDER}\big[ (\ANNEALEDDISORDER^2 - 1) \exp u(\APPROXCONSTRAINT)\big] }{\E_{\ANNEALEDDISORDER}\big[ \exp u(\APPROXCONSTRAINT)\big]}\,,\label{eq:interpolationrrbardef2}
\end{align}
where the second equality is only valid when $\rho>q$.
Here the equality follows by Gaussian Integration by Parts. 
As we will justify later in Lemma \upref{lem:realRSeqsolsboundedtechnicalcheck}, we have $r \ge \bar{r}$ (irrespective of whether $\rho=q$ or $\rho>q$). As $|u'|, |u''| \le D$ in both the interpolators and GMM case and as $\alpha \le \DENSITYBOUND \le 2$, we have $r, \bar{r} \le K(D)$, which we use repeatedly in the following. We also note $r, \bar{r}$ are independent of $\eps''$.

\paragraph{Cavity in $\SAMPLES$:} We next consider another interpolation parameter $v \in [0,1]$.
To set up the cavity in $\SAMPLES$, we need to introduce the notion of \textit{replicas}. 
Given a probability distribution $\nu$, we consider $\REPLICAS$ i.i.d. draws $\theta^{\REPLICAIDX}, 1 \le \REPLICAIDX \le \REPLICAS$ from $\nu$.
These i.i.d. draws are known $\theta^{\REPLICAIDX}$ are known as \textit{replicas}.
That is, $(\theta^{\REPLICAIDX})_{\ell=1}^L \sim \nu^{\otimes \REPLICAS}$.
It is often useful to consider functions of $\REPLICAS > 1$ replicas; one example is the overlap $R_{1,2} = \frac1d \langle \bar\theta^1, \bar\theta^2 \rangle$.
In particular, functions of $\REPLICAS > 1$ replicas naturally arise when differentiating quantities of interest such as $R_{1,1}$ and $R_{1,2}$ with respect to $t$ and $v$ in the cavity in $d$ and $n$.

Proceeding with the setup, we let $\QUENCHEDDISORDER$ and $\ANNEALEDDISORDER^{\REPLICAIDX}$ denote i.i.d. $N(0, 1)$ scalars independent of everything else for each replica $\REPLICAIDX$, $1 \le \REPLICAIDX \le \REPLICAS$. 
For $v \in [0,1]$, define for each replica $\REPLICAIDX$, $1 \le \REPLICAIDX \le \REPLICAS$,
\begin{align}
S_{\SAMPLES, t, v}^\REPLICAIDX = S_v^\REPLICAIDX &:= \sqrt{v} \Big(\frac1{\sqrt{\DIMENSION}} \sum_{2 \le \DIMENSIONIDX \le \DIMENSION-1} g_{\SAMPLES, \DIMENSIONIDX} \theta^{\REPLICAIDX}_{\DIMENSIONIDX} + \sqrt{\frac{t}{\DIMENSION}} g_{\SAMPLES, \DIMENSION} \theta^{\REPLICAIDX}_{\DIMENSION} \Big) + \sqrt{1-v} \big( \sqrt{q_{\SAMPLES, \DIMENSION}} \QUENCHEDDISORDER + \sqrt{\rho_{\SAMPLES, \DIMENSION} -q_{\SAMPLES, \DIMENSION}} \ANNEALEDDISORDER^\REPLICAIDX \big) \notag \\
&\qquad \qquad + s_\SAMPLES \Big(v\frac{\theta^\REPLICAIDX_1}{\sqrt{\DIMENSION}} + (1-v) x \Big) \,.\notag
\end{align}
Note that in the following when we consider $S_{\SAMPLES, t, v}^\REPLICAIDX$, $t$ will be always held fixed. When it is clear from context we will just write $S_{v}^\REPLICAIDX$. Note $S_{\SAMPLES, t, 1}^\REPLICAIDX = S_{\SAMPLES, t}^\REPLICAIDX$. 

Next, we let $\langle\, \cdot\, \rangle_{t, \sim}$ denote a Gibbs measure defined by 
\begin{align}
\langle f \rangle_{t, \sim} &:= \frac1{\NORMALIZEINTERPOLATION^L} \int \E_{\ANNEALEDDISORDER} \Bigg[ f\, \prod_{1 \le \REPLICAIDX \le L} \one\{\theta^{\REPLICAIDX}_1/\sqrt{\DIMENSION} \in I\} \exp\big( H_{t, \SAMPLES-1, \DIMENSION}(\btheta^{\REPLICAIDX}) - \GAUSSIANMEASURE \| \btheta^{\REPLICAIDX} \|^2 \big) \Bigg]\, \rmd \btheta^1 \cdots \rmd \btheta^L\label{eq:interpolationtildemeasure}
\end{align}
for $f$ a function of $\REPLICAS$ replicas and the corresponding annealed disorder $\ANNEALEDDISORDER^{\REPLICAIDX}$, where we define:
\begin{align}
\NORMALIZEINTERPOLATION &:= \int \one\{\theta_1/\sqrt{\DIMENSION} \in I\} \exp\big( H_{t, \SAMPLES-1, \DIMENSION}(\btheta) - \GAUSSIANMEASURE \| \btheta \|^2 \big) \, \rmd \btheta\,,\notag \\
H_{t, \SAMPLES-1, \DIMENSION}(\btheta) &:= \sum_{1 \le \SAMPLEIDX \le \SAMPLES-1} u(S_{\SAMPLEIDX, t}(\btheta)) + \theta_\DIMENSION\sqrt{(1-t)r} Y - \frac{(1-t)(r-\bar{r})}2 \theta_{\DIMENSION}^2\,.
\end{align}
(Note for a fixed value of $\REPLICAS$, if $f$ is a function of $\REPLICAS+1$ replicas, the Gibbs measure \equpref{eq:interpolationtildemeasure} will use $\REPLICAIDX+1$ replicas $\REPLICAIDX=1, \ldots, \REPLICAIDX, \REPLICAIDX+1$.)
Now let 
\begin{align}\label{eq:nu_t_v_def}
\nu_{t, v}(f) := \E \frac{\big\langle f \exp \big( \sum_{1 \le \REPLICAIDX \le \REPLICAS} u(S_v^\REPLICAIDX) \big) \big\rangle_{t, \sim} }{\big\langle \exp u(S_v^1)\big\rangle_{t, \sim}^\REPLICAS } = \E \frac{\big\langle f \exp \big( \sum_{1 \le \REPLICAIDX \le \REPLICAS} u(S_{\SAMPLES, t, v}^\REPLICAIDX) \big) \big\rangle_{t, \sim} }{ \big\langle \exp u(S_{\SAMPLES, t, v}^{\REPLICAS+1})\big\rangle_{t, \sim}^\REPLICAS }\,.
\end{align}
Observe that $\nu_{t} \equiv \nu_{t,1}$. 
We will establish that as one datapoint is removed, the resulting difference in the free energy approximately equals the following expression:
\begin{proposition}\label{prop:interpolationperconstraint}
For all $d \ge d_{\upref{prop:interpolationperconstraint}}(\eps'') \ge \frac1{\eps''^2}$, we have 
\begin{align}
&\Big| \phi_{\beta, \SAMPLES, \DIMENSION}(I) - \phi_{\beta, \SAMPLES-1, \DIMENSION}(I) - \frac1{\DIMENSION} \E \Big[ \log \E_{\ANNEALEDDISORDER} \exp u\big( x s_\SAMPLES + \sqrt{q_{\SAMPLES, \DIMENSION}} \QUENCHEDDISORDER + \sqrt{\rho_{\SAMPLES, \DIMENSION}-q_{\SAMPLES, \DIMENSION}}\ANNEALEDDISORDER \big) \Big] \Big| \le \frac{K_{\upref{prop:interpolationperconstraint}}\eps''}{\DIMENSION}\,,\label{eq:interpolationperconstrainteq}
\end{align}
where $K_{\upref{prop:interpolationperconstraint}}$ depends on $D, \GAUSSIANMEASURE_0, \GAUSSIANMEASURE_1, \SIGNALMEDIUMDEPENDENCE$, and also on $\MARGIN$ for the interpolators.
\end{proposition}
To prove Proposition \upref{prop:interpolationperconstraint}, we will need to establish the following.
First, in Appendix \upref{subsec:pf_interpolationkeyconcentration}, we will show the following critical Lemma, which lets us transfer expectations from $\nu_{t, v}$ to $\nu_{t, 1}$.
\begin{lemma}\label{lem:interpolationkeyconcentration}
For any function $f:(\R^{\DIMENSION-1})^\REPLICAS \rightarrow \R$ with $f \ge 0$,
\begin{align}
\nu_{t, v}(f) \le K \Big( \nu_{t, 1}(f) + \Big( \eps'' + \frac1d \Big) \sup_{t, v \in [0, 1]^2} \nu_{t, v}(f^2)^{1/2} \Big)\,,
\end{align}
where $K$ depends on $D, \GAUSSIANMEASURE_0, \GAUSSIANMEASURE_1, \SIGNALSHORTDEPENDENCE, \REPLICAS$, and also on $\MARGIN$ for the interpolators.
\end{lemma}  
We will apply the above for $f=|R_{1, 1}-\rho_{\SAMPLES, \DIMENSION}|, |R_{1,2}-q_{\SAMPLES, \DIMENSION}|$, and thus are left with the task of controlling $\nu_{t, 1}\big( |R_{1, 1}-\rho_{\SAMPLES, \DIMENSION}| \big), \nu_{t, 1}\big( |R_{1, 2}-q_{\SAMPLES, \DIMENSION}| \big)$ that appear on the right hand side of the above Lemma. To this end, we will also establish the following via concentration of Lipschitz functions of log-concave measures in Appendix \upref{subsec:logconcaveconcentration}.
\begin{proposition}\label{prop:logconcaveconcentration}
For all $d \ge d_{\upref{prop:logconcaveconcentration}}(\eps'')$, all $t, v \in [0,1]$ and all $k \le \DIMENSION/4$, for $f \equiv R_{1,1}$ or $R_{1,2}$, 
\begin{align*}
\nu_{t,v}\Big( \big(f - \nu_{t,v}(f) \big)^{2k} \Big) \le \Big( \frac{K k}{\DIMENSION} \Big)^k\,.
\end{align*}
and 
\begin{align*}
\log \nu_{t,v}\Big( \exp \Big( \frac{\| \bar{\btheta} \|^2}{K} \Big) \Big) &\le K \DIMENSION\quad \,, \quad \log \nu_{t,v}\Big( \exp \Big( \frac{\theta_{\DIMENSION}^2}{K} \Big) \Big) \le K\,,
\end{align*}
where $K$ depends on $D, \GAUSSIANMEASURE_0, \GAUSSIANMEASURE_1, \SIGNALSHORTDEPENDENCE$, and also on $\MARGIN$ for the interpolators.
\end{proposition}
We also cite a simple fact to bound higher order moments of relevant random variables:
\begin{lemma}[Lemma 3.1.8, \cite{talagrand2010mean}]\label{lem:mgftomomentslemma}
Consider a random variable $X \ge 0$ and let $C = \log \E \big[ \exp X \big]$. Then for each $k \ge 1$ we have
\[ \E\big[ X^k \big] \le 2^k(k^k+C^k)\,.\]
\end{lemma}
\vspace{1pt}
Now we have the necessary tools to prove Proposition \upref{prop:interpolationperconstraint}.
\begin{proof}[Proof of Proposition \upref{prop:interpolationperconstraint}]
We let $d_{\upref{prop:interpolationperconstraint}}(\eps'') = d_{\upref{prop:logconcaveconcentration}}(\eps'')$ and consider any $d \ge d_{\upref{prop:interpolationperconstraint}}(\eps'')$. Define for the proof of Proposition \upref{prop:interpolationperconstraint} the following function of $v$:
\begin{align*}
\varphi(v) &:= \E \Big[ \log \big\langle \exp u(S_{\SAMPLES, 1, v}) \big\rangle_{1,\sim} \Big]\,.
\end{align*}
(The above function is not to be confused with the link function $\varphi$, which does not explicitly appear in the proof of Proposition \upref{prop:interpolationperconstraint}.)
Notice that 
\begin{align*}
\frac1{\DIMENSION} \varphi(1) &= \phi_{\beta, \SAMPLES, \DIMENSION}(I) - \phi_{\beta, \SAMPLES-1,\DIMENSION}(I) \, ,\\
\varphi(0) &= \E \Big[ \log \E_{\ANNEALEDDISORDER} \exp u\big( x s_\SAMPLES  + \sqrt{q_{\SAMPLES, \DIMENSION}}\QUENCHEDDISORDER + \sqrt{\rho_{\SAMPLES, \DIMENSION}-q_{\SAMPLES, \DIMENSION}}\ANNEALEDDISORDER \big) \Big] \,.
\end{align*}
Consequently it remains to prove that 
\begin{align*}
\Big| \frac{\rmd}{\rmd v} \varphi(v) \Big| \le K(D, \GAUSSIANMEASURE_0, \GAUSSIANMEASURE_1, \MARGIN, \SIGNALMEDIUMDEPENDENCE) \Big( \eps''+\frac1{\sqrt{\DIMENSION}} \Big)\,,
\end{align*}
where here and in the rest of the proof of Proposition \upref{prop:interpolationperconstraint}, there is no dependence in $K$ on $\MARGIN$ for the posterior. To this end, we compute 
\begin{align*}
\frac{\rmd}{\rmd v} \varphi(v) 
&= \E \frac{\frac{\rmd}{\rmd v}\big\langle \exp u(S_{\SAMPLES, 1, v}) \big\rangle_{1,\sim}}{\big\langle \exp u(S_{\SAMPLES, 1, v}) \big\rangle_{1,\sim}} \\
&= \E \frac{\big\langle \exp u(S_{\SAMPLES, 1, v}) u'(S_{\SAMPLES, 1, v}) \cdot s_{\SAMPLES}\big(\frac{\theta_1}{\sqrt{\DIMENSION}} - x\big) \big\rangle_{1,\sim}}{\big\langle \exp u(S_{\SAMPLES, 1, v}) \big\rangle_{1,\sim}} \\
&\qquad +\E \frac{\Big\langle \exp u(S_{\SAMPLES, 1, v}) u'(S_{\SAMPLES, 1, v})\cdot  \big( \frac1{2\sqrt{v\DIMENSION}} \big\langle g_M, \bar{\btheta} \rangle - \frac1{2\sqrt{1-v}} \big( \sqrt{q_{\SAMPLES, \DIMENSION}}\QUENCHEDDISORDER + \sqrt{\rho_{\SAMPLES, \DIMENSION}-q_{\SAMPLES, \DIMENSION}}\ANNEALEDDISORDER \big) \big) \Big\rangle_{1,\sim}}{\big\langle \exp u(S_{\SAMPLES, 1, v}) \big\rangle_{1,\sim}}\\
&:= \expressionI + \expressionII\,.
\end{align*}
Next observe as $\exp(\cdot) \ge 0$, we have
\begin{align*}
\big| \expressionI \big| &\le  \E  \frac{\Big|\big\langle \exp u(S_{\SAMPLES, 1, v}) u'(S_{\SAMPLES, 1, v}) \cdot s_{\SAMPLES}\big(\frac{\theta_1}{\sqrt{\DIMENSION}} - x\big) \big\rangle_{1,\sim} \Big|}{\big\langle \exp u(S_{\SAMPLES, 1, v}) \big\rangle_{1,\sim}} \\
&\le \E  \frac{\Big\langle \big| \exp u(S_{\SAMPLES, 1, v}) u'(S_{\SAMPLES, 1, v}) \cdot s_{\SAMPLES}\big(\frac{\theta_1}{\sqrt{\DIMENSION}} - x\big) \big| \Big\rangle_{1,\sim} }{\big\langle \exp u(S_{\SAMPLES, 1, v}) \big\rangle_{1,\sim}}\,.
\end{align*}
Next remark, with the quenched disorder fixed, we have
\begin{align*}
&\Big\langle \bigl| \exp u(S_{\SAMPLES, 1, v}) u'(S_{\SAMPLES, 1, v}) \cdot s_{\SAMPLES}\bigl(\tfrac{\theta_1}{\sqrt{\DIMENSION}} - x\bigr) \bigr| \Big\rangle_{1,\sim} \\
&\quad = \frac{1}{Z_{t,n-1,d}} \int \E_{\ANNEALEDDISORDER}\Big[ \one\{\theta_1/\sqrt{\DIMENSION} \in I \}  \exp\Big( H_{1,\SAMPLES-1}(\btheta) - \GAUSSIANMEASURE \| \btheta \|^2 \Big) \\
&\qquad\qquad\qquad\qquad\qquad\qquad\qquad
\cdot \exp u(S_{\SAMPLES, 1, v}) \cdot \bigl|u'(S_{\SAMPLES, 1, v})\bigr| \cdot \bigl|\tfrac{\theta_1}{\sqrt{\DIMENSION}} - x \bigr| \cdot |s_{\SAMPLES}|\Big]\, \rmd\btheta \\
&\quad \le D \eps'' |s_{\SAMPLES}| \cdot \big\langle \exp u(S_{\SAMPLES, 1, v}) \big\rangle_{1,\sim}\, .
\end{align*}
Here, the last inequality follows as $s_{\SAMPLES}$ is independent of $\ANNEALEDDISORDER$ and $\btheta$, and as the numerator of the last integral is only nonzero when $\theta_1/\sqrt{\DIMENSION} \in I$, which implies $\big| \frac{\theta_1}{\sqrt{\DIMENSION}} - x \big| \le \eps''$. Combining with the earlier display implies
\begin{align*}
\big| \expressionI \big| \le \E \Bigg[ \frac{ D \eps'' |s_{\SAMPLES}| \cdot \big\langle \exp u(S_{\SAMPLES, 1, v}) \big\rangle_{1,\sim} }{ \big\langle \exp u(S_{\SAMPLES, 1, v}) \big\rangle_{1,\sim} } \Bigg] = D \eps'' \cdot \E[ |s_{\SAMPLES}| ] \le K(D, \SIGNALSHORTDEPENDENCE) \eps''\,.
\end{align*}
Moreover, applying Gaussian Integration by Parts on $g_n$, $\QUENCHEDDISORDER$, and $\ANNEALEDDISORDER$ yields that
\begin{align*}
\big| \expressionII \big| &\le K(D) \nu_{1, v}\Big( \big|R_{1,1} - \rho_{\SAMPLES, \DIMENSION} \big| + \big| R_{1,2}-q_{\SAMPLES, \DIMENSION}\big| \Big)\,. 
\end{align*}
The calculation in the above equality is analogous as in the proof of Lemma 3.3.5 of \cite{talagrand2010mean}. The reason why the argument still applies is that the Gaussian random variables we integrate by parts with respect to are the exact same as in \cite{talagrand2010mean}; there is no integration by parts w.r.t. the disorder from the signal part $\expressionI$, which we instead bounded directly above. See the proof of Lemma \upref{lem:cavityinMintegrationbyparts} for the details of a similar computation.

We thus have, applying Lemma \upref{lem:interpolationkeyconcentration} to the non-negative function $\big|R_{1,1} - \rho_{\SAMPLES, \DIMENSION} \big| + \big|R_{1,2} - q_{\SAMPLES, \DIMENSION} \big|$, 
\begin{align*}
&\Big| \frac{\rmd}{\rmd v} \varphi(v) \Big| \\
&\quad\le \big|\expressionI \big| + \big| \expressionII \big| \\ 
&\quad \le K(D, \SIGNALSHORTDEPENDENCE) \eps'' + K(D) \Big( \nu_{1,v} \big( \big|R_{1,1} - \rho_{\SAMPLES, \DIMENSION} \big| + \big|R_{1,2} - q_{\SAMPLES, \DIMENSION} \big| \big) \Big) \\
&\quad\le K(D, \GAUSSIANMEASURE_0, \GAUSSIANMEASURE_1, \MARGIN, \SIGNALSHORTDEPENDENCE) \Bigg( \nu_{1,1} \big( \big|R_{1,1} - \rho_{\SAMPLES, \DIMENSION} \big| + \big|R_{1,2} - q_{\SAMPLES, \DIMENSION} \big| \big) \\
&\qquad\qquad\qquad\qquad\qquad\qquad+ \Big( \eps'' + \frac1{\DIMENSION} \Big) \sup_{t,v\in [0,1]^2} \Big( \nu_{t,v} \big( (R_{1,1} - \rho_{\SAMPLES, \DIMENSION})^2 \big)^{1/2} + \nu_{t,v} \big( (R_{1,2} - q_{\SAMPLES, \DIMENSION})^2 \big)^{1/2} \Big) \Bigg) \\
&\qquad\qquad\qquad+ K(D, \SIGNALSHORTDEPENDENCE) \eps'' \,.
\end{align*}
As $\nu_{1,1}(R_{1,1}) = \rho_{\SAMPLES, \DIMENSION}$, $\nu_{1,1}(R_{1,2}) = q_{\SAMPLES, \DIMENSION}$, we thus obtain from the first part of Proposition \upref{prop:logconcaveconcentration} that
\begin{align*}
\nu_{1,1} \big( \big|R_{1,1} - \rho_{\SAMPLES, \DIMENSION} \big| + \big|R_{1,2} - q_{\SAMPLES, \DIMENSION} \big| \big) \le \frac{K(D, \GAUSSIANMEASURE_0, \GAUSSIANMEASURE_1, \MARGIN, \SIGNALMEDIUMDEPENDENCE)}{\sqrt{\DIMENSION}}\,.
\end{align*}
The second part of Proposition \upref{prop:logconcaveconcentration} and Lemma \upref{lem:mgftomomentslemma} implies $\nu_{t,v}\big( R_{1,2}^2 \big) \le K(D, \GAUSSIANMEASURE_0, \GAUSSIANMEASURE_1, \MARGIN, \SIGNALSHORTDEPENDENCE)$. This gives
\begin{align*}
\sup_{t,v\in [0,1]^2} \Big( \nu_{t,v} \big( (R_{1,1} - \rho_{\SAMPLES, \DIMENSION})^2 \big)^{1/2} + \nu_{t,v} \big( (R_{1,2} - q_{\SAMPLES, \DIMENSION})^2 \big)^{1/2} \Big) \le K(D, \GAUSSIANMEASURE_0, \GAUSSIANMEASURE_1, \MARGIN, \SIGNALSHORTDEPENDENCE)\, .
\end{align*}
We thus can bound 
\begin{align*}
\Big| \frac{\rmd}{\rmd v} \varphi(v) \Big| &\le K(D, \SIGNALSHORTDEPENDENCE) \eps'' + K(D, \GAUSSIANMEASURE_0, \GAUSSIANMEASURE_1, \MARGIN, \SIGNALMEDIUMDEPENDENCE) \Big( \eps'' + \frac1{\sqrt{\DIMENSION}} + \frac1{\DIMENSION} \Big)\,,
\end{align*}
which as remarked earlier is sufficient to conclude Proposition \upref{prop:interpolationperconstraint}, since $d \ge \frac1{\eps''^2}$.
\end{proof}

\subsection{Simplifying the expression from cavity in $n$}\label{subsec:RSequationseliminatemollifier}
Again let $\alpha = n/d$, and let $\psi_{\alpha}, \bar{\psi}_{\alpha}$ be defined in terms of $\alpha$ as per \equpref{eq:psidef}, \equpref{eq:psibardef}.
Our goal here in Appendix \upref{subsec:RSequationseliminatemollifier} is to show that the expression from Proposition \upref{prop:interpolationperconstraint}, specifically
\begin{align}\label{eq:interpolate_oneattime_expression}
\frac1{\DIMENSION} \E \Big[ \log \E_{\ANNEALEDDISORDER} \exp u\big( x s_\SAMPLES + \sqrt{q_{\SAMPLES, \DIMENSION}} \QUENCHEDDISORDER + \sqrt{\rho_{\SAMPLES, \DIMENSION}-q_{\SAMPLES, \DIMENSION}}\ANNEALEDDISORDER \big) \Big]\,,
\end{align}
is very close to $\RSG_I$. The difference between \equpref{eq:interpolate_oneattime_expression} and $\RSG_I$ is that $\RSG_I$ instead depends on $\rho_0, q_0$, and additionally for the interpolators, the mollifier $u(\cdot)$ in \equpref{eq:interpolate_oneattime_expression} is replaced by the hard step function $\one\{x \ge \MARGIN\}$ in $\RSG_I$.

Here we will do the work in this argument that is common to the proof for the interpolators and posterior, specifically Propositions \upref{prop:approxrealRSequations}, \upref{prop:realRSeqsolsbounded}, \upref{prop:boundedapproxRSeq}.
The additional work for the interpolators that is needed to replace the mollifier $u(\cdot)$ in \equpref{eq:interpolate_oneattime_expression} by the hard step function $\one\{x \ge \MARGIN\}$ in $\RSG_I$ is done in Propositions
\upref{prop:diffvsuniquesol} and \upref{prop:mollifiererror}, presented later in Appendix \upref{subsec:finishinterpolation}.

Beginning with the argument, consider the following system of four variables and four equations:
\begin{align}
r = \psi_{\alpha}(q, \rho)\, ,\, \bar{r} = \bar{\psi}_{\alpha}(q, \rho)\, ,\, \rho = \frac1{2\GAUSSIANMEASURE+r-\bar{r}} + \frac{r}{(2\GAUSSIANMEASURE+r-\bar{r})^2}\, ,\, q = \frac{r}{(2\GAUSSIANMEASURE+r-\bar{r})^2}\, .\label{eq:realRSequations}
\end{align}
We will establish that $\rho_{\SAMPLES, \DIMENSION}, q_{\SAMPLES, \DIMENSION}$ approximately satisfy \equpref{eq:realRSequations}. This is where we use the careful choice of $r, \bar{r}$ in terms of $\psi_{\alpha}, \bar{\psi}_{\alpha}$ from \equpref{eq:interpolationrrbardef1}, \equpref{eq:interpolationrrbardef2}, which lets us complete the cavity in $\SAMPLES$ argument. In particular, in Appendix \upref{subsec:pf_derivativeintbounded}, we will complete the cavity in $\SAMPLES$ and also execute the cavity in $d$ to establish the following.
\begin{lemma}\label{lem:derivativeintbounded}
For all $t \in [0,1]$, we have
\begin{align}
\Big| \frac{\rmd}{\rmd t} \nu_t\big( \theta_{\DIMENSION}^1 \theta_{\DIMENSION}^2 \big) \Big|\,,\,\Big| \frac{\rmd}{\rmd t} \nu_t\big( (\theta_{\DIMENSION}^1)^2 \big) \Big| \le K\Big(\eps''+\frac1{\sqrt{\DIMENSION}} \Big)\,,\label{eq:approxrealRSequationsboundderiv}
\end{align}
where $K$ depends on $D,\GAUSSIANMEASURE_0,\GAUSSIANMEASURE_1, \SIGNALMEDIUMDEPENDENCE$.
\end{lemma}
With Lemma \upref{lem:derivativeintbounded} in hand, we now establish:
\begin{proposition}\label{prop:approxrealRSequations}
The following system of equations is satisfied by $\rho_{\SAMPLES, \DIMENSION}, q_{\SAMPLES, \DIMENSION}$:
\begin{align}
r = \psi_{\alpha} (q, \rho)\, , \,
\bar{r} = \bar{\psi}_{\alpha} (q, \rho)\, ,\,
\rho = \frac1{2\GAUSSIANMEASURE+r-\bar{r}} + \frac{r}{(2\GAUSSIANMEASURE+r-\bar{r})^2} + \delta_1 \, ,\,
q = \frac{r}{(2\GAUSSIANMEASURE+r-\bar{r})^2} + \delta_2\, ,\label{eq:approxRSequations}
\end{align}
where 
\begin{align}
\big| \delta_1 \big|\, ,\,  \big| \delta_2 \big| \le K_{\upref{prop:approxrealRSequations}}\Big( \eps'' + \frac1{\sqrt{\DIMENSION}} \Big)\,,\label{eq:approxRSequationserror}
\end{align}
for $K_{\upref{prop:approxrealRSequations}}$ depending on $D, \GAUSSIANMEASURE_0, \GAUSSIANMEASURE_1, \SIGNALMEDIUMDEPENDENCE$. Here the errors $|\delta_1|, |\delta_2|$ are upper bounded by a quantity independent of $\alpha = \frac{n}{d}$.\footnote{The upper bound only depends on $\DENSITYBOUND$ and is monotonically increasing, so we can apply the bound $\DENSITYBOUND \le 2$.} 
Moreover, as noted earlier, we have $r \ge \bar{r}$.
\end{proposition}

\begin{proof}[Proof of Proposition \upref{prop:approxrealRSequations}]
Consider $\nu_1(R_{1,1}), \nu_1(R_{1,2})$. By symmetry between coordinates, 
\[ \nu_1(R_{1,1}) = \nu_1\big( (\theta_{\DIMENSION}^1)^2 \big)\,,\,\nu_1(R_{1,1}) = \nu_1\big( \theta_{\DIMENSION}^1 \theta_{\DIMENSION}^2 \big)\,.\]
Next observe that for the measure $\nu_0$, the $\DIMENSION$-th coordinate completely decouples from the others. Moreover, the restriction $\one\{ \theta_1/\sqrt{\DIMENSION} \in I\}$ only concerns the first coordinate $\theta_1$. Thus, letting $Z_t$ denote the relevant normalizing constant, we have 
\begin{align*}
\nu_0\big( (\theta_{\DIMENSION}^1)^2 \big) &= \E \frac{\int \theta_{\DIMENSION}^2 \one\{ \theta_1/\sqrt{\DIMENSION} \in I\} \exp\Big( \sum_{\SAMPLEIDX \le \SAMPLES} u\big(S_{\SAMPLEIDX,0}(\btheta)\big) - \GAUSSIANMEASURE \| \btheta\|^2 + \theta_{\DIMENSION}\sqrt{r} Y - \frac{r-\bar{r}}2 \theta_{\DIMENSION}^2 \Big) \rmd \btheta }{Z_t} \\
&= \E \frac{\int \theta_{\DIMENSION}^2 \exp\Big( \theta_{\DIMENSION} \sqrt{r} Y - \frac{2\GAUSSIANMEASURE + r-\bar{r}}2 \theta_{\DIMENSION}^2 \Big) \rmd \theta_{\DIMENSION} }{ \int \exp\Big( \theta_{\DIMENSION} \sqrt{r} Y - \frac{2\GAUSSIANMEASURE + r-\bar{r}}2 \theta_{\DIMENSION}^2 \Big) \rmd \theta_{\DIMENSION} }\,.
\end{align*}
Similarly, we have 
\begin{align*}
\nu_0\big( \theta_{\DIMENSION}^1 \theta_{\DIMENSION}^2 \big) &= \E \frac{\int \theta^1_{\DIMENSION} \theta^2_{\DIMENSION} \exp\Big( (\theta^1_{\DIMENSION} + \theta^2_{\DIMENSION}) \sqrt{r} Y - \frac{2\GAUSSIANMEASURE + r-\bar{r}}2 \Big( \big(\theta^1_{\DIMENSION}\big)^2+\big(\theta^2_{\DIMENSION}\big)^2 \Big) \Big) \rmd \theta^1_{\DIMENSION} \theta^2_{\DIMENSION} }{ \int \exp\Big( (\theta^1_{\DIMENSION} + \theta^2_{\DIMENSION}) \sqrt{r} Y - \frac{2\GAUSSIANMEASURE + r-\bar{r}}2 \Big( \big(\theta^1_{\DIMENSION}\big)^2+\big(\theta^2_{\DIMENSION}\big)^2 \Big) \Big) \rmd \theta^1_{\DIMENSION} \theta^2_{\DIMENSION} }  \\
&= \E  \Bigg( \frac{\int \theta_{\DIMENSION} \exp\Big( \theta_{\DIMENSION} \sqrt{r} Y - \frac{2\GAUSSIANMEASURE + r-\bar{r}}2 \theta_{\DIMENSION}^2 \Big) \rmd \theta_{\DIMENSION} }{ \int \exp\Big( \theta_{\DIMENSION} \sqrt{r} Y - \frac{2\GAUSSIANMEASURE + r-\bar{r}}2 \theta_{\DIMENSION}^2 \Big) \rmd \theta_{\DIMENSION} } \Bigg)^2 \,.
\end{align*}
The above are univariate Gaussian integrals w.r.t. $\theta_{\DIMENSION}$. As detailed on p. 219 of \cite{talagrand2010mean}, letting $Z$ be a centered univariate Gaussian and $s \in \R$, we obtain from Gaussian Integration by Parts that
\begin{align*}
\E\big[ Z e^{sZ} ] = s \E\big[ Z^2 \big] \E\big[ e^{sZ} \big]\,,\,\E\big[ Z^2 e^{sZ} ] = \E\big[Z^2 \big] \Big( \E\big[ e^{sZ} \big] + s^2 \E\big[ Z^2 \big] \E\big[ e^{sZ} \big] \Big)\,,
\end{align*}
where expectation is w.r.t. the law of $Z$. We obtain
\begin{align*}
\Bigg( \frac{\E\big[ Z e^{sZ} ]}{\E\big[ e^{sZ} \big]} \Bigg)^2 = s^2 \E\big[ Z^2 \big]^2\,,\, \frac{\E\big[ Z^2 e^{sZ} ]}{\E\big[ e^{sZ} \big]} = \E\big[Z^2] + s^2 \E\big[ Z^2 \big]^2\,. 
\end{align*}
We apply the above display with $s = \sqrt{r}Y$, $Z \sim N\big(0, \frac1{2\GAUSSIANMEASURE+r-\bar{r}}\big)$ to compute the quantity inside the expectation in the integrals equaling  $\nu_0\big( \theta_{\DIMENSION}^1 \theta_{\DIMENSION}^2 \big)$, $\nu_0\big( (\theta_{\DIMENSION}^1)^2 \big)$. This yields
\begin{align*}
\nu_0\big( \theta_{\DIMENSION}^1 \theta_{\DIMENSION}^2 \big) &= \E  \Bigg( \frac{\int \theta_{\DIMENSION} \exp\Big( \theta_{\DIMENSION} \sqrt{r} Y - \frac{2\GAUSSIANMEASURE + r-\bar{r}}2 \theta_{\DIMENSION}^2 \Big) \rmd \theta_{\DIMENSION} }{ \int \exp\Big( \theta_{\DIMENSION} \sqrt{r} Y - \frac{2\GAUSSIANMEASURE + r-\bar{r}}2 \theta_{\DIMENSION}^2 \Big) \rmd \theta_{\DIMENSION} } \Bigg)^2 = \frac{r Y^2}{(2\GAUSSIANMEASURE+r-\bar{r})^2}\,, \\
\nu_0\big( (\theta_{\DIMENSION}^1)^2 \big) &= \E \frac{\int \theta_{\DIMENSION}^2 \exp\Big( \theta_{\DIMENSION} \sqrt{r} Y - \frac{2\GAUSSIANMEASURE + r-\bar{r}}2 \theta_{\DIMENSION}^2 \Big) \rmd \theta_{\DIMENSION} }{ \int \exp\Big( \theta_{\DIMENSION} \sqrt{r} Y - \frac{2\GAUSSIANMEASURE + r-\bar{r}}2 \theta_{\DIMENSION}^2 \Big) \rmd \theta_{\DIMENSION} } = \frac{1}{2\GAUSSIANMEASURE+r-\bar{r}} + \frac{r Y^2}{(2\GAUSSIANMEASURE+r-\bar{r})^2}\,.
\end{align*}
Finally taking expectation over $Y \sim N(0, 1)$, we establish that 
\begin{align*}
\rho &= \nu_1\big( (\theta_{\DIMENSION}^1)^2 \big) = \frac1{2\GAUSSIANMEASURE+r-\bar{r}} + \frac{r}{(2\GAUSSIANMEASURE+r-\bar{r})^2} + \Big( \nu_1\big( (\theta_{\DIMENSION}^1)^2 \big) - \nu_0\big( (\theta_{\DIMENSION}^1)^2 \big)\Big) \, , \notag\\
q &=  \nu_1\big( \theta_{\DIMENSION}^1 \theta_{\DIMENSION}^2 \big) = \frac{r}{(2\GAUSSIANMEASURE+r-\bar{r})^2} + \Big( \nu_1\big( \theta_{\DIMENSION}^1 \theta_{\DIMENSION}^2 \big) - \nu_0\big( \theta_{\DIMENSION}^1 \theta_{\DIMENSION}^2 \big)\Big)\,.
\end{align*}
The result follows by applying Lemma \upref{lem:derivativeintbounded} to control $\nu_1\big( (\theta_{\DIMENSION}^1)^2 \big) - \nu_0\big( (\theta_{\DIMENSION}^1)^2 \big)$, $\nu_1\big( \theta_{\DIMENSION}^1 \theta_{\DIMENSION}^2 \big) - \nu_0\big( \theta_{\DIMENSION}^1 \theta_{\DIMENSION}^2 \big)$.
\end{proof}
In Appendix \upref{subsec:real_rs_bounded}, we will establish:
\begin{proposition}\label{prop:realRSeqsolsbounded}
For all $(\alpha, \GAUSSIANMEASURE, x) \in [0, \DENSITYBOUND] \times [\GAUSSIANMEASURE_0, \GAUSSIANMEASURE_1] \times [-1, 1]$, any solution $(q, \rho, r, \bar{r})$ of \equpref{eq:realRSequations} must satisfy $(q, \rho) \in (0, \Csol) \times \big(\frac1{\Csol}, \Csol \big)$ and $\frac{1}{\rho-q} < \Csol$. 
\end{proposition}

We will use Proposition \upref{prop:approxrealRSequations} and \upref{prop:realRSeqsolsbounded} to establish the following.
\begin{proposition}\label{prop:boundedapproxRSeq}
Consider any $(\alpha, \GAUSSIANMEASURE, x) \in [0, \DENSITYBOUND] \times [\GAUSSIANMEASURE_0, \GAUSSIANMEASURE_1] \times [-1, 1]$. Then there is a $ 0 < \eps''_{\upref{prop:boundedapproxRSeq}} \le \eps'$ depending on $u, \beta_0, \beta_1, \SIGNALMEDIUMDEPENDENCE$ satisfying the following. For all $\eps'' \le \eps''_{\upref{prop:boundedapproxRSeq}}$ and all $d \ge d_{\upref{prop:boundedapproxRSeq}} \ge \frac1{\eps''^2}$, for any $n \ge 1$ such that $\frac{n}{d} = \alpha \in [0, \DENSITYBOUND]$, we have 
\begin{align*}
(q_{n, d}, \rho_{n, d}) \in \big[0, \Csol\big] \times \Big[\frac1{\Csol}, \Csol\Big]\quad\text{ and }\quad \frac1{\rho_{n,d} - q_{n, d}} \le \Csol\,,
\end{align*}
Furthermore, in the case of the interpolators, supposing the mollifier $u$ satisfies the conditions of Theorem \upref{thm:mollifiedhamiltonianfreeenergy} in terms of some parameter $\eps'$, we have $|\delta_1|, |\delta_2| \le \eps'$ where $\delta_1, \delta_2$ are the errors from Proposition \upref{prop:approxrealRSequations} for $\rho_{n,d}, q_{n,d}$.
\end{proposition}
\begin{proof}
By definition of \equpref{eq:approxRSequations}, Proposition \upref{prop:approxrealRSequations}, and as $\alpha \le 2$, for suitable $\eps''_{\upref{prop:boundedapproxRSeq}}$ depending on $u, \beta_0, \beta_1, \SIGNALMEDIUMDEPENDENCE$ and $\eps'' \le \eps''_{\upref{prop:boundedapproxRSeq}}$ and $d \ge d_{\upref{prop:boundedapproxRSeq}} \ge \frac1{\eps''^2}$, we have the bounds $\big| \rho_{n, d} \big|, \big| q_{n, d} \big| \le K_{\textsc{bound}}(D, \beta_0, \beta_1, \SIGNALSHORTDEPENDENCE)$ for some $K_{\textsc{bound}}(D, \beta_0, \beta_1, \SIGNALSHORTDEPENDENCE)>0$ (dependence on $D$ comes through $u$). Without loss of generality, suppose $K_{\textsc{bound}}(D, \beta_0, \beta_1, \SIGNALSHORTDEPENDENCE) > \Csol$. Consider 
\begin{align*}
\cD_u &:= \Big[-K_{\textsc{bound}}(D, \beta_0, \beta_1, \SIGNALSHORTDEPENDENCE), K_{\textsc{bound}}(D, \beta_0, \beta_1, \SIGNALSHORTDEPENDENCE) \Big] \\
&\qquad \qquad \qquad \qquad \times \Big[-K_{\textsc{bound}}(D, \beta_0, \beta_1, \SIGNALSHORTDEPENDENCE), K_{\textsc{bound}}(D, \beta_0, \beta_1, \SIGNALSHORTDEPENDENCE) \Big]\,,\\
\cD &:= \Big\{ (q, \rho) \in \big(0, \Csol\big) \times \Big(\frac1{\Csol}, \Csol\Big), \frac1{\rho-q} < \Csol\Big\} \,.
\end{align*}
Let 
\begin{align*}
R(q, \rho, \alpha, \beta) &:= \Bigg( \rho - \frac1{2\beta + \psi_{\alpha}(q, \rho) - \bar{\psi}_{\alpha}(q, \rho)} - \frac{\psi_{\alpha}(q, \rho)}{\big( 2\beta + \psi_{\alpha}(q, \rho) - \bar{\psi}_{\alpha}(q, \rho) \big)^2}, \\
&\qquad \qquad \qquad \qquad q -  \frac{\psi_{\alpha}(q, \rho)}{\big( 2\beta + \psi_{\alpha}(q, \rho) - \bar{\psi}_{\alpha}(q, \rho) \big)^2}\Bigg)\,,\\
\underline{c} &:= \inf_{(q, \rho) \in \cD_u - \cD, \alpha \in [0, 2], \beta \in [\beta_0, \beta_1]} \big\| R(q, \rho, \alpha, \beta) \big\| \le \inf_{(q, \rho) \in \cD_u - \cD, \alpha \in [0, \DENSITYBOUND], \beta \in [\beta_0, \beta_1]} \big\| R(q, \rho, \alpha, \beta) \big\|\,. 
\end{align*}
Note $R(q, \rho, \alpha, \beta)$ is continuous in $q, \rho, \alpha, \beta$. Moreover in $\cD_u - \cD$, by Proposition \upref{prop:realRSeqsolsbounded}, we have that $R(q, \rho, \alpha, \beta)>0$ pointwise; else we obtain a solution to \equpref{eq:realRSequations} outside $\cD$, contradicting Proposition \upref{prop:realRSeqsolsbounded}. Thus as $\big\{ (q, \rho, \alpha, \beta):(q, \rho) \in \cD_u - \cD, \alpha \in [0, 2], \beta \in [\beta_0, \beta_1]\big\}$ is compact, $\underline{c}>0$. 

Hence in the case of interpolators, we may take $\eps''_{\upref{prop:boundedapproxRSeq}}=\eps''_{\upref{prop:boundedapproxRSeq}}(u, \beta_0, \beta_1, \SIGNALMEDIUMDEPENDENCE)$ small enough, and $d_{\upref{prop:boundedapproxRSeq}}(\eps'') \ge \frac1{\eps''^2}$ large enough in terms of $K_{\upref{prop:approxrealRSequations}}(D, \beta_0, \beta_1, \SIGNALMEDIUMDEPENDENCE)$ and $\eps'$ (recall that $D$ depends solely on $u$) so that $|\delta_1|, |\delta_2|$ from \equpref{eq:approxRSequationserror} are each strictly less than $\min\big\{ \frac{\underline{c}}2, \frac{\eps'}2 \big\}$. 
In the case of the posterior, we ignore dependence in $\eps'$ in the above definitions of $\eps''_{\upref{prop:boundedapproxRSeq}}=\eps''_{\upref{prop:boundedapproxRSeq}}(u, \beta_0, \beta_1, \SIGNALMEDIUMDEPENDENCE)$ and $d_{\upref{prop:boundedapproxRSeq}}(\eps'')$, now guaranteeing $|\delta_1|, |\delta_2| < \underline{c}/2$.

By Proposition \upref{prop:approxrealRSequations}, it follows that $\| R(q_{n,d}, \rho_{n,d}, \alpha, \beta) \| < \underline{c}$. 
Since $(q_{n,d}, \rho_{n,d}) \in \cD_u$, by definition of $\underline{c}$, it follows that $(q_{n,d}, \rho_{n,d}) \in \cD$. The conclusion of this Proposition follows.
\end{proof}
Note this argument is where we use that $|\delta_1|, |\delta_2|$ from \equpref{eq:approxRSequationserror} are upper bounded independently of $\alpha$. 
In particular this lets us take $\eps''_{\upref{prop:boundedapproxRSeq}}=\eps''_{\upref{prop:boundedapproxRSeq}}(u, \beta_0, \beta_1, \SIGNALMEDIUMDEPENDENCE)$ independently of $\alpha$, which is crucial as we will use the same $\eps''$ across $n=1, 2, \ldots, d$ when computing the expected free energy. 
The compactness supplied from the above Proposition is crucial in the arguments to follow in Propositions
\upref{prop:diffvsuniquesol} and \upref{prop:mollifiererror}.

The above work lets us interpolate one datapoint at a time and compute the difference in free energy, as we do in Appendix \upref{subsec:finishinterpolation}. 
The last step we need is to compute $\phi_{\beta, 0, d}(I)$, which is the base case arising when all $n$ datapoints have been accounted for.
\begin{lemma}\label{lem:interpolationlastcalc}
For any $\eps>0$, for some $\eps''_{\upref{lem:interpolationlastcalc}}, d_{\upref{lem:interpolationlastcalc}} > 0$ depending on $\eps$, we have for all $\eps'' \le \eps''_{\upref{lem:interpolationlastcalc}}$ and $d \ge d_{\upref{lem:interpolationlastcalc}}$ that 
\begin{align*}
\Big| \phi_{\beta, 0, d}(I) - \Big( \RS_{0,I}(0) + \frac12 \log (2\pi e) \Big) \Big| \le \eps\,.
\end{align*}
\end{lemma}
\begin{proof}
A direct computation gives that
\begin{align*}
\phi_{\beta, 0, d}(I) = \frac12 \log \frac{\pi}{\beta} - \beta r(x) + O\Big( \frac{\log d}d + \eps'' \Big)\,,
\end{align*}
where we recall the definition of $r(x)$ from \equpref{eq:intervalrxdef}. Here $O(\cdot)$ only hides universal constants. 
Now note when $\alpha=0$, $\frac{\partial F}{\partial q} = 0$ implies $q=0$, hence $q_0(0)=0$. Also, considering $\frac{\partial F}{\partial \rho} = 0$ when $\alpha=0$ implies that $\rho_0(0) = \frac1{2\beta}$.
Therefore 
\begin{align*}
\RS_{0,I}(0) + \frac12 \log (2\pi e) = \frac12 \log\big(2\pi e \rho_0(0)\big) - \beta \big( \rho_0(0) + r(x) \big) = \frac12 \log \frac{\pi}{\beta} - \beta r(x)\,.
\end{align*}
For suitable $\eps'' \le \eps''_{\upref{lem:interpolationlastcalc}}(\eps)$ and $d \ge d_{\upref{lem:interpolationlastcalc}}(\eps)$, this implies the Lemma.
\end{proof}

\subsection{Proof of Theorem \upref{thm:mollifiedhamiltonianfreeenergy} for interpolators}\label{subsec:finishinterpolation}
We now complete the proof of Theorem \upref{thm:mollifiedhamiltonianfreeenergy} by summing Proposition \upref{prop:interpolationperconstraint} and the results of Appendix \upref{subsec:RSequationseliminatemollifier} across all datapoints $i$, $1 \le i \le n$.
We first handle the necessary steps to replace the mollifier $u(\cdot)$ in \equpref{eq:interpolate_oneattime_expression} by the hard step function $\one\{x \ge \MARGIN\}$ in $\RSG_I$.
Specifically, this is done in the following Propositions \upref{prop:mollifiererror} and \upref{prop:diffvsuniquesol}, which establish that the impact of the mollifier $u(\cdot)$ is insignificant. At a high level, this is because the conditions of Theorem \upref{thm:mollifiedhamiltonianfreeenergy} imply $u(\cdot)$ is an accurate approximation to $\one\{x \ge \MARGIN\}$. 
\begin{proposition}\label{prop:mollifiererror}
Consider any $\delta>0$ and let $\MARGIN^+ = \max\{\MARGIN, 0\}$. Then there exists $\eps'_{\upref{prop:mollifiererror}}>0$ depending on $\delta, \Csol, \MARGIN^+, \lambda$ such that the following holds. 
Suppose $u \le 0$ satisfies for some $\eps' \le \eps'_{\upref{prop:mollifiererror}}$ that:
\begin{enumerate}
    \item $u(x) = 0$ for $x \ge \MARGIN$,
    \item $\exp u(x) \le \eps'$ for all $x \le \MARGIN - \eps'$.
\end{enumerate}
Then for any $(q, \rho)$ such that $(q, \rho) \in [0, \Csol] \times [\frac1{\Csol}, \Csol]$ and $\frac{1}{\rho-q} \le \Csol$, and any $\SIGNALDISORDER \sim \cD$, we have 
\begin{align*}
&\Big| \E \Big[ \log \E_{\ANNEALEDDISORDER} \exp u\big( x\SIGNALDISORDER + \sqrt{q}\QUENCHEDDISORDER + \sqrt{\rho-q}\ANNEALEDDISORDER \big) \Big]  - \E \Big[ \log \P_{\ANNEALEDDISORDER} \big( x\SIGNALDISORDER + \sqrt{q}\QUENCHEDDISORDER + \sqrt{\rho-q}\ANNEALEDDISORDER \ge \MARGIN \big) \Big] \Big| \le \delta\, .
\end{align*}
\end{proposition}
\begin{proof}
Note $\E_{\ANNEALEDDISORDER} \exp u\big( x\SIGNALDISORDER + \sqrt{q}\QUENCHEDDISORDER + \sqrt{\rho-q}\ANNEALEDDISORDER \big) \ge \P_{\ANNEALEDDISORDER} \big( x\SIGNALDISORDER + \sqrt{q}\QUENCHEDDISORDER + \sqrt{\rho-q}\ANNEALEDDISORDER \ge \MARGIN \big)$.
Breaking into cases on the value of $x\SIGNALDISORDER + \sqrt{q}\QUENCHEDDISORDER + \sqrt{\rho-q}\ANNEALEDDISORDER$ and using that $\exp u(x) \le \eps'$ for all $x \le \MARGIN-\eps'$ and $u \le 0$, we thus may rewrite
\begin{align*}
0 &\le \E \Big[ \log \E_{\ANNEALEDDISORDER} \exp u\big( x\SIGNALDISORDER + \sqrt{q}\QUENCHEDDISORDER + \sqrt{\rho-q}\ANNEALEDDISORDER \big) \Big]  - \E \Big[ \log \P_{\ANNEALEDDISORDER} \big( x\SIGNALDISORDER + \sqrt{q}\QUENCHEDDISORDER + \sqrt{\rho-q}\ANNEALEDDISORDER \ge \MARGIN \big)  \Big] \\
&= \E \Bigg[ \log \frac{ \P_{\ANNEALEDDISORDER} \big( x\SIGNALDISORDER + \sqrt{q}\QUENCHEDDISORDER + \sqrt{\rho-q}\ANNEALEDDISORDER \ge \MARGIN \big) + \eps' + \P_{\ANNEALEDDISORDER}\big(x\SIGNALDISORDER + \sqrt{q}\QUENCHEDDISORDER + \sqrt{\rho-q}\ANNEALEDDISORDER \in (\MARGIN - \eps' , \MARGIN) \big) }{ \P_{\ANNEALEDDISORDER} \big( x\SIGNALDISORDER + \sqrt{q}\QUENCHEDDISORDER + \sqrt{\rho-q}\ANNEALEDDISORDER \ge \MARGIN \big) }\Bigg] \\
&= \E \Bigg[ \log \Bigg(1 + \frac{ \eps' + \P_{\ANNEALEDDISORDER}\big(x\SIGNALDISORDER + \sqrt{q}\QUENCHEDDISORDER + \sqrt{\rho-q}\ANNEALEDDISORDER \in (\MARGIN - \eps' , \MARGIN) \big) }{ \P_{\ANNEALEDDISORDER} \big( x\SIGNALDISORDER + \sqrt{q}\QUENCHEDDISORDER + \sqrt{\rho-q}\ANNEALEDDISORDER \ge \MARGIN \big) } \Bigg)\Bigg] \\
&\le \E \Bigg[ \frac{\eps' + \P_{\ANNEALEDDISORDER}\big(x\SIGNALDISORDER + \sqrt{q}\QUENCHEDDISORDER + \sqrt{\rho-q}\ANNEALEDDISORDER \in (\MARGIN - \eps' , \MARGIN) \big)}{\P_{\ANNEALEDDISORDER} \big( x\SIGNALDISORDER + \sqrt{q}\QUENCHEDDISORDER + \sqrt{\rho-q}\ANNEALEDDISORDER \ge \MARGIN^+ \big)} \Bigg]\,,
\end{align*}
where we used the inequality $\log(1+x) \le x$ and replaced $\MARGIN$ by $\MARGIN^+$ in the denominator by monotonicity. We next bound
\begin{align*}
\P_{\ANNEALEDDISORDER}\Big(x\SIGNALDISORDER + \sqrt{q}\QUENCHEDDISORDER + \sqrt{\rho-q}\ANNEALEDDISORDER \in (\MARGIN - \eps' , \MARGIN) \Big) &= \P_{\ANNEALEDDISORDER} \Big( \ANNEALEDDISORDER \in \Big[\frac{\MARGIN-\eps'-sx-z\sqrt{q}}{\sqrt{\rho-q}}, \frac{\MARGIN-sx-z\sqrt{q}}{\sqrt{\rho-q}} \Big] \Big) \\
&\le \frac{\eps'}{\sqrt{\rho-q}} \cdot \sup_{w \in \R} \frac1{\sqrt{2\pi}} \exp(-w^2/2) \\
&\le \Csol^{1/2} \eps'\,.
\end{align*}
Also note that as $\rho-q \ge \frac1{\Csol}$, $\E\Big[ \frac1{\P_{\ANNEALEDDISORDER} \big( x\SIGNALDISORDER + \sqrt{q}\QUENCHEDDISORDER + \sqrt{\rho-q}\ANNEALEDDISORDER \ge \MARGIN^+ \big)} \Big]$ is a continuous function of $q, \rho$ and only depends on $q, \rho, \MARGIN^+, \lambda$. Hence we may upper bound its value in the compact region $(q, \rho) \in [0, \Csol] \times [\frac1{\Csol}, \Csol], \frac1{\rho-q} \le \Csol$ as follows:
\begin{align*}
\sup_{(q, \rho) \in [0, \Csol] \times [\frac1{\Csol}, \Csol], \frac1{\rho-q} \le \Csol} \E\Bigg[ \frac1{\P_{\ANNEALEDDISORDER} \big( x\SIGNALDISORDER + \sqrt{q}\QUENCHEDDISORDER + \sqrt{\rho-q}\ANNEALEDDISORDER \ge \MARGIN^+ \big)} \Bigg] \le K(\Csol, \MARGIN^+, \lambda)\,.
\end{align*}
It follows that 
\begin{align*}
0 &\le \E \Big[ \log \E_{\ANNEALEDDISORDER} \exp u\big( x\SIGNALDISORDER + \sqrt{q}\QUENCHEDDISORDER + \sqrt{\rho-q}\ANNEALEDDISORDER \big) \Big]  - \E \Big[ \log \P_{\ANNEALEDDISORDER} \big( x\SIGNALDISORDER + \sqrt{q}\QUENCHEDDISORDER + \sqrt{\rho-q}\ANNEALEDDISORDER \ge \MARGIN \big)  \Big] \\
&\le \eps' (\Csol^{1/2}+1) \E\Bigg[ \frac1{\P_{\ANNEALEDDISORDER} \big( x\SIGNALDISORDER + \sqrt{q}\QUENCHEDDISORDER + \sqrt{\rho-q}\ANNEALEDDISORDER \ge \MARGIN^+ \big)} \Bigg] \\
&\le K(\Csol, \MARGIN^+, \lambda) \big(\Csol^{1/2}+1 \big) \eps' \,.
\end{align*}
Thus we may take $\eps'_{\upref{prop:mollifiererror}} = \eps'_{\upref{prop:mollifiererror}}(\delta, \Csol, \MARGIN^+, \lambda) \le \frac{\delta}{K(\Csol, \MARGIN^+, \lambda) (\Csol^{1/2}+1)}$.
\end{proof}

We next use the uniqueness of the solution $(q_0, \rho_0)$ to \equpref{eq:niceRSeqs}, given by Proposition \upref{prop:nicersuniquesol}, to control the difference of $(q_{n,d}, \rho_{n,d})$ vs. $(q_0, \rho_0)$.
\begin{proposition}\label{prop:diffvsuniquesol}
Consider any $(\beta, x) \in [\beta_0, \beta_1] \times [-1, 1]$. Given any $\delta>0$, there exists $\eps'_{\upref{prop:diffvsuniquesol}}>0$ depending on $\delta$ such that the following holds. 
Consider any $u \le 0$ satisfies the conditions of Theorem \upref{thm:mollifiedhamiltonianfreeenergy} with $\eps' = \eps'_{\upref{prop:diffvsuniquesol}}$. 
Then for all $\eps'' \le \eps''_{\upref{prop:boundedapproxRSeq}}$ and all $d \ge d_{\upref{prop:boundedapproxRSeq}}$, and for any $n \ge 1$ such that $\frac{n}{d} = \alpha \in [0, \DENSITYBOUND]$, the following holds.
Defining $\big( q_0(\alpha, \beta, x), \rho_0(\alpha, \beta, x) \big)$ as per Proposition \upref{prop:nicersuniquesol},
\begin{align*}
\big| q_{n,d} - q_0(\alpha, \beta, x) \big|\quad \text{and}\quad \big|\rho_{n,d} - \rho_0(\alpha, \beta, x) \big| \le \delta\,.
\end{align*}
\end{proposition}
\begin{proof}
Suppose for the sake of contradiction that there exists a $\delta>0$ and a sequence $\eps'_k \rightarrow 0$, a sequence of functions $u_k$ such that $u_k(x)=0$ for $x \ge \MARGIN$ and $\exp u_k(x) \le \eps'_k$ for all $x \le \MARGIN-\eps'_k$, $\eps''_k \le \eps''_{\upref{prop:boundedapproxRSeq}}(u_k, \beta_0, \beta_1, \SIGNALMEDIUMDEPENDENCE)$, $d_k \ge d_{\upref{prop:boundedapproxRSeq}}(\eps''_k)$, $n_k$ with $\frac{n_k}{d_k} = \alpha_k \in [0, \DENSITYBOUND]$, such that
\begin{align}\label{eq:mollifiererror_contradiction}
\max\Big\{ \big|q_{n_k, d_k} -q_0(\alpha_k, \beta, x)\big|\, ,\, \big|\rho_{n_k, d_k} -\rho_0(\alpha_k, \beta, x)\big| \Big\} \ge \delta\, .    
\end{align}
By Proposition \upref{prop:boundedapproxRSeq}, since $\eps''_k \le \eps''_{\upref{prop:boundedapproxRSeq}}(u_k, \beta_0, \beta_1, \SIGNALMEDIUMDEPENDENCE)$, $d_k \ge d_{\upref{prop:boundedapproxRSeq}}(\eps''_k)$, we have 
\begin{align}
(q_{n_k, d_k}, \rho_{n_k, d_k}) \in \big[0, \Csol\big] \times \Big[\frac1{\Csol}, \Csol\Big]\, ,\, \frac1{\rho_{n_k, d_k} - q_{n_k, d_k}} \le \Csol\,.\label{eq:mollifiererroruniquenesssolbounded}
\end{align}
Also, letting $\delta_{1,k}, \delta_{2,k}$ be the corresponding $\delta_1, \delta_2$ from the system \equpref{eq:approxRSequations} for $u_k, \alpha_k, q_{n_k, d_k}, \rho_{n_k, d_k}$, we have by the second part of Proposition \upref{prop:boundedapproxRSeq} that $|\delta_{1,k}|, |\delta_{2,k}| \le \eps'_k$.

Thus $(q_{n_k, d_k}, \rho_{n_k, d_k})$ lie in a set bounded independently of $k$, so we may extract a bounded subsequence $(\alpha_k, q_{n_k, d_k}, \rho_{n_k, d_k})$ converging to some limit $(\alpha^{\star}, q^{\star}, \rho^{\star})$. 
Upon reindexing, as $\eps'_k \rightarrow 0$, we may suppose that
\begin{align}\label{eq:mollifiererror_limit_subseq}
(\eps'_k, \alpha_k, q_{n_k, d_k}, \rho_{n_k, d_k}, \delta_{1,k}, \delta_{2,k}) \rightarrow (0, \alpha^{\star}, q^{\star}, \rho^{\star}, 0, 0)\,.
\end{align}
Thus, we have 
\begin{align}
(q^{\star}, \rho^{\star}) \in \big[0, \Csol\big] \times \Big[\frac1{\Csol}, \Csol\Big]\, ,\, \frac1{\rho^{\star} - q^{\star}} \le \Csol\,.\label{eq:mollifiererror_limit_uniquenesssolbounded}
\end{align}
By \equpref{eq:mollifiererror_contradiction}, and as $\lim_{k\rightarrow\infty} q_0(\alpha_k, \beta, x) = q_0(\alpha^{\star}, \beta, x)$, $\lim_{k\rightarrow\infty} \rho_0(\alpha_k, \beta, x) = \rho_0(\alpha^{\star}, \beta, x)$ by Proposition \upref{prop:nicersuniquesol}, it follows that 
\begin{align}
\max\Big\{ \big| q^{\star} - q_0(\alpha^{\star}, \beta, x) \big|, \big| \rho^{\star} - \rho_0(\alpha^{\star}, \beta, x) \big| \Big\} \ge \delta\,.\label{eq:mollifiererroruniquenesssolcontradict}
\end{align}
By \equpref{eq:mollifiererroruniquenesssolbounded}, \equpref{eq:mollifiererror_limit_subseq}, \equpref{eq:mollifiererror_limit_uniquenesssolbounded}, in particular as $\rho_{n_k, d_k} - q_{n_k, d_k}, \rho^\star - q^\star \ge \frac1{\Csol}$, we also have the existence of the following limit and the following equality:
\begin{equation}\label{eq:mollifier_error_limit1}
\begin{aligned}
&\lim_{k \rightarrow \infty} \E_{\ANNEALEDDISORDER}\Big[ f(\ANNEALEDDISORDER) \one\big\{ x\SIGNALDISORDER + \sqrt{q_{n_k, d_k}}\QUENCHEDDISORDER + \sqrt{\rho_{n_k, d_k} - q_{n_k, d_k}}\ANNEALEDDISORDER \ge \MARGIN \big\} \Big] \\
&\quad = \E_{\ANNEALEDDISORDER}\Big[ f(\ANNEALEDDISORDER) \one\big\{ x\SIGNALDISORDER + \sqrt{q^\star}\QUENCHEDDISORDER + \sqrt{\rho^\star - q^\star}\ANNEALEDDISORDER \ge \MARGIN \big\} \Big] \,.
\end{aligned}
\end{equation}
Next, letting $f(\ANNEALEDDISORDER) \in \big\{1, \ANNEALEDDISORDER, \ANNEALEDDISORDER^2-1\big\}$ and using our conditions on $u_k$,
\begin{align}
&\Big| \E_{\ANNEALEDDISORDER}\Big[ f(\ANNEALEDDISORDER) \exp u_k\big( x\SIGNALDISORDER + \sqrt{q_{n_k, d_k}}\QUENCHEDDISORDER + \sqrt{\rho_{n_k, d_k} - q_{n_k, d_k}}\ANNEALEDDISORDER \big) \Big] \notag \\
&\qquad \qquad \qquad - \E_{\ANNEALEDDISORDER}\Big[ f(\ANNEALEDDISORDER) \one\big\{ x\SIGNALDISORDER + \sqrt{q_{n_k, d_k}}\QUENCHEDDISORDER + \sqrt{\rho_{n_k, d_k} - q_{n_k, d_k}}\ANNEALEDDISORDER \ge \MARGIN \big\} \Big] \Big| \notag \\
&\quad \le \eps'_k \Big| \E_{\ANNEALEDDISORDER}\big[ f(\ANNEALEDDISORDER) \big] \Big| + \Big| \E_{\ANNEALEDDISORDER}\Big[ f(\ANNEALEDDISORDER) \one\big\{ x\SIGNALDISORDER + \sqrt{q_{n_k, d_k}}\QUENCHEDDISORDER + \sqrt{\rho_{n_k, d_k} - q_{n_k, d_k}}\ANNEALEDDISORDER \in \big(\MARGIN-\eps'_k, \MARGIN\big) \big\} \Big] \Big| \notag \\
&\quad = \eps'_k \Big| \E_{\ANNEALEDDISORDER}\big[ f(\ANNEALEDDISORDER) \big] \Big| + \int_{\MARGIN_1}^{\MARGIN_2} |f(\ANNEALEDDISORDER)| \pdfnormal(\ANNEALEDDISORDER)\, \rmd \ANNEALEDDISORDER \notag \\
&\quad \le \eps'_k \Big| \E_{\ANNEALEDDISORDER}\big[ f(\ANNEALEDDISORDER) \big] \Big| + \frac{K \eps'_k}{\sqrt{\rho_{n_k, d_k} - q_{n_k, d_k}}} \notag \\
&\quad \le K \Csol^{1/2} \eps'_k\,,\label{eq:mollifiererroruniquenesssolineq}
\end{align}
where we let $\MARGIN_1 = \frac{\MARGIN - \eps'_k - x\SIGNALDISORDER - \sqrt{q_{n_k, d_k}}\QUENCHEDDISORDER}{\sqrt{\rho_{n_k, d_k}-q_{n_k, d_k}}}$, $\MARGIN_2 = \frac{\MARGIN - sX - \sqrt{q_{n_k, d_k}}\QUENCHEDDISORDER}{\sqrt{\rho_{n_k, d_k}-q_{n_k, d_k}}} = \MARGIN_1 + \frac{\eps'_k}{\sqrt{\rho_{n_k, d_k} - q_{n_k, d_k}}}$.
Here in the last step we used that $|f(y)| \pdfnormal(y) \le K$ uniformly for all $y \in \R$.

In light of \equpref{eq:mollifier_error_limit1}, \equpref{eq:mollifiererroruniquenesssolineq}, since $\Csol$ is independent of the mollifier and as $\eps'_k \rightarrow 0$, it follows that the following limit exists and can be computed as:
\begin{equation}\label{eq:mollifier_error_limit}
\begin{aligned}  
&\lim_{k \rightarrow \infty} \E_{\ANNEALEDDISORDER}\Big[ f(\ANNEALEDDISORDER) \exp u_k\big( x\SIGNALDISORDER + \sqrt{q_{n_k, d_k}}\QUENCHEDDISORDER + \sqrt{\rho_{n_k, d_k} - q_{n_k, d_k}}\ANNEALEDDISORDER \big) \Big] \\
&\quad = \E_{\ANNEALEDDISORDER}\Big[ f(\ANNEALEDDISORDER) \one\big\{ x\SIGNALDISORDER + \sqrt{q^{\star}}\QUENCHEDDISORDER + \sqrt{\rho^{\star} - q^{\star}}\ANNEALEDDISORDER  \ge \MARGIN \big\} \Big]\,.
\end{aligned}
\end{equation}
(Note this argument does not assume anything a-priori about the limit of the $u_k$'s; the existence of $q^{\star}, \rho^{\star}$ and \equpref{eq:mollifiererroruniquenesssolineq} shows that the limit exists.) 

Now consider the system \equpref{eq:approxRSequations}, which we know each $(u_k, \alpha_k, q_{n_k, d_k}, \rho_{n_k, d_k}, \delta_{1,k}, \delta_{2,k})$ solves by Proposition \upref{prop:approxrealRSequations}. 
Since $\rho^{\star} - q^{\star} > \frac1{\Csol}$, we then may rewrite $\psi_{\alpha}, \bar{\psi}_{\alpha}$ and hence $r, \bar{r}$ as per \equpref{eq:interpolationrrbardef1}, \equpref{eq:interpolationrrbardef2}. 

Since $\delta_{1,k}, \delta_{2,k} \rightarrow 0$, we may now take limits $k \rightarrow \infty$ on both sides of \equpref{eq:approxRSequations} (where now $r, \bar{r}$ are given as per \equpref{eq:interpolationrrbardef1}, \equpref{eq:interpolationrrbardef2}), since the denominators of \equpref{eq:approxRSequations} are bounded away from $0$ independently of $k$ as $r \ge \bar{r}$ by Lemma \upref{lem:realRSeqsolsboundedtechnicalcheck} and as $\beta \le \beta_1$. 
Using \equpref{eq:mollifier_error_limit} to compute the limit of the expressions involving the $u_k$'s and letting $\MARGIN' = \frac{\MARGIN-x\SIGNALDISORDER-\sqrt{q}\QUENCHEDDISORDER}{\sqrt{\rho-q}}$, it follows that $(q^{\star}, \rho^{\star})$ solves the system
\begin{align*}
r = \frac{\alpha^{\star}}{\rho-q} \E\Big[ \Big(\frac{\pdfnormal( \MARGIN' )}{\cN( \MARGIN')}\Big)^2 \Big] \quad \,&, \quad 
\bar{r} = \frac{\alpha^{\star}}{\rho-q} \E\Big[ \frac{\MARGIN' \pdfnormal( \MARGIN')}{\cN( \MARGIN')} \Big]\,, \\
\rho = \frac1{2\beta+r-\bar{r}} + \frac{r}{(2\beta+r-\bar{r})^2}\quad \,&, \quad q = \frac{r}{(2\beta+r-\bar{r})^2}\,.
\end{align*}
Manipulating the above system and using the identity $\cA'(x) = \cA(x)^2 - x\cA(x)$ where $\cA(x) = \frac{\pdfnormal(x)}{\cN(x)}$ denotes the Inverse Mills' Ratio (see Lemma \upref{lem:millsratioproperties}), it follows that $(q^{\star}, \rho^{\star})$ solves \equpref{eq:niceRSeqs}. 

Since $(q^{\star}, \rho^{\star}) \in [0, \Csol] \times \big[ \frac1{\Csol}, \Csol\big]$ and $\frac1{\rho^{\star}-q^{\star}} \le \Csol$, and as $(q^{\star}, \rho^{\star})$ solves \equpref{eq:niceRSeqs}, uniqueness given by Proposition \upref{prop:nicersuniquesol} now implies that $\rho^{\star} = \rho_0(\alpha^{\star}, \beta, x)$ and $q^{\star} = q_0(\alpha^{\star}, \beta, x)$. 
This contradicts \equpref{eq:mollifiererroruniquenesssolcontradict}.
\end{proof}
We now complete the proof of Theorem \upref{thm:mollifiedhamiltonianfreeenergy}.
As $\Csol$ only depends on $\beta_0, \beta_1$, we may let $H_1=H_1(\beta_0, \beta_1, \MARGIN, \lambda)$ denote an upper bound on the norm of the gradient of the $C^\infty$ function $\E\big[ \log \P_{\ANNEALEDDISORDER}\big(x\SIGNALDISORDER + \sqrt{q}\QUENCHEDDISORDER + \sqrt{\rho-q}\ANNEALEDDISORDER \ge \MARGIN\big) \big]$ over the compact set $\big\{ (x, \rho, q) : -1 \le x \le 1, 0 \le q \le \Csol, \frac1{\Csol} \le \rho \le \Csol, \frac1{\rho-q} \le \Csol\big\}$. Here the function is $C^\infty$ in this compact domain as $\rho-q \ge \frac1{\Csol}$. Dependence on $\lambda$ arises as the law of $S$ depends on $\lambda$. Without loss of generality suppose $H_1 \ge 1$.

By Proposition \upref{prop:nicersuniquesol}, we may also let $H_2=H_2(\beta_0, \beta_1, \MARGIN, \lambda)$ be an upper bound on the gradient of the $C^\infty$ function $\E \big[ \log \P_{\ANNEALEDDISORDER} \big( x\SIGNALDISORDER + \sqrt{q_0(\alpha, \beta, x)}\QUENCHEDDISORDER + \sqrt{\rho_0(\alpha, \beta, x) - q_0(\alpha, \beta, x) }\ANNEALEDDISORDER \ge \MARGIN \big) \big]$ on $(\alpha, \beta, x) \in [0,\DENSITYBOUND]\times [\beta_0, \beta_1] \times [-1,1]$. Here the function is $C^\infty$ in this compact domain as $\rho_0(\alpha, \beta, x) - q_0(\alpha, \beta, x) \ge \frac1{2\Csol}$ by Proposition \upref{prop:nicersuniquesol}.

Now, we let 
\begin{align*}
\delta = \delta(\eps, \beta_0, \beta_1, \MARGIN, \lambda) & := \frac{\eps}{12H_1}\,, \\
\eps' = \eps'(\eps, \beta_0, \beta_1, \MARGIN, \lambda) &:= \min\Big\{ \eps'_{\upref{prop:mollifiererror}}(\delta, \Csol, \MARGIN^+, \lambda ), \eps'_{\upref{prop:diffvsuniquesol}}(\delta) \Big\}\,, \\
\eps_{\RS_0} = \eps_{\RS_0}(\eps, u, \beta_0, \beta_1, \MARGIN, \SIGNALMEDIUMDEPENDENCE) &:= \min\Big\{ \eps''_{\upref{prop:boundedapproxRSeq}}(u, \beta_0, \beta_1, \SIGNALMEDIUMDEPENDENCE), \frac{\eps}{12 K_{\upref{prop:interpolationperconstraint}}(D, \beta_0, \beta_1, \MARGIN, \lambda)}, \eps''_{\upref{lem:interpolationlastcalc}}(\eps/4) \Big\}\,.
\end{align*}
Consider any $0 < \eps'' \le \eps_{\RS_0}$ as in the statement of this Theorem. (Note dependence on $\eps'$ in $\eps_{\RS_0}$ comes implicitly through $u$.) Finally we let 
\begin{align*}
d(\eps, \eps'', \beta_0, \beta_1, \MARGIN, \lambda) &:= \max\Big\{ d_{\upref{prop:boundedapproxRSeq}}(\eps''), \frac{12H_2}{\eps}, d_{\upref{lem:interpolationlastcalc}}(\eps/4), d_{\upref{prop:interpolationperconstraint}}(\eps'') \Big\}\,.
\end{align*}
Consider any $d \ge d(\eps, \eps'', \beta_0, \beta_1, \MARGIN, \lambda)$ as in the statement of this Theorem. We now consider any $n$ such that $\frac{n}{d} \in [0, \DENSITYBOUND]$. 

Next, we consider any $1 \le i \le n$ and let $\alpha = \frac{i}{d}$.
In what follows, we apply the results from Appendices \upref{subsec:interpolation1}, \upref{subsec:RSequationseliminatemollifier} with $i$ in place of $n$.
By Proposition \upref{prop:interpolationperconstraint} and as $d \ge d_{\upref{prop:interpolationperconstraint}}(\eps'')$ and $\eps'' \le \frac{\eps}{12K_{\upref{prop:interpolationperconstraint}}(D, \beta_0, \beta_1, \MARGIN, \SIGNALSHORTDEPENDENCE)}$, we have 
\begin{equation}\label{eq:interpolate_free_energy_finish1}
\begin{aligned}
&\Big| \phi_{\beta, i, \DIMENSION}(I) - \phi_{\beta, i-1, \DIMENSION}(I) - \frac1{\DIMENSION} \E \Big[ \log \E_{\ANNEALEDDISORDER} \exp u\big( x s_i + \sqrt{q_{i, \DIMENSION}} \QUENCHEDDISORDER + \sqrt{\rho_{i, \DIMENSION}-q_{i, \DIMENSION}}\ANNEALEDDISORDER \big) \Big] \Big| \\
&\quad \le \frac{K_{\upref{prop:interpolationperconstraint}}(D, \GAUSSIANMEASURE_0, \GAUSSIANMEASURE_1, \MARGIN, \SIGNALSHORTDEPENDENCE)}{\DIMENSION} \eps'' \le \frac{\eps}{12d}\,.
\end{aligned}
\end{equation}
By Proposition \upref{prop:boundedapproxRSeq}, as $\eps'' \le \eps''_{\upref{prop:boundedapproxRSeq}}(u, \beta_0, \beta_1, \SIGNALMEDIUMDEPENDENCE)$ and $d \ge d_{\upref{prop:boundedapproxRSeq}}(\eps'')$, we have 
\begin{align*}
\big( q_{i, d}, \rho_{i, d} \big) \in [0, \Csol] \times \Big[ \frac1{\Csol}, \Csol \Big]\,,\, \frac1{\rho_{i,d} - q_{i,d}} \le \Csol\,.
\end{align*}
Consequently by Proposition \upref{prop:mollifiererror}, as $\eps' \le \eps'_{\upref{prop:mollifiererror}}(\delta, \Csol, \MARGIN^+, \lambda)$ and $s_i \sim \cD$,
\begin{equation}\label{eq:interpolate_free_energy_finish2}
\begin{aligned}
&\Big| \E \Big[ \log \E_{\ANNEALEDDISORDER} \exp u\big( xs_i + \sqrt{q_{i, d}}\QUENCHEDDISORDER + \sqrt{\rho_{i, d}-q_{i, d}}\ANNEALEDDISORDER \big) \Big] \\
&\qquad \qquad \qquad - \E \Big[ \log \P_{\ANNEALEDDISORDER} \big( xs_i + \sqrt{q_{i, d}}\QUENCHEDDISORDER + \sqrt{\rho_{i, d}-q_{i, d}}\ANNEALEDDISORDER \ge \MARGIN \big) \Big] \Big| \le \delta \le \frac{\eps}{12}\, .
\end{aligned}
\end{equation}
Next, by Proposition \upref{prop:diffvsuniquesol}, we have 
\begin{align*}
\big| q_{i,d} - q_0(\alpha, \beta, x) \big|\,,\,\big|\rho_{i,d} - \rho_0(\alpha, \beta, x)\big| \le \delta = \frac{\eps}{12H_1}\,.
\end{align*}
Consequently
\begin{equation}\label{eq:interpolate_free_energy_finish3}
\begin{aligned}
&\Big| \E \Big[ \log \P_{\ANNEALEDDISORDER} \big( xs_i + \sqrt{q_{i, d}}\QUENCHEDDISORDER + \sqrt{\rho_{i, d}-q_{i, d}}\ANNEALEDDISORDER \ge \MARGIN \big) \Big] \\
&\qquad \qquad - \E \Big[ \log \P_{\ANNEALEDDISORDER} \big( xs_i + \sqrt{q_0(\alpha, \beta, x)}\QUENCHEDDISORDER + \sqrt{\rho_0(\alpha, \beta, x) - q_0(\alpha, \beta, x) }\ANNEALEDDISORDER \ge \MARGIN \big) \Big] \Big| \\
&\quad \le \frac{\eps \sqrt{2}}{12H_1} \cdot H_1 < \frac{\eps}8\,.
\end{aligned}
\end{equation}
Finally, observe that for any $\alpha' \in \big[\frac{i-1}{d}, \frac{i}{d} \big]$, $\big|\alpha - \alpha'\big| \le \frac1d$, so as $d \ge 12H_2/\eps$ we have
\begin{equation}\label{eq:interpolate_free_energy_finish4}
\begin{aligned}
&\Big| \E \Big[ \log \P_{\ANNEALEDDISORDER} \big( xs_i + \sqrt{q_0(\alpha, \beta, x)}\QUENCHEDDISORDER + \sqrt{\rho_0(\alpha, \beta, x) - q_0(\alpha, \beta, x) }\ANNEALEDDISORDER \ge \MARGIN \big) \Big] \\
&\qquad \qquad - \E \Big[ \log \P_{\ANNEALEDDISORDER} \big( xs_i + \sqrt{q_0(\alpha', \beta, x)}\QUENCHEDDISORDER + \sqrt{\rho_0(\alpha', \beta, x) - q_0(\alpha', \beta, x) }\ANNEALEDDISORDER \ge \MARGIN \big) \Big] \Big| \\
&\quad \le H_2 \cdot \frac{1}{d} \le \frac{\eps}{12}\,.
\end{aligned}
\end{equation}
By Proposition \upref{prop:nicersuniquesol}, as $\beta, x, I$ are fixed, we may consider $\frac{\rmd}{\rmd \alpha} \RS_{0,I}(\alpha)$. Since $\big(q_0(\alpha, \beta, x), \rho_0(\alpha, \beta, x)\big)$ solves $\frac{\partial F_I}{\partial q} = \frac{\partial F_I}{\partial \rho}=0$, we have
\begin{align*}
\frac{\rmd}{\rmd \alpha} \RS_{0,I}(\alpha) &= \frac{\partial F_I}{\partial \alpha} \big(x, q_0(\alpha, \beta, x), \rho_0(\alpha, \beta, x)\big) \\
&\qquad \qquad + \frac{\partial F_I}{\partial q} \big(x, q_0(\alpha, \beta, x), \rho_0(\alpha, \beta, x)\big) \cdot \frac{\partial}{\partial \alpha} q_0(\alpha, \beta, x) \\
&\qquad \qquad + \frac{\partial F_I}{\partial \rho} \big(x, q_0(\alpha, \beta, x), \rho_0(\alpha, \beta, x)\big) \cdot \frac{\partial}{\partial \alpha} \rho_0(\alpha, \beta, x) \\
&=\E \Big[ \log \P_{\ANNEALEDDISORDER} \big( xs_i + \sqrt{q_0(\alpha, \beta, x)}\QUENCHEDDISORDER + \sqrt{\rho_0(\alpha) - q_0(\alpha, \beta, x) }\ANNEALEDDISORDER \ge \MARGIN \big) \Big]\,,
\end{align*}
as $s_i \sim \cD$ and as independent of $Z, W$. The Mean Value Theorem now implies that for some $\alpha' \in \big[\frac{i-1}{d}, \frac{i}{d} \big]$,
\begin{equation}\label{eq:interpolate_free_energy_finish5}
\begin{aligned}
&\RS_0\Big(\frac{i}{d}\Big) - \RS_{0,I}\Big(\frac{i-1}{d}\Big) \\
&\quad = \frac1{d} \E \Big[ \log \P_{\ANNEALEDDISORDER} \big( xs_i + \sqrt{q_0(\alpha', \beta, x)}\QUENCHEDDISORDER + \sqrt{\rho_0(\alpha') - q_0(\alpha', \beta, x) }\ANNEALEDDISORDER \ge \MARGIN \big) \Big]\,.
\end{aligned}
\end{equation}
Combining \equpref{eq:interpolate_free_energy_finish1}, \equpref{eq:interpolate_free_energy_finish2}, \equpref{eq:interpolate_free_energy_finish3}, \equpref{eq:interpolate_free_energy_finish4}, \equpref{eq:interpolate_free_energy_finish5} implies for all $1 \le i \le n$,
\begin{align*}
\Big| \phi_{\beta, i, d}(I) - \phi_{\beta, i-1, d}(I) - \Big( \RS_{0,I}\Big(\frac{i}{d}\Big) - \RS_{0,I}\Big(\frac{i-1}{d}\Big) \Big) \Big| < \frac{3\eps}{8d}\,.
\end{align*}
Summing across all such $i$ and noting there are at most $n \le \DENSITYBOUND d \le 2d$ such $i$ yields 
\begin{align*}
\Big| \phi_{\beta, n, d}(I) - \phi_{\beta, 0, d}(I) - \Big( \RS_{0,I}\Big(\frac{n}{d}\Big) - \RS_{0,I}(0) \Big) \Big| < \frac{3\eps}4\,.
\end{align*}
Last, we combine with Lemma \upref{lem:interpolationlastcalc}, yielding
\begin{align*}
\Big| \phi_{\beta, n, d}(I) - \Big( \RS_{0,I}\Big(\frac{n}{d}\Big) + \frac12 \log(2\pi e) \Big) \Big| < \eps\,.
\end{align*}
This completes the proof of Theorem \upref{thm:mollifiedhamiltonianfreeenergy}.

\subsection{Proof of Theorem \upref{thm:posterior_mollifiedhamiltonianfreeenergy} for posterior}\label{subsec:posterior_finishinterpolation}
Here we prove Theorem \upref{thm:posterior_mollifiedhamiltonianfreeenergy} with the same argument as the proof of Theorem \upref{thm:mollifiedhamiltonianfreeenergy}. 
Note $u(x) \ge -D(\lambda)(1+|x|)$ for $u(x)$ from \equpref{eq:defuSgmm}, \equpref{eq:defuSlogistic}.
Next, by the same argument as the proof of Proposition \upref{prop:diffvsuniquesol}, we have the following for any $(\beta, x) \in [\beta_0, \beta_1] \times [-1, 1]$. For all $\eps'' \le \eps''_{\upref{prop:boundedapproxRSeq}}$, all $d \ge d_{\upref{prop:boundedapproxRSeq}}$, and all $n \ge 1$ with $n/d = \alpha$,
\begin{align}\label{eq:posterior_diffvsuniquesol}
\big| q_{n,d} - q_0(\alpha, \beta, x) \big|\, ,\, \big|\rho_{n,d} - \rho_0(\alpha, \beta, x) \big| \le \delta\,.
\end{align}
Here the proof follows by the same argument as the proof of Proposition \upref{prop:diffvsuniquesol}, using Propositions \upref{prop:nicersuniquesol}, \upref{prop:approxrealRSequations}, and \upref{prop:boundedapproxRSeq}, but now we do not consider the $u_k$ or replace $u_k$ by $\one\{x \ge \MARGIN\}$ as $k \rightarrow \infty$. Here we use that a solution to \equpref{eq:realRSequations} solves \equpref{eq:niceRSeqs} where $u(\cdot)$ is identical in both systems; this can be seen via Gaussian integration by parts.

Now to prove Theorem \upref{thm:mollifiedhamiltonianfreeenergy}, we let $H_1=H_1(\beta_0, \beta_1, \lambda)$ be an upper bound on the gradient of the continuously differentiable function $\E\big[ \log \E_{\ANNEALEDDISORDER} \exp u\big(x\SIGNALDISORDER + \sqrt{q}\QUENCHEDDISORDER + \sqrt{\rho-q}\ANNEALEDDISORDER \big) \big]$ over the compact set $\big\{ (x, \rho, q) : -1 \le x \le 1, 0 \le q \le \Csol, \frac1{\Csol} \le \rho \le \Csol, \frac1{\rho-q} \le \Csol\big\}$. Here we can verify that the function is continuously differentiable as $\rho-q \ge \frac1{\Csol}$, as $u$ is twice differentiable, and as $u(x) \ge -D(\lambda)(1+|x|)$. Without loss of generality suppose $H_1 \ge 1$.

Similarly by Proposition \upref{prop:nicersuniquesol}, we may let $H_2=H_2(\beta_0, \beta_1, \lambda)$ be an upper bound on the gradient of the continuously differentiable $\E \big[ \log \E_{\ANNEALEDDISORDER} \exp u\big( x\SIGNALDISORDER + \sqrt{q_0(\alpha, \beta, x)}\QUENCHEDDISORDER + \sqrt{\rho_0(\alpha, \beta, x) - q_0(\alpha, \beta, x) }\ANNEALEDDISORDER \big) \big]$ on $(\alpha, \beta, x) \in [0,\DENSITYBOUND]\times [\beta_0, \beta_1] \times [-1,1]$. Here we use that $\rho_0(\alpha, \beta, x) - q_0(\alpha, \beta, x) \ge \frac1{2\Csol}$ by Proposition \upref{prop:nicersuniquesol}.
Now, we let 
\begin{align*}
\delta = \delta(\eps, \beta_0, \beta_1, \lambda) &:= \frac{\eps}{12H_1}\,, \\
\eps_{\RS_0} = \eps_{\RS_0}(u, \beta_0, \beta_1, \SIGNALMEDIUMDEPENDENCE) &:= \min\Big\{ \eps''_{\upref{prop:boundedapproxRSeq}}(u, \beta_0, \beta_1, \SIGNALMEDIUMDEPENDENCE), \frac{\eps}{12 K_{\upref{prop:interpolationperconstraint}}(D, \beta_0, \beta_1, \lambda)}, \eps''_{\upref{lem:interpolationlastcalc}}(\eps/4) \Big\}\,.
\end{align*}
Considering any $0 < \eps'' \le \eps_{\RS_0}$ as in the statement of this Theorem, we let
\begin{align*}
d(\eps, \eps'', \beta_0, \beta_1, \lambda) &:= \max\Big\{ d_{\upref{prop:boundedapproxRSeq}}(\eps''), \frac{12H_2}{\eps}, d_{\upref{lem:interpolationlastcalc}}(\eps/4), d_{\upref{prop:interpolationperconstraint}}(\eps'') \Big\}\,.
\end{align*}
Consider any $d \ge d(\eps, \eps'', \beta_0, \beta_1, \lambda)$ as in the statement of this Theorem. We now consider any $n$ such that $\frac{n}{d} \in [0, \DENSITYBOUND]$. 
For any $1 \le i \le n$, letting $\alpha = \frac{i}{d}$, we obtain by Proposition \upref{prop:interpolationperconstraint} that 
\begin{align} 
&\Big| \phi_{\beta, i, \DIMENSION}(I) - \phi_{\beta, i-1, \DIMENSION}(I) - \frac1{\DIMENSION} \E \Big[ \log \E_{\ANNEALEDDISORDER} \exp u\big( x s_i + \sqrt{q_{i, \DIMENSION}} \QUENCHEDDISORDER + \sqrt{\rho_{i, \DIMENSION}-q_{i, \DIMENSION}}\ANNEALEDDISORDER \big) \Big] \Big| \le \frac{\eps}{12d}\,.\notag
\end{align}
By Proposition \upref{prop:boundedapproxRSeq}, we have 
\begin{align*}
\big( q_{i, d}, \rho_{i, d} \big) \in [0, \Csol] \times \Big[ \frac1{\Csol}, \Csol \Big]\,,\, \frac1{\rho_{i,d} - q_{i,d}} \le \Csol\,.
\end{align*}
Next by \equpref{eq:posterior_diffvsuniquesol}, we have 
\begin{align*}
\big| q_{i,d} - q_0(\alpha, \beta, x) \big|\,,\,\big|\rho_{i,d} - \rho_0(\alpha, \beta, x)\big| \le \delta = \frac{\eps}{12H_1}\,.
\end{align*}
Consequently 
\begin{align}
&\Big| \E \Big[ \log \E_{\ANNEALEDDISORDER} \exp u\big( xs_i + \sqrt{q_{i, d}}\QUENCHEDDISORDER + \sqrt{\rho_{i, d}-q_{i, d}}\ANNEALEDDISORDER\big) \Big] \notag \\
&\qquad \qquad - \E \Big[ \log \E_{\ANNEALEDDISORDER} \exp u\big( xs_i + \sqrt{q_0(\alpha, \beta, x)}\QUENCHEDDISORDER + \sqrt{\rho_0(\alpha, \beta, x) - q_0(\alpha, \beta, x) }\ANNEALEDDISORDER\big) \Big] \Big| < \frac{\eps}8\,. \notag
\end{align}
Finally, we observe that for any $\alpha' \in \big[\frac{i-1}{d}, \frac{i}{d} \big]$, we have $\big|\alpha - \alpha'\big| \le \frac1d$, so as $d \ge 12H_2/\eps$, 
\begin{align}
&\Big| \E \Big[ \log \E_{\ANNEALEDDISORDER} \exp u\big( xs_i + \sqrt{q_0(\alpha, \beta, x)}\QUENCHEDDISORDER + \sqrt{\rho_0(\alpha, \beta, x) - q_0(\alpha, \beta, x) }\ANNEALEDDISORDER \big) \Big] \notag \\
&\qquad \qquad - \E \Big[ \log \E_{\ANNEALEDDISORDER} \exp u\big( xs_i + \sqrt{q_0(\alpha', \beta, x)}\QUENCHEDDISORDER + \sqrt{\rho_0(\alpha', \beta, x) - q_0(\alpha', \beta, x) }\ANNEALEDDISORDER \big) \Big] \Big| \le \frac{\eps}{12}\,. \notag 
\end{align}
With the above bounds in hand, we finish identically to the proof of Theorem \upref{thm:mollifiedhamiltonianfreeenergy}, using Proposition \upref{prop:nicersuniquesol} and Lemma \upref{lem:interpolationlastcalc}. This proves Theorem \upref{thm:posterior_mollifiedhamiltonianfreeenergy}.

\section{Concentration with respect to Gaussian Measure}\label{sec:gaussianconversion}
The goal of this section is to show that the log-partition function concentrates with respect to the Gaussian measure, as stated in Theorem \upref{thm:gaussianmeasurersformula} for interpolators and Theorem \upref{thm:posteriorgaussianmeasurersformula} for the posterior below. 
Recall that in Theorems \upref{thm:mollifiedhamiltonianfreeenergy}, \upref{thm:posterior_mollifiedhamiltonianfreeenergy}, we computed the value of the free energy, that is the expected normalized log-partition function, corresponding to the interpolators and the posterior respectively. Here, we upgrade these statements to exponentially high-probability concentration of the normalized log-partition function about this same value.

\paragraph{Theorem \upref{thm:gaussianmeasurersformula} for interpolators:} Recall the definition of the set $U_i$ interpolating the $i$-th datapoint $i=1, \ldots, n$ (equivalently, satisfying each constraint) defined in (\ref{eq:Ui_init_def}), and consider $C$, the intersection of all the $U_i$. Recalling the definition of $S_i$ in \equpref{eq:Si_def}, these sets are defined by
\begin{align}
U_i &:= \bigl\{ \theta \in \R^d \, :\, S_i \ge \MARGIN \bigr\}\quad\text{and}\quad C := \bigcap_{i=1}^n U_i \,.\label{eq:hyperplane_intersection_def} 
\end{align}
In Theorem \upref{thm:mollifiedhamiltonianfreeenergy}, we showed $\phi_{\beta, n, d}(I)$ from \equpref{eq:restricted_free_energy_mollified} equals $\RSG_I(n/d, \beta, x)$, where $\phi_{\beta, n, d}(I)$ is the normalized logarithm of \equpref{eq:restricted_partition_function_mollified}. 
We now show after replacing the mollifier $\sum_{1 \le i \le n} u\big(\big( s_i \theta_1 + \langle g_i, \bar\theta\rangle \big)/\sqrt{d} \big)$ in \equpref{eq:restricted_partition_function_mollified} by the indicator $\one\{\bigcap_{1 \le i \le n} U_i\}$, the corresponding normalized log-partition-function concentrates exponentially around the same value. Note this hard indicator arises as we want to study the set of interpolators.
We now are in a position to state Theorem \upref{thm:gaussianmeasurersformula}:
\begin{theorem}[For interpolators]\label{thm:gaussianmeasurersformula}
Consider $n$ such that $1 \le n \le \DENSITYBOUND d$. Then for any $\eps>0$, there is $\eps_{\RSG}>0$ depending on $\eps, \DENSITYBOUND, \beta_0, \beta_1, \MARGIN, \lambda$ such that the following holds. For any $\GAUSSIANMEASURE \in [\GAUSSIANMEASURE_0, \GAUSSIANMEASURE_1]$, recalling the definition of $C$ in \equpref{eq:hyperplane_intersection_def}, define for any interval $I \subseteq \R$,
\begin{align}
Z_{\beta}(I) &:= \int \one\big\{\{\theta_1/\sqrt{d} \in I \} \cap C\big\}\, \exp\big( -\GAUSSIANMEASURE \|\btheta\|^2 \bigr)\, \rmd \btheta\, ,\label{eq:free_energy_gaussian_def}
\end{align}
Then for any $x \in [-1,1]$ and $0 < \eps'' \le \eps_{\RSG}$, letting $I=[x-\eps'', x+\eps'']$, we have for $d$ large enough in terms of $\eps, \eps'', \DENSITYBOUND, \beta_0, \beta_1, \MARGIN, \lambda$ that
\begin{align}\label{eq:gaussian_conversion_convergence}
&\P\Bigg( \Big| \frac1d \log Z_{\beta}(I) - \RSG_I(n/d, \GAUSSIANMEASURE, x) \Big| \ge \eps \Bigg) \le \exp\big(-d/K\big)\,,
\end{align}
where $\RSG_I$ is defined as in \equpref{eq:RSGdef} with $\exp u(x) = \one\{x \ge \MARGIN\}$, and where $K$ is large enough in terms of $\eps, \DENSITYBOUND, \beta_0, \beta_1, \MARGIN, \lambda$.
\end{theorem}
We prove Theorem \upref{thm:gaussianmeasurersformula} in Appendix \upref{subsec:gaussianconversioninterpolants}. 
Specifically, we must handle the non-Gaussian signal component to establish several key technical results on stochastic processes and concentration of measure, which establishes the key `add one constraint' style estimate Theorem \upref{thm:gaussiansetup}. This result upper bounds the probability that a Gibbs measure places small mass on one side of a given constraint. This is done in Appendix \upref{subsec:gaussian_measure_preliminaries}. To handle the non-Gaussianity, it turns out we need the width $\eps''$ of $I$ to be small enough in terms of $\lambda$ but independent of $d$. See the discussion at the start of Appendix \upref{subsec:gaussian_measure_preliminaries} for more details. 

\paragraph{Theorem \upref{thm:posteriorgaussianmeasurersformula} for posterior:} 
We have the following similar result for the posterior, which establishes exponential concentration of the log-partition function about $\RSG_I$. 
Note this Theorem does not involve hard indicators as we are interested in studying the posterior $p_{X,y}$ in \equpref{eq:posterior}, which is defined in terms of $u(x)$ from \equpref{eq:defuSgmm}, \equpref{eq:defuSlogistic} in the GMM and logistic cases respectively. 
\begin{theorem}[For posterior]\label{thm:posteriorgaussianmeasurersformula}
Consider $n$ such that $1 \le n \le \DENSITYBOUND d$. Consider $u(x)$ from \equpref{eq:defuSgmm}, \equpref{eq:defuSlogistic} in the GMM and logistic cases respectively. Then for any $\eps>0$, there is $\eps_{\RSG}>0$ depending on $\eps, \DENSITYBOUND, \beta_0, \beta_1, \lambda$ such that the following holds. For any $\GAUSSIANMEASURE \in [\GAUSSIANMEASURE_0, \GAUSSIANMEASURE_1]$, define for any interval $I \subseteq \R$,
\begin{align}
Z_{\beta}(I) &:= \int \one\big\{  \theta_1/\sqrt{\DIMENSION} \in I \big\} \exp\Big(\sum_{i=1}^n u( S_i ) - \GAUSSIANMEASURE \| \theta \|^2 \Big)\, \rmd \btheta\, ,\label{eq:free_energy_posteriorgaussian_def}
\end{align}
where $S_i$ is defined as in \equpref{eq:Si_def}. Then for any $x \in [-1,1]$ and $0 < \eps'' \le \eps_{\RSG}$, letting $I=[x-\eps'', x+\eps'']$, we have for $d$ large enough in terms of $\eps, \eps'', \DENSITYBOUND, \beta_0, \beta_1, \lambda$ that
\begin{align}\label{eq:posteriorgaussian_conversion_convergence}
&\P\Bigg( \Big| \frac1d \log Z_{\beta}(I) - \RSG_I(n/d, \GAUSSIANMEASURE, x) \Big| \ge \eps \Bigg) \le \exp\big(-d/K\big)\,,
\end{align}
where $\RSG_I$ is defined as in \equpref{eq:RSGdef} with this $\exp u(x)$, and where $K$ is large enough in terms of $\eps, \DENSITYBOUND, \beta_0, \beta_1, \lambda$.
\end{theorem}
We prove Theorem \upref{thm:posteriorgaussianmeasurersformula} in Appendix \upref{subsec:posteriorgaussianconversion}. Unlike Theorem \upref{thm:gaussianmeasurersformula}, the proof of Theorem \upref{thm:posteriorgaussianmeasurersformula} is relatively direct as a consequence of concentration of Lipschitz functions of log-concave measures, as stated below. Note that in contrast, such tools are not available to prove Theorem \upref{thm:gaussianmeasurersformula} due to the hard indicators in the Hamiltonian therein.

\paragraph{Concentration of Lipschitz functions log-concave measures.} We next introduce concentration of Lipschitz functions of strongly log-concave measures in the form given in Theorem \upref{thm:logconcaveconcentrationbasics} below, originally credited to B. Maurey. See also the text of Bakry, Gentil, and Ledoux \cite{bakry2013analysis}, Chapter 5. This fact allows us to establish the high-probability concentration bounds of both Theorem \upref{thm:gaussianmeasurersformula} and Theorem \upref{thm:posteriorgaussianmeasurersformula}. We will also use this fact crucially in proving Proposition \upref{prop:logconcaveconcentration} in Appendix \upref{subsec:logconcaveconcentration}. 

Note that this Theorem holds when $\psi$ can take on the value $-\infty$, that is when $\psi$ is supported on a convex set that is a strict subset of $\R^d$, as noted in (3.21) of \cite{talagrand2010mean}. Technically, this result is written in \cite{talagrand2010mean} when $f$ is $H$-Lipschitz on all of $\R^d$; the result writte in \cite{talagrand2010mean} immediately implies the following by the McShane-Whitney Extension Theorem \cite{mcshane1934extension, whitney1934analytic}.
\begin{theorem}[Theorem 3.1.4, \cite{talagrand2010mean}]\label{thm:logconcaveconcentrationbasics}
Consider a measure $\mu(\btheta) \propto \exp\Big( \psi(\btheta) \Big)$ defined over a convex set $C$ such that for some $\GAUSSIANMEASURE>0$, we have 
\[ \frac12 \Big( \psi(\btheta_1) + \psi(\btheta_2) \Big) - \psi\Big( \frac{\btheta_1+\btheta_2}2 \Big) \le -\beta \Big\| \frac{\btheta_1-\btheta_2}2 \Big\|^2~. \]
(Note $\psi$ can take on the value $-\infty$ in the above.)
Then for any set $C \subseteq \R^d$,
\begin{equation*}
\int \exp \Big( \frac{\beta}{2} d^2(x, C) \Big)\, \rmd\mu(x) \leq \frac{1}{\mu(C)}\, ,
\end{equation*}
where $d(x, C) = \inf\{d(x, y) : y \in C\}$ is the distance from $x$ to $C$. Moreover, if $f$ is a function on $\R^d$ with Lipschitz constant $H$ on the support of $\mu$, i.e., for all $x, y \in \text{supp}(\mu)$ we have $\| f(x) - f(y) \| \leq H \| x-y \|$, then
\begin{equation*}
\int \exp\Big( \frac{\GAUSSIANMEASURE}{8H^2} \big( f(x) - \int f\, \rmd\mu \big)^2 \Big) ~\rmd\mu(x) \leq 4~,
\end{equation*}
and
\begin{equation*}
\text{ for all } k \geq 1\, ,\, \quad
\int \left( f(x) - \int f\rmd\mu \right)^{2k} \rmd\mu(x)
\leq 4 \left( \frac{8kH^2}{\beta} \right)^k\, .
\end{equation*}
\end{theorem}

\subsection{Proof of Theorem \upref{thm:gaussianmeasurersformula} for interpolators}\label{subsec:gaussianconversioninterpolants}
Here throughout Appendix \upref{subsec:gaussianconversioninterpolants}, we let $Z_{\beta}(I)$ be as per \equpref{eq:free_energy_gaussian_def}. The goal of this section is to prove Theorem \upref{thm:gaussianmeasurersformula}, and to this end there are two main steps.
\begin{enumerate}
    \item First, we establish high-probability concentration of the normalized logarithm of \equpref{eq:restricted_partition_function_mollified} around its expectation $\RSG_I(n/d, \beta, x)$, and more generally, concentration of the normalized free energy around its expectation. In particular, this also applies for \equpref{eq:free_energy_gaussian_def}. 
This is done via a technically involved application of Bernstein's Inequality in Lemma \upref{lem:controlmollifierbernsteinspiked}.
\item Second, we show that the mollifier $u(\cdot)$ in \equpref{eq:restricted_partition_function_mollified} can be eliminated without significantly changing the expectation of its logarithm, as discussed in the proof ideas in Section \upref{sec:proofideas}. This is a delicate step done in Lemma \upref{lem:controlmollifiers}. 
\end{enumerate}
With both of these steps in hand, we then prove Theorem \upref{thm:gaussianmeasurersformula}. All these steps are done in Appendix \upref{subsec:gaussian_measure_rs_pf}. Preliminary steps on stochastic processes and concentration of measure necessary to prove these Lemmas -- specifically, the key `add one constraint' style estimate, Theorem \upref{thm:gaussiansetup} -- are presented in Appendix \upref{subsec:gaussian_measure_preliminaries}. These initial steps are important but highly technical.

\subsubsection{Preliminaries}\label{subsec:gaussian_measure_preliminaries}
Central to the subsequent proof in Appendix \upref{subsec:gaussian_measure_rs_pf} will be the combination of the following Theorem \upref{thm:gaussiansetup} and Lemma \upref{lem:momentsofRV838tal}. These are important results on concentration of measure that let us control hard-to-bound quantities which involve the addition of the $i$-th datapoint or constraint, which behave at an exponential scale. In turn, these results hinge on the key Lemma \upref{lem:keyconcentrationgaussianconversion} on stochastic processes, which is proved by several comparison inequalities. 

Although similar in spirit, the statements and proofs differ from those in Section 8.2 of \cite{talagrand2011advanced}, as handling the non-Gaussian `signal' component $s_i \sim \cD$ of the disorder requires more technical ingredients. Notably, these results hinge on being able to choose the width $\eps''$ of $I$ small enough in terms of $\lambda$, but still independently of $d$, so that for $d \ge d(\eps'')$, various partition functions are lower bounded by $\exp(-O(d))$ where $O(\cdot)$ is independent of $\eps''$. In particular, see the assumption $Z \ge \exp(-ad)$ in Theorem \upref{thm:gaussiansetup}; this assumption applies in the proof of Theorem \upref{thm:gaussianmeasurersformula} with exponentially high probability by Lemma \upref{lem:controlmollifierbernsteinspiked}.

Recall as per Appendix \upref{sec:interpolationfreeenergy} that $\signalsubexp, \SIGNALGAUSSIAN$ depend on $\lambda$. However the following bounds depend delicately on these parameters, so here in Appendix \upref{subsec:gaussian_measure_preliminaries} and only here, we will make the dependence on $\signalsubexp, \SIGNALGAUSSIAN$ fully explicit.

We begin with several preliminary results on stochastic processes.
\begin{lemma}\label{lem:gaussiancomparisonsoftmax}
Let $\xi_1$ and $(X'_\ell)_{1 \le \ell \le L}$ be a collection of real valued random variables where $\xi_1$ has zero mean, and $\xi_1$ and the collection $(X'_{\ell})_{1 \le \ell \le L}$ are independent (but the elements of $(X'_{\ell})_{1 \le \ell \le L}$ are not necessarily independent of each other). 
Let $u_1,\ldots,u_L \in \R$ be fixed such that $|u_\ell - x| \le \eps''$ for all $1 \le \ell \le L$ for some $x \in \R$. 
Define $X=(X_\ell)_{\ell \le L}$ by $X_\ell := u_\ell \xi_1 + X'_\ell$ for all $1\le \ell \le L$. 
Then for all $t \ge 0$,
\[\E\Big[ \log \sum_{\ell=1}^L \exp(t X_\ell) \Big] \ge \E\Big[ \log \sum_{\ell=1}^L \exp(t X'_\ell) \Big] \,.\]
\end{lemma}
\begin{proof}
Fix $s\in\R$. Since the function $F(x) := \log \sum_{\ell=1}^L \exp(t x_\ell)$ for $x\in\R^L$ is convex, we have 
\[F(X) \ge F(X') + \langle \nabla F(X'), X-X'\rangle\,.\]
Let $u = (u_\ell)_{\ell=1}^L$. The expectation of the second term on the right-hand side is $\E \big[\xi_1\langle \nabla F(X'), u\rangle \big] = 0$ as $\xi_1$ is zero mean and independent of $X'$, proving the Lemma. 
\end{proof}
\begin{lemma}[Slepian's Inequality, in the form of Proposition 8.2.2, \cite{talagrand2011advanced}]\label{lem:gaussiansoftmaxcomparetoindependent}
Consider two jointly Gaussian families $ (U_{\ell})_{\ell \le L}$, $(V_{\ell})_{\ell \le L}$, and suppose that for all $1 \le \ell \le L$ we have $\E[ U_{\ell}^2 ] \ge \E[ V_{\ell}^2 ]$, and for all $1 \le \ell_1 \neq \ell_2 \le L$ we have $\E[ U_{\ell_1} U_{\ell_2} ] \le \E[ V_{\ell_1} V_{\ell_2} ]$. Then 
\begin{align*}
\E\Big[ \log \sum_{1 \le \ell \le L} \exp(tU_{\ell}) \Big] \ge \E\Big[ \log \sum_{1 \le \ell \le L} \exp(tV_{\ell}) \Big]\,.
\end{align*}
\end{lemma}
\begin{lemma}[Proposition 8.2.3, \cite{talagrand2011advanced}]\label{lem:lowerboundsoftmaxindependent}
There exists a $K_{\upref{lem:lowerboundsoftmaxindependent}}>0$ with the following property. Letting $(\ANNEALEDDISORDER_{\ell})_{\ell \le L} \sim N(0, 1)$ be i.i.d., for $K_{\upref{lem:lowerboundsoftmaxindependent}} \le t \le \sqrt{\log L}/K_{\upref{lem:lowerboundsoftmaxindependent}}$, we have
\begin{align*}
\E\Big[ \log \sum_{\ell \le L} \exp(t \ANNEALEDDISORDER_{\ell}) \Big] \ge \log L + \frac{t^2}5\,.
\end{align*}
\end{lemma}

We will now prove the crucial Lemma \upref{lem:keyconcentrationgaussianconversion} that combines the above results. Its aim is to upper bound the probability of the following event: that for $L$ replicas, the number of $\ell, 1 \le \ell \le L$ such that $\frac1{\sqrt{d}}\big( \theta_1^{\ell} s + \langle \bar\theta^{\ell}, \bar g \rangle \big) \ge \MARGIN$ is small. The proof is involved and spans the next several pages. We will then leverage Lemma \upref{lem:keyconcentrationgaussianconversion} to prove the important Theorem \upref{thm:gaussiansetup}.
\begin{lemma}\label{lem:keyconcentrationgaussianconversion}
Consider any $x \in [-1,1]$ and $0 \le c_1 < c_2 < c_3$ with
\begin{align}\label{eq:keyconcentration_gaussianconversion_ineq}
c_3 \ge \frac1{c_2 - c_1}\quad\text{and}\quad c_1 \ge x^2 \E[s^2]\,.
\end{align}
Then there exists $K_{\upref{lem:keyconcentrationgaussianconversion}} \ge 1$ depending on  $c_3, \signalsubexp, \SIGNALGAUSSIAN$ and $K'_{\upref{lem:keyconcentrationgaussianconversion}} \ge \max\big\{ 12 c_3 \signalsubexp, 1\big\}$ depending on $c_3, \signalsubexp$ satisfying the following. 
Consider $0 < \eps'' \le \frac1{K'_{\upref{lem:keyconcentrationgaussianconversion}}}$ and 
\begin{align}
\btheta^1, \ldots, \btheta^L \in \R^d\quad\text{such that}\quad\theta_1^{\ell}/\sqrt{d} \in [x-\eps'', x+\eps'']\, \forall \, 1 \le \ell \le L\,.\label{eq:keyconcentrationgaussianconversion_narrowinterval}
\end{align}
Consider independent $s \sim \cD$, $\bar{g} \sim N(0, I_{d-1})$, and a family of random variables $(X_{\ell})_{\ell \le L}$ defined by
\begin{align*}
X_{\ell} = \frac{\theta_1^{\ell} s + \langle \bar\theta^{\ell}, \bar{g} \rangle}{\sqrt{d}}\,,
\end{align*}
such that
\begin{align}
&c_2 \le \E\big[ X^2_{\ell} \big] \le c_3\text{ for all }1 \le \ell \le L\quad\text{ and }\quad \E\big[ X_{\ell_1} X_{\ell_2} \big] \le c_1\text{ for all }1 \le \ell_1 \neq \ell_2 \le L\,.
\end{align}
Then for $K_{\upref{lem:keyconcentrationgaussianconversion}} \le t \le \sqrt{\log L}/K_{\upref{lem:keyconcentrationgaussianconversion}}$, we have 
\begin{align}
&\P\Big( \#\Big\{ \ell \le L\,:\,X_{\ell} \ge \frac{t}{K_{\upref{lem:keyconcentrationgaussianconversion}}} \Big\} \le L\exp\big(-K_{\upref{lem:keyconcentrationgaussianconversion}}t^2\big) \Big) \le K_{\upref{lem:keyconcentrationgaussianconversion}} \exp\Big( -\frac{t^2}{K_{\upref{lem:keyconcentrationgaussianconversion}}}\Big)\,.\label{eq:keyconcentrationgaussianconversion_pt1}
\end{align}
Consequently, for all $L^{-1/K_{\upref{lem:keyconcentrationgaussianconversion}}} \le t' \le \exp\big( -K_{\upref{lem:keyconcentrationgaussianconversion}} \max\{1, \MARGIN\}^2 \big)$, 
\begin{align}
\P\Big( \#\Big\{ \ell \le L\,:\,X_{\ell} \ge \MARGIN\Big\} \le Lt' \Big) \le K_{\upref{lem:keyconcentrationgaussianconversion}}\, t'^{1/K_{\upref{lem:keyconcentrationgaussianconversion}}^2}\,.\label{eq:keyconcentrationgaussianconversion_pt2}
\end{align}
\end{lemma}
\begin{proof}
We will focus on proving \equpref{eq:keyconcentrationgaussianconversion_pt1}. The proof of \equpref{eq:keyconcentrationgaussianconversion_pt2} then directly follows from \equpref{eq:keyconcentrationgaussianconversion_pt1}.

\paragraph{Proof of \equpref{eq:keyconcentrationgaussianconversion_pt1}.} The idea behind proving \equpref{eq:keyconcentrationgaussianconversion_pt1} is to apply the second moment method to $\exp(t X_{\ell})$. Specifically, let us consider $t$ such that $K_{\upref{lem:keyconcentrationgaussianconversion}} \le t \le \sqrt{\log L}/K_{\upref{lem:keyconcentrationgaussianconversion}}$ and define the following event $E$:
\begin{align}
E &:= \Big\{ \sum_{1\le \ell \le L} \exp(t X_{\ell}) \ge L \exp\Big( \frac{t^2}{60 c_3} \Big) \Big\} \label{eq:keyconcentration_eventE_1}\\
&\qquad\qquad \cap \Big\{ \sum_{1 \le \ell \le L} \exp(2t X_{\ell}) \le L \exp\big( (1+K(c_3, \signalsubexp)) t^2 \big)\, \Big\}\,. \label{eq:keyconcentration_eventE_2}
\end{align}
The main claim is that 
\begin{align}\label{eq:keyconcentration_lowerbd_PE}
\P(E) \ge 1 - 2 \exp\Big( -\frac{t^2}{K(c_3, \SIGNALGAUSSIAN) }\Big) - \exp\big( -t^2 \big)\,.
\end{align}
We will prove \equpref{eq:keyconcentration_lowerbd_PE} at the end of the proof of this Lemma. Assuming \equpref{eq:keyconcentration_lowerbd_PE} for now, let us prove \equpref{eq:keyconcentrationgaussianconversion_pt1}.
Let $\cU[L]$ denote the uniform measure on $[L]$.
Consider any $(X_1, \ldots, X_L) \in E$. 
By definition of $E$, we have from the Paley-Zygmund Inequality
\begin{align}
\P_{\ell \sim \cU[L]}\Big( \exp(t X_{\ell}) \ge \frac12 \E_{\ell \sim \cU[L]}\Big[ \exp(t X_{\ell}) \Big] \Big) \ge \frac{\E_{\ell \sim \cU[L]}\Big[ \exp(t X_{\ell}) \Big]^2}{4 \E_{\ell \sim \cU[L]}\Big[ \exp(2t X_{\ell}) \Big]}\,.\label{eq:keyconcentration_paley_zygmund}
\end{align}
Now, we note $\E_{\ell \sim \cU[L]}\Big[ \exp( t X_{\ell}) \Big] = \frac1L \sum_{1 \le \ell \le L} \exp( t X_{\ell})$ for any $  \in \R$. 
By definition of $E$, specifically \equpref{eq:keyconcentration_eventE_1}, we have for $t \ge K_{\upref{lem:keyconcentrationgaussianconversion}}$ that
\begin{align*}
\E_{\ell \sim \cU[L]}\Big[ \exp(t X_{\ell}) \Big] \ge \exp\Big( \frac{t^2}{60 c_3} \Big) \ge 2 \exp\Big( \frac{t^2}{80 c_3} \Big) \ge 4\,.
\end{align*}
Consequently, $(X_1, \ldots, X_L) \in E$ implies the following by \equpref{eq:keyconcentration_paley_zygmund} and \equpref{eq:keyconcentration_eventE_2}:
\begin{align*}
\frac1L \#\Big\{ \ell \le L\,:\,X_{\ell} \ge \frac{t}{80 c_3} \Big\} &= \P_{\ell \sim \cU[L]}\Big( \exp(t X_{\ell}) \ge \exp\Big( \frac{t^2}{80 c_3} \Big) \Big) \\
&\ge \P_{\ell \sim \cU[L]}\Big( \exp(t X_{\ell}) \ge \frac12 \E_{\ell \sim \cU[L]}\Big[ \exp(t X_{\ell}) \Big] \Big) \\
&\ge \exp\big( -(1+K(c_3, \signalsubexp)) t^2 \big)\,.
\end{align*}
Combining the above display with \equpref{eq:keyconcentration_lowerbd_PE}, we obtain
\begin{align*}
&\P\Big( \#\Big\{ \ell \le L\,:\,X_{\ell} \ge \frac{t}{80 c_3} \Big\} \ge L\exp\big( -(1+K(c_3, \signalsubexp)) t^2 \big) \Big) \\
&\quad \ge \P(E) \\
&\quad \ge 1 - 2 \exp\Big( -\frac{t^2}{K(c_3, \SIGNALGAUSSIAN) }\Big) - \exp\big( -t^2 \big) \\
&\quad \ge 1 - 3\exp\Big( -\frac{t^2}{K(c_3, \SIGNALGAUSSIAN)}\Big)\,.
\end{align*}
Thus for $K_{\upref{lem:keyconcentrationgaussianconversion}} = K_{\upref{lem:keyconcentrationgaussianconversion}}(c_3, \signalsubexp, \SIGNALGAUSSIAN)$ large enough, we obtain from the above that
\begin{align*}
&\P\Big( \#\Big\{ \ell \le L\,:\,X_{\ell} \ge \frac{t}{K_{\upref{lem:keyconcentrationgaussianconversion}}} \Big\} \le L\exp\big( -K_{\upref{lem:keyconcentrationgaussianconversion}} t^2 \big) \Big) \\
&\quad \le \P\Big( \#\Big\{ \ell \le L\,:\,X_{\ell} \ge \frac{t}{80 c_3} \Big\} \le L\exp\big( -K_{\upref{lem:keyconcentrationgaussianconversion}} t^2 \big) \Big) \\
&\quad \le \P\Big( \#\Big\{ \ell \le L\,:\,X_{\ell} \ge \frac{t}{80 c_3} \Big\} \le L\exp\big( -(1+K(c_3, \signalsubexp)) t^2 \big) \Big) \\
&\quad \le 3\exp\Big( -\frac{t^2}{K(c_3, \SIGNALGAUSSIAN)}\Big) \\
&\quad \le K_{\upref{lem:keyconcentrationgaussianconversion}} \exp\Big( -\frac{t^2}{K_{\upref{lem:keyconcentrationgaussianconversion}}} \Big)\,,
\end{align*}
proving \equpref{eq:keyconcentrationgaussianconversion_pt1}. 

\paragraph{Proof of \equpref{eq:keyconcentrationgaussianconversion_pt2}.} Note by monotonicity in $\MARGIN$, it suffices to prove the result for $\MARGIN \ge 0$. 
First suppose $\MARGIN \ge 1$. Consider any $t$ such that $K_{\upref{lem:keyconcentrationgaussianconversion}} \MARGIN \le t \le \sqrt{\log L}/K_{\upref{lem:keyconcentrationgaussianconversion}}$. As $\MARGIN \ge 1$, we may apply \equpref{eq:keyconcentrationgaussianconversion_pt1}, which yields
\begin{align*}
\P\Big( \#\Big\{ \ell \le L\,:\,X_{\ell} \ge \MARGIN \Big\} \le L\exp\big(-K_{\upref{lem:keyconcentrationgaussianconversion}}t^2\big) \Big) &\le \P\Big( \#\Big\{ \ell \le L\,:\,X_{\ell} \ge \frac{t}{K_{\upref{lem:keyconcentrationgaussianconversion}}} \Big\} \le L\exp\big(-K_{\upref{lem:keyconcentrationgaussianconversion}}t^2\big) \Big) \\
&\le K_{\upref{lem:keyconcentrationgaussianconversion}} \exp\Big( -\frac{t^2}{K_{\upref{lem:keyconcentrationgaussianconversion}}}\Big)\,.
\end{align*}
Letting $t' = \exp\big(-K_{\upref{lem:keyconcentrationgaussianconversion}}t^2\big)$, we obtain for $L^{-1/K_{\upref{lem:keyconcentrationgaussianconversion}}} \le t' \le \exp\big( -K_{\upref{lem:keyconcentrationgaussianconversion}}^3 \MARGIN^2 \big) \le \exp\big( -K_{\upref{lem:keyconcentrationgaussianconversion}} \MARGIN^2 \big)$,
\begin{align*}
\P\Big( \#\Big\{ \ell \le L\,:\,X_{\ell} \ge \MARGIN \Big\} \le Lt' \Big) \le K_{\upref{lem:keyconcentrationgaussianconversion}} t'^{1/K_{\upref{lem:keyconcentrationgaussianconversion}}^2}
\end{align*}
Now suppose $\MARGIN \in [0,1]$. Consider any $t$ such that $K_{\upref{lem:keyconcentrationgaussianconversion}} \le t \le \sqrt{\log L}/K_{\upref{lem:keyconcentrationgaussianconversion}}$. Since $\MARGIN \le 1 \le \frac{t}{K_{\upref{lem:keyconcentrationgaussianconversion}}}$, we obtain from \equpref{eq:keyconcentrationgaussianconversion_pt1} that
\begin{align*}
\P\Big( \#\Big\{ \ell \le L\,:\,X_{\ell} \ge \MARGIN \Big\} \le L\exp\big(-K_{\upref{lem:keyconcentrationgaussianconversion}}t^2\big) \Big) & \le \P\Big( \#\Big\{ \ell \le L\,:\,X_{\ell} \ge 1 \Big\} \le L\exp\big(-K_{\upref{lem:keyconcentrationgaussianconversion}}t^2\big) \Big) \\
& \le \P\Big( \#\Big\{ \ell \le L\,:\,X_{\ell} \ge \frac{t}{K_{\upref{lem:keyconcentrationgaussianconversion}}} \Big\} \le L\exp\big(-K_{\upref{lem:keyconcentrationgaussianconversion}}t^2\big) \Big) \\
& \le K_{\upref{lem:keyconcentrationgaussianconversion}} \exp\Big( -\frac{t^2}{K_{\upref{lem:keyconcentrationgaussianconversion}}}\Big)\,.
\end{align*}
The result now follows from setting $t' = \exp\big(-K_{\upref{lem:keyconcentrationgaussianconversion}}t^2\big)$ as in the $\MARGIN \ge 1$ case above. As remarked earlier, the result for $\MARGIN < 0$ follows from the $\MARGIN=0$ case. This proves \equpref{eq:keyconcentrationgaussianconversioneq2}.

\paragraph{Proof of \equpref{eq:keyconcentration_lowerbd_PE}.}
We will lower bound the probabilities of each of the events \equpref{eq:keyconcentration_eventE_1}, \equpref{eq:keyconcentration_eventE_2} defining $E$ separately and then take a Union Bound. For brevity, consider the vector $X=(X_1, \ldots, X_L) \in \R^L$ and let 
\begin{align*}
F(X) := \log \sum_{1 \le \ell \le L} \exp(t X_{\ell}) = \log \sum_{1 \le \ell \le L} \exp\Big( t \cdot \frac{\theta_1^{\ell} s + \langle \bar\theta^{\ell}, \bar{g} \rangle}{\sqrt{d}}\Big)\,.
\end{align*}

First we lower bound the probability of \equpref{eq:keyconcentration_eventE_1}. Note $F(X)$ can be viewed a function of $g'=(s, \bar{g})^T$, which has $\min\{1, \SIGNALGAUSSIAN\}$ strongly log-concave law. 
Letting $w_{\ell}(g') \propto \exp\big( t \langle \theta^{\ell}, g' \rangle / \sqrt{d}\big)$ be weights summing to 1 for $\ell, 1 \le \ell \le L$, a direct calculation shows that the gradient $\grad_{g'} F$ of $F$ w.r.t. $g'$ is $t \sum_{\ell \le L} w_{\ell}(g') \btheta^{\ell} / \sqrt{d}$. Also note as $\bar g \sim N(0, I_{d-1})$ independent of $s$,
\begin{align*}
\E[X_{\ell}^2] = \frac{(\theta_1^{\ell})^2 \E[s^2] + \| \bar\theta^{\ell} \|^2}{d} = \frac{\| \theta^{\ell}\|^2 + (\E[s^2]-1) (\theta_1^{\ell})^2}d\,.
\end{align*}
Consequently since $\theta_1^{\ell}/\sqrt{d} \in I$ and $I \subseteq [-2, 2]$,
\begin{align*}
\| \grad_{g'} F\|^2 \le t^2 \max_{\ell} \| \theta^{\ell} \|^2 / d = t^2 \max_{\ell}\Big( \E[X_{\ell}^2] + (1-\E[s^2]) \cdot \frac{(\theta_1^{\ell})^2}{d} \Big) \le t^2 ( c_3 + 4)\,.
\end{align*}
Applying Theorem \upref{thm:logconcaveconcentrationbasics} on the concentration of Lipschitz functions of strongly log-concave measures now gives for all $t'>0$,
\begin{align}
\P\Big( \big| F(X) - \E\big[ F(X) \big] \big| \ge t' \Big) \le 2 \exp\Big( -\frac{t'^2}{K(c_3, \SIGNALGAUSSIAN) t^2 }\Big)\,.\label{eq:basic_stochastic_process_concentrationlogconcave}
\end{align}

We now aim to lower bound $\E\big[ F(X) \big]$. Define
\begin{align*}
u_{\ell} := \frac{\theta_1^{\ell}}{\sqrt{d}}\,,\,\xi_1 := s - \E[s]\,,\, X'_{\ell} := \frac{\theta_1^{\ell} \E[s]}{\sqrt{d}} + \frac{\langle \bar\theta^{\ell}, \bar{g} \rangle}{\sqrt{d}}\,.
\end{align*}
Note that $X_{\ell} = u_{\ell} \xi_1 + X'_{\ell}$, that $\xi_1$ has mean 0, and that $\xi_1, X'_{\ell}$ are independent as the $\theta^{\ell}$ here are fixed and as $\E[s], s$ are independent. Lemma \upref{lem:gaussiancomparisonsoftmax} now gives 
\begin{align}
\E\big[ F(X) \big] &\ge \E\Big[ \log \sum_{1 \le \ell \le L } \exp (t X'_{\ell}) \Big] \,.\label{eq:basic_stochastic_process_bound1}
\end{align}
Recall $\big|\E[s]\big| \le \E[s^2]^{1/2} \le \signalsubexp^{1/2}$, and note $\big| \theta_1^{\ell} / \sqrt{d}\big| \le 2$ for all $\ell$ as $I \subseteq [-2, 2]$. Thus
\begin{align*}
X'_{\ell} \ge \bar X_{\ell} - \big|\E[s]\big| \cdot \max_{1 \le \ell \le L} \big| \theta_1^{\ell} / \sqrt{d}\big| \ge \bar X_{\ell} - 2 \signalsubexp^{1/2} \quad\text{where}\quad \bar X_{\ell} := \frac{1}{\sqrt{d}}\langle \bar\theta^{\ell}, \bar{g} \rangle\,.
\end{align*}
Note $\bar X_{\ell}$ is a centered Gaussian family. Combining the above display with \equpref{eq:basic_stochastic_process_bound1} gives
\begin{align}
\E\big[ F(X) \big] &\ge \E\Big[ \log \sum_{1 \le \ell \le L } \exp (t \bar X_{\ell}) \Big] - 2t \signalsubexp^{1/2} \,.\label{eq:basic_stochastic_process_after_bias}
\end{align}
We now observe that as $\bar g \sim N(0, I_{d-1})$ independent of $s$,
\begin{align*}
\E[X_{\ell}^2] &= \frac{(\theta_1^{\ell})^2 \E[s^2] + \| \bar\theta^{\ell} \|^2}{d} = \E[\bar X_{\ell}^2] + \frac{(\theta_1^{\ell})^2}d \E[s^2]\,, \\
\E[X_{\ell_1} X_{\ell_2} ] &= \frac{\theta_1^{\ell_1} \theta_1^{\ell_2} \E[s^2] + \langle \bar \theta^{\ell_1}, \bar \theta^{\ell_2} \rangle}{d} = \E[\bar X_{\ell_1} \bar X_{\ell_2} ] + \frac{\theta_1^{\ell_1} \theta_1^{\ell_2}}d \E[s^2] \,.
\end{align*}
As $\theta_1^{\ell}/\sqrt{d} \in [x-\eps'', x+\eps'']$, it follows that $\big| \frac{(\theta_1^{\ell})^2}d - x^2 \big|, \big| \frac{\theta_1^{\ell_1} \theta_1^{\ell_2}}d - x^2 \big| \le 3\eps''$ and therefore
\begin{align}
\big| \E[\bar X_{\ell}^2] - \big(\E[ X_{\ell}^2]-x^2 \E[s^2] ) \big|\,,\, \big| \E[\bar X_{\ell_1} \bar X_{\ell_2} ] - \big(\E[ X_{\ell_1} X_{\ell_2} ]-x^2 \E[s^2] \big) \big| \le 3\eps'' \signalsubexp\,.\label{eq:barXvsXvarcovar}
\end{align}
Let 
\begin{align*}
\bar c_1 := c_1 - x^2 \E[s^2] + 3\eps'' \signalsubexp\,,\, \bar c_2 := c_2 - x^2 \E[s^2] - 3\eps'' \signalsubexp\,.
\end{align*}
It follows by \equpref{eq:barXvsXvarcovar} that for $\eps'' \le \frac1{K'_{\upref{lem:keyconcentrationgaussianconversion}}(c_3, \signalsubexp)} \le \frac1{12 \signalsubexp c_3}$, because $c_2-c_1 \ge \frac1{c_3}$, we have
\begin{align}
\bar c_2 - \bar c_1 &\ge (c_2-x^2 \E[s^2]) - (c_1-x^2 \E[s^2]) - 6 \eps'' \signalsubexp \ge c_2 - c_1 - \frac{1}{2c_3} \ge \frac1{2 c_3}\,, \label{eq:bar_c_vs_normal1} \\
\bar c_1 &> c_1 - x^2 \E[s^2] \ge 0\,,\label{eq:bar_c_vs_normal2}
\end{align}
where the last step uses the condition \equpref{eq:keyconcentration_gaussianconversion_ineq}. Additionally by \equpref{eq:barXvsXvarcovar}, we have
\begin{align*}
\E[\bar X_{\ell_1} \bar X_{\ell_2} ] &\le \E[ X_{\ell_1} X_{\ell_2} ] - x^2 \E[s^2] + 3\eps'' \signalsubexp = \bar c_1\,,\\
\E[\bar X_{\ell}^2] &\ge \E[X_{\ell}^2 ] - x^2 \E[s^2] - 3\eps''\signalsubexp = \bar c_2\,.
\end{align*}

We now consider i.i.d. $\QUENCHEDDISORDER, (\ANNEALEDDISORDER_{\ell})_{\ell \le L} \sim N(0,1)$. By \equpref{eq:bar_c_vs_normal1}, \equpref{eq:bar_c_vs_normal2}, we may define the centered Gaussian family $(V_{\ell})_{\ell=1}^L$ by
\begin{align*}
V_{\ell} := \QUENCHEDDISORDER \sqrt{\bar c_1} + \ANNEALEDDISORDER_{\ell}\sqrt{\bar c_2 - \bar c_1}\,.
\end{align*}
Thus, $\E[V_{\ell}^2] = \bar c_2 \le \E[\bar X_{\ell}^2]$ and $\E[V_{\ell_1} V_{\ell_2}] = \bar c_1 \ge \E[\bar X_{\ell_1} \bar X_{\ell_2}]$ for $\ell_1 \neq \ell_2$. 
Consequently as $\bar X_{\ell}$ is a centered Gaussian family, the conditions of Lemma \upref{lem:gaussiansoftmaxcomparetoindependent} apply. Combining it with \equpref{eq:basic_stochastic_process_after_bias} gives
\begin{align*}
\E\big[ F(X) \big] &\ge \E\Big[ \log \sum_{1 \le \ell \le L } \exp (t \bar X_{\ell}) \Big] - 2t \signalsubexp^{1/2}  \\ 
&\ge \E\Big[ \log \sum_{1 \le \ell \le L} \exp(tV_{\ell}) \Big] - 2t \signalsubexp^{1/2} \\ 
&= \E\Big[ \log \sum_{1 \le \ell \le L} \exp(t\sqrt{\bar c_2 - \bar c_1}\ANNEALEDDISORDER_{\ell}) \Big] - 2t \signalsubexp^{1/2} \,. 
\end{align*}
Take $K_{\upref{lem:keyconcentrationgaussianconversion}} = K_{\upref{lem:keyconcentrationgaussianconversion}}(c_3, \signalsubexp, \SIGNALGAUSSIAN) \ge K_{\upref{lem:lowerboundsoftmaxindependent}} \sqrt{2c_3} \ge \frac{K_{\upref{lem:lowerboundsoftmaxindependent}}}{\sqrt{\bar c_2- \bar c_1}}$, where the last inequality follows from \equpref{eq:bar_c_vs_normal1}.
Thus we may apply Lemma \upref{lem:lowerboundsoftmaxindependent} with $t\sqrt{\bar c_2 - \bar c_1}$ in place of $t$, which gives that for $t \ge K_{\upref{lem:keyconcentrationgaussianconversion}} = K_{\upref{lem:keyconcentrationgaussianconversion}}(c_3, \signalsubexp, \SIGNALGAUSSIAN)$,
\begin{align*}
\E\big[ F(X) \big] &\ge \E\Big[ \log \sum_{1 \le \ell \le L} \exp(t\sqrt{\bar c_2 - \bar c_1}\ANNEALEDDISORDER_{\ell}) \Big] - 2t \signalsubexp^{1/2}  \\ 
&\ge \log L + \frac{(\bar c_2 - \bar c_1)t^2}5 - 2t \signalsubexp^{1/2} \\ 
&\ge \log L + \frac{t^2}{20 c_3} - 2t \signalsubexp^{1/2} \\
&\ge \log L + \frac{t^2}{30 c_3}\,,
\end{align*}
where we use \equpref{eq:bar_c_vs_normal1} to lower bound $\bar c_2 - \bar c_1$. 
Combining with \equpref{eq:basic_stochastic_process_concentrationlogconcave} gives
\begin{equation}\label{eq:keyconcentrationgaussianconversioneq1}
\begin{aligned}
\P\Big( \sum_{1\le \ell \le L} \exp(t X_{\ell}) \ge L \exp\Big( \frac{t^2}{60 c_3} \Big) \Big) &\ge 1 - \P\Big( F(X) \le \E\big[ F(X) \big] - \frac{t^2}{60 c_3} \Big) \\
&\ge 1 - 2 \exp\Big( -\frac{t^2}{K(c_3, \SIGNALGAUSSIAN) }\Big)\,.
\end{aligned}
\end{equation}
This lower bounds the probability of \equpref{eq:keyconcentration_eventE_1}.

Next, we lower bound the probability of \equpref{eq:keyconcentration_eventE_2} as follows. As each $X_{\ell} = \frac{\theta^{\ell}_1 s + \langle \bar\theta^{\ell}, \bar{g} \rangle}{\sqrt{d}}$ has identical law and as $s, \bar{g}$ are independent,
\begin{align}\label{eq:keyconcentration_2ndmoment_1}
\E\Big[ \sum_{1 \le \ell \le L} \exp(2t X_{\ell})\Big] = L \E\big[ \exp(2tX_1) \big] = L\E\Big[ \exp\Big(\frac{2t\theta^1_1}{\sqrt{d} s} \Big) \Big] \E\Big[ \exp\Big( \frac{2t \langle \bar\theta^1, \bar{g} \rangle}{\sqrt{d}} \Big)\Big]\,.
\end{align}
Standard calculations on the MGF of a Gaussian yield
\begin{align}\label{eq:keyconcentration_2ndmoment_2}
\E\Big[ \exp\big( \frac{2t \langle \bar\theta^1, \bar{g} \rangle}{\sqrt{d}} \big)\Big] = \exp\Big( \frac{2t^2 \|\bar\btheta^1\|^2}d \Big)\,.    
\end{align}
Next we have $\frac{2t\theta^1_1}{\sqrt{d}} s \le \frac{t^2 (\theta^1_1)^2 \signalsubexp}d + \frac{s^2}{\signalsubexp}$, therefore as $s^2$ is $\signalsubexp$-sub-Exponential,
\begin{align}\label{eq:keyconcentration_2ndmoment_3}
\E\Big[ \exp\Big(\frac{2t\theta^1_1 s}{\sqrt{d}} \Big) \Big] \le \exp\Big( \frac{t^2 (\theta^1_1)^2 \signalsubexp}d\Big) \E\Big[ \exp\Big( \frac{s^2}{\signalsubexp} \Big)\Big] \le 2\exp\Big( \frac{t^2 (\theta^1_1)^2 \signalsubexp}d\Big)\,.
\end{align}
Since $c_3 \ge \E\big[ X_{\ell}^2 \big] = \frac{(\theta^{\ell}_1)^2 \E[s^2] + \|\bar\theta^{\ell}\|^2}d$ and as $\theta_1^1 / \sqrt{d} \in I \subseteq [-2,2]$, it follows that
\begin{align*}
\frac1d \|\theta^1\|^2 \le c_3 + \frac{(\theta_1^1)^2}d \le c_3 + 4\,.
\end{align*}
Combining the above bound with \equpref{eq:keyconcentration_2ndmoment_1}, \equpref{eq:keyconcentration_2ndmoment_2}, \equpref{eq:keyconcentration_2ndmoment_3} gives
\begin{align*}
\E\Big[ \sum_{1 \le \ell \le L} \exp(2t X_{\ell})\Big] \le 2L\exp\Big( \frac{t^2 K(\signalsubexp) \| \theta^1 \|^2}{d} \Big) \le L \exp\big(  K(c_3, \signalsubexp) t^2 \big)\,,
\end{align*}
where we used that $t \ge K_{\upref{lem:keyconcentrationgaussianconversion}} = K_{\upref{lem:keyconcentrationgaussianconversion}}(c_3, \signalsubexp, \SIGNALGAUSSIAN)$. 
Thus by Markov's Inequality,
\begin{align}
\P\Big( \sum_{1 \le \ell \le L} \exp(2t X_{\ell}) \le L \exp\big( (1+K(c_3, \signalsubexp) ) t^2 \big) \Big) \ge 1 - \exp( -t^2 )\,.\label{eq:keyconcentrationgaussianconversioneq2}
\end{align}
This lower bounds the probability of \equpref{eq:keyconcentration_eventE_2}. Combining \equpref{eq:keyconcentrationgaussianconversioneq1}, \equpref{eq:keyconcentrationgaussianconversioneq2} proves \equpref{eq:keyconcentration_lowerbd_PE}, thus completing the proof of Lemma \upref{lem:keyconcentrationgaussianconversion}.
\end{proof}

We also cite the following fact from \cite{talagrand2011advanced}. Technically, in \cite{talagrand2011advanced}, Theorem \upref{thm:gaussiansetuptalagrand} is stated with $\theta$ in place of $\bar \theta$ in the conclusion \equpref{eq:gibbsmeasurenice1_prelim}, \equpref{eq:gibbsmeasurenice2_prelim}, however the proof is identical.
\begin{theorem}[Corollary of the proof of Theorem 8.2.7, \cite{talagrand2011advanced}]\label{thm:gaussiansetuptalagrand}
Consider a concave $U \le 0$ defined on $\R^d$, and a convex set $C \subseteq \R^d$. For $\GAUSSIANMEASURE \in [\GAUSSIANMEASURE_0, \GAUSSIANMEASURE_1]$, define a measure $G_C$ on $\R^d$ as follows, where $Z'$ denotes the appropriate normalization constant:
\[ \text{ for all } B \subseteq \R^d\,,\, G_C(B) = \frac1{Z'} \int \one\{B \cap C\} \exp \Big(U(\btheta)-\GAUSSIANMEASURE \| \btheta \|^2 \Big) \rmd \btheta\, .\]
Assume that for some $a>0$ we have 
\[ Z' \ge \exp(-ad)\, .\]
Then for some $0 \le c'_1 < c'_2 < c'_3$ with $c'_3 \ge \frac1{c'_2-c'_1}$, where $c'_3$ only depends on $a, \beta_0, \beta_1$, we have that $G_C$ satisfies:
\begin{align}
G_C\Bigl(\big\{\btheta : c'_2 d \le \|\bar \btheta\|^2 \le c'_3 d\big\}\Bigr) &\ge 1 - \exp\Bigl(-\frac{d}{c'_3}\Bigr)\, , \label{eq:gibbsmeasurenice1_prelim}\\
G_C^{\otimes 2}\Bigl(\big\{(\btheta^1,\btheta^2) : | \langle \bar \btheta^1, \bar \btheta^2 \rangle | \le c'_1 d\big\}\Bigr)
&\ge 1 - \exp\Bigl(-\frac{d}{c'_3}\Bigr)\, .\label{eq:gibbsmeasurenice2_prelim}
\end{align}
\end{theorem}
\begin{proof}
This follows from the exact same proof of Theorem 8.2.7, \cite{talagrand2011advanced}, except we now consider $R_{1,1} = \frac1d \| \bar \theta \|^2, R_{1,2} = \frac1d \langle \bar \theta_1, \bar \theta_2 \rangle$ as we do throughout this paper, as per \equpref{eq:overlap_defs_thispaper}. Following the same proof as therein, it is shown that \equpref{eq:gibbsmeasurenice1_prelim}, \equpref{eq:gibbsmeasurenice2_prelim} hold where we can take $c'_3 = \frac4{d'^2}$ where $d' = \exp\big( -2(a+a^2/4\beta+3)\big)$. Note now rather than following the proof of (3.26) and Theorem 3.1.11 in \cite{talagrand2010mean}, we follow the proof of Lemma \upref{lem:logconcaveboundmomentslemma} and Proposition \upref{prop:concentrationoverlapfixedgibbs}, directly using the bound $Z' \ge \exp(-ad)$ rather than Lemma \upref{lem:lem316}. (There is no dependence on $D$ here as dependence on $D$ in the proof of Lemma \upref{lem:logconcaveboundmomentslemma}, Proposition \upref{prop:concentrationoverlapfixedgibbs} only came due to Lemma \upref{lem:lem316}.) The rest of the proof is identical to the proof of Theorem 8.2.7 in \cite{talagrand2011advanced}. Note that $c'_1 \ge 0$ as $\langle R_{1,2} \rangle \ge 0$.
\end{proof}

Now we prove the following very useful result that upper bounds the probability that a Gibbs measure places small mass on one side of a random hyperplane in $\R^d$, where the hyperplane has law identical to the constraints (the probability is over the hyperplane). We will use this result repeatedly next in Appendix \upref{subsec:gaussian_measure_rs_pf}.
\begin{theorem}\label{thm:gaussiansetup}
Consider a concave $U \le 0$ defined on $\R^d$, and a convex set $C \subseteq \R^d$. Consider any $x \in [-1,1]$ and let $I = [x-\eps'', x+\eps'']$ for some $\eps''>0$. For $\GAUSSIANMEASURE \in [\GAUSSIANMEASURE_0, \GAUSSIANMEASURE_1]$, define a measure $G_C$ on $\R^d$ as follows, where $Z'$ denotes the appropriate normalization constant:
\[ \text{ for all } B \subseteq \R^d\,,\, G_C(B) = \frac1{Z'} \int \one\big\{ B \cap C \cap \{ \theta_1/\sqrt{d} \in I\} \big\} \exp \Big(U(\btheta)-\GAUSSIANMEASURE \| \btheta \|^2 \Big) \rmd \btheta\, .\]
Assume that for some $a>0$ we have 
\[ Z' \ge \exp(-ad)\, .\]
Then for $K_{\upref{thm:gaussiansetup}} \ge 2$ depending on $a, \beta_0, \beta_1, \signalsubexp, \SIGNALGAUSSIAN$, $K'_{\upref{thm:gaussiansetup}} \ge 2$ depending on $a, \beta_0, \beta_1, \signalsubexp$, and $K_{\upref{lem:keyconcentrationgaussianconversion}} \ge 1$ depending on $a, \beta_0, \beta_1, \signalsubexp, \SIGNALGAUSSIAN$, we have the following.
Consider $s \sim \cD$ and $\bar g \sim N(0, I_{d-1})$ such that $s, \bar g$ are independent.
If $\eps'' \le \frac1{K'_{\upref{thm:gaussiansetup}}}$, then $G_C$ satisfies the following implication: 
\begin{equation}\label{eq:hyperplanesizekeyprop}
\begin{aligned}
&\text{If } \exp\Big( -\frac{d}{K_{\upref{thm:gaussiansetup}}}\Big) \le t \le \frac14 \exp\Big( -K_{\upref{lem:keyconcentrationgaussianconversion}} \max\{1,\MARGIN\}^2 \Big)\,, \\
&\text{ then }\P\Bigg( G_C\Big( \Big\{ \theta\,:\, \frac{\theta_1 s + \langle \bar\theta, \bar{g} \rangle}{\sqrt{d}} \ge \MARGIN \Big\}\Big) \le t \Bigg) \le 16 K_{\upref{lem:keyconcentrationgaussianconversion}}  t^{1/K^2_{\upref{lem:keyconcentrationgaussianconversion}}}\,.
\end{aligned}
\end{equation}
\end{theorem}
\begin{proof}
We define for all $\theta$,
\begin{align*}
X(\theta) := \frac{s \theta_1 + \langle \bar\theta, \bar{g} \rangle}{\sqrt{d}}\,,
\end{align*}
and we let $X_{\ell} := X(\theta^{\ell})$ for all replicas $1 \le \ell \le L$. 

First, we prove that the conditions of Lemma \upref{lem:keyconcentrationgaussianconversion} apply with high enough probability over $G_C$. By Theorem \upref{thm:gaussiansetuptalagrand} applied with the convex set $\{\theta_1/\sqrt{d}\in I\} \cap C$, 
for some $0 \le c'_1 < c'_2 < c'_3 = c'_3(a, \beta_0, \beta_1)$ with $c'_3 \ge \frac{1}{c'_2-c'_1}$, we have \equpref{eq:gibbsmeasurenice1_prelim}, \equpref{eq:gibbsmeasurenice2_prelim}. (Crucially $c'_3$ does not explicitly depend on $I$, in particular it does not depend on $\eps''$, which will be very important in the following arguments.) Since $G_C\big(\{ \theta\,:\,\theta_1/\sqrt{d} \not\in I \}\big)=0$ (recall the constraint $\theta_1/\sqrt{d} \in I$ is present in the definition of $G_C$ from this Theorem), it follows that 
\begin{align}
G_C\Bigl(\big\{\btheta : c'_2 d \le \|\bar \btheta\|^2 \le c'_3 d\,,\,\theta_1/\sqrt{d} \in I \big\}\Bigr) &\ge 1 - \exp\Bigl(-\frac{d}{c'_3}\Bigr)\, , \label{eq:gibbs_measure_nice_conditions1} \\
G_C^{\otimes 2}\Bigl(\big\{(\btheta^1,\btheta^2) : | \langle \bar \btheta^1, \bar \btheta^2 \rangle | \le c'_1 d\,,\,\theta_1/\sqrt{d} \in I \big\}\Bigr)
&\ge 1 - \exp\Bigl(-\frac{d}{c'_3}\Bigr)\,. \label{eq:gibbs_measure_nice_conditions2}
\end{align}
Now note for any $\theta^{\ell}, \theta^{\ell_1}, \theta^{\ell_2}$, we have as $\bar g \sim N(0, I_{d-1})$ independent of $s$,
\begin{align*}
\E\big[X^2_{\ell}\big] = \frac1d \Big( (\theta_1^{\ell})^2 \E[s^2] + \| \bar\theta^{\ell} \|^2 \Big)\quad,\quad
\E\big[X_{\ell_1} X_{\ell_2}\big] = \frac1d \Big( \theta_1^{\ell_1} \theta_1^{\ell_2} \E[s^2] + \langle \bar\theta^{\ell_1}, \bar\theta^{\ell_2} \rangle \Big)\,. 
\end{align*}
Consequently for $\theta^{\ell}$ such that $c'_2 d \le \|\bar \btheta^{\ell}\|^2 \le c'_3 d\,,\,(\theta^{\ell})_1/\sqrt{d} \in I$, we have 
\begin{align*}
\E\big[X^2_{\ell}\big] &\le c'_3 + \max\big\{ |x + \eps''|, |x - \eps''| \big\}^2 \signalsubexp) \le c'_3 + 4\signalsubexp\,, \\
\E\big[X^2_{\ell}\big] &\ge c'_2 + (x^2 - 3\eps'') \E[s^2]\,.
\end{align*}
Similarly for $(\theta^{\ell_1}, \theta^{\ell_2})$ such that $| \langle \btheta^{\ell_1}, \btheta^{\ell_2} \rangle | \le c'_1 d\,,\,(\theta^{\ell_1})_1/\sqrt{d} \in I\,,\,(\theta^{\ell_2})_1/\sqrt{d} \in I$, we have
\begin{align*}
\E\big[ X_{\ell_1} X_{\ell_2} \big] &\le c'_1 + (x^2+3\eps'') \E[s^2]\,.
\end{align*}
Define
\begin{align}\label{eq:gibbs_measure_nice_choices}
c_3 = 2c'_3 + 4\signalsubexp\,,\, c_1 := c'_1 + x^2 \E[s^2] + \frac1{4c'_3}\,,\,c_2 := c'_2 + x^2 \E[s^2] - \frac1{4c'_3}\,.
\end{align}
Let 
\begin{align}
K'_{\upref{thm:gaussiansetup}}(a, \beta_0, \beta_1, \signalsubexp) = K'_{\upref{lem:keyconcentrationgaussianconversion}}(c_3, \signalsubexp) \ge 12 c_3 \signalsubexp\,.
\end{align}
Hence since $\eps'' \le \frac1{K'_{\upref{thm:gaussiansetup}}(a, \beta_0, \beta_1, \signalsubexp)} \le \frac1{12 c_3 \signalsubexp}$, it follows that for such $\theta$ satisfying the conditions from \equpref{eq:gibbs_measure_nice_conditions1} we have $c_3 \ge \E[X^2_{\ell}] \ge c_2$, and similarly for $\theta^1, \theta^2$ satisfying the conditions from \equpref{eq:gibbs_measure_nice_conditions2} we have $\E[X_{\ell_1}X_{\ell_2}] \le c_1$. 
As $c_3 \ge c'_3$, and as $G_C\big(\{ \theta\,:\,\theta_1/\sqrt{d} \not\in I \}\big)=0$, it follows that
\begin{align}
G_C\Bigl(\big\{\btheta : c_2 \le \E[X^2_{\ell}] \le c_3 \big\}\Bigr) &\ge 1 - \exp\Bigl(-\frac{d}{c_3}\Bigr)\,,\label{eq:gibbsmeasurenice1}\\
G_C^{\otimes 2}\Bigl(\big\{(\btheta^1,\btheta^2) : \E[X_{\ell_1} X_{\ell_2}] \le c_1 \big\}\Bigr)
&\ge 1 - \exp\Bigl(-\frac{d}{c_3}\Bigr)\,.\label{eq:gibbsmeasurenice2}
\end{align}

Consider $L$ such that $2L^2 \le \exp(d/c_3)$ and define 
\begin{align}
E^L &:= \Big\{ (\theta^1, \ldots, \theta^L)\,:\, c_2 \le \E[X^2_{\ell}] \le c_3\, \forall \, \ell\,,\, \E[X_{\ell_1}X_{\ell_2}] \le c_1\,\forall\, \ell_1 \neq \ell_2\,,\, \theta^{\ell}/\sqrt{d} \in I\, \forall\, \ell\Big\}\,.\label{eq:stochastic_process_ELdef}
\end{align}
Now observe that the conditions of Lemma \upref{lem:keyconcentrationgaussianconversion} apply with $c_1, c_2, c_3$ chosen in \equpref{eq:gibbs_measure_nice_choices} because we have $\eps'' \le \frac1{K'_{\upref{lem:keyconcentrationgaussianconversion}}(c_3, \signalsubexp)}$, and because $c'_1 \ge 0$ and therefore
\begin{align*}
c_2-c_1 = c'_2-c'_1 - \frac1{2c'_3} \ge \frac1{2c'_3} \ge \frac1{c_3}\quad,\quad c_1 > c'_1 + x^2 \E[s^2] \ge x^2 \E[s^2]\,.
\end{align*}
Applying Lemma \upref{lem:keyconcentrationgaussianconversion}, we obtain the following. For any $(\theta^1, \ldots, \theta^L) \in E^L$ and $t$ satisfying
\begin{align}\label{eq:gaussiansetup_condition_init}
L^{-1/K_{\upref{lem:keyconcentrationgaussianconversion}}} \le t \le \exp\big( -K_{\upref{lem:keyconcentrationgaussianconversion}} \max\{1,\MARGIN\}^2 \big)\,,
\end{align}
where $K_{\upref{lem:keyconcentrationgaussianconversion}} = K_{\upref{lem:keyconcentrationgaussianconversion}}(c_3, \signalsubexp, \SIGNALGAUSSIAN) \ge 1$ comes from Lemma \upref{lem:keyconcentrationgaussianconversion} with $c_3$ defined in \equpref{eq:gibbs_measure_nice_choices}, we have 
\begin{align}\label{eq:keyconcentration_consequence}
\P\Big( \#\big\{ \ell \le L\,:\,X_{\ell} \ge \MARGIN\big\} \le t L\Big) \le K_{\upref{lem:keyconcentrationgaussianconversion}}\, t^{1/K_{\upref{lem:keyconcentrationgaussianconversion}}^2}\,,
\end{align}
where probability in the above is over $\bar g, s$.

Consider the event $G_C^{\otimes L}\Big( (\theta^1, \ldots, \theta^L) \in E^L\,:\, \#\big\{ \ell \ge L\,:\, X_{\ell} \ge \MARGIN \big\} \ge t L \Big)$. Observe that when $G_C^{\otimes L}\Big( (\theta^1, \ldots, \theta^L) \in E^L\,:\, \#\big\{ \ell \le L\,:\, X_{\ell} \ge \MARGIN \big\} \ge t L \Big) \ge \frac14$, we have by Linearity of Expectation,
\begin{align}\label{eq:keyconcentration_implication1}
G_C\big( \big\{ \theta\,:\, X(\theta) \ge \MARGIN \big\}\big) &= \frac1L \int \#\big\{ \ell \le L\,:\, X(\theta^{\ell}) \ge \MARGIN \big\}\, \rmd G^{\otimes L}_C(\theta) \notag \\
&\ge \frac1L \cdot tL \cdot G_C^{\otimes L}\Big( (\theta^1, \ldots, \theta^L) \in E^L\,:\, \#\big\{ \ell \le L\,:\, X_{\ell} \ge \MARGIN \big\} \ge t L \Big) \notag \\
&\ge \frac{t}4\,.
\end{align}
Next observe that 
\begin{align}\label{eq:keyconcentration_implication2}
&G_C^{\otimes L}\Big( (\theta^1, \ldots, \theta^L) \in E^L\,:\, \#\big\{ \ell \le L\,:\, X'_{\ell} \ge \MARGIN \big\} \ge t L \Big) \notag \\
&\quad = G_C^{\otimes L}(E^L) - G_C^{\otimes L}\Big( (\theta^1, \ldots, \theta^L) \in E^L\,:\, \#\big\{ \ell \le L\,:\, X'_{\ell} \ge \MARGIN \big\} \le t L \Big) \notag \\
&\quad \ge 1 - L^2 \exp\Big(-\frac{d}{c_3}\Big) - G_C^{\otimes L}\Big( (\theta^1, \ldots, \theta^L) \in E^L\,:\, \#\big\{ \ell \le L\,:\, X'_{\ell} \ge \MARGIN \big\} \le t L \Big) \notag \\
&\quad \ge \frac12 - G_C^{\otimes L}\Big( (\theta^1, \ldots, \theta^L) \in E^L\,:\, \#\big\{ \ell \le L\,:\, X'_{\ell} \ge \MARGIN \big\} \le t L \Big)\,, 
\end{align}
where we apply the condition on $L$ and use \equpref{eq:gibbsmeasurenice1}, \equpref{eq:gibbsmeasurenice2} and that $G_C\big( \{\theta:\theta_1/\sqrt{d}\not\in I\}\big)=0$.

Combining \equpref{eq:keyconcentration_implication1}, \equpref{eq:keyconcentration_implication2}, it follows that as events we have
\begin{align}
&\Bigg\{ G_C^{\otimes L}\Big( (\theta^1, \ldots, \theta^L) \in E^L\,:\, \#\big\{ \ell \le L\,:\, X'_{\ell} \ge \MARGIN \big\} \le t L \Big) \le \frac14 \Bigg\} \notag \\
&\quad \subseteq \Bigg\{ G_C^{\otimes L}\Big( (\theta^1, \ldots, \theta^L) \in E^L\,:\, \#\big\{ \ell \le L\,:\, X'_{\ell} \ge \MARGIN \big\} \ge t L \Big) \ge \frac14 \Bigg\} \notag \\
&\quad \subseteq \Big\{ G_C\big( \big\{ \theta\,:\, X(\theta) \ge \MARGIN \big\}\big) \ge \frac{t}4\Big\}\,.\label{eq:keyconcentration_implication3}
\end{align}
Now by Markov's Inequality and \equpref{eq:keyconcentration_consequence}, we obtain
\begin{align*}
&\P\Bigg( G_C^{\otimes L}\Big( (\theta^1, \ldots, \theta^L) \in E^L\,:\, \#\big\{ \ell \le L\,:\, X'_{\ell} \ge \MARGIN \big\} \le t L \Big) \ge \frac14 \Bigg) \\
&\quad \le 4\E\Bigg[ G_C^{\otimes L}\Big( (\theta^1, \ldots, \theta^L) \in E^L\,:\, \#\big\{ \ell \le L\,:\, X'_{\ell} \ge \MARGIN \big\} \le t L \Big) \Bigg] \\
&\quad = 4 \int \one\big\{(\theta^1, \ldots, \theta^L) \in E^L \big\}\, \P\Big( \one\Big\{ \#\big\{ \ell \le L\,:\, X'_{\ell} \ge \MARGIN \big\} \le t L\Big\} \Big)\, \rmd G_C^{\otimes L}(\theta) \\
&\quad \le 4K_{\upref{lem:keyconcentrationgaussianconversion}} t^{1/K^2_{\upref{lem:keyconcentrationgaussianconversion}}}\,.
\end{align*}
Therefore by the contrapositive of \equpref{eq:keyconcentration_implication3}, combining with the above display yields
\begin{align}
\P\Big( G_C\big( \big\{ \theta\,:\, X(\theta) \ge \MARGIN \big\}\big) \le \frac{t}4 \Big) &\le \P\Bigg( G_C^{\otimes L}\Big( (\theta^1, \ldots, \theta^L) \in E^L\,:\, \#\big\{ \ell \le L\,:\, X'_{\ell} \ge \MARGIN \big\} \le t L \Big) \ge \frac14 \Bigg) \notag \\
&\le 4K_{\upref{lem:keyconcentrationgaussianconversion}} t^{1/K^2_{\upref{lem:keyconcentrationgaussianconversion}}}\,.\label{eq:keyconcentration_implication4}
\end{align}
Hence, \equpref{eq:gaussiansetup_condition_init} implies \equpref{eq:keyconcentration_implication4}. The result follows by taking $4t$ in place of $t$ in this implication, noting $4^{1/K^2_{\upref{lem:keyconcentrationgaussianconversion}}} \le 4$ as $K_{\upref{lem:keyconcentrationgaussianconversion}} \ge 1$, and recalling this holds for any $L$ with $2L^2 \le \exp(d/c_3)$.
\end{proof}

Now we find a generic way to upper bound random variables satisfying the implication in Theorem \upref{thm:gaussiansetup}. This is analogous to Lemma 8.3.8 of \cite{talagrand2011advanced}, but we require some more generality.
\begin{lemma}\label{lem:momentsofRV838tal}
Consider a random variable $V \ge 0$ and assume that for certain numbers $C_0, C_1, C_2, D_0 \ge 1$ and some $a \in \R$ we have 
\begin{align*}
&\P(V \ge \exp(ad)) = 0\, ,\\
&D_0 \le t \le \exp(d/C_0) \implies \P(V \ge t) \le C_1 t^{-1/C_2}\, .
\end{align*}
Then we have
\begin{align*}
\E [ V^{\EXPONENTEXPSCALE} ] \le K_{\upref{lem:momentsofRV838tal}} = D_0+2C_1 C_2+C_1 \quad\text{ for }\quad \EXPONENTEXPSCALE = \max\Big\{1, \frac1{2C_2}, \frac1{aC_0 C_2}\Big\} \in (0, 1)\,.
\end{align*}
\end{lemma}
\begin{proof}
As $\P(V \ge \exp(ad))=0$, we may write
\begin{align*}
\E[ V^{\EXPONENTEXPSCALE}] &= \EXPONENTEXPSCALE \int_0^{\infty} t^{\EXPONENTEXPSCALE-1} \P(V \ge t)\,\rmd t \\
&= \EXPONENTEXPSCALE \int_0^{D_0} t^{\EXPONENTEXPSCALE-1} \P(V \ge t)\,\rmd t + \EXPONENTEXPSCALE \int_{D_0}^{\exp(d/C_0)} t^{\EXPONENTEXPSCALE-1} \P(V \ge t)\,\rmd t + \EXPONENTEXPSCALE \int_{\exp(d/C_0)}^{\exp(ad)} t^{\EXPONENTEXPSCALE-1} \P(V \ge t)\,\rmd t\,.
\end{align*}
We upper bound these three integrals as follows:
\begin{itemize}
\item First, we have that
\begin{align*}
\EXPONENTEXPSCALE \int_0^{D_0} t^{\EXPONENTEXPSCALE-1} \P(V \ge t)\,\rmd t \le t^{\EXPONENTEXPSCALE} \Big|_0^{D_0} = D_0^{\EXPONENTEXPSCALE} \le D_0\,.
\end{align*}

\item Second, as $\EXPONENTEXPSCALE \le \frac1{2C_2} \le 1$ we have by the upper bound on $\P(V \ge t)$ that
\begin{align*}
\EXPONENTEXPSCALE \int_{D_0}^{\exp(d/C_0)} t^{\EXPONENTEXPSCALE-1} \P(V \ge t)\,\rmd t &\le \EXPONENTEXPSCALE C_1 \int_{D_0}^{\exp(d/C_0)} t^{\EXPONENTEXPSCALE - 1 - 1/C_2}\,\rmd t \\
&\le C_1 \cdot 2C_2 D_0^{-1/2C_2} \\
&\le 2C_1 C_2\,,
\end{align*}
where we used that $D_0 \ge 1$. 

\item Third, we have by the upper bound on $\P(V \ge t)$ that
\begin{align*}
\EXPONENTEXPSCALE \int_{\exp(d/C_0)}^{\exp(ad)} t^{\EXPONENTEXPSCALE-1} \P(V \ge t)\,\rmd t &\le \EXPONENTEXPSCALE C_1 \int_{\exp(d/C_0)}^{\exp(ad)} t^{\EXPONENTEXPSCALE - 1 - 1/C_2}\,\rmd t \\
&\le \EXPONENTEXPSCALE C_1 \exp\Big(-\frac{d}{C_0 C_2} \Big) \int_{\exp(d/C_0)}^{\exp(ad)} t^{\EXPONENTEXPSCALE-1}\,\rmd t \\
&\le \frac{\EXPONENTEXPSCALE C_1 \exp\big(-\frac{d}{C_0 C_2} \big)}{\EXPONENTEXPSCALE} \cdot \exp(a \EXPONENTEXPSCALE d) \le C_1\,,
\end{align*}
as $t \ge \exp(d/C_0)$ and $\EXPONENTEXPSCALE \le \frac1{aC_0 C_2}$.
\end{itemize}
The Lemma is proved.
\end{proof}

We will also need the following useful fact:
\begin{lemma}\label{lem:upperboundintegral}
Consider any $U \le 0$ defined on $\R^d$, an interval $I \subseteq \R$, and a convex set $C \subseteq \R^d$. For $\GAUSSIANMEASURE \in [\GAUSSIANMEASURE_0, \GAUSSIANMEASURE_1]$, define a measure $G_C$ on $\R^d$ as follows, where $Z'$ denotes the appropriate normalization constant:
\[ \text{ for all } B \subseteq \R^d\,,\, G_C(B) = \frac1{Z'} \int_{B \cap C} \one\{ \theta_1/\sqrt{d} \in I \}\, \exp \Big(U(\btheta)-\GAUSSIANMEASURE \| \btheta \|^2 \Big) \rmd \btheta\, .\]
Assume that for some $a>0$ we have 
\[ Z' \ge \exp(-ad)\, .\]
Then for $r \le r_{\upref{lem:upperboundintegral}}$ where $r_{\upref{lem:upperboundintegral}}$ depends on $a, \GAUSSIANMEASURE_0, \GAUSSIANMEASURE_1$, we have 
\[ G_C\Big( \big\{ \| \bar{\btheta} \| \le r \sqrt{d} \big\} \Big) \le \exp(-d)\, .\]
\end{lemma}
\begin{proof}
By the condition on $Z'$ and as $U \le 0$, we have 
\begin{align*}
G_C\big( \big\{ \| \bar{\btheta} \| \le r \sqrt{d} \big\} \big) &\le \exp(ad) \int_{\| \bar{\btheta} \| \le r \sqrt{d}} \one\{ \theta_1/\sqrt{d} \in I \}\, \exp\big(-\GAUSSIANMEASURE \| \btheta \|^2 \big)\, \rmd \btheta \\
&= \exp(ad) \Big(\frac{\beta}{\pi}\Big)^{d/2} \int_\R \one\{ \theta_1/\sqrt{d} \in I \}\,\exp(-\GAUSSIANMEASURE \theta_1^2)\,\rmd\btheta \cdot \int_{\| \bar{\btheta} \| \le r \sqrt{d}} \exp\big(-\GAUSSIANMEASURE \| \bar\btheta \|^2 \big)\, \rmd \bar\btheta \\
&\le \exp\big(d(a+K(\beta) )\big) \P_{Z}\big(\| Z \| \le r\sqrt{d}\big)\,,
\end{align*}
where $Z \sim N\big(0, \frac1{\beta} I_{d-1} \big)$. Standard bounds on the concentration of the norm of multivariate Guassians (see e.g. Theorem 3.1.1, \cite{vershynin2018high}) now give that for suitable $r_{\upref{lem:upperboundintegral}}(a, \beta_0, \beta_1)$ and for all $r \le r_{\upref{lem:upperboundintegral}}(a, \beta_0, \beta_1)$, we have $\P_{Z}\big(\| Z \| \le r\sqrt{d}\big) \le \exp\big(-d(a+K(\beta)+1)\big)$. This yields the desired conclusion.
\end{proof}

\subsubsection{Proof of Theorem \upref{thm:gaussianmeasurersformula}}\label{subsec:gaussian_measure_rs_pf}
We now have established the technical ingredients necessary to prove Theorem \upref{thm:gaussianmeasurersformula}. As discussed earlier, the proof proceeds through two key Lemmas. In Lemma \upref{lem:controlmollifierbernsteinspiked} we show concentration of the normalized free energy around its expectation for quite general mollified Hamiltonians. This generality is needed as we apply this result for both $u(\cdot)$ from Theorem \upref{thm:mollifiedhamiltonianfreeenergy} and with $\exp u(x)=\one\{x \ge \MARGIN\}$. In Lemma \upref{lem:controlmollifiers} we show the impact of the mollifier on the expectation is small.

Since these arguments will involve changing the mollifier $u(\cdot)$, overloading notation, we make the following definitions for only Appendix \upref{subsec:gaussianconversioninterpolants}.
For any function $u(\cdot)$ and interval $I \subseteq \R$, let 
\begin{align}\label{eq:gaussian_mollifier_freeenergy}
H_u(\btheta) := \sum_{1 \le i \le n} u(S_i)\quad\text{and}\quad Z_u(I) := \int \one\{\theta_1/\sqrt{d} \in I\} \exp \Big( H_u(\btheta) - \GAUSSIANMEASURE \| \btheta \|^2 \Big)\, \rmd \btheta\, .
\end{align}
We do not vary $\beta$ in this proof, so there will be no ambiguity in the above definitions.

Next, for $A \ge 0$, we define the truncated logarithm
\begin{align}\label{eq:truncated_log_def}
\log_A x := \max(\log x, -A)\, .
\end{align}
Several useful properties of the truncated logarithm, Lemmas \upref{lem:truncatedlog1correctproof}, \upref{lem:truncatedlog2}, \upref{lem:truncatedlog3}, are presented in Appendix \upref{subsec:appendixtruncatedlog}. The reason we introduce the truncated logarithm is to truncate various partition functions $Z$ by $\log_{ad}(Z)$, allowing us to use Theorem \upref{thm:gaussiansetup}. Later in the proof when we conclude Theorem \upref{thm:gaussianmeasurersformula}, we will choose $a$ large enough so the truncation is of no impact with high probability using the concentration result Lemma \upref{lem:controlmollifierbernsteinspiked}.

Our first step is to show strong concentration properties of $\log_{ad}Z_{u}(I)$. From here on out, following the convention of this paper excluding Appendix \upref{subsec:gaussian_measure_preliminaries}, we express dependence on the quantities $\signalsubexp, \SIGNALGAUSSIAN$ through the signal-to-noise ratio $\lambda$ governing them. 
\begin{lemma}\label{lem:controlmollifierbernsteinspiked}
Consider any $\MARGIN \in \R$, $\GAUSSIANMEASURE \in [\GAUSSIANMEASURE_0, \GAUSSIANMEASURE_1]$, and $a>0$. Then whenever $u \le 0$ is concave with $u(x)=0$ for all $x \ge \MARGIN$, we have the following. Consider any interval $I = [x-\eps'', x+\eps'']$ with $\eps'' \le \frac1{K'_{\upref{thm:gaussiansetup}}}$ where $K'_{\upref{thm:gaussiansetup}}$ depends on $a, \beta_0, \beta_1, \lambda$. Then for all $d$ large enough in terms of $a, \GAUSSIANMEASURE_0, \GAUSSIANMEASURE_1, \MARGIN, \lambda$ and all $n \le 2d$, we have for all $t > 0$,
\begin{align}
&\P\Bigl( \Big| \frac1d \log_{ad} Z_{u}(I) - \E \Big[\frac1d \log_{ad} Z_{u}(I) \Big] \Big| \ge t \Bigr) \le 2\exp\Bigl( -\frac{d}{K(a, \beta_0, \beta_1, \MARGIN, \lambda)}\min\{ t^2, t\} \Bigr)\, .\label{eq:mollifierbernstein}
\end{align}
\end{lemma}
Note Lemma \upref{lem:controlmollifierbernsteinspiked} applies for any concave $u \le 0$ with $u(x)=0$ for all $x \ge \MARGIN$, and does not impose further regularity assumptions on $u$. In particular, it applies for $u$ given by $\exp u(x) = \one\{x \ge \MARGIN\}$.

\begin{proof}
Let $\mathfrak{F}_i$ be the $\sigma$-algebra generated by the randomness of the first $i$ samples. Let $\E^i$ denote conditional expectation w.r.t. $\mathfrak{F}_i$. Let
\begin{align}\label{eq:mollifier_bernstein_Xi_def}
X_i(I) := \E^i \Big[ \frac1d \log_{ad} Z_{u}(I) \Big]-\E^{i-1} \Big[ \frac1d \log_{ad} Z_{u}(I) \Big]\,.
\end{align}
Note $X_i(I)$ is $\mathfrak{F}_i$-measurable and $\E^{i-1}[X_i(I)]=0$. 
First, we claim it suffices to prove for all $i$ that:
\begin{align}
\E^{i-1}\Big[ \exp \Big( \frac{d\,|X_i(I)|}{K(a, \beta_0, \beta_1, \MARGIN, \lambda)} \Big) \Big] \le 2\, ,\label{eq:firstgoal}
\end{align}
To this end, note that upon establishing \equpref{eq:firstgoal} for all $1 \le i \le n$, applying Bernstein's Inequality on the $X_i(I)$ and recalling $n \le 2d$ yields
\begin{align*}
\P\Biggl( \Big| \frac1d \log_{ad} Z_{\beta}(I) - \E \Big[\frac1d \log_{ad} Z_{\beta}(I) \Big] \Big| \ge t \Biggr) &= \P\Biggl( \Big|\sum_{i=1}^n X_i(I) \Big| \ge t  \le 2\exp\Bigl( -\frac{d \min\{ t^2, t\}}{K(a, \beta_0, \beta_1, \MARGIN, \lambda)} \Bigr) \Bigg)\, .
\end{align*}
The rest of the proof of Lemma \upref{lem:controlmollifierbernsteinspiked} is now devoted to showing \equpref{eq:firstgoal}.
We define
\begin{align}
W_i(I) &:= \int \one\{\theta_1/\sqrt{d} \in I\} \exp \Big(\sum_{i' \neq i, 1 \le i' \le n} u(S_{i'})-\GAUSSIANMEASURE \| \btheta \|^2 \Big)\, \rmd \btheta\,, \label{eq:bernstein_W_i_normalizingconst} \\
Y_1(I) &:= \one \Big\{W_i(I) \ge \exp(-ad) \Big\} \Big| \log_{ad} \Big(  \frac{Z_{u}(I)}{W_i(I)} \Big) \Big|\, , \label{eq:bernstein_Y1} \\
Y_2(I) &:= \one \Big\{ \frac{Z_{u}(I)}{W_i(I)} < \exp(-ad), W_i(I) > 1 \Big\} \Big| \log_{ad} Z_{u}(I) - \log_{ad} W_i(I) \Big|\,, \label{eq:bernstein_Y2} \\
Y(I) &:= Y_1(I) + Y_2(I)\,. \label{eq:bernstein_Y} 
\end{align}
Now, we begin with some preparatory work to show \equpref{eq:firstgoal}.
\begin{lemma}\label{lem:truncatedlogfix1}
We have the pointwise upper bound
\begin{align*}
\Big| \log_{ad} Z_{u}(I) - \log_{ad}W_i(I) \Big| &\le Y(I)\, .
\end{align*}
\end{lemma}
\begin{proof}[Proof of Lemma \upref{lem:truncatedlogfix1}]
Notice $Z_{u}(I) \le W_i(I)$ as $u \le 0$. We break into the three possible cases:
\begin{itemize}
    \item If $W_i(I) \le \exp(-ad)$, then the left hand side is 0 so the above holds.
    \item If $Z_{u}(I) \le \exp(-ad) < W_i(I)$: If $W_i(I) \le 1$, the above follows as a consequence of Lemma \upref{lem:truncatedlog1correctproof} on truncated logarithms. Else if $W_i(I) > 1$, we must have $Z_{u}(I) / W_i(I) \le \exp(-ad)$ in this case, and therefore the left hand side above is upper bounded by $Y_2(I)$. 
    \item If $Z_{u}(I), W_i(I) > \exp(-ad)$: First suppose $W_i(I) \le 1$. Then this is again a consequence of Lemma \upref{lem:truncatedlog1correctproof}. Else if $W_i(I) > 1$ and $\frac{Z_{u}(I)}{W_i(I)} \ge \exp(-ad)$, then $Z_{u}(I) \ge \exp(-ad)$, and the Lemma follows by definition of logarithms. Finally if $W_i(I) > 1$ and $\frac{Z_{u}(I)}{W_i(I)} < \exp(-ad)$, the left hand side above is upper bounded by $Y_2(I)$. 
\end{itemize}
This proves Lemma \upref{lem:truncatedlogfix1}.
\end{proof}
Notice $W_i(I)$ does not depend on the randomness defining the $i$-th sample. Therefore we have that $\E^i \Big[ \log_{ad} W_i(I) \Big] = \E^{i-1} \Big[ \log_{ad} W_i(I) \Big]$. 
Lemma \upref{lem:truncatedlogfix1} thus gives
\begin{align*}
d \big|X_i(I) \big| &= \Big| \E^i \big[ \log_{ad} Z_{u}(I) \big] - \E^{i-1} \big[  \log_{ad} Z_{\beta}(I) \big] \Big| \\
&= \Big| \E^i \big[ \log_{ad} Z_{u}(I) - \log_{ad} W_i(I) \big] -  \E^{i-1} \big[ \log_{ad} Z_{u}(I) - \log_{ad} W_i(I) \big] \Big| \\
&\le \E^i  \Big[ \big| \log_{ad} Z_{u}(I) - \log_{ad} W_i(I) \big| \Big] + \E^{i-1} \Big[ \big| \log_{ad} Z_{u}(I) - \log_{ad} W_i(I) \big| \Big] \\
&\le \E^i \big[ |Y(I)| \big] + \E^{i-1} \big[ |Y(I)| \big]\, .
\end{align*}
Thus for any $\EXPONENTEXPSCALE>0$, we have
\begin{align}
\E^{i-1}\Big[ \exp\big( \EXPONENTEXPSCALE d \big|X_i(I) \big| \big) \Big] &\le \E^{i-1} \Big[ \exp\big(\EXPONENTEXPSCALE \E^i |Y(I)|\big) \exp\big(\EXPONENTEXPSCALE \E^{i-1} |Y(I)|\big) \Big] \notag \\
&= \exp\Big(\EXPONENTEXPSCALE \E^{i-1}\big[ |Y(I)| \big] \Big) \E^{i-1} \Big[ \exp\big(\EXPONENTEXPSCALE \E^i |Y(I)|\big) \Big] \notag \\
&\le \E^{i-1} \Big[ \exp\big(\EXPONENTEXPSCALE |Y(I)|\big) \Big] \cdot \E^{i-1} \E^i \Big[ \exp\big( \EXPONENTEXPSCALE |Y(I)| \big) \Big] \notag \\
&= \E^{i-1} \Big[ \exp\big(\EXPONENTEXPSCALE |Y(I)|\big) \Big]^2\, .\label{eq:controlmollifierbernsteinintermediate}
\end{align}
Let $\E_i$ denote expectation in the randomness defining the $i$-th sample, thus $\E^{i-1}=\E_i \E^i = \E^i \E_i$. Hence 
\[ \E^{i-1} \Big[ \exp\big(\EXPONENTEXPSCALE Y(I) \big) \Big] =\E^i\E_i \Big[ \exp\big(\EXPONENTEXPSCALE Y(I) \big) \Big]\, .\]
Now, note it suffices to show that for $\EXPONENTEXPSCALE = \EXPONENTEXPSCALE(a, \GAUSSIANMEASURE_0, \GAUSSIANMEASURE_1, \lambda)>0$, we have 
\begin{align}
\E_i \Big[ \exp\big(\EXPONENTEXPSCALE Y(I) \big) \Big] \le K_{\upref{eq:secondgoal}}(a, \beta_0, \beta_1, \MARGIN, \SIGNALMEDIUMDEPENDENCE)\, . \label{eq:secondgoal}
\end{align}
This is because Jensen's Inequality then implies 
\begin{align*}
\E_i \Big[\exp\Big(\frac{\EXPONENTEXPSCALE}{K_{\upref{eq:secondgoal}}(a, \beta_0, \beta_1, \MARGIN, \SIGNALMEDIUMDEPENDENCE)} Y(I) \Big) \Big] &\le \E_i \Big[ \exp\big(\EXPONENTEXPSCALE Y(I) \big) \Big]^{\frac{1}{K_{\upref{eq:secondgoal}}(a, \beta_0, \beta_1, \MARGIN, \SIGNALMEDIUMDEPENDENCE)}} \\
&\le K_{\upref{eq:secondgoal}}(a, \beta_0, \beta_1, \MARGIN, \SIGNALMEDIUMDEPENDENCE)^{1/K_{\upref{eq:secondgoal}}(a, \beta_0, \beta_1, \MARGIN, \SIGNALMEDIUMDEPENDENCE)} \le 2\, ,
\end{align*}
hence implying \equpref{eq:firstgoal} when combined with the above.

If $W_i(I) \le \exp(-ad)$ then $Y(I)=0$, and \equpref{eq:secondgoal} immediately follows. Thus suppose from now on that $W_i(I) \ge \exp(-ad)$; we will establish \equpref{eq:secondgoal}. Note this step is permitted because $\E_i$ is only with respect to the randomness defining the $i$-th sample, while $W_i(I)$ does not depend on this randomness, and hence the randomness defining $W_i(I)$ and $\E_i$ are independent.

By AM-GM, to show \equpref{eq:secondgoal}, it suffices to show for some $\EXPONENTEXPSCALE = \EXPONENTEXPSCALE(a, \GAUSSIANMEASURE_0, \GAUSSIANMEASURE_1, \MARGIN)>0$ that
\begin{align}
\E_i \Big[ \exp\big(\EXPONENTEXPSCALE Y_1(I) \big) \Big]\, ,\, \E_i \Big[ \exp\big(\EXPONENTEXPSCALE Y_2(I) \big) \Big]\le K(a, \beta_0, \beta_1, \MARGIN, \SIGNALMEDIUMDEPENDENCE)\, . \label{eq:realgoal}
\end{align}
To this end, define the measure $G_i$ on $\R^d$ by 
\[ G_i(\btheta) \propto \one\{ \theta_1/\sqrt{d} \in I\}\, \exp\Big( \sum_{1 \le i' \le n, i' \neq i} u(S_{i'})-\GAUSSIANMEASURE \| \btheta \|^2 \Big)\, .\]
Note as the normalizing constant of $G_i(\btheta)$ is exactly $W_i(I)$ by \equpref{eq:bernstein_W_i_normalizingconst}, and as $u(x) = 0$ for all $x \ge \MARGIN$ and $u(x) \le 0$, we have
\begin{align}
\frac{Z_{u}(I)}{W_i(I)} = \int \exp u(S_i)\, \rmd G_i(\btheta)\quad,\quad 1 \ge \frac{Z_{u}(I)}{W_i(I)}= \int \exp u(S_i)\, \rmd G_i(\btheta) \ge G_i(S_i \ge \MARGIN)\, . \label{eq:desiredimplicationintermediate}
\end{align}
We apply Theorem \upref{thm:gaussiansetup} with $U(\btheta)$ defined by $\exp( U(\btheta) )=\one\{ \theta_1/\sqrt{d} \in I \} \exp\big( \sum_{i' \neq i} u(S_{i'}) \big)$. 
Note $U \le 0$ is concave as $u \le 0$ is concave and the $S_i$ are affine.
Since $G_i(S_i \ge \MARGIN) \le Z_{u}(I) / W_i(I)$ from \equpref{eq:desiredimplicationintermediate}, and as $\eps'' \le \frac1{K'_{\upref{thm:gaussiansetup}}}$ in the statement of this Lemma and the condition $W_i(I) \ge \exp(-ad)$, we obtain
\begin{align}
\exp\big( -d / K_{\upref{thm:gaussiansetup}}\big) \le t \le \frac14 \exp\big( -K_{\upref{lem:keyconcentrationgaussianconversion}} \max\{1,\MARGIN\}^2 \big) \implies \P\Bigg( \frac{Z_{u}(I)}{W_i(I)} \le t \Bigg) \le 16 K_{\upref{lem:keyconcentrationgaussianconversion}}  t^{1/K^2_{\upref{lem:keyconcentrationgaussianconversion}}}\,,\label{eq:setupforlemma}
\end{align}
for $K_{\upref{thm:gaussiansetup}}=K_{\upref{thm:gaussiansetup}}(a, \beta_0, \beta_1, \SIGNALMEDIUMDEPENDENCE), K'_{\upref{thm:gaussiansetup}}=K'_{\upref{thm:gaussiansetup}}(a, \beta_0, \beta_1, \SIGNALSHORTDEPENDENCE) \ge 2$ and $K_{\upref{lem:keyconcentrationgaussianconversion}} = K_{\upref{lem:keyconcentrationgaussianconversion}}(a, \beta_0, \beta_1, \SIGNALMEDIUMDEPENDENCE) \ge 1$.

\paragraph{Upper bounding $\E_i\big[ \exp(\EXPONENTEXPSCALE Y_1(I) ) \big]$.}
We apply Lemma \upref{lem:momentsofRV838tal} with
\begin{align*}
V = \exp\big( Y_1(I) \big)\,,\,C_0 = K_{\upref{thm:gaussiansetup}}\,,\,C_1 = 16K_{\upref{lem:keyconcentrationgaussianconversion}}\,,\,C_2 = K^2_{\upref{lem:keyconcentrationgaussianconversion}}\,,\,D_0=4 \exp\big(K_{\upref{lem:keyconcentrationgaussianconversion}} \max\{1,\MARGIN\}^2 \big)\,.
\end{align*}
As $Z_{u}(I) \le W_i(I)$ by \equpref{eq:desiredimplicationintermediate}, we have $Y_1(I) = \Big|\log_{ad}\Big( \frac{Z_{u}(I)}{W_i(I)} \Big)\Big|$ and $\exp Y_1(I) \le \min\Big\{ \frac{ W_i(I) }{ Z_u(I)}, \exp(ad)\Big\}$.
Thus by \equpref{eq:setupforlemma}, it follows that the conditions for Lemma \upref{lem:momentsofRV838tal} are satisfied. Thus there is such a $\EXPONENTEXPSCALE_1=\EXPONENTEXPSCALE_1(a, \beta_0, \beta_1, \SIGNALMEDIUMDEPENDENCE)>0$ for which
\begin{align}
\E_i\Big[ \exp\big(\EXPONENTEXPSCALE_1 \, Y_1(I) \big) \Big] \le K(a, \beta_0, \beta_1, \MARGIN, \SIGNALMEDIUMDEPENDENCE)\, .\label{eq:expY_1bounded}
\end{align}

\paragraph{Upper bounding $\E_i\big[ \exp(\EXPONENTEXPSCALE Y_2(I) ) \big]$.} As $u \le 0$, we have the pointwise bound
\[ Z_{u}(I)\,,\, W_i(I) \le \int_{\R^d} \exp\big( -\GAUSSIANMEASURE \| \btheta \|^2 \big)\, \rmd \btheta = \Big( \frac{2\pi}{\GAUSSIANMEASURE} \Big)^{d/2}\,.\]
By definition of truncated logarithm, this implies
\[ \big| \log_{ad} Z_{u}(I) - \log_{ad} W_i(I) \big| \le K(a, \GAUSSIANMEASURE_0) d\,.\]
Thus for all $\EXPONENTEXPSCALE>0$, we have from the definition of $Y_2(I)$ in \equpref{eq:bernstein_Y2} that
\begin{align}
\E_i \big[ \exp(\EXPONENTEXPSCALE Y_2) \big] &\le \P_i\Big( \frac{Z_{u}(I)}{W_i(I)} < \exp(-ad) \Big) \exp\big( \EXPONENTEXPSCALE K(a, \GAUSSIANMEASURE_0) d \big) + 1 \,.\label{eq:boundexponential}
\end{align} 

Our strategy now will be to show $\P_i \big( Z_{u}(I) / W_i(I) < \exp(-ad) \big)$ is exponentially small in $d$, and then choose $\EXPONENTEXPSCALE$ small enough to upper bound \equpref{eq:boundexponential}.
First, note as $d > K(a, \beta_0, \beta_1, \MARGIN, \SIGNALMEDIUMDEPENDENCE)$, we have 
\begin{align}
\exp\big( -d / K_{\upref{thm:gaussiansetup}} \big)\, ,\, \exp(-ad) \le \frac14 \exp\big( -K'_{\upref{lem:keyconcentrationgaussianconversion}} \max\{1,\MARGIN\}^2 \big)\, .\label{eq:setupforlemmaconditions}
\end{align} 
Suppose $\exp(-ad) \ge \exp\big( -d / K_{\upref{thm:gaussiansetup}} \big)$. By \equpref{eq:setupforlemmaconditions}, the condition for \equpref{eq:setupforlemma} with $t = \exp(-ad)$ thus holds, and \equpref{eq:setupforlemma} yields
\[ \P_i\Big( \frac{Z_{u}(I)}{W_i(I)} < \exp(-ad) \Big) \le 16 K_{\upref{lem:keyconcentrationgaussianconversion}} \exp\big( -ad / K^2_{\upref{lem:keyconcentrationgaussianconversion}} \big)\,.\]
Else, say $\exp(-ad) < \exp\big( -d / K_{\upref{thm:gaussiansetup}} \big)$. By \equpref{eq:setupforlemmaconditions}, the conditions for \equpref{eq:setupforlemma} with $t = \exp\big( -d / K_{\upref{thm:gaussiansetup}} \big)$ thus holds, and \equpref{eq:setupforlemma} yields
\begin{align*}
\P_i\Big( \frac{Z_{u}(I)}{W_i(I)} < \exp(-ad) \Big) \le \P_i\Big( \frac{Z_{u}(I)}{W_i(I)} < \exp\big( -d / K_{\upref{thm:gaussiansetup}} \big) \Big) \le 16 K_{\upref{lem:keyconcentrationgaussianconversion}} \exp\Big( -\frac{d}{K_{\upref{thm:gaussiansetup}} K^2_{\upref{lem:keyconcentrationgaussianconversion}}}\Big)\,.
\end{align*}
Combining with \equpref{eq:boundexponential} implies that for all $\EXPONENTEXPSCALE>0$, we have
\begin{align*}
\E_i \big[ \exp(\EXPONENTEXPSCALE Y_2) \big] \le 16 K_{\upref{lem:keyconcentrationgaussianconversion}} \exp\Big( -d \min \Big\{\frac{a}{K^2_{\upref{lem:keyconcentrationgaussianconversion}}}, \frac{1}{K_{\upref{thm:gaussiansetup}} K^2_{\upref{lem:keyconcentrationgaussianconversion}}} \Big\} + \EXPONENTEXPSCALE K(a, \GAUSSIANMEASURE_0) d \Big)\, .    
\end{align*}
Taking 
\begin{align*}
\EXPONENTEXPSCALE_2=\EXPONENTEXPSCALE_2(a, \beta_0, \beta_1, \SIGNALMEDIUMDEPENDENCE)=\frac{1}{2K(a, \GAUSSIANMEASURE_0)} \min \Big\{\frac{a}{K^2_{\upref{lem:keyconcentrationgaussianconversion}}}, \frac{1}{K_{\upref{thm:gaussiansetup}} K^2_{\upref{lem:keyconcentrationgaussianconversion}}} \Big\} > 0    
\end{align*}
to be small enough, we obtain from the above that
\begin{align}
\E_i \big[ \exp(\EXPONENTEXPSCALE_2 Y_2) \big] \le 16K_{\upref{lem:keyconcentrationgaussianconversion}}+1\,.\label{eq:expY_2bounded}
\end{align}

Taking $\EXPONENTEXPSCALE = \min\{ \EXPONENTEXPSCALE_1, \EXPONENTEXPSCALE_2 \}>0$, as $Y_1(I), Y_2(I) \ge 0$, we obtain from \equpref{eq:expY_1bounded}, \equpref{eq:expY_2bounded} that 
\[ \E_i\big[ \exp(\lambda Y_1(I) ) \big]\, ,\, \E_i \big[ \exp(\lambda Y_2(I) )\big] \le K(a, \beta_0, \beta_1, \MARGIN, \SIGNALMEDIUMDEPENDENCE)\,,\]
which establishes \equpref{eq:realgoal} and thus finishes the proof of Lemma \upref{lem:controlmollifierbernsteinspiked}.
\end{proof}

Our next goal is to show that when the $\eps'$ defining the mollifiers $u(\cdot)$ are chosen suitably, the mollification has only a minor impact on the free energy.
\begin{lemma}\label{lem:controlmollifiers}
For any $a>0$ and $\eps \in (0, 1)$, there is $\eps' > 0$ depending on $\eps, a, \beta_0, \beta_1, \MARGIN, \SIGNALMEDIUMDEPENDENCE$ such that for any concave $u \le 0$ satisfying
\[ u(x)=0\text{ for all }x \ge \MARGIN \text{ and } \exp u\big( \MARGIN-\eps' \big) \le \eps'\,,\]
the following holds. For all $\GAUSSIANMEASURE \in [\GAUSSIANMEASURE_0, \GAUSSIANMEASURE_1]$, any $x \in [-1,1]$, and any interval $I = [x-\eps'', x+\eps'']$ with $\eps'' < \frac1{K'_{\upref{thm:gaussiansetup}}}$ where $K'_{\upref{thm:gaussiansetup}}$ depends on $a, \beta_0, \beta_1, \SIGNALSHORTDEPENDENCE$, we have for all $d \ge d(\eps, a, \beta_0, \beta_1, \MARGIN, \SIGNALMEDIUMDEPENDENCE)$ and all $n \le 2d$,
\begin{align}
&\Big| \frac1d\E \big[ \log_{ad} Z_{\beta}(I) \big] - \frac1d \E \big[ \log_{ad} Z_u(I) \big] \Big| \le \eps\, ,\label{eq:contromollifierstobound}
\end{align}
where we recall the definition of $Z_{\beta}(I)$ from \equpref{eq:free_energy_gaussian_def} and $Z_u(I)$ from \equpref{eq:gaussian_mollifier_freeenergy}.
\end{lemma}
\begin{proof}
We replace each $\one\{U_i\}$ with the corresponding mollified quantity $\exp\big( u(S_i) \big)$ and show this leads to only a small difference per $i$. For $i \le n$, define
\begin{align}
C_i :=\bigcap_{i'<i}U_{i'}\quad ,\quad V_i(I) := \int \one\big\{ C_i \cap \{  \theta_1/\sqrt{d} \in I \} \big\} \exp\Big( \sum_{i \le i' \le n} u(S_{i'})-\GAUSSIANMEASURE \| \btheta \|^2 \Big)\, \rmd \btheta\, .\label{eq:gaussian_conversion_Cidef}    
\end{align}
The quantity of interest \equpref{eq:contromollifierstobound} in the statement of the Lemma is thus 
\[ \Big| \frac1d \E \big[ \log_{ad}V_1(I) \big] - \frac1d \E\big[ \log_{ad}V_{n+1}(I) \big] \Big| \le \frac1d \sum_{1 \le i \le n} \Big| \E \big[ \log_{ad} V_{i+1}(I) \big] - \E \big[ \log_{ad}V_{i}(I) \big] \Big|\, .\]
We will bound each of the terms above. To do so, consider any $i, 1 \le i \le n$. We use $\E_i$ to denote expectation and $\P_i$ to denote probability in the randomness defining the $i$-th sample. Notice
\begin{align*}
\Big| \E \big[ \log_{ad} V_{i+1}(I) \big] - \E \big[ \log_{ad}V_{i}(I) \big] \Big| &\le \E \E_i \Big[ \big| \log_{ad}V_{i+1}(I) - \log_{ad}V_{i}(I) \big|\Big]\, .
\end{align*}
We will prove for any $i, 1 \le i \le n$ that
\begin{align}
\E_i \Big[ \big| \log_{ad}V_{i+1}(I) - \log_{ad}V_{i}(I) \big|\Big] \le K\eps\, ,\label{eq:gaussianconversionmollifierclosehighlevel}
\end{align}
where $K>0$ is a universal constant independent of everything else. Upon scaling $\eps$ by a universal constant, \equpref{eq:gaussianconversionmollifierclosehighlevel} proves the desired result. 

To this end, define the Gibbs measure $G_i$ on $\R^d$ as follows: for all subsets $B \subseteq \R^d$, we let
\[ G_i(B) := \frac1{Z_i(I)} \int \one\big\{C_i \cap B \cap \{ \theta_1/\sqrt{d} \in I\} \big\} \exp\Big( \sum_{i < i' \le n} u(S_{i'})-\GAUSSIANMEASURE \| \btheta \|^2 \Big)\, \rmd \btheta\,,\]
where 
\[ Z_i(I) := \int_{C_i} \one\big\{C_i \cap \{ \theta_1/\sqrt{d} \in I\} \big\} \exp\Big( \sum_{i < i' \le n} u(S_{i'})-\GAUSSIANMEASURE \| \btheta \|^2 \Big)\, \rmd \btheta \]
is the normalizing constant. Recalling the definition of $C_i$ from \equpref{eq:gaussian_conversion_Cidef}, note that $Z_i(I)$ does not depend on the randomness defining the $i$-th sample.

Let $\langle \cdot \rangle$ denote expectation w.r.t. $G_i$. Recalling the definition of $C_i$ from \equpref{eq:gaussian_conversion_Cidef} once more, we notice for analogous reasons as above, $\langle \cdot \rangle$ does not depend on the randomness defining the $i$-th sample. Next, notice 
\begin{align}
V_i(I) &= Z_i(I)\, \big\langle \exp u(S_i) \big\rangle \, , \label{eq:gaussian_conversion_Vi}\\
V_{i+1}(I) &= Z_i(I)\, \big\langle \one\{U_i\} \big\rangle = Z_i(I)\, G_i(U_i)=Z_i(I)\, G_i\big(\{S_i \ge \MARGIN\}\big)\, .\label{eq:gaussian_conversion_Vi+1}
\end{align}
As $u \le 0$, we have $V_i(I), V_{i+1}(I) \le Z_i(I)$.

Now, we break into cases based on $Z_i(I)$. This step is permitted because $\E_i$ is only w.r.t. the randomness defining the $i$-th sample, while as noted above, $Z_i(I)$ does not depend on this randomness. 
First, suppose $Z_i(I) \le \exp(-ad)$. Here we have
\[ \E_i \Big[ \big| \log_{ad}V_{i+1}(I) - \log_{ad}V_{i}(I) \big| \Big]=0\, ,\]
establishing \equpref{eq:gaussianconversionmollifierclosehighlevel}.

Now, suppose $Z_i(I) \ge \exp(-ad)$ from here on out, which is permitted as justified above. Note as $u(x)=0$ for $x \ge \MARGIN$, observe that
\begin{align}
0 \le Y  \le X \le 1\qquad\text{ where } \qquad X := \big\langle \exp u(S_i) \big\rangle \,,\, Y := G_i\big(\{S_i \ge \MARGIN \}\big)\, .\label{eq:gaussianconversionmollifierXYineq}
\end{align}
The idea is that $X$ and $Y$ should be close if $u(\cdot)$ is a mollifier that accurately represents the indicator; formalizing this intuition is how we make quantitative bounds. 
Before breaking into further cases on $Z_i(I)$, we first present two useful Lemmas.
\begin{lemma}\label{lem:gaussianconversionmollifierprobbound}
Let 
\begin{align}
c := \Big( \frac{\eps^2}{4 K_{\upref{lem:keyconcentrationgaussianconversion}} K_1} \Big)^{K^2_{\upref{lem:keyconcentrationgaussianconversion}}} > 0\, ,\label{eq:gaussianconversionmollifierchoiceofc}
\end{align}
where $K_{\upref{lem:keyconcentrationgaussianconversion}}>0$ is from Theorem \upref{thm:gaussiansetup} and depends on $a, \beta_0, \beta_1, \SIGNALMEDIUMDEPENDENCE$, and where $K_1>0$ is large enough in terms of $a, \beta_0, \beta_1, \MARGIN, \SIGNALMEDIUMDEPENDENCE$. Then if $Z_i(I) \ge \exp(-ad)$, for $d$ large enough in terms of $\eps, a, \beta_0, \beta_1, \MARGIN, \SIGNALMEDIUMDEPENDENCE$, we have 
\[ \E_i \Big[ \big|\log_{ad} Y \big|\, \one\big\{Y \le c \big\} \Big] \le \eps\, ,\]
where $Y$ is defined as per \equpref{eq:gaussianconversionmollifierXYineq}.
\end{lemma}
\begin{proof}[Proof of Lemma \upref{lem:gaussianconversionmollifierprobbound}]
Since $Z_i(I) \ge \exp(-ad)$ and as $\eps'' \le \frac1{K'_{\upref{thm:gaussiansetup}}}$ in the statement of Lemma \upref{lem:controlmollifiers}, we may apply Theorem \upref{thm:gaussiansetup} with $U(\btheta) = \sum_{i < i' \le n} u(S_{i'})$ and $C = C_i$. Thus $G_C$ defined in Theorem \upref{thm:gaussiansetup} is precisely $G_i$. Note $U \le 0$ and is concave as $u \le 0$ is concave and the $S_i$ are affine. We obtain
\begin{align}
\exp( -d / K_{\upref{thm:gaussiansetup}}) \le t \le \frac14 \exp\Big( -K_{\upref{lem:keyconcentrationgaussianconversion}} \max\{1,\MARGIN\}^2 \Big) \implies \P( Y \le t ) \le 16 K_{\upref{lem:keyconcentrationgaussianconversion}}  t^{1/K^2_{\upref{lem:keyconcentrationgaussianconversion}}}\label{eq:annoyinggaussianargumentimplication}\,.
\end{align}
for $K_{\upref{thm:gaussiansetup}}=K_{\upref{thm:gaussiansetup}}(a, \beta_0, \beta_1, \SIGNALMEDIUMDEPENDENCE)$, $K'_{\upref{thm:gaussiansetup}}=K'_{\upref{thm:gaussiansetup}}(a, \beta_0, \beta_1, \SIGNALSHORTDEPENDENCE)$, and $K_{\upref{lem:keyconcentrationgaussianconversion}} = K_{\upref{lem:keyconcentrationgaussianconversion}}(a, \beta_0, \beta_1, \SIGNALMEDIUMDEPENDENCE)$.

Next, note $Y \le 1$ by \equpref{eq:gaussianconversionmollifierXYineq}, so
\[ \big|\log_{ad} Y \big| = \log V\text{ where } V := \min\big\{ \exp (ad), 1/Y \big\} \ge 1\, .\]
Thus we may apply Lemma \upref{lem:momentsofRV838tal} with $V$ defined as above, and $C_0 = K_{\upref{thm:gaussiansetup}}$, $C_1 = 16K_{\upref{lem:keyconcentrationgaussianconversion}}$, $C_2 = K^2_{\upref{lem:keyconcentrationgaussianconversion}}$, and $D_0=4 \exp\big(K_{\upref{lem:keyconcentrationgaussianconversion}} \max\{1,\MARGIN\}^2 \big)$. 
We thus obtain that for $\EXPONENTEXPSCALE=\EXPONENTEXPSCALE(a, \beta_0, \beta_1, \SIGNALMEDIUMDEPENDENCE)>0$,
\[ \E_i \Big[ V^{\EXPONENTEXPSCALE(a, \beta_0, \beta_1, \SIGNALMEDIUMDEPENDENCE)} \Big] \le K(a, \beta_0, \beta_1, \MARGIN, \SIGNALMEDIUMDEPENDENCE)\,.\]
As $V \ge 1$, we have $\log V < \tilde{K}(\EXPONENTEXPSCALE) V^{\EXPONENTEXPSCALE}$ for all $\EXPONENTEXPSCALE>0$, where $\tilde{K}(\EXPONENTEXPSCALE)$ is an explicit univariate function. Applying this for $\EXPONENTEXPSCALE = \EXPONENTEXPSCALE(a, \beta_0, \beta_1, \SIGNALMEDIUMDEPENDENCE)/2$ and combining with the above display gives, by taking $K_1$ larger  if needed,
\begin{align*}
\E_i \Big[ \big(\log V \big)^2 \Big] &\le \tilde{K}(\EXPONENTEXPSCALE(a, \beta_0, \beta_1, \SIGNALMEDIUMDEPENDENCE)/2)^2 \cdot \E_i \Big[ V^{\EXPONENTEXPSCALE(a, \beta_0, \beta_1, \SIGNALMEDIUMDEPENDENCE)} \Big] \\
&\le \tilde{K}(\EXPONENTEXPSCALE(a, \beta_0, \beta_1, \SIGNALMEDIUMDEPENDENCE)/2)^2 \cdot K(a, \beta_0, \beta_1, \MARGIN, \SIGNALMEDIUMDEPENDENCE) \,.
\end{align*} 
Thus, we may take $K_1$ large enough in terms of $a, \beta_0, \beta_1, \MARGIN, \SIGNALMEDIUMDEPENDENCE$ so that 
\begin{align}
\E_i \Big[ \big(\log V \big)^2 \Big] &\le K_1\,, \label{eq:gaussianconversion_K1_choice1}\\
\Big( \frac1{4 K_{\upref{lem:keyconcentrationgaussianconversion}} K_1} \Big)^{K^2_{\upref{lem:keyconcentrationgaussianconversion}}} &\le \frac14 \exp\big( -K_{\upref{lem:keyconcentrationgaussianconversion}} \max\{1,\MARGIN\}^2 \big)\,.\label{eq:gaussianconversion_K1_choice2}
\end{align}
Since $\eps \le 1$, for $d \ge d(\eps, a, \beta_0, \beta_1, \MARGIN, \SIGNALMEDIUMDEPENDENCE)$, our choice of $c>0$ implies that the conditions of \equpref{eq:annoyinggaussianargumentimplication} hold with $t=c$.
Now \equpref{eq:annoyinggaussianargumentimplication} gives $\P_i(Y \le c) \le 4 K'_{\upref{lem:keyconcentrationgaussianconversion}}  c^{1/K'^2_{\upref{lem:keyconcentrationgaussianconversion}}}$. Cauchy-Schwarz and our choice of $c$ now implies
\begin{align*}
\E_i \Big[ \big|\log_{ad} Y \big|\, \one\big\{Y \le c\big\} \Big] &\le \E_i \Big[ \big|\log_{ad} Y \big|^2 \Big]^{1/2} \E_i \Big[ \one\big\{Y \le c\big\} \Big]^{1/2} \\
&=  \E_i \Big[ \big(\log V \big)^2 \Big]^{1/2} \, \P_i \big(Y \le c \big)^{1/2} \\
&\le K_1^{1/2} \cdot \big( 4 K'_{\upref{lem:keyconcentrationgaussianconversion}}  c^{1/K'^2_{\upref{lem:keyconcentrationgaussianconversion}}} \big)^{1/2} \\
&\le \eps\, ,
\end{align*}
by our choice of $c>0$. This proves Lemma \upref{lem:gaussianconversionmollifierprobbound}.
\end{proof}

\begin{lemma}\label{lem:gaussianconversionmollifierXYbound}
Define $X, Y$ as per \equpref{eq:gaussianconversionmollifierXYineq}. Then if $Z_i(I) \ge \exp(-ad)$, for $\eps' \le \min\{ 1, r_{\upref{lem:upperboundintegral}}^2\}$ where $r_{\upref{lem:upperboundintegral}}$ depends on $a, \GAUSSIANMEASURE_0, \GAUSSIANMEASURE_1$, we have 
\begin{align}
\E_i \Big[ \big|X-Y \big| \Big] \le \exp(-d) + K\sqrt{\eps'}\,,\label{eq:gaussianconversionXminusYbound}
\end{align}
where $K>0$ is a universal constant.
\end{lemma}
\begin{proof}[Proof of Lemma \upref{lem:gaussianconversionmollifierXYbound}]
Since $\exp u(\MARGIN-\eps') \le \eps'$, we have the following pointwise:
\[ 0 \le X-Y \le \eps'+\Big\langle \one\big\{\MARGIN-\eps' \le S_i \le \MARGIN \big\} \Big\rangle\, ,\]
where we recall that $\langle \cdot \rangle$ denotes expectation w.r.t. $G_i$, which does not depend on the randomness defining the $i$-th sample. Since the event $\big\{ \| \bar{\btheta} \| > \sqrt{\eps' d} \big\}$ is independent of the randomness defining the $i$-th sample, we obtain
\begin{align}
0 &\le \E_i\Big[ \big|X-Y \big| \Big] \notag \\
&\le \eps'+\Big \langle \E_i \Big[ \one\big\{\MARGIN-\eps' \le S_i \le \MARGIN \big\} \Big] \Big\rangle \notag \\
&= \eps' + \Big\langle \P_i \Big( \{\MARGIN-\eps' \le S_i \le \MARGIN \} \, \one\{ \big\| \bar{\btheta} \| > \sqrt{\eps' d} \big\} \Big) \Big\rangle  + \Big\langle \P_i \Big( \{\MARGIN-\eps' \le S_i \le \MARGIN \} \, \one\big\{ \| \bar{\btheta} \| \le \sqrt{\eps' d} \big\} \Big) \Big\rangle  \notag \\
&\le \eps' + \Big\langle \P_i \Big( \MARGIN-\eps' \le S_i \le \MARGIN \Big)\, \one\big\{ \| \bar{\btheta} \| > \sqrt{\eps' d} \big\} \Big\rangle + \Big\langle \one\big\{ \| \bar{\btheta} \| \le \sqrt{\eps' d} \big\} \Big\rangle\, .\label{eq:gaussianconversionXminusYboundintermediate}
\end{align}
Recall $Z_i(I) \ge \exp(-ad)$. Thus for $\eps' \le r_{\upref{lem:upperboundintegral}}(a, \beta_0, \beta_1)^2$, noting the conditions apply for the Gibbs measure $\langle \cdot \rangle$ as $u \le 0$, we obtain from Lemma \upref{lem:upperboundintegral} that
\[ \Big\langle \one\big\{ \| \bar{\btheta} \| \le \sqrt{\eps' d} \big\} \Big\rangle \le \exp(-d)\, .\]

Now consider any $\MARGIN_1<\MARGIN_2$. Letting $\bar{\P}, \bar{\E}$ denote probability and expectation in $g_i$, observe that
\begin{align*}
\P_i \big( \MARGIN_1 < S_i < \MARGIN_2 \big) &= \int_{\R} \phi(s)\, \bar{\P} \Big(\frac{\langle g_i, \bar{\btheta} \rangle}{\sqrt{d}} \in \Big[ \MARGIN_1 - \frac{s \theta_1}{\sqrt{d}}, \MARGIN_2 - \frac{s \theta_1}{\sqrt{N}} \Big] \Big)\, \rmd s\, ,
\end{align*}
where $\phi(\cdot)$ denotes the density of $s_i$. For a given $s \in \R$, let $\MARGIN_1' := \MARGIN_1 - \frac{s \theta_1}{\sqrt{d}}$, $\MARGIN_2' := \MARGIN_2 - \frac{s \theta_1}{\sqrt{d}}$ and $Z = \frac{ \langle g_i, \bar{\btheta} \rangle }{\sqrt{d}}$. Note as $Z \sim N\Big(0, \frac{\|\bar{ \btheta}\|^2}d\Big)$ is Gaussian,
\begin{align*}
\bar{\P} \Big(\frac{ \langle g_i, \bar{\btheta} \rangle }{\sqrt{d}}\in \Big[ \MARGIN_1 - \frac{s \theta_1}{\sqrt{d}}, \MARGIN_2 - \frac{s \theta_1}{\sqrt{d}} \Big] \Big) &= \frac{1}{\sqrt{2\pi \E [Z^2] }} \int_{\MARGIN'_1}^{\MARGIN'_2} \exp\Big( -\frac{t^2}{2 \E [Z^2]} \Big) \, \rmd t \\
&\le \frac{K}{\sqrt{\E [Z^2]}} (\MARGIN'_2-\MARGIN'_1) \\
&= \frac{K\sqrt{d}}{\| \bar{\btheta} \|} (\MARGIN_2-\MARGIN_1)\,,
\end{align*}
for $K>0$ a universal constant. Thus
\begin{align*}
\P_i \big( \MARGIN_1 < S_i < \MARGIN_2 \big) &\le \int_{\R} \phi(s)\, \frac{K\sqrt{d}}{\| \bar{\btheta} \|} (\MARGIN_2-\MARGIN_1)\, \rmd s = \frac{K\sqrt{d}}{\| \bar{\btheta} \|} (\MARGIN_2-\MARGIN_1)\, ,
\end{align*}
yielding that for any $\theta$,
\[ \P_i \big( \MARGIN-\eps' \le S_i \le \MARGIN \big)\, \one\big\{ \| \bar{\btheta} \| > \sqrt{\eps' d} \big\} \le \frac{K\eps' \sqrt{d}}{\| \bar{\btheta} \|} \, \one\{ \| \bar{\btheta} \| > \sqrt{\eps' d} \} \le K \sqrt{\eps'}\, .\]
Combining these steps with \equpref{eq:gaussianconversionXminusYboundintermediate} yields 
\begin{align*}
\E_i \Big[ \big|X-Y\big| \Big] \le \eps' + \exp(-d) + K\sqrt{\eps'}\, ,
\end{align*}
and hence \equpref{eq:gaussianconversionXminusYbound}.
\end{proof}

We now return to the proof of Lemma \upref{lem:controlmollifiers}. We break into cases on the value of $Z_i(I)$, recalling that it suffices to prove \equpref{eq:gaussianconversionmollifierclosehighlevel} when $Z_i(I) \ge \exp(-ad)$.

\paragraph{If $\exp(-ad) \le Z_i(I) \le 1$.}
Applying Lemma \upref{lem:truncatedlog2} and Lemma \upref{lem:truncatedlog3} on truncated logarithms in the first and second inequality below respectively, as $X \ge Y$ and combining \equpref{eq:gaussian_conversion_Vi}, \equpref{eq:gaussian_conversion_Vi+1} and \equpref{eq:gaussianconversionmollifierXYineq}, we obtain for any $c>0$ that
\begin{align}
\big| \log_{ad}V_{i+1}-\log_{ad}V_i \big| &= \Big| \log_{ad} \big( Z_i(I) Y \big) - \log_{ad}\big( Z_i(I) X \big) \Big| \notag \\
&\le \big|\log_{ad} Y - \log_{ad} X \big| \notag \\
&= \big|\log_{ad} X - \log_{ad} Y \big| \notag \\
&\le \big|\log_{ad} Y \big|\, \one\{Y \le c\} + \frac{|X-Y|}c\, .\notag
\end{align}
Taking $c = c(\eps, a, \beta_0, \beta_1, \MARGIN, \SIGNALMEDIUMDEPENDENCE) > 0$ as per \equpref{eq:gaussianconversionmollifierchoiceofc}, Lemmas \upref{lem:gaussianconversionmollifierprobbound} and \upref{lem:gaussianconversionmollifierXYbound} give
\begin{align*}
\E_i \Big[ \big| \log_{ad}V_{i+1}-\log_{ad}V_i \big| \Big] &\le \E_i \Big[ \big|\log_{ad} Y \big|\, \one\{Y \le c\} \Big] + \frac1c \E_i \Big[ \big|X-Y\big| \Big] \\
&\le \eps + \frac1c \big( \exp(-d) + K\sqrt{\eps'} \big)\, .
\end{align*}
Since $c = c(\eps, a, \beta_0, \beta_1, \MARGIN, \SIGNALMEDIUMDEPENDENCE)$ in \equpref{eq:gaussianconversionmollifierchoiceofc} does not depend on $\eps'$, we have for $\eps' \le \eps'(\eps, a, \beta_0, \beta_1, \MARGIN, \SIGNALMEDIUMDEPENDENCE)$ and $d \ge d(\eps, a, \beta_0, \beta_1, \MARGIN, \SIGNALMEDIUMDEPENDENCE)$ that
\begin{align*}
\E_i \Big[ \big| \log_{ad}V_{i+1}-\log_{ad}V_i \big| \Big] \le K\eps\, ,
\end{align*}
establishing \equpref{eq:gaussianconversionmollifierclosehighlevel} in this case.

\paragraph{If $Z_i(I) > 1$. }Let $c = c(\eps, a, \beta_0, \beta_1, \MARGIN, \SIGNALMEDIUMDEPENDENCE)>0$ be from \equpref{eq:gaussianconversionmollifierchoiceofc}. Note for $d \ge d(\eps, a, \beta_0, \beta_1, \MARGIN, \SIGNALMEDIUMDEPENDENCE)$ that $\exp(-ad) \le c$. Thus we may define the following events:
\begin{align*}
 E_1 := \Big\{Y \le \exp(-ad) \Big\}\, ,\, E_2 := \Big\{ \exp(-ad) < Y \le c\Big\}\, ,\, E_3 := \Big\{Y \ge c \Big\}\,.
\end{align*}
Recall $1 \ge X = \langle \exp u(S_i) \rangle \ge Y$ from \equpref{eq:gaussianconversionmollifierXYineq} and that here $Z_i(I) > 1$. Thus as $Y \ge c > 0$ on $E_3$,
\begin{align*}
Z_i(I) X \one\{E_3\}\, ,\, Z_i(I) Y \one\{E_3\} \ge c \one\{E_3\} \ge \exp(-ad) \one\{E_3\}\,,    
\end{align*}
and therefore
\begin{align*}
\Big| \log_{ad} \big( Z_i(I) X \big) - \log_{ad} \big(Z_i(I) Y \big) \Big|\, \one\{E_3\} &= \Big| \log \big( Z_i(I) X \big) - \log \big( Z_i(I) Y \big) \Big|\, \one\{E_3\}  \\
&= \Big| \log \frac{X}{Y} \Big|\, \one\{E_3\} \\
&\le \Big( \frac{X}{Y}-1 \Big)\, \one\{E_3\} \\
&\le \frac1{c} |X-Y|\,  \one\{E_3\}\, .
\end{align*}
By \equpref{eq:gaussian_conversion_Vi}, \equpref{eq:gaussian_conversion_Vi+1}, we thus may write
\begin{align}
\E_i \Big[ \big| \log_{ad}V_{i+1}-\log_{ad}V_i \big| \Big] &= \E_i \Big[ \big|\log_{ad} \big( Z_i(I) X \big) - \log_{ad} \big( Z_i(I) Y \big) \big| \Big] \notag \\
&= \E_i \Big[ \one\{E_1\} \big| \log_{ad} \big( Z_i(I) X \big) - \log_{ad} \big( Z_i(I) Y \big) \big| \Big] \notag \\
&\qquad + \E_i \Big[ \one\{E_2\} \big| \log_{ad} \big( Z_i(I) X \big) - \log_{ad} \big( Z_i(I) Y \big) \big| \Big] \notag \\
&\qquad + \E_i \Big[ \one\{E_3\} \cdot \frac1c |X-Y| \Big] \notag \\
&:= \expressionI + \expressionII + \expressionIII\,. \label{eq:annoyingfixdecomposition}
\end{align}
Now we upper bound each of $\expressionI, \expressionII, \expressionIII$:
\begin{itemize}
\item \textbf{Upper bounding $\expressionI$:} Since $Z_i(I) \ge \exp(-ad)$ and $\eps'' \le \frac1{K'_{\upref{thm:gaussiansetup}}}$, we may apply Theorem \upref{thm:gaussiansetup} identically as in the proof of Lemma \upref{lem:gaussianconversionmollifierprobbound}, which establishes that the implication \equpref{eq:annoyinggaussianargumentimplication} holds.
First, suppose $\exp(-ad) \ge \exp\big( -d / K_{\upref{thm:gaussiansetup}} \big)$. Taking $t = \exp(-ad)$ in \equpref{eq:annoyinggaussianargumentimplication} gives
\[ \P_i\big( Y < \exp(-ad) \big) \le 16 K_{\upref{lem:keyconcentrationgaussianconversion}}  \exp\Big( -\frac{ad}{K^2_{\upref{lem:keyconcentrationgaussianconversion}}}\Big)\, .\]
Else, say $\exp(-ad) < \exp\big( -d / K_{\upref{thm:gaussiansetup}} \big)$. Taking $t = \exp\big( -d / K_{\upref{thm:gaussiansetup}} \big)$ in \equpref{eq:annoyinggaussianargumentimplication} gives
\begin{align*}
\P_i\big( Y < \exp(-ad) \big) &\le \P_i\Big( Y < \exp\big( -d / K_{\upref{thm:gaussiansetup}} \big)\Big) \le 16 K_{\upref{lem:keyconcentrationgaussianconversion}}  \exp\Big( -\frac{d}{K_{\upref{thm:gaussiansetup}} K^2_{\upref{lem:keyconcentrationgaussianconversion}}}\Big) \, .
\end{align*}
Therefore, it follows that
\begin{align*}
\P_i(E_1) &= \P_i\big( Y < \exp(-ad) \big) \le 16 K_{\upref{lem:keyconcentrationgaussianconversion}} \exp\Big( -\frac{d}{K^2_{\upref{lem:keyconcentrationgaussianconversion}}\max\{1/a, K_{\upref{thm:gaussiansetup}} \}}  \Big)\, .
\end{align*}
Now as $u \le 0$, we have
\begin{align*}
Z_i(I) &\le \int_{\R^d} \exp(-\GAUSSIANMEASURE \| \btheta \|^2)\, \rmd \btheta \le \Big( \frac{2\pi}{\GAUSSIANMEASURE_0} \Big)^{d/2}\, .
\end{align*}
Thus as $X, Y \le 1$ by \equpref{eq:gaussianconversionmollifierXYineq} and by definition of the truncated logarithm, we have 
\[ \Big| \log_{ad} \big( Z_i(I) X \big) - \log_{ad} \big( Z_i(I) Y\big) \Big| \le K(a, \GAUSSIANMEASURE_0) d\, .\]
Thus for $d \ge d(\eps, a, \beta_0, \beta_1, \SIGNALMEDIUMDEPENDENCE)$, we obtain
\begin{align}
&\E_i \Big[ \one\{E_1\} \Big| \log_{ad} \big( Z_i(I) X \big) - \log_{ad} \big( Z_i(I) Y\big) \Big| \Big] \notag \\
&\quad\le K(a, \GAUSSIANMEASURE_0) d \cdot \P_i(E_1) \notag \\
&\quad\le K(a, \GAUSSIANMEASURE_0) d \cdot 16 K_{\upref{lem:keyconcentrationgaussianconversion}} \exp\Big( -\frac{d}{K^2_{\upref{lem:keyconcentrationgaussianconversion}}\max\{1/a, K_{\upref{thm:gaussiansetup}} \}}  \Big) \notag \\
&\quad\le \eps\,.\label{eq:gaussianconversion_mollifier_bound1}
\end{align}
    
\item \textbf{Upper bounding $\expressionII$:} 
On the event $E_2$, we have $c \ge Y > \exp(-ad)$. Thus by \equpref{eq:gaussianconversionmollifierXYineq}, on the event $E_2$, we have $1 \ge X \ge Y \ge \exp(-ad)$ as well. As $Z_i(I) \ge 1$ in this case, we obtain
\begin{align*}
\one\{E_2\} \Big| \log_{ad} \big(Z_i(I) X \big) - \log_{ad} \big( Z_i(I) Y \big) \Big| &= \one\{E_2\} \Big| \log \big( Z_i(I) X \big) - \log \big( Z_i(I) Y \big) \Big| \\
&= \one\{E_2\} | \log X - \log Y | \\
&\le \one\{E_2\} |\log Y| \\
&= \one\{E_2\} \big|\log_{ad} Y \big| \\
&\le \one\{Y \le c\} \big|\log_{ad} Y \big|\,.
\end{align*}
Therefore as $Z_i(I) \ge 1 \ge \exp(-ad)$, Lemma \upref{lem:gaussianconversionmollifierprobbound} implies for $d \ge d(\eps, a, \beta_0, \beta_1, \MARGIN, \SIGNALMEDIUMDEPENDENCE)$,
\begin{align}\label{eq:gaussianconversion_mollifier_bound2}
\E_i \Big[ \one\{E_2\}\, \big| \log_{ad} \big(Z_i(I) X \big) - \log_{ad} \big( Z_i(I) Y \big) \big| \Big] &\le \E_i \Big[  \one\{Y \le c\}\, \big|\log_{ad} Y \big| \Big] \le \eps\, .
\end{align}

\item \textbf{Upper bounding \expressionIII:}
By Lemma \upref{lem:gaussianconversionmollifierXYbound}, as $Z_i(I) \ge 1$ in this case, for $\eps' \le \eps'(a, \GAUSSIANMEASURE_0, \GAUSSIANMEASURE_1)$ we have 
\[ \E_i\Big[ \one\{E_3\} \cdot \frac1c |X-Y| \Big] \le \exp(-d) + K\sqrt{\eps'}\, .\]
Hence, for $\eps' \le \eps'(\eps, a, \beta_0, \beta_1, \MARGIN, \SIGNALMEDIUMDEPENDENCE)$ and for $d \ge d(\eps, a, \beta_0, \beta_1, \MARGIN, \SIGNALMEDIUMDEPENDENCE)$, we have
\begin{align}\label{eq:gaussianconversion_mollifier_bound3}
\E_i\Big[ \one\{E_3\} \cdot \frac1c |X-Y| \Big] \le \eps\, .
\end{align}
\end{itemize}
Combining \equpref{eq:gaussianconversion_mollifier_bound1}, \equpref{eq:gaussianconversion_mollifier_bound2}, \equpref{eq:gaussianconversion_mollifier_bound3} with \equpref{eq:annoyingfixdecomposition}, we obtain that
\[ \E_i \Big[ \big| \log_{ad}V_{i+1}-\log_{ad}V_i \big| \Big] \le 3\eps\,,\]
hence establishing \equpref{eq:gaussianconversionmollifierclosehighlevel} when $Z_i(I) \ge 1$. This completes the proof of Lemma \upref{lem:controlmollifiers}.
\end{proof}

Finally, we put all these steps together to prove Theorem \upref{thm:gaussianmeasurersformula}.
\begin{proof}[Proof of Theorem \upref{thm:gaussianmeasurersformula}]
Without loss of generality, we assume $\eps \le \frac14$ in the following. 
Consider $\bar\Phi(x,q,\rho)$ as defined in \equpref{eq:phi_master_rho} with $\exp u(x)=\one\{x \ge \MARGIN\}$, and note this function is a continuous function of its three arguments. 
By continuity supplied by the second part of Proposition \upref{prop:nicersuniquesol} and compactness, we may let $a = a(\beta_0, \beta_1, \MARGIN, \lambda ) > 0$ be large enough so that the following holds. For all $(\alpha,\GAUSSIANMEASURE, x) \in [0, \DENSITYBOUND] \times [\GAUSSIANMEASURE_0, \GAUSSIANMEASURE_1] \times [-1, 1]$, letting $q_0=q_0(\alpha, \beta, x), \rho_0=\rho_0(\alpha, \beta, x)$ in the display below for brevity, we have 
\[ \bar\Phi (x, q_0, \rho_0 ) - \beta \rho + \frac12 \log(2\pi e) \,,\, \bar\Phi (x, q_0, \rho_0 ) - \beta (\rho+x^2) + \frac12 \log(2\pi e) \ge -a+1\,.\]
Here we note $\DENSITYBOUND$ depends only on $\beta_0, \beta_1, \lambda$ as $\Csol$ depends on $\beta_0, \beta_1$.
Note $a$ does not depend on $\eps'', \eps'$, since $\bar\Phi(x, q, \rho), q_0(\alpha, \beta, x), \rho_0(\alpha, \beta, x)$ do not depend on $I$ as per Remark \upref{rem:sol_nointervaldependence}. 
As $\RSG_I$ equals $\bar\Phi (x, q_0, \rho_0 ) - \beta \rho + \frac12 \log(2\pi e)$ or $\bar\Phi (x, q_0, \rho_0 ) - \beta (\rho+x^2) + \frac12 \log(2\pi e)$ pointwise, we obtain for all $(\alpha, \beta, x) \in [0, \DENSITYBOUND] \times [\beta_0, \beta_1] \times [-1,1]$ that
\begin{align}\label{eq:gaussian_conversion_choiceofa_use}
\RSG_I(\alpha, \beta, x) \ge -a+1\,.
\end{align}

By Lemma \upref{lem:controlmollifiers}, there is $\eps'_1
= \eps'_1(\eps, a, \beta_0, \beta_1, \MARGIN, \SIGNALMEDIUMDEPENDENCE)>0$ such that for all $d \ge d(\eps, a, \beta_0, \beta_1, \MARGIN, \SIGNALMEDIUMDEPENDENCE)$ and all $n \le 2d$, the following holds if $\eps'' \le \frac1{K'_{\upref{thm:gaussiansetup}}}$, where $K'_{\upref{thm:gaussiansetup}}=K'_{\upref{thm:gaussiansetup}}(a, \beta_0, \beta_1, \lambda)$. For any concave $u \le 0$ such that 
\[ u(x)=0\text{ for all }x \ge \MARGIN \text{ and } \exp u\big( \MARGIN-\eps'_1 \big) \le \eps'_1\, ,\]
then 
\begin{align}
\Big| \frac1d\E \big[ \log_{ad} Z_{\beta}(I) \big] - \frac1d \E \big[ \log_{ad} Z_u(I) \big] \Big| \le \eps\, .\label{eq:mollifierclose}
\end{align}
Likewise by Theorem \upref{thm:mollifiedhamiltonianfreeenergy}, there exists $\eps'_2 = \eps'_2(\eps, \beta_0, \beta_1, \MARGIN^+)$ such that for any concave $u \le 0$ such that 
\begin{enumerate}
    \item $u' \ge 0$ with $u' \neq 0$ on a set of positive Lebesgue measure,
    \item $u(x) = 0$ for $x \ge \MARGIN$,
    \item $\exp u(\MARGIN-\eps'_2) \le \eps'_2$,
    \item $u$ is four times differentiable and $|u^{(l)}| \le D$ for $l = 1,2,3,4$,
\end{enumerate}
there is a $\eps_{\RS_0}=\eps_{\RS_0}(\eps, u, \beta_0, \beta_1, \lambda)>0$ such that if $0 < \eps'' \le \eps_{\RS_0}$, the following holds. 
For $\DIMENSION \ge d(\eps, \eps'', \beta_0, \beta_1, \MARGIN, \lambda)$, we have 
\begin{align}
\Big| \frac1d \E\big[ \log Z_u(I) \big] - \RSG_I(n/d, \GAUSSIANMEASURE, x) \Big| \le \eps\, .\label{eq:expectationhamiltonian} 
\end{align}
Let 
\[ \eps'=\eps'(\eps, a, \beta_0, \beta_1, \MARGIN, \SIGNALMEDIUMDEPENDENCE) := \min \{\eps'_1, \eps'_2 \}\, .\]
Given $\eps'$, it is straightforward to construct a mollifier $\hat{u} \le 0$ that is concave, satisfying $\hat{u}' \ge 0$ and $\hat{u}' \neq 0$ on a set of positive Lebesgue measure, satisfying $|\hat{u}^{(l)}| \le D(\eps', \MARGIN)$ for $l=1,2,3,4$, and satisfying  
\[ \hat{u}(x)=0\text{ for all }x \ge \MARGIN\, ,\, \exp \hat{u}(\MARGIN-\eps') \le \eps'\, .\]
Hence as $\hat{u}' \ge 0$ and $\eps' \le \eps'_1, \eps'_2$, it follows that $\exp \hat{u}(\MARGIN-\eps'_1) \le \eps'_1$ and $\exp \hat{u}(\MARGIN-\eps'_2) \le \eps'_2$. 

Since $\eps'=\eps'(\eps, a, \beta_0, \beta_1, \MARGIN, \SIGNALMEDIUMDEPENDENCE)$, we may let 
\[ \eps_{\RSG} = \eps_{\RSG}\big(\eps, \DENSITYBOUND, \beta_0, \beta_1, \MARGIN, \lambda \big) := \min\big\{ 1/K'_{\upref{thm:gaussiansetup}}, \eps_{\RS_0}\big\}\,,\]
where now $\eps_{\RS_0}=\eps_{\RS_0}(\hat u, \beta_0, \beta_1, \lambda)$. This is well-defined as $\hat{u}$ is constructed using $\MARGIN, \eps'$, and as $\eps'$ in turn depends on $\eps, a=a(\DENSITYBOUND, \beta_0, \beta_1, \MARGIN, \lambda), \beta_0, \beta_1, \MARGIN$, and $\lambda$.

Consider any $0 \le \eps'' < \eps_{\RSG}$ and let $I = [x-\eps'', x+\eps'']$. 
Note as $\exp \hat{u}(\MARGIN-\eps'_1) \le \eps'_1$ and $\exp \hat{u}(\MARGIN-\eps'_2) \le \eps'_2$, for $d \ge d(\eps, \eps'', a, \beta_0, \beta_1, \MARGIN, \SIGNALMEDIUMDEPENDENCE)$, we have both \equpref{eq:mollifierclose} and \equpref{eq:expectationhamiltonian} for $Z_{\hat u}(I)$. 

With these steps in hand, the finish is as follows. We need one last argument about controlling the error of the truncated logarithm:
\begin{lemma}\label{lem:gaussianconversionmollifiercloseexpectation}
For $d$ large enough in terms of $\eps, \eps'', a, \beta_0, \beta_1, \MARGIN, \SIGNALMEDIUMDEPENDENCE$, we have 
\begin{align}
\Big| \frac1d \E \big[\log_{ad} Z_{\hat{u}}(I)\big] - \frac1d \E \big[ \log Z_{\hat{u}}(I)\big] \Big| \le \eps\, .\label{eq:mollifiercloseexpectation}
\end{align}
\end{lemma}
\begin{proof}[Proof of Lemma \upref{lem:gaussianconversionmollifiercloseexpectation}]
Let 
\[ \Omega := \big\{Z_{\hat{u}}(I) \le \exp(-ad) \big\}\, . \]
Thus by Cauchy-Schwarz,
\begin{align*}
\Big| \frac1d \E \big[\log_{ad} Z_{\hat{u}}(I)\big] - \frac1d \E \big[ \log Z_{\hat{u}}(I)\big] \Big| &= \Big| \E\Big[ \one\{\Omega\} \Big( a+\frac1d \log Z_{\hat{u}}(I) \Big) \Big] \Big| \\
&\le \P(\Omega) \E\Big[ \Big( a+\frac1d \log Z_{\hat{u}}(I) \Big)^2 \Big] \\
&\le K \P(\Omega) \E\Big[ a^2 + \frac1{d^2} \log^2 Z_{\hat{u}}(I) \Big]\, .
\end{align*}
First we upper bound $\E\Big[ a^2 + \frac1{d^2} \log^2 Z_{\hat{u}}(I) \Big]$. Notice as $\hat{u} \le 0$,
\[ Z_{\hat{u}}(I) \le \int \exp \Big(H_{\hat{u}}(\btheta) - \GAUSSIANMEASURE \| \btheta \|^2 \Big)\, \rmd \btheta \le \int \exp (-\GAUSSIANMEASURE \| \btheta \|^2)\, \rmd \btheta \le \Big(\frac{2\pi}{\GAUSSIANMEASURE_0}\Big)^{d/2}\, , \]
thus 
\[ \E\Big[ a^2 + \frac1{d^2} \log^2 Z_{\hat{u}}(I) \Big] \le K(a, \beta_0)\, .\]

Now we upper bound $\P(\Omega)$. As $\log_A x \ge \log x$ for all $x \in \R$ and as $\eps \le \frac14$, \equpref{eq:expectationhamiltonian} and \equpref{eq:gaussian_conversion_choiceofa_use} yields that for $d \ge d(\eps, \eps'', a, \beta_0, \beta_1, \MARGIN, \SIGNALMEDIUMDEPENDENCE)$,
\[ \frac1d \E\big[ \log_{ad} Z_{\hat{u}}(I) \big] \ge  \frac1d \E \big[ \log Z_{\hat{u}}(I) \big] \ge \RSG_I(n/d, \GAUSSIANMEASURE, x) - \eps \ge -a+1-\eps \ge -a+3/4\, .\]
Therefore, 
\[ \Omega \subseteq \Big\{ \Big| \frac1d \log_{ad} Z_{\hat{u}}(I) - \frac1d \E\big[ \log_{ad} Z_{\hat{u}}(I)\big] \Big| \ge \frac34 \Big\}\, . \]
Now as $\eps'' \le 1/K'_{\upref{thm:gaussiansetup}}$, we may apply Lemma \upref{lem:controlmollifierbernsteinspiked} with $u = \hat{u}$. Thus for $d \ge d(\eps, \eps'', a, \beta_0, \beta_1, \MARGIN, \SIGNALMEDIUMDEPENDENCE)$,
\begin{align*}
\P(\Omega) &\le \P\Big( \Big| \frac1d \log_{ad} Z_{\hat{u}}(I) - \frac1d \E\big[ \log_{ad} Z_{\hat{u}}(I)\big] \Big| \ge \frac34 \Big) \le 2\exp\Big( -\frac{d}{K(a, \beta_0, \beta_1, \MARGIN, \SIGNALMEDIUMDEPENDENCE)} \Big)\,.
\end{align*}
Combining the above bounds, it follows for such $d$ that 
\begin{align*}
\Big| \frac1d \E \big[\log_{ad} Z_{\hat{u}}(I)\big] - \frac1d \E \big[ \log Z_{\hat{u}}(I)\big] \Big| \le K(a, \beta_0) \cdot 2\exp\Big( -\frac{d}{K(a, \beta_0, \beta_1, \MARGIN, \SIGNALMEDIUMDEPENDENCE)} \Big) \le \eps\, .
\end{align*}
This proves Lemma \upref{lem:gaussianconversionmollifiercloseexpectation}.
\end{proof}

Now summing together \equpref{eq:mollifierclose}, \equpref{eq:expectationhamiltonian}, and \equpref{eq:mollifiercloseexpectation}, we obtain for $d \ge d(\eps, \eps'', a, \beta_0, \beta_1, \MARGIN, \SIGNALMEDIUMDEPENDENCE)$ that
\begin{align*}
\Big| \frac1d \E\big[ \log_{ad} Z_{\beta}(I) \big]  - \RSG_I(n/d, \GAUSSIANMEASURE, x) \Big| &\le \Big| \frac1d \E\big[ \log_{ad} Z_{\beta}(I) \big] - \frac1d \E\big[ \log_{ad} Z_{\hat{u}}(I) \big] \Big| \\
&\qquad + \Big| \frac1d \E\big[ \log_{ad} Z_{\hat{u}}(I) \big] - \frac1d \E\big[ \log Z_{\hat{u}}(I)\big]\Big| \\
&\qquad + \Big| \frac1d \E\big[ \log Z_{\hat{u}}(I)\big] - \RSG_I(n/d, \GAUSSIANMEASURE, x) \Big| \\
&\le 3\eps\, .
\end{align*}
By Lemma \upref{lem:controlmollifierbernsteinspiked} with $\exp u(x) = \one\{x \ge \MARGIN\}$ and $t=\eps$, we obtain for such $d$,
\begin{align*}
&\P \Big( \Big| \frac1d \log_{ad} Z_{\beta}(I) -  \frac1d \E\big[ \log_{ad} Z_{\beta}(I) \big] \Big| \ge \eps \Big) \le 2 \exp \Big( - \frac{d}{K(\eps, a, \beta_0, \beta_1, \MARGIN, \SIGNALMEDIUMDEPENDENCE)} \Big)\,,
\end{align*} 
where we obtain dependence on $\eps$ in $K$ as we took $t=\eps$ in Lemma \upref{lem:controlmollifierbernsteinspiked}. Combining the above two displays, we obtain for $d \ge d(\eps, \eps'', a, \beta_0, \beta_1, \MARGIN, \SIGNALMEDIUMDEPENDENCE)$ that
\begin{align}\label{eq:gaussian_conversion_last}
&\P\Big( \Big| \frac1d \log_{ad} Z_{\beta}(I) - \RSG_I(n/d, \GAUSSIANMEASURE, x) \Big|  \ge 4\eps \Big) \le \exp \Big( - \frac{d}{K(\eps, a, \beta_0, \beta_1, \MARGIN, \SIGNALMEDIUMDEPENDENCE)} \Big)\,.
\end{align} 
However, by \equpref{eq:gaussian_conversion_choiceofa_use}, we have $\RSG_I(n/d, \GAUSSIANMEASURE, x) \ge -a+1\ge-a+4\eps$. Thus as events, 
\begin{align*}
&\Big\{ \Big| \frac1d \log_{ad} Z_{\beta}(I) - \RSG_I(n/d, \GAUSSIANMEASURE, x) \Big| \le 4\eps \Big\} \equiv \Big\{ \Big| \frac1d \log Z_{\beta}(I)  - \RSG_I(n/d, \GAUSSIANMEASURE, x) \Big| \le 4\eps \Big\}\, ,
\end{align*}
since under this event $\log_{ad} \equiv \log$ as functions. Combining the above display with \equpref{eq:gaussian_conversion_last} and applying this result for $\eps/4$ and recalling $a = a(\beta_0, \beta_1, \MARGIN, \lambda )$ proves Theorem \upref{thm:gaussianmeasurersformula}.
\end{proof}

\subsection{Proof of Theorem \upref{thm:posteriorgaussianmeasurersformula} for posterior} \label{subsec:posteriorgaussianconversion}
Our strategy will be to apply Theorem \upref{thm:logconcaveconcentrationbasics} on concentration of Lipschitz functions of strongly log-concave measures. Specifically, we consider $\frac1d \log Z_{\beta}(I)$ as a function of the disorder $(s_i, \bar g_i)_{i=1}^n$, and bound its resulting Lipschitz constant; recall as in Appendix \upref{sec:interpolationfreeenergy} that the first four derivatives of $u$ are upper bounded by some quantity $D$, which only depends on $\lambda$. Here in this proof we write dependence on $D$, and in fact we only need the bound on the first derivative of $u$.

Let us begin by considering the space $\mathcal{S} = \R^{nd}$ where the first $n$ coordinates correspond to the $s_j$ and the next $(n-1)d$ coordinates correspond to $g_{i,j}$.
We endow $\mathcal{S}$ with the product measure $\gamma$ where the measure on each coordinate is given by law of each corresponding part of the disorder (that is, $s_j$ or $g_{i,j}$). 
Define
\begin{align}
\label{eq:posteriorgaussian_Bstar_bd}
&\begin{aligned}
B^{\star} &:= \Big( \sum_{i=1}^n \| \bar{g}_i \|^2 \Big)^{1/2} + \frac1{16d} \sum_{j=2}^d \Big(\sum_{i=1}^n g_{ij}\Big)^2  \\
&\qquad \qquad + \Big(n \sum_{i=1}^n s_i^2 \Big)^{1/2}  + \frac1{8d\signalsubexp} \Big( \sum_{i=1}^n s_i \Big)^2 + Y^2 + \QUENCHEDDISORDER^2 + \sqrt{2\pi/\beta_0} d \,, 
\end{aligned} \\
&\cC := \{ B^{\star} \le 3K_{\upref{lem:posteriorgaussian_Bstar_control}}(\beta_0, \lambda)\, d \}\,, \label{eq:posteriorgaussian_C_set}
\end{align}
where $K_{\upref{lem:posteriorgaussian_Bstar_control}}(\beta_0, \lambda)$ is from Lemma \upref{lem:posteriorgaussian_Bstar_control}. Note $B^{\star}$ is a convex function of the disorder and therefore $\cC$ is a convex subset of $\R^d$. Also recall $\signalsubexp$ depends solely on $\lambda$.

Next, let $\gamma'$ be a probability measure defined on $\cC$ with density proportional to that of $\gamma$. By Assumption \upref{ass:disorder_technical_assumption}, $\gamma'$ is $\min\{1, \SIGNALGAUSSIAN\}$ strongly log-concave. 

Now, let us differentiate $\frac1d \log Z_{\beta}(I)$ w.r.t. the disorder, letting $\langle \cdot \rangle$ below and in what follows here in Appendix \upref{subsec:posteriorgaussianconversion} denote the Gibbs average w.r.t $\one\{\theta_1/\sqrt{d} \in I\} \exp\big( \sum_{i=1}^n u(S_i) - \beta \| \theta \|^2 \big)$:
\begin{align*}
\frac{\partial}{\partial s_i} \Big( \frac1d \log Z_{\beta}(I) \Big) &= \frac1d \Bigg( \displaystyle\int_{\R^d} \one\{\theta_1/\sqrt{d} \in I\} \exp\Big( \sum_{i=1}^n u(S_i) - \beta \| \theta \|^2 \Big) \cdot u'(S_i) \cdot \theta_1/\sqrt{d}\, \rmd \theta \Bigg) / Z_{\beta}(I) \\
&= \frac1d \big\langle u'(S_i) \theta_1/\sqrt{d} \big\rangle\,, \\
\frac{\partial}{\partial g_{ij}} \Big( \frac1d \log Z_{\beta}(I) \Big) &= \frac1d \Bigg( \displaystyle\int_{\R^d} \one\{\theta_1/\sqrt{d} \in I\} \exp\Big( \sum_{i=1}^n u(S_i) - \beta \| \theta \|^2 \Big) \cdot u'(S_i) \cdot \theta_j/\sqrt{d}\, \rmd \theta \Bigg) / Z_{\beta}(I) \\
&= \frac1d \big\langle u'(S_i) \theta_j/\sqrt{d} \big\rangle\,.
\end{align*}
Therefore, letting $\grad$ denote the gradient of $\frac1d \log Z_{\beta}(I)$ w.r.t. the disorder, and defining $R_{1,1}$ as per \equpref{eq:overlap_defs_thispaper}, we have
\begin{align}\label{eq:posteriorgaussian_lipschitzbd}
\| \grad \|^2 &= \frac1{d^2} \Big( \sum_{i=1}^n \big\langle u'(S_i) \theta_1/\sqrt{d} \big\rangle^2 + \sum_{i=1}^n \sum_{j=1}^d \big\langle u'(S_i) \theta_j/\sqrt{d} \big\rangle \Big) \notag \\
&\le \frac1{d^2} \Big( \sum_{i=1}^n \big\langle u'(S_i)^2 \theta_1^2/ d \big\rangle + \sum_{i=1}^n \sum_{j=1}^d \big\langle u'(S_i)^2 \theta_j^2/d \big\rangle^2 \Big) \notag \\
&\le \frac{D^2 n}{d^2} \Big( 4 + \sum_{j=1}^d \big\langle \theta_j^2/d \big\rangle \Big) \notag \\
&\le \frac{2D^2}{d} \big( 4 + R_{1,1} \big)\,.
\end{align}
Here we used that $n \le \DENSITYBOUND d \le 2d$, and that $\big\langle \theta_1^2/d \big\rangle \le 4$ due to the indicator $\one\{\theta_1/\sqrt{d} \in I\}$ in the Gibbs average $\langle \cdot \rangle$.

\paragraph{Upper bounding $R_{1,1}$.} We now upper bound $R_{1,1}$ in terms of $B^{\star}$ defined in \equpref{eq:posteriorgaussian_Bstar_bd}. To this end we introduce the following Lemmas.
\begin{lemma}[Lemma 3.1.5, \cite{talagrand2010mean}]\label{lem:lem315}
Consider a measure $\mu$ on $\R^d$ defined as follows: for some $C \subseteq \R^d$, we have for all $B \subseteq \R^d$,
\begin{align}\label{eq:lem315_measureform}
\mu(B) = \frac{1}{Z} \int \one\{B \cap C\} \E_{W}\Big[ \exp\Big( U(\btheta) - \GAUSSIANMEASURE \|\btheta\|^2 + \sum_{1 \le j \le d} a_j \btheta_j + a_0 \Big) \Big]\, \rmd\btheta\,,    
\end{align}
where $U \le 0$, $Z$ here denotes the appropriate normalizing constant, and $a_0$ potentially depends on the annealed disorder $W$. Then
\begin{align}\label{eq:lem315_conclusion}
\int \exp\Big( \frac\GAUSSIANMEASURE2 \| \btheta\|^2\Big)\, \rmd\mu(\btheta) \le \frac{1}{Z} \Big(\frac{2\pi}{\GAUSSIANMEASURE}\Big)^{d/2} \exp\Big( \frac1{2\GAUSSIANMEASURE} \sum_{1 \le j \le d} a_j^2 \Big) \E_W\big[ \exp a_0 \big]\,.
\end{align}
\end{lemma}
Technically the above Lemma is stated in \cite{talagrand2010mean} for $C = \R^d$ and with $a_0=0$, but the proof is identical. Note we will need this Lemma with $a_0$ depending on the annealed disorder $W$ from the interpolation argument in Appendix \upref{subsec:logconcaveconcentration}.

We next find a way to upper bound the moments w.r.t. a measure in the form of Lemma \upref{lem:lem315}, provided that the normalizing constant is not too small.
\begin{lemma}\label{lem:logconcaveboundmomentslemma}
Consider a measure $\mu$ in the form \equpref{eq:lem315_measureform} given in Lemma \upref{lem:lem315}, let $a_j$ be defined as in Lemma \upref{lem:lem315}, and let $Z$ be the corresponding normalizing constant. Suppose $Z \ge \exp\big( -K' B)$ for some $K'$ depending on $D, \beta_0, \beta_1, \SIGNALSHORTDEPENDENCE$ and some $B \ge d$. Then we have:
\begin{enumerate}
\item We have
\[ \Big\langle \exp\Big( \frac{\GAUSSIANMEASURE}2 \| \bar\btheta \|^2 \Big) \Big\rangle \le \Big\langle \exp\Big( \frac{\GAUSSIANMEASURE}2 \| \btheta \|^2 \Big) \Big\rangle \le \exp\Big( K_{\upref{lem:logconcaveboundmomentslemma}} B + \frac1{2\beta_0} \sum_{1 \le j \le d} a_j^2 + \log \E_W\big[ \exp a_0 \big] \Big)\,,\]
where $K_{\upref{lem:logconcaveboundmomentslemma}} \ge 1$ depends on $K', D, \GAUSSIANMEASURE_0, \GAUSSIANMEASURE_1, \SIGNALSHORTDEPENDENCE$.
\item For $k \le 4d$, we have 
\[ \big\langle \| \bar\btheta \|^{2k}\big\rangle \le \big\langle \| \btheta \|^{2k}\big\rangle \le \Big( K B + \frac1{2\beta_0} \sum_{1 \le j \le d} a_j^2 + \log \E_W\big[ \exp a_0 \big]\Big)^k\,,\]
where $K$ depends on $K', D, \beta_0, \beta_1, \SIGNALSHORTDEPENDENCE$.
\end{enumerate}
\end{lemma}
\begin{proof}[Proof of Lemma \upref{lem:logconcaveboundmomentslemma}]
For the first part of this Lemma, we directly combine Lemma \upref{lem:lem315} with the given hypothesis $Z' \ge \exp\big( -K'B)$. For the second part of this Lemma, apply Lemma \upref{lem:mgftomomentslemma} with $X = \frac{\GAUSSIANMEASURE}2 \| \btheta \|^2$ and expectation w.r.t. $\langle \cdot \rangle$. By the first part of this Lemma, we have 
\[ \log \E \big[ \exp X \big] \le K_{\upref{lem:logconcaveboundmomentslemma}} B + \frac1{2\beta_0} \sum_{1 \le j \le d} a_j^2 + \log \E_W\big[ \exp a_0 \big]\,, \]
where expectation here is w.r.t. $\langle \cdot \rangle$. The proof is complete upon noting $k \le 4d \le 4B$ and multiplying by $(2/\GAUSSIANMEASURE)^k$ on both sides.
\end{proof}

Now, note that $u(x) \ge -D(1+|x|)$ where $u(x)$ is from \equpref{eq:defuSgmm}, \equpref{eq:defuSlogistic} in the GMM and logistic cases respectively. Thus Lemma \upref{lem:lem316} implies that for $d \ge d(\eps'')$, we have 
\begin{align}\label{eq:posteriorgaussian_partitionfunc}
Z_{\beta}(I) \ge \exp\big( -K(D, \beta_0, \beta_1, \lambda) B^{\star} \big)\,,
\end{align}
where $B^{\star}$ is given in \equpref{eq:posteriorgaussian_Bstar_bd}. (The proof of Lemma \upref{lem:lem316} is deferred to Appendix \upref{subsec:logconcaveconcentration} for simplicity of notation, since the Lemma is stated for the interpolating Gibbs measures from Appendix \upref{sec:interpolationfreeenergy}, of which the measure considered here is a special case with $t=v=1$.) 

We now show in both the GMM and logistic cases that
\begin{align}\label{eq:posteriorgaussian_R11bd}
R_{1,1} \le \frac{K(D, \beta_0, \beta_1, \lambda)\,B^{\star} }d\,.
\end{align}
\begin{itemize}
    \item In the logistic case where $u(x)$ is from \equpref{eq:defuSlogistic}, we have $\varphi \le 1$ and therefore $u \le 0$. Hence by \equpref{eq:posteriorgaussian_partitionfunc}, we may apply Lemma \upref{lem:logconcaveboundmomentslemma} with the $a_j=0$, which gives \equpref{eq:posteriorgaussian_R11bd}.
    \item In the GMM case where $u(x)$ is from \equpref{eq:defuSgmm}, as $u(x) = \sqrt{\lambda} x$, we can write
\begin{align}\label{eq:posteriorgaussian_GMM_aj}
\sum_{i=1}^n u(S_i) = \sum_{i=1}^n \sqrt{\frac{\lambda}d } \Big( s_i \theta_1 + \sum_{j=2}^d g_{ij} \theta_j \Big) = \sqrt{\frac{\lambda}d} \Big(\sum_{i=1}^n s_i \Big) \theta_1 + \sum_{j=2}^d \sqrt{\frac{\lambda}d} \Big( \sum_{i=1}^n g_{ij} \Big) \theta_j\,.
\end{align}
Thus applying Lemma \upref{lem:logconcaveboundmomentslemma} with \equpref{eq:posteriorgaussian_partitionfunc}, we obtain that 
\begin{align*}
R_{1,1} \le \frac{1}{d} \Bigg( K(D, \beta_0, \beta_1, \lambda) B^{\star} + \frac{\lambda}{2\beta_0 d} \Big( \sum_{i=1}^n s_i \Big)^2 + \frac{\lambda}{2\beta_0 d} \sum_{j=2}^d \Big( \sum_{i=1}^n g_{ij} \Big)^2 \Bigg) \le \frac{K(D, \beta_0, \beta_1, \lambda)\,B^{\star} }d\,,
\end{align*}
where the last step follows from the definition of $B^{\star}$ in \equpref{eq:posteriorgaussian_Bstar_bd}. Again this proves \equpref{eq:posteriorgaussian_R11bd}.
\end{itemize}

\paragraph{Completing the argument.} Now we apply Theorem \upref{thm:logconcaveconcentrationbasics} with $f=\frac1d \log Z_{\beta}(I)$ and the strongly log-concave measure $\gamma'$. By \equpref{eq:posteriorgaussian_lipschitzbd} and \equpref{eq:posteriorgaussian_R11bd}, $f$ is Lipschitz with constant $K(D, \beta_0, \beta_1, \lambda) / \sqrt{d}$ on $\cC$.
Since $\gamma'$ is supported on $\cC$ and $\cC$ is convex, Theorem \upref{thm:logconcaveconcentrationbasics} yields for $d \ge d(\eps'')$ that
\begin{align}\label{eq:posteriorgaussian_logconcaveconcentration}
\P_{\gamma'}\Big( \Big| \frac1d \log Z_{\beta}(I) - \frac1d \E_{\gamma'}\big[ \log Z_{\beta}(I) \big] \Big| \ge \eps/2 \Big) \le \exp\big( -d / K(\eps, D, \beta_0, \beta_1, \lambda)\big)\,,
\end{align}
where $\eps>0$ is as in the statement of Theorem \upref{thm:posteriorgaussianmeasurersformula}.

We next state the following Lemma that essentially lets us reduce the desired tail bound in Theorem \upref{thm:posteriorgaussianmeasurersformula} to \equpref{eq:posteriorgaussian_logconcaveconcentration}. Its proof will be presented later.
\begin{lemma}\label{lem:posteriorgaussian_controlmeandiff}
For $d \ge d(\eps, \eps'', D, \beta_0, \beta_1, \lambda)$, we have
\begin{align}
\P_{\gamma}\big( \cC^c \big) &\le \exp(-d)\,, \label{eq:posteriorgaussian_complementprobability} \\
\Big| \frac1d \E_{\gamma'}\big[ \log Z_{\beta}(I) \big] - \frac1d \E_{\gamma}\big[ \log Z_{\beta}(I) \big] \Big| &\le \eps/4\,.\label{eq:posteriorgaussian_meandiffcontrol}
\end{align}
\end{lemma}

We will now prove Theorem \upref{thm:posteriorgaussianmeasurersformula} given Lemma \upref{lem:posteriorgaussian_controlmeandiff}. 
Let
\begin{align}
E_{\eps}' &:= \Big\{ \Big| \frac1d \log Z_{\beta}(I) - \frac1d \E_{\gamma'}\big[ \log Z_{\beta}(I) \big] \Big| \ge \eps/4 \Big\}\,. \label{eq:posteriorgaussian_ezeventdef}
\end{align}
By Lemma \upref{lem:posteriorgaussian_controlmeandiff}, for $d \ge d(\eps, \eps'', D, \beta_0, \beta_1, \lambda)$,
\begin{align}
&\P_{\gamma}\Big( \Big| \frac1d \log Z_{\beta}(I) - \frac1d \E_{\gamma}\big[ \log Z_{\beta}(I) \big] \Big| \ge \eps/2 \Big) \notag \\
&\quad= \P_{\gamma}\Big( \Big| \frac1d \log Z_{\beta}(I) - \frac1d \E_{\gamma'}\big[ \log Z_{\beta}(I) \big] + \frac1d \E_{\gamma'}\big[ \log Z_{\beta}(I) \big] - \frac1d \E_{\gamma}\big[ \log Z_{\beta}(I) \big] \Big| \ge \eps/2 \Big) \notag \\
&\quad \le \P_{\gamma}\Big( \Big| \frac1d \log Z_{\beta}(I) - \frac1d \E_{\gamma'}\big[ \log Z_{\beta}(I) \big] \Big| \ge \eps/4 \Big) \notag \\
&\quad = \E_{\gamma}\Big[ \one\{E'_{\eps} \} \Big] \notag \\
&\quad = \E_{\gamma}\Big[ \one\{ E'_{\eps}\} \one\{\cC\} \Big] + \E_{\gamma}\Big[ \one\{ E'_{\eps} \} \one\{\cC^c\} \Big] \notag \\
&\quad \le \gamma(\cC) \E_{\gamma'}\Big[ \one\{ E'_{\eps}\} \one\{\cC\} \Big] + \E_{\gamma}\Big[ \one\{\cC^c\} \Big] \notag \\
&\quad \le \exp\big( -d / K(\eps, D, \beta_0, \beta_1, \lambda)\big)\,. \label{eq:posteriorgaussian_concentrationfinish}
\end{align}
Here in the last inequality, we applied \equpref{eq:posteriorgaussian_logconcaveconcentration}. 

Finally, recall by Theorem \upref{thm:posterior_mollifiedhamiltonianfreeenergy} that for $d \ge d(\eps, \eps'', D, \beta_0, \beta_1, \lambda)$, we have 
\begin{align*}
\Big| \frac1d \E_{\gamma}\big[ \log Z_{\beta}(I) \big] - \RSG_I(n/d, \GAUSSIANMEASURE, x) \Big| \le \eps/2\,.
\end{align*}
Therefore combining with \equpref{eq:posteriorgaussian_concentrationfinish} gives 
\begin{align*}
\P_{\gamma}\Big( \Big| \frac1d \log Z_{\beta}(I) - \RSG_I(n/d, \GAUSSIANMEASURE, x) \big] \Big| \ge \eps \Big) &\le \P_{\gamma}\Big( \Big| \frac1d \log Z_{\beta}(I) - \frac1d \E_{\gamma}\big[ \log Z_{\beta}(I) \big] \Big| \ge \eps/2 \Big) \\
&\le \exp\big( -d / K(\eps, D, \beta_0, \beta_1, \lambda)\big)\,. 
\end{align*}
This proves Theorem \upref{thm:posteriorgaussianmeasurersformula} given Lemma \upref{lem:posteriorgaussian_controlmeandiff}.

\paragraph{Proof of Lemma \upref{lem:posteriorgaussian_controlmeandiff}.} Finally, we prove Lemma \upref{lem:posteriorgaussian_controlmeandiff}. We first introduce the following Lemma which gives us a tail bound on $B^{\star}$. This Lemma follows as a corollary of the proof of Lemma \upref{lem:lem310}, which is proven later in Appendix \upref{subsec:logconcaveconcentration}.
\begin{lemma}\label{lem:posteriorgaussian_Bstar_control}
For all $k \le d$ we have for $K_{\upref{lem:posteriorgaussian_Bstar_control}}$ depending on $\beta_0, \SIGNALMEDIUMDEPENDENCE$,
\begin{align*}
\E\big[ B^{\star k} \big] \le \big( K_{\upref{lem:posteriorgaussian_Bstar_control}} d \big)^k\,.
\end{align*}
\end{lemma}
Given Lemma \upref{lem:posteriorgaussian_Bstar_control}, by the definition of $\cC$, we directly obtain \equpref{eq:posteriorgaussian_complementprobability}. 
To prove \equpref{eq:posteriorgaussian_meandiffcontrol}, we introduce an additional Lemma:
\begin{lemma}\label{lem:posteriorgaussian_controlpartitionfunc}
For $d \ge d(\eps'')$, we have
\begin{align*}
\big|\log Z_{\beta}(I) \big| \le K(D, \beta_0, \beta_1, \lambda) B^{\star}\,.
\end{align*}
\end{lemma}
\begin{proof}[Proof of Lemma \upref{lem:posteriorgaussian_controlpartitionfunc}]
The lower bound on $\log Z_{\beta}(I)$ follows from \upref{eq:posteriorgaussian_partitionfunc}. For the upper bound, note 
\begin{align*}
Z_{\beta}(I) \le \int_{\R^d} \exp\Big( -\beta\| \theta\|^2 + \sum_{1 \le j \le d} a_j \theta_j \Big)\, \rmd\theta \le (2\pi / \beta_0)^{d/2} \,,
\end{align*}
for some $(a_j)_{j=1}^d \in \R$. (Note the integral is invariant to coordinate translation.) Specifically, in the logistic case we use that $u \le 0$, and in the GMM case we use \equpref{eq:posteriorgaussian_GMM_aj}. This proves Lemma \upref{lem:posteriorgaussian_controlpartitionfunc}.
\end{proof}

Now to prove \equpref{eq:posteriorgaussian_meandiffcontrol}, we write
\begin{align*}
\frac1d \E_{\gamma}\big[ \log Z_{\beta}(I) \big] &= \frac1d \E_{\gamma}\big[ \log Z_{\beta}(I) \one\{\cC\} \big] + \frac1d \E_{\gamma}\big[ \log Z_{\beta}(I) \one\{\cC^c\} \big] \\
&= \frac1d \gamma(\cC) \E_{\gamma'}\big[ \log Z_{\beta}(I) \big] + \frac1d \E_{\gamma}\big[ \log Z_{\beta}(I) \one\{\cC^c\} \big]\,,
\end{align*}
where we note $\E_{\gamma'}\big[ \log Z_{\beta}(I) \one\{\cC\} \big] = \E_{\gamma'}\big[ \log Z_{\beta}(I) \big]$ as the support of $\gamma'$ equals $\cC$.
Therefore
\begin{align}
&\frac1d \E_{\gamma'}\big[ \log Z_{\beta}(I) \big] - \frac1d \E_{\gamma}\big[ \log Z_{\beta}(I) \big] \notag \\
&\quad = \frac1d \E_{\gamma'}\big[ \log Z_{\beta}(I) \big] - \Big( \frac1d \E_{\gamma}\big[ \log Z_{\beta}(I) \one\{\cC\} \big] + \frac1d \E_{\gamma}\big[ \log Z_{\beta}(I) \one\{\cC^c\} \big] \Big) \notag \\
&\quad = \frac{1-\gamma(\cC)}d \E_{\gamma'}\big[ \log Z_{\beta}(I) \big] - \frac1d \E_{\gamma}\big[ \log Z_{\beta}(I) \one\{\cC^c\} \big]\,.\label{eq:posteriorgaussian_rewritediff}
\end{align}
By Lemma \upref{lem:posteriorgaussian_controlpartitionfunc}, and since $\gamma'$ is supported on $\cC$ on which $B^{\star} \le K(D, \beta_0, \beta_1, \lambda) d$, we have 
\begin{align}\label{eq:posteriorgaussian_bounddiff1}
\Big| \frac1d \E_{\gamma'}\big[ \log Z_{\beta}(I) \big]\Big| \le \frac1d \E_{\gamma'}\big[ \big|\log Z_{\beta}(I) \big| \big] \le K(D, \beta_0, \beta_1, \lambda)\,.
\end{align}
Next by Lemma \upref{lem:posteriorgaussian_Bstar_control}, note
\begin{align}\label{eq:posteriorgaussian_tailbd}
\P_{\gamma}\big(B^{\star} \ge t\big) \le \big(K_{\upref{lem:posteriorgaussian_Bstar_control}}(D, \beta_0, \beta_1, \lambda) / t \big)^d \text{ for all }t \ge 0\,.
\end{align}
Thus by Lemma \upref{lem:posteriorgaussian_controlpartitionfunc} and integrating the tail bound \equpref{eq:posteriorgaussian_tailbd}, we have
\begin{align}
\Big| \E_{\gamma}\big[ \log Z_{\beta}(I) \one\{\cC^c\} \big] \Big| &\le K(D, \beta_0, \beta_1, \lambda) \E_{\gamma}\big[ B^{\star} \one\{\cC^c\} \big] \notag \\
&\le K(D, \beta_0, \beta_1, \lambda) \int_{ 3K_{\upref{lem:posteriorgaussian_Bstar_control}}(D, \beta_0, \beta_1, \lambda) d }^\infty \P_{\gamma}\big(B^{\star} \ge t\big)\, \rmd t \notag \\
&\le K(D, \beta_0, \beta_1, \lambda) \exp(-d)\,.\label{eq:posteriorgaussian_boundwtailbd}
\end{align}
Combining \equpref{eq:posteriorgaussian_complementprobability}, \equpref{eq:posteriorgaussian_rewritediff}, \equpref{eq:posteriorgaussian_bounddiff1}, \equpref{eq:posteriorgaussian_boundwtailbd} now proves \equpref{eq:posteriorgaussian_meandiffcontrol} for $d \ge d(\eps, \eps'', D, \beta_0, \beta_1, \lambda)$, finishing the proof of Lemma \upref{lem:posteriorgaussian_controlmeandiff}.

\section{Conversion from Gaussian to spherical measure: proof of Theorems \ref{thm:mainresult}, \ref{thm:mainresultposterior}}\label{sec:sphericalconversion}
Here we complete the proof of our main results Theorems \upref{thm:mainresult}, \upref{thm:mainresultposterior}.
Theorem \upref{thm:mainresult} on interpolators and Theorem \upref{thm:mainresultposterior} on the posterior follow as a corollary of the following Theorems \upref{thm:sphericalmeasurersformula}, \upref{thm:posteriorrsformula} respectively.

\paragraph{Theorem \upref{thm:sphericalmeasurersformula} for interpolators:} This result states that for a narrow interval $I$ near $x \in [-1,1]$ (whose width is independent of $d$), the normalized logarithm of the spherical measure of the interpolators constrained to this interval, approximately equals $\inf_{q \in [0,1)}\Phi(x, q)$ for $d$ large enough:
\begin{theorem}[For interpolators]\label{thm:sphericalmeasurersformula}
Let $\mu_d(\cdot)$ denote the uniform measure on $\cS^{d-1}(\sqrt{d})$.
Consider any $\alpha \le \DENSITYBOUND$ and $\MARGIN \in \R$, and let $n = \lfloor \alpha d \rfloor$. 
Define $\delta=\delta(\alpha, \lambda, \MARGIN_+)$ as follows, where the inequality uses $\alpha \le \DENSITYBOUND$:
\begin{align}
\delta=\delta(\alpha, \lambda, \MARGIN_+) := \frac{1/2 + \alpha K_{\upref{lem:controlexpectationofAW}} \Csol^{1/2}/\sqrt{2}}{9} \le \frac1{10} \,.\label{eq:deltaalphafulldef}
\end{align}
Recalling the definition of $C$ in \equpref{eq:hyperplane_intersection_def} and letting $\Phi(x, q)$ be as in \equpref{eq:phi_master} with $\exp u(x) = \one\{ x \ge \MARGIN\}$, we have for all $x \in [-1+\delta, 1-\delta]$ and $\eps>0$:
\begin{enumerate}
\item Upper bound: 
Let $\eps_{\RSG} = \eps_{\RSG}(\eps, \DENSITYBOUND, \beta_0, \beta_1, \MARGIN, \lambda)$ be from Theorem \upref{thm:gaussianmeasurersformula} and consider any $0 < \eps'' \le \min\{ \eps, \eps_{\RSG} \}$. Define the interval
\begin{align}\label{eq:sphericalconversion_I_upperbound}
I := \begin{cases} \big[\frac{x-\eps''}{1-\eps''}, x+\eps'' \big] &: x \ge 0\,,\, \MARGIN < 0 \,, \\
\big[x-\eps'', \frac{x+\eps''}{1-\eps''} \big] &: x<0 \,,\, \MARGIN < 0 \,, \\
\big[x-\eps'', \frac{x+\eps''}{1+\eps''} \big] &: x \ge 0 \,,\, \MARGIN \ge 0 \,, \\
\big[\frac{x-\eps''}{1+\eps''}, x+\eps'' \big] &: x<0 \,,\, \MARGIN \ge 0 \,.\end{cases}
\end{align}
If $|x| > \eps''$, for $d$ large enough in terms of $\eps, \eps'', \alpha, \beta_0, \beta_1, \MARGIN, \lambda$ and with probability at least $1-\exp\bigl(-d / K \bigr)$ where $K$ depends on $\eps, \alpha, \beta_0, \beta_1, \MARGIN, \lambda$, we have
\begin{align}\label{eq:spherical_conversion_upperbound}
\frac1d \log \mu_d \Bigl(\one\{\theta_1/\sqrt{\DIMENSION} \in I\}  \cap C \Bigr) \le \inf_{q \in [0,1)}\Phi(x, q) + \eps\, .
\end{align}
(Note this bound is non-vacuous; we have $I \neq \emptyset$.)

\item Lower bound: Suppose $\eps < \frac1{K'}$ where $K'$ is large enough in terms of $\alpha, \beta_0, \beta_1, \MARGIN, \lambda$. Then there exists $\eps_{\RS}>0$ depending on $\eps, \alpha, \beta_0, \beta_1, \MARGIN, \lambda$ such that for all $0 < \eps'' \le \eps_{\RS}$, the following holds. 
Define the interval
\begin{align}\label{eq:sphericalconversion_I_lowerbound}
I := \begin{cases} \big[\frac{x-\eps''}{1+\eps}, x+\eps''\big] &: x \ge 0\,,\, \MARGIN < 0 \,, \\
\big[x-\eps'', \frac{x+\eps''}{1+\eps} \big] &: x < 0\,,\, \MARGIN < 0 \,, \\
\big[x-\eps'', \frac{x+\eps''}{1-\eps}\big] &: x \ge 0\,,\, \MARGIN \ge 0 \,, \\
\big[\frac{x-\eps''}{1-\eps}, x+\eps''\big] &: x < 0\,,\, \MARGIN \ge 0 \,.\end{cases}
\end{align}
If $|x| > \eps''$, for $d$ large enough in terms of $\eps, \eps'', \alpha, \beta_0, \beta_1, \MARGIN, \lambda$ and with probability at least $1-\exp\bigl(-d/K \bigr)$ where $K$ depends on $\eps, \alpha, \beta_0, \beta_1, \MARGIN, \lambda$, we have
\begin{align}\label{eq:spherical_conversion_lowerbound}
\frac1d \log \mu_d \Bigl(\one\{\theta_1/\sqrt{\DIMENSION} \in I\}  \cap C\Bigr) \ge \inf_{q \in [0,1)}\Phi(x, q) - \eps\, .
\end{align}
\end{enumerate}
\end{theorem}

\paragraph{Theorem \upref{thm:posteriorrsformula} for posterior:} This result is similar to Theorem \upref{thm:sphericalmeasurersformula}, but now is for the normalized logarithm of the integral of the posterior constrained to the interval $I$, rather than the spherical measure of the interpolators:
\begin{theorem}[For posterior]\label{thm:posteriorrsformula}
Consider any $\alpha \le \DENSITYBOUND$ and let $n = \lfloor \alpha d \rfloor$. 
Consider $u(x)$ from \equpref{eq:defuSgmm}, \equpref{eq:defuSlogistic} in the GMM and logistic cases respectively and let $\Phi(x, q)$ be as in \equpref{eq:phi_master} with this $u(x)$. Define $\delta=\delta(\alpha, \lambda)$ as follows, where the inequality uses $\alpha \le \alpha_0$:
\begin{align}\label{eq:posteriordeltadef}
\delta = \delta(\alpha, \lambda) := \frac1{18 - \alpha (2D^2+D)} \le \frac1{10}\,.
\end{align}
For any interval $I \subseteq \R$, let
\begin{align}\label{eq:posterior_integral_def}
Z(I) := \int_{\cS^{d-1}(\sqrt{d})} \one\big\{  \theta_1/\sqrt{\DIMENSION} \in I \big\} \exp\Big(\sum_{i=1}^n u( S_i ) \Big)\, \rmd \theta\,,
\end{align}
where $S_i=S_i(\theta)$ is as in \equpref{eq:Si_def}. Then for all $x \in [-1+\delta, 1-\delta]$ and $\eps>0$, we have:
\begin{enumerate}
\item Upper bound: 
Let $\eps_{\RSG} = \eps_{\RSG}(\eps, \DENSITYBOUND, \beta_0, \beta_1, \lambda)$ be from Theorem \upref{thm:posteriorgaussianmeasurersformula}, consider any $0 < \eps'' \le \min\{ \eps / K', \eps_{\RSG} \}$ where $K'$ depends on $\lambda$. Define the interval $I$ by 
\begin{align}\label{eq:posterior_sphericalconversion_I_upperbound}
I := \begin{cases} 
\big[x-\eps'', \frac{x+\eps''}{1+\eps''} \big] &: x \ge 0 \,, \\
\big[\frac{x-\eps''}{1+\eps''}, x+\eps'' \big] &: x<0 \,.\end{cases}
\end{align}
If $|x| > \eps''$, for $d$ large enough in terms of $\eps, \eps'', \alpha, \beta_0, \beta_1, \lambda$ and with probability at least $1-\exp\bigl(-d / K \bigr)$ where $K$ depends on $\eps, \alpha, \beta_0, \beta_1, \lambda$, we have
\begin{align}\label{eq:posterior_spherical_conversion_upperbound}
\frac1d \log Z(I) \le \inf_{q \in [0,1)}\Phi(x, q) + \frac12 \log(2 \pi e) + \eps\, .
\end{align}
(As in Theorem \upref{thm:sphericalmeasurersformula}, note this bound is non-vacuous; we have $I \neq \emptyset$.)

\item Lower bound: Suppose $\eps < \frac1{K'}$ where $K'$ is large enough in terms of $\alpha, \beta_0, \beta_1, \lambda$. Then there exists $\eps_{\RS}>0$ depending on $\eps, \alpha, \beta_0, \beta_1, \lambda$ such that for all $0 < \eps'' \le \eps_{\RS}$, the following holds. 
Define the interval
\begin{align}\label{eq:posterior_sphericalconversion_I_lowerbound}
I := \begin{cases} 
\big[x-\eps'', \frac{x+\eps''}{1-\eps}\big] &: x \ge 0\,, \\
\big[\frac{x-\eps''}{1-\eps}, x+\eps''\big] &: x < 0 \,.\end{cases}
\end{align}
If $|x| > \eps''$, for $d$ large enough in terms of $\eps, \eps'', \alpha, \beta_0, \beta_1, \lambda$ and with probability at least $1-\exp\bigl(-d/K \bigr)$ where $K$ depends on $\eps, \alpha, \beta_0, \beta_1, \lambda$, we have
\begin{align}\label{eq:posterior_spherical_conversion_lowerbound}
\frac1d \log Z(I) \ge \inf_{q \in [0,1)}\Phi(x, q) + \frac12 \log(2 \pi e) - \eps\,.
\end{align}
\end{enumerate}
\end{theorem}

We will prove Theorem \upref{thm:sphericalmeasurersformula} in the $\MARGIN < 0$ case in Appendix \upref{subsec:sphericalpf_negativemargin}, prove Theorem \upref{thm:sphericalmeasurersformula} in the $\MARGIN \ge 0$ case in Appendix \upref{subsec:sphericalpf_positivemargin}, and prove Theorem \upref{thm:posteriorrsformula} in Appendix \upref{subsec:sphericalpf_posterior}. All three proofs follow the same unified strategy that leverages our work in Appendix \upref{sec:gaussianconversion}, specifically Theorem \upref{thm:gaussianmeasurersformula} for the interpolators and Theorem \upref{thm:posteriorgaussianmeasurersformula} for the posterior. We present a full proof for $\MARGIN < 0$ in Appendix \upref{subsec:sphericalpf_negativemargin} and highlight the necessary modifications to complete the proof when $\MARGIN>0$ or for the posterior in Appendices \upref{subsec:sphericalpf_positivemargin}, \upref{subsec:sphericalpf_posterior}.

\begin{remark}
We emphasize that the interval $I$ from Theorems \upref{thm:sphericalmeasurersformula}, \upref{thm:posteriorrsformula} does \emph{not} equal $[x-\eps'', x+\eps'']$. This is because our analysis relies on several comparison inequalities where we scale $[x-\eps'', x+\eps'']$ by a factor close to 1 independent of $d$; see e.g. \equpref{eq:sphericalconversioncontainment_larger1}, \equpref{eq:sphericalconversioncontainment_smaller1}.
$I$ is then a larger interval containing the scaled intervals for the lower bound, or a smaller interval contained within the scaled intervals for the upper bound. 
For the upper bound, the parameter $w$ governing the range of these scalings equals $1 \pm \eps''$ (the sign depending on whether $\MARGIN \ge 0, \MARGIN < 0$), and we can check that $I \neq \emptyset$. For the lower bound, $\eps''$ delicately depends on the parameter $w=1 \pm \eps$ governing these scalings, but now $I \supseteq [x-\eps'', x+\eps'']$ and is therefore nonempty.
\end{remark}

Before we present the proofs of Theorems \upref{thm:sphericalmeasurersformula}, \upref{thm:posteriorrsformula}, we combine them with an application of the Laplace Method over $x \in [-1,1]$ to prove our main results: Theorem \upref{thm:mainresult} on interpolators and Theorem \upref{thm:mainresultposterior} on the posterior respectively. Note the Laplace Method gives the `$\sup_x \inf_q$' structure therein.
In preparation, we will present the following properties of the RS equations, which follow from the convex-concave structure established in Appendix \upref{sec:RSequationsunique}.
The proof of this Lemma is in Appendix \upref{subsec:spherical_conversion_RS_pfs}.
\begin{lemma}\label{lem:rs_formula_properties}
Let $\delta=\delta(\alpha, \lambda, \MARGIN_+)$ from \equpref{eq:deltaalphafulldef} for the interpolators, or $\delta=\delta(\alpha, \lambda)$ from \equpref{eq:posteriordeltadef} for the posterior. We have the following:
\begin{enumerate}
\item For all $(\alpha, \beta, x) \in [0, \DENSITYBOUND] \times [\beta_0, \beta_1] \times [-1, 1]$, we have $\frac{\partial \rho_0(\alpha, \GAUSSIANMEASURE, x)}{\partial \GAUSSIANMEASURE} < 0$.

\item For all $(\alpha, x) \in [0, \DENSITYBOUND] \times [-1+\delta, 1-\delta]$, there is a unique $\GAUSSIANMEASURE=\GAUSSIANMEASURE(\alpha, x) \in [\GAUSSIANMEASURE_0, \GAUSSIANMEASURE_1]$ such that $\rho_0(\alpha, \GAUSSIANMEASURE, x)=1-x^2$, and moreover $\GAUSSIANMEASURE(\alpha, x) \in (1/4, 9)$. 

\item For all $(\alpha, x) \in [0, \DENSITYBOUND] \times [-1+\delta, 1-\delta]$, recalling the definition of $\bar\Phi$ in \equpref{eq:phi_master_rho},
\begin{align*}
f(\GAUSSIANMEASURE) := \GAUSSIANMEASURE + \bar\Phi\big(x, q_0(\alpha, \beta, x), \rho_0(\alpha, \beta, x)\big) - \beta \big(\rho_0(\alpha, \beta, x) + x^2\big)
\end{align*}
attains a unique minimum on $\GAUSSIANMEASURE \in [\GAUSSIANMEASURE_0, \GAUSSIANMEASURE_1]$ at $\GAUSSIANMEASURE=\GAUSSIANMEASURE(\alpha, x)$. Thus for any interval $I$, if $0 \not\in I$, then $\GAUSSIANMEASURE + \RSG_{I}(\alpha, \GAUSSIANMEASURE, x)$ attains a unique minimum on $\GAUSSIANMEASURE \in [\GAUSSIANMEASURE_0, \GAUSSIANMEASURE_1]$ at $\GAUSSIANMEASURE=\GAUSSIANMEASURE(\alpha, x)$.

\item For all $(\alpha, x) \in [0, \DENSITYBOUND] \times [-1+\delta, 1-\delta]$, we have 
\begin{align*}
\inf_{q \in [0,1)}\Phi(x, q) = f\big(\beta(\alpha, x)\big)\,.
\end{align*}

\item For all $(\alpha, x) \in [0, \DENSITYBOUND] \times [-1+\delta, 1-\delta]$, $\inf_{q \in [0,1)}\Phi(x, q)$ is continuous in $(\alpha, x)$ (note $\alpha$ is implicitly an argument in $\Phi(x, q)$).
\end{enumerate}
\end{lemma}

We now prove Theorems \upref{thm:mainresult}, \upref{thm:mainresultposterior} using Theorems \upref{thm:sphericalmeasurersformula}, \upref{thm:posteriorrsformula} and the above Lemma \upref{lem:rs_formula_properties}.

\subsection{Proof of Theorems \upref{thm:mainresult}, \upref{thm:mainresultposterior}}
\label{sec:proof_main_results}
We provide the proof of Theorem \upref{thm:mainresult} using Theorem \upref{thm:sphericalmeasurersformula} when $\MARGIN < 0$; the proof of Theorem \upref{thm:mainresult} in the $\MARGIN \ge 0$ case and the proof of Theorem \upref{thm:mainresultposterior} using Theorem \upref{thm:posteriorrsformula} are analogous.
Recall that $y_i X_i \equiv (s_i, g_i)^T$ in law as random variables and that without loss of generality we can set $\theta_\star = \sqrt{d} e_1$.

First, we show \equpref{eq:spherical_conversion_upperbound}, \equpref{eq:spherical_conversion_lowerbound} with the interval $I$ replaced by any interval $I \subseteq [-1+\delta, 1-\delta]$. Recall as per \equpref{eq:csoldef} that $\beta_0, \beta_1$ are universal constants, therefore in this proof we will not explicitly write dependence on them.

\paragraph{Proof of \equpref{eq:spherical_conversion_upperbound}:} Note for $d \ge d(\eps)$, $\frac1d \log \mu_d \bigl(\one\{\theta_1/\sqrt{\DIMENSION} \in [-1/\sqrt{d}, 1/\sqrt{d}] \}\bigr) \le \eps$. 
Now if $0 \in I$, we subdivide $I$ into the two closed intervals $I_1 = I \cap [-1, -1/\sqrt{d}]$, $I_2 = I \cap [1/\sqrt{d}, 1]$, and the remaining part $I-I_1-I_2$ (if $0 \not\in I$ this step is not necessary, and we can just apply the same argument as in the following upper bound for $I_1$).

We now upper bound $\frac1d \log \mu_d \bigl(\one\{\theta_1/\sqrt{\DIMENSION} \in I_1 \} \cap C\bigr)$, the upper bound for $I_2$ being analogous. 
We divide the closed interval $I_1$ into a set of consecutive closed intervals $\bar I_k$ each of width $\eps_{\RSG}=\eps_{\RSG}(\eps, \DENSITYBOUND, \MARGIN, \lambda)$ (recall $\beta_0, \beta_1$ are universal constants), $1 \le k \le K(\eps_{\RSG})$, that overlap at exactly their endpoints. 
Here $\eps_{\RSG}$ comes from Theorem \upref{thm:gaussianmeasurersformula}.

For each interval $\bar I_k=[a,b]$, note $0<a< b \le 1$, $b-a \le \eps_{\RSG}$. 
Letting $\eps'' = \frac{b-a}{2-a} < b-a$ and $x = b-\eps''$, we have $\bar I_k = [a,b] = \big[\frac{x-\eps''}{1-\eps''}, x+\eps'' \big]$, $x>\eps''$, $\eps'' \le b-a \le \eps_{\RSG}$. 
Since $\eps'' \ge b-a = \eps_{\RSG}=\eps_{\RSG}(\eps, \DENSITYBOUND, \MARGIN, \lambda)$, Theorem \upref{thm:sphericalmeasurersformula} yields that for $d \ge d(\eps, \alpha, \MARGIN, \lambda)$, with probability at least $1-\exp\bigl(- d / K(\eps, \alpha, \MARGIN, \lambda) \bigr)$,
\begin{align*}
\frac1d \log \mu_d \Bigl(\one\{\theta_1/\sqrt{\DIMENSION} \in \bar I_k\}  \cap C \Bigr) &\le \inf_{q \in [0,1)}\Phi(x, q) + \eps \le \sup_{x \in I} \inf_{q \in [0,1)}\Phi(x, q) + \eps\,.
\end{align*}
Here the second inequality uses that $\bar I_k \subseteq I_1$.
As we have $K(\eps_\RSG)$ such intervals $\bar I_k$, combining the above bound over all $1 \le k \le K(\eps_\RSG)$ and applying a Union Bound (note we sum the bounds on $\mu_d \bigl(\one\{\theta_1/\sqrt{\DIMENSION} \in I_1\}  \cap C \bigr)$, which behaves at the exponential scale) yields that for $d \ge d(\eps, \alpha, \MARGIN, \lambda)$, with probability at least $1-\exp\big(-d / K(\eps, \alpha, \MARGIN, \lambda) \bigr)$,
\begin{align*}
\frac1d \log \mu_d \Bigl(\one\{\theta_1/\sqrt{\DIMENSION} \in I_1\}  \cap C \Bigr) \le \sup_{x \in I} \inf_{q \in [0,1)}\Phi(x, q) + 2\eps\,.
\end{align*}
An analogous bound applies for $I_2$ via the same proof, and for $I-I_1-I_2$ we recall the bound $\frac1d \log \mu_d \bigl(\one\{\theta_1/\sqrt{\DIMENSION} \in [-1/\sqrt{d}, 1/\sqrt{d}] \}\bigr) \le \eps$. Combining these three upper bounds and applying a Union Bound yields for $d \ge d(\eps, \alpha, \MARGIN, \lambda)$, with probability at least $1-\exp\bigl(-d / K(\eps, \alpha, \MARGIN, \lambda) \bigr)$,
\begin{align*}
\frac1d \log \mu_d \Bigl(\one\{\theta_1/\sqrt{\DIMENSION} \in I\}  \cap C\Bigr) \le \sup_{x \in I} \inf_{q \in [0,1)}\Phi(x, q) + 4\eps\,.
\end{align*}
Upon taking $\eps \leftarrow \eps/4$, this proves \equpref{eq:spherical_conversion_upperbound} for this $I$.

\paragraph{Proof of \equpref{eq:spherical_conversion_lowerbound}:}  
In this proof, we write $I = [a,b]$ for $-1 \le a < b \le 1$.
Note it suffices to prove the result when $\eps < \frac1{K'(\alpha, \MARGIN, \lambda)}$.
First, continuity from Lemma \upref{lem:rs_formula_properties} implies $\sup_{x \in I} \inf_{q \in [0,1)}\Phi(x, q)$ is attained at some $x^\star \in I$. 
Noting $\inf_{q \in [0,1)}\Phi(x, q)$ only depends on $\alpha, \MARGIN, \lambda$, by continuity supplied by Lemma \upref{lem:rs_formula_properties} and compactness, we let $\bar\delta = \bar\delta(\eps, \alpha, \MARGIN, \lambda)$ be small enough so that
\begin{align}\label{eq:sphericalconversion_laplace_uniformcontinuity}
\big| \inf_{q \in [0,1)}\Phi(x_1, q) - \inf_{q \in [0,1)}\Phi(x_2, q) \big| \le \eps\,, \qquad \forall\, x_1,x_2 \in [-1+\delta, 1-\delta]\,,\, |x_1-x_2| \le \bar\delta\,.
\end{align}
Since $\eps \le \eps_0$, taking $\eps_0$ small enough in terms of $\alpha, \MARGIN, \lambda, |I|$, we may suppose that $|I| \ge \bar\delta$.

We first prove the result when $I \subseteq [0, \infty)$; the proof when $I \subseteq (-\infty, 0]$ is analogous, and later we will provide the argument when neither case holds. 
\begin{enumerate}
    \item Suppose $x^{\star} > 0$ is distance at least $\bar\delta/4$ from $a, b$. 
Let $\bar\eps = \bar\eps(\eps, \bar\delta, \alpha, \lambda) \le \min\{ \bar\delta/8, \eps \}$ and $\eps'' < \min\big\{ \eps_{\RS}(\bar\eps, \alpha, \MARGIN, \lambda), \bar\delta/8 \big\}$ be small enough, where $\eps_{\RS}$ comes from Theorem \upref{thm:sphericalmeasurersformula}. Let
\begin{align}\label{eq:sphericalconversion_laplace_barIdef}
\bar I := \Big[\frac{x^\star-\eps''}{1+\bar\eps}, x^\star+\eps'' \Big]\,.
\end{align}
As $|a| \le 1$, $\bar\eps + \eps'' < \bar\delta/4$, $|x^{\star}| \le 1$ and $x^\star$ is distance at least $\bar\delta/4$ from $a,b$, we have $\bar I \subseteq I$.
Also, clearly $x^\star \ge \bar\delta/4 > \eps''$.
Applying Theorem \upref{thm:sphericalmeasurersformula}, recalling that $\bar\delta$ and $\bar\eps$ depend only on $\eps, \alpha, \MARGIN, \lambda$, we obtain for $d \ge d(\bar\eps, \eps'', \alpha, \MARGIN, \lambda) = d(\eps, \alpha, \MARGIN, \lambda)$ that with probability at least $1-\exp\bigl(-d / K(\eps, \alpha, \MARGIN, \lambda) \bigr)$,
\begin{align}
\frac1d \log \mu_d \Bigl(\one\{\theta_1/\sqrt{\DIMENSION} \in I\} \cap C \Bigr) &\ge \frac1d \log \mu_d \Bigl(\one\{\theta_1/\sqrt{\DIMENSION} \in \bar I\} \cap C \Bigr) \notag \\
&\ge \inf_{q \in [0,1)}\Phi(x^\star, q) - \eps\,, \label{eq:sphericalconversion_laplace_basiclowerbd}
\end{align}
where we use that $\bar\eps \le \eps$ in the last inequality.

\item Suppose $x^\star$ is within distance $\bar\delta/4$ from $a,b$. Now, we simply apply the above argument with $x^\star$ replaced by $\bar x^\star = a+\bar\delta/4$ if $|a-x^\star| \le \bar\delta/4$ or $\bar x^\star = b-\bar\delta/4$ if $|b-x^\star| \le \bar\delta/4$. 
Since $|I| \ge \bar\delta$, the resulting $\bar x^\star$ is always distance at least $\bar\delta/4$ from $a,b$, and hence $\bar x^\star > \eps''$ as well.
Define $\bar \eps$, $\eps''$ identically to case 1) above and $\bar I$ now in terms of $\bar x^\star$ as per \equpref{eq:sphericalconversion_laplace_barIdef}. Applying Theorem \upref{thm:sphericalmeasurersformula} analogously, by $|x^\star - \bar x^\star | \le \bar\delta$ and \equpref{eq:sphericalconversion_laplace_uniformcontinuity}, we obtain for the same $d$ and with the same probability as \equpref{eq:sphericalconversion_laplace_basiclowerbd} that
\begin{align}
\frac1d \log \mu_d \Bigl(\one\{\theta_1/\sqrt{\DIMENSION} \in I\}  \cap C \Bigr) &\ge \frac1d \log \mu_d \Bigl(\one\{\theta_1/\sqrt{\DIMENSION} \in \bar I\} \cap C \Bigr) \notag \\ 
&\ge \inf_{q \in [0,1)}\Phi(\bar x^\star, q) - \bar\eps \notag \\
&\ge \inf_{q \in [0,1)}\Phi(x^\star, q) - 2\eps\,.\label{eq:sphericalconversion_laplace_lowerbd2}
\end{align}
\end{enumerate}
Finally suppose $I$ is not contained within $[0,\infty)$ or $(-\infty, 0]$. 
We suppose $x^\star \ge 0$ in what follows, the argument when $x^\star < 0$ is analogous. 
\begin{enumerate}
    \item Suppose $\big| I \cap [0, \infty) \big| \ge \bar\delta/2$.
    If $x^\star \ge \bar\delta/4$ we simply apply the argument in either cases 1) or 2) above, defining $\bar\eps, \eps'',\bar I$ the same way as there.
In particular, we still have $\bar I \subseteq I$ as $\big| I \cap [0, \infty)\big| \ge \bar\delta/2$; if we need to define $\bar x^\star$ as in case 2) above (note this only arises when $x^\star$ is within $\bar\delta/4$ of the right endpoint of $I$), since $\big| I \cap [0, \infty)\big| \ge \bar\delta/2$, the resulting $\bar x^\star$ is always distance at least $\bar\delta/4 > \eps''$ from the endpoints of $I$ and 0.

Next suppose $0 \le x^\star \le \bar\delta/4$. 
We apply the argument in case 2) above with $\bar x^\star = \bar\delta/4$. Note that defining $\bar I$ analogously, we have $\bar I \subseteq I$ as $\big| I \cap [0, \infty)\big| \ge \bar\delta/2$. Combining with the fact that $|\bar x^\star - x^\star| \le \bar\delta/4$ and \equpref{eq:sphericalconversion_laplace_uniformcontinuity}, we obtain \equpref{eq:sphericalconversion_laplace_lowerbd2}.

\item Suppose $\big| I \cap [0, \infty)\big| < \bar\delta/2$.
Therefore as $|I| \ge \bar\delta$, $\big| I \cap (-\infty, 0] \big| \ge \bar\delta/2$.
We now apply the argument in either cases 1) or 2) above now with $\bar x^{\star} = -\bar\delta/4$, defining $\bar \eps, \eps''$ analogously and now defining $\bar I := \big[ \bar x^{\star} -\eps'', \frac{\bar x^{\star}+\eps''}{1+\bar\eps} \big]$.

Since the interval $I$ is not fully contained in $[0, \infty)$ or $(-\infty, 0]$ and $\big| I \cap (-\infty, 0] \big| \ge \bar\delta/2$, we have $\bar x^\star \in I$, and furthermore that $\bar I \subseteq I \cap (-\infty, 0]$. Also, $|\bar x^\star| \ge \bar\delta/4 > \eps''$. Finally, as $x^\star \ge 0$ and $\big| I \cap [0, \infty)\big| < \bar\delta/2$, we have $|\bar x^\star - x^\star| \le 3\bar\delta/4$. Combining with \equpref{eq:sphericalconversion_laplace_uniformcontinuity} gives \equpref{eq:sphericalconversion_laplace_lowerbd2}.
\end{enumerate}
This yields the desired lower bound in all cases. Now, \equpref{eq:spherical_conversion_lowerbound} follows by taking $\eps \leftarrow \eps/2$. \newline

To complete the proof, we must upper bound $\frac1d \log \mu_d\big( \one\{|\theta_1/\sqrt{d}| \ge 1-\delta\} \cap C \big)$. As $\delta \le \frac1{10}$, it remains to show that with probability at least $1 - \exp\big(-d/K(\eps, \alpha, \MARGIN, \lambda)\big)$ for $d \ge d(\eps, \alpha, \MARGIN, \lambda)$, we have the following for a universal constant $c>0$. For the interpolators,
\begin{align}\label{eq:interpolants_finish_desired}
\frac1d \log \mu_d\Big( \one\big\{|\theta_1/\sqrt{d}| \ge 9/10 \big\} \cap C \Big) < \sup_{x \in [-1 + \delta, 1-\delta]} \inf_{q \in [0,1)} \Phi(x,q) - c\,,
\end{align}
and for the posterior, defining $Z(I)$ as per \equpref{eq:posterior_integral_def},
\begin{align}\label{eq:posterior_finish_desired}
\frac1d \log \Big( Z\big( [-1, -9/10] \big) + Z\big( [9/10, 1] \big) \Big) - \frac12 \log(2\pi e) < \sup_{x \in [-1 + \delta, 1-\delta]} \inf_{q \in [0,1)} \Phi(x,q) - c\,.
\end{align}
We now upper bound the $\LHS$ in \equpref{eq:interpolants_finish_desired}, \equpref{eq:posterior_finish_desired}.
By following the same proof as that of Lemma \upref{lem:posterior_spherical_continuity}, we have with probability at least $1 - \exp\big(-d/K(\eps, \alpha, \MARGIN, \lambda)\big)$,
\begin{align}\label{eq:posterior_spherical_finish_gmmbd}
\sup_{\theta \in \cS^{d-1}(\sqrt{d})} \sum_{i=1}^n u(S_i) \le K(\lambda) \sqrt{n d}\,.
\end{align}
Standard bounds on the surface area of a spherical cap gives
\begin{align*}
\frac1d \log \mu_d\Big( \one\{ |\theta_1/\sqrt{d}| \ge 9/10 \Big) \le -c'\,,
\end{align*}
where $c'>0$ is a universal constant.
Letting $\LHS$ denote the $\LHS$ in \equpref{eq:interpolants_finish_desired}, \equpref{eq:posterior_finish_desired}, combined with the bounds $u = \log \varphi \le 0$ for the logistic case for the posterior or \equpref{eq:posterior_spherical_finish_gmmbd} for the GMM case for the posterior, we obtain with probability at least $1 - \exp\big(-d/K(\eps, \alpha, \MARGIN, \lambda)\big)$ and for $d \ge K$,
\begin{align}\label{eq:posterior_spherical_finish_lhsbd}
\LHS \le -c' + K(\lambda) \sqrt{\alpha}\,.
\end{align}
On the other hand, by the proof of 1) in Appendix \upref{subsec:univariate_maximizer_pf}, we have $\sup_{x \in [-1,1]} \inf_{q \in [0,1)} \Phi(x,q) = \Phi(x^\star, q^\star)$ for a unique $(x^\star, q^\star) \in [0, K(\lambda) \alpha]^2$. Since $\alpha \le \alpha_0(\lambda, \MARGIN_+)$ and $\delta \le 1/10$, the bound \equpref{eq:rs_bivariate_lowerbd_q} implies that 
\begin{align}\label{eq:posterior_spherical_finish_rhsbd}
\sup_{x \in [-1 + \delta, 1-\delta]} \inf_{q \in [0,1)} \Phi(x,q) = \Phi(x^\star, q^\star) \ge -\alpha K(\lambda, \MARGIN_+)\,.
\end{align}
Since $\alpha \le \alpha_0(\lambda, \MARGIN_+)$ and $c'>0$ is a universal constant, combining \equpref{eq:posterior_spherical_finish_lhsbd}, \equpref{eq:posterior_spherical_finish_rhsbd} now implies \equpref{eq:interpolants_finish_desired} for the interpolators or \equpref{eq:posterior_finish_desired} for the posterior. This completes the proof of Theorems \upref{thm:mainresult}, \upref{thm:mainresultposterior}.

\subsection{Proof of Theorem \upref{thm:sphericalmeasurersformula} for interpolators: $\MARGIN < 0$}\label{subsec:sphericalpf_negativemargin}
Throughout this proof, we let $\bar I := [x-\eps'', x+\eps'']$. 
Note $0 \not\in \bar I$ as $|x|>\eps''$. 
Consequently recalling \equpref{eq:Fqrhodef}, we note $\RSG_{\bar I}$ does not depend on $\eps''$ as $r_{\bar I}(x)$, the only part depending on $\bar I$, equals $x^2$. 
In the following Appendices \upref{subsec:sphericalpf_negativemargin}, \upref{subsec:sphericalpf_positivemargin}, \upref{subsec:sphericalpf_posterior}, we let $\AREA(A)$ refer to the surface area any Lebesgue measurable subset of $v \cS^{d-1}(\sqrt{d})$ for any $v>0$.

\paragraph{Preliminaries:} 
Due to the lack of continuity caused by the hard indicator functions, we instead use the monotonicity from the following comparison inequalities.
For any $v>0$, consider
\begin{align}
vU_i = \Bigl\{ v\theta \in \R^{\DIMENSION} \, :\, \frac1{\sqrt{\DIMENSION}}\Bigl( s_i \theta_1 + \langle g_i, \bar{\btheta} \rangle \Bigr) \ge \MARGIN \Bigr\} = \Bigl\{ \theta \in \R^{\DIMENSION} \, :\, \frac1{\sqrt{\DIMENSION}}\Bigl( s_i \theta_1 + \langle g_i, \bar{\btheta} \rangle \Bigr) \ge v\MARGIN \Bigr\}\, , \label{eq:v_plane_logic}
\end{align} 
where the last equality follows from the transformation $\btheta \leftarrow v\btheta$.
Therefore as $\MARGIN < 0$, we have 
\begin{align}
&\text{For $v \ge 1$:}\quad vU_i \supseteq \Bigl\{ \theta \in \R^{\DIMENSION} \, :\, \frac1{\sqrt{\DIMENSION}}\Bigl( s_i \theta_1 + \langle g_i, \bar{\btheta} \rangle \Bigr) \ge \MARGIN \Bigr\} = U_i, \text{ hence } vC \supseteq C\,, \label{eq:sphericalvlarger1containment}\\
&\text{For $v \le 1$:}\quad vU_i \subseteq \Bigl\{ \theta \in \R^{\DIMENSION} \, :\, \frac1{\sqrt{\DIMENSION}}\Bigl( s_i \theta_1 + \langle g_i, \bar{\btheta} \rangle \Bigr) \ge \MARGIN \Bigr\} = U_i, \text{ hence } vC \subseteq C\, .\label{eq:sphericalvsmaller1containment}
\end{align}
We next set up several comparison inequalities letting us to relate the indicator $\one\{\theta_1/\sqrt{d} \in I\}$ to $\one\{\theta_1/\sqrt{d} \in \bar I\}$ (the difference is between $I, \bar I$).
For any interval $[a,b]$, remark
\begin{align}
&\text{For $1 \le v \le w$}:\quad v[a, b] \supseteq [wa, b]\text{ if $b \ge a \ge 0$}\quad , \quad  v[a,b] \supseteq [a, wb]\text{ if $a \le b \le 0$}\,,\label{eq:sphericalconversioncontainment_larger1} \\
&\text{For $w \le v \le 1$}:\quad v[a, b] \subseteq [wa, b]\text{ if $b \ge a \ge 0$}\quad , \quad v[a,b] \subseteq [a, wb]\text{ if $a \le b \le 0$}\,.\label{eq:sphericalconversioncontainment_smaller1}
\end{align} 

\paragraph{Proof of \equpref{eq:spherical_conversion_upperbound}:}
Let $w=1-\eps''$ here in the proof of \equpref{eq:spherical_conversion_upperbound} and let $I$ be as per \equpref{eq:sphericalconversion_I_upperbound}. 
Consider any $\GAUSSIANMEASURE \in [\GAUSSIANMEASURE_0, \GAUSSIANMEASURE_1]$. By \equpref{eq:sphericalvsmaller1containment} and \equpref{eq:sphericalconversioncontainment_smaller1}, we have the following. If $x \ge \eps'' >0$,
\begin{align*}
\LHS := &\int \one\big\{\{\theta_1/\sqrt{d} \in \bar I \} \cap C \big\} \exp(-\GAUSSIANMEASURE \| \btheta \|^2)\, \rmd \btheta \\
&\ge \VOL\Bigl(\one \big\{\theta_1/\sqrt{\DIMENSION} \in [x-\eps'', x+\eps''] \big\} \cap C \cap \one \big\{w^2 d \le \|\btheta \|^2 \le d \big\} \Bigr) \exp(-\GAUSSIANMEASURE d) \\
&=  \exp(-\GAUSSIANMEASURE \DIMENSION)\sqrt{\DIMENSION}\int_w^1 \AREA\Big( \one\big\{\theta_1/\sqrt{\DIMENSION} \in [x-\eps'', x+\eps'']\big\} \cap C \cap v\cS^{d-1}(\sqrt{d}) \Big) \, \rmd v \\
&\ge \exp(-\GAUSSIANMEASURE \DIMENSION)\sqrt{\DIMENSION}\int_w^1 \AREA \Big(\one\Big\{\theta_1/\sqrt{\DIMENSION} \in v \Big[ \frac{x-\eps''}{w}, x+\eps'' \Big]\Big\}  \cap vC \cap v\cS^{d-1}(\sqrt{d}) \Big)\, \rmd v \\
&= \exp(-\GAUSSIANMEASURE \DIMENSION) \sqrt{\DIMENSION} \int_w^1 \AREA\Big( \one\Big\{\theta_1/\sqrt{\DIMENSION} \in \Big[\frac{x-\eps''}{w}, x+\eps'' \Big]\Big\} \cap C \cap \cS^{d-1}(\sqrt{d}) \Big) v^{d-1}\, \rmd v \\
&= \exp(-\GAUSSIANMEASURE \DIMENSION)\, \AREA \Big(\one\{\theta_1/\sqrt{\DIMENSION} \in I\} \cap C \cap \cS^{d-1}(\sqrt{d}) \Big)\,   \frac1{\sqrt{\DIMENSION}}\, (1-w^d)\, \\
&\ge \exp(-\GAUSSIANMEASURE \DIMENSION)\, \AREA \Big(\one\{\theta_1/\sqrt{\DIMENSION} \in I\} \cap C \cap \cS^{d-1}(\sqrt{d}) \Big)\,   \frac1{2\sqrt{\DIMENSION}}\,.
\end{align*}
Here we used the change of variable $\theta \leftarrow v \theta$, and that $w^d = (1-\eps'')^d \le \frac12$ for $d \ge d(\eps'')$.
Similarly if $x \le -\eps''<0$, by an analogous derivation as above, again using \equpref{eq:sphericalvsmaller1containment} and \equpref{eq:sphericalconversioncontainment_smaller1},
\begin{align*}
\LHS 
&\ge \exp(-\GAUSSIANMEASURE \DIMENSION)\sqrt{\DIMENSION}\int_w^1 \AREA \Big(\one\big\{\theta_1/\sqrt{\DIMENSION} \in [x-\eps'', x+\eps'']\big\} \cap C \cap v\cS^{d-1}(\sqrt{d}) \Big)\, \rmd v \\
&\ge \exp(-\GAUSSIANMEASURE \DIMENSION)\sqrt{\DIMENSION}\int_w^1 \AREA \Big(\one\Big\{\theta_1/\sqrt{\DIMENSION} \in v \Big[ x-\eps'', \frac{x+\eps''}{w} \Big]\Big\}  \cap vC \cap v\cS^{d-1}(\sqrt{d}) \Big)\, \rmd v \\
&\ge \exp(-\GAUSSIANMEASURE \DIMENSION)\, \AREA \Big(\one\{\theta_1/\sqrt{\DIMENSION} \in I\} \cap C \cap \cS^{d-1}(\sqrt{d}) \Big)\,   \frac1{2\sqrt{\DIMENSION}}\,.
\end{align*}
For both cases, taking logarithms gives
\begin{align}\label{eq:sphericalconversion_upperbd_comparisonresult}
\frac1d \log\big( \LHS \big) 
&\ge -\GAUSSIANMEASURE + \frac1d \log \AREA \Big(\one\{\theta_1/\sqrt{\DIMENSION} \in I\} \cap C \cap \cS^{d-1}(\sqrt{d}) \Big) + \frac1d \log\Big( \frac1{2\sqrt{\DIMENSION}}\Big)\, .
\end{align}
By Theorem \upref{thm:gaussianmeasurersformula} applied with the interval $\bar I$, since $\eps'' \le \eps_{\RSG}(\eps, \DENSITYBOUND, \beta_0, \beta_1, \MARGIN, \lambda)$, we have for $d \ge d(\eps, \eps'', \DENSITYBOUND, \beta_0, \beta_1, \MARGIN, \lambda)$ with probability at least $1-\exp\big(-d / K(\eps, \DENSITYBOUND, \beta_0, \beta_1, \MARGIN, \lambda)\big)$,
\begin{align}
&\frac1d \log\big( \LHS \big) \in \Big[ \RSG_{\bar I}(n/d, \GAUSSIANMEASURE, x) - \eps, \RSG_{\bar I}(n/d, \GAUSSIANMEASURE, x) + \eps \Big]\, .\label{eq:sphericalconversion_upperbd_gaussian1}
\end{align}
By Proposition \upref{prop:nicersuniquesol}, as $0 \not\in \bar I$, the Lipschitz constant of $\RSG_{\bar I}(\alpha, \beta, x)$ over $(\alpha, \beta, x) \in [0, \DENSITYBOUND] \times [\beta_0, \beta_1] \in [-1, 1]$ is upper bounded in terms of $\DENSITYBOUND, \beta_0, \beta_1, \MARGIN, \lambda$. 
Therefore, we obtain that for $d \ge d(\eps, \DENSITYBOUND, \beta_0, \beta_1, \MARGIN, \lambda)$, we have $\big| \RSG_{\bar I}(n/d, \beta, x) - \RSG_{\bar I}(\alpha, \beta, x) \big| \le \eps$.
Noting 
\begin{align*}
\AREA\Big(\one\{\theta_1/\sqrt{\DIMENSION} \in I\} \cap C \cap \cS^{d-1}(\sqrt{d}) \Big) = \mu_d\Big(\one\{\theta_1/\sqrt{\DIMENSION} \in I\} \cap C \Big) \AREA\big( \cS^{d-1}(\sqrt{d}) \big)\, ,
\end{align*}
combining with \equpref{eq:sphericalconversion_upperbd_comparisonresult}, \equpref{eq:sphericalconversion_upperbd_gaussian1} yields for the same $d$ and the same probability as \equpref{eq:sphericalconversion_upperbd_gaussian1} that
\begin{equation}\label{eq:sphericalconversionupperboundfinal}
\begin{aligned}
&\frac1d \log \mu_d\Big(\one\{\theta_1/\sqrt{\DIMENSION} \in I\} \cap C \Big) \\
&\quad \le \GAUSSIANMEASURE + \RSG_{\bar I}(\alpha, \GAUSSIANMEASURE, x) + 2\eps + \frac{\log(2\sqrt{\DIMENSION})}d - \frac1d \log \AREA(\cS^{d-1}(\sqrt{d}))\, .
\end{aligned}
\end{equation}
By Lemma \upref{lem:rs_formula_properties}, as $\alpha \le \DENSITYBOUND$ and as $0 \not\in \bar I$, $\inf_{q \in [0,1)} \Phi(x, q)=\GAUSSIANMEASURE+\RSG_{\bar I}(\alpha,\GAUSSIANMEASURE, x)-\frac12 \log(2\pi e)$ for $\beta=\beta(\alpha, x)$.
Moreover, recall $\lim_{d\rightarrow\infty} \frac1d \log \AREA(\cS^{d-1}(\sqrt{d})) = \frac12 \log(2e\pi)$.

Taking $\GAUSSIANMEASURE=\GAUSSIANMEASURE(\alpha, x)$ in \equpref{eq:sphericalconversionupperboundfinal} thus gives for $d \ge d(\eps, \eps'', \alpha, \beta_0, \beta_1, \MARGIN, \lambda)$ and with probability at least $1-\exp\big(-d / K(\eps, \alpha, \beta_0, \beta_1, \MARGIN, \lambda)\big)$,
\[ \frac1d \log \mu_d\Big(\one\{\theta_1/\sqrt{\DIMENSION} \in I\} \cap C \Big) \le \inf_{q \in [0,1)} \Phi(x, q) + 3\eps\, ,\]
where we note $\lim_{d\rightarrow\infty} \frac1d \log \AREA(\cS^{d-1}(\sqrt{d}))-\frac12 \log(2e\pi) = 0$. The result follows taking $\eps \leftarrow \eps/3$.

\paragraph{Proof of \equpref{eq:spherical_conversion_lowerbound}:} We proceed with a similar strategy to the proof of \equpref{eq:spherical_conversion_upperbound} above, though this direction is more involved. 
First, by the second part of Proposition \upref{prop:nicersuniquesol}, $\rho_0(\alpha, \beta, x)$ is infinitely differentiable in $\GAUSSIANMEASURE$ and $x$ for $(\GAUSSIANMEASURE, x) \in [\GAUSSIANMEASURE_0, \GAUSSIANMEASURE_1] \times [-1, 1]$. We thus may define
\begin{align}
\underH = \underH(\alpha, \MARGIN, \lambda) &:= \inf_{(\GAUSSIANMEASURE, x) \in [\GAUSSIANMEASURE_0, \GAUSSIANMEASURE_1] \times [-1, 1]} -\frac{\partial}{\partial \GAUSSIANMEASURE} \rho_0(\alpha,\GAUSSIANMEASURE, x)\,, \label{eq:underLdef}\\
\overH = \overH(\alpha, \MARGIN, \lambda) &:= \sup_{(\GAUSSIANMEASURE, x) \in [\GAUSSIANMEASURE_0, \GAUSSIANMEASURE_1] \times [-1, 1]} -\frac{\partial}{\partial \GAUSSIANMEASURE} \rho_0(\alpha,\GAUSSIANMEASURE, x)\,.\label{eq:overLdef}
\end{align}
By Lemma \upref{lem:rs_formula_properties}, Proposition \upref{prop:nicersuniquesol}, and compactness, $\underH > 0$. $\underH > 0$ yields strong convexity, which will be crucial in our subsequent arguments. 
For the rest of the proof of \equpref{eq:spherical_conversion_lowerbound}, we consider $\eps < \min\big\{1, \frac{1}{12} \underH \big\}$.

\paragraph{Choice of parameters and properties of variational problem.}
Let $w = 1 + \eps$ and define
\begin{align}\label{eq:sphericalconversion_barf_def}
\bar f(\alpha, \beta, x) := \bar\Phi\big(x, q_0(\alpha, \beta, x), \rho_0(\alpha, \beta, x)\big) - \beta\big(\rho_0(\alpha, \beta, x) + x^2)\,,
\end{align}
where consider $\bar\Phi(x,q,\rho)$ with $u(x)$ therein such that $\exp u(x) = \one\{x \ge \MARGIN\}$. (Here as $0 \not\in \bar I$, we have $\bar f = \RSG_{\bar I} - \frac12 \log(2\pi e)$, however since $\bar I$ is technically defined in terms of $\eps''$, we will define $\eps_{\RS}$ in terms of $\bar f$ and then consider $0 < \eps'' \le \eps_{\RS}$.) We now define the following positive parameters:
\begin{align}
\VARPROBLEMSPHERICAL>0 &\text{ as the unique solution to } \rho_0\big(\alpha, \GAUSSIANMEASURE(\alpha, x) - \VARPROBLEMSPHERICAL, x\big) = w^2 - x^2\, ,\label{eq:sphericaldeltachoice} \\
\VARPROBLEMSPHERICAL_1, \VARPROBLEMSPHERICAL_2 &:= \VARPROBLEMSPHERICAL / 2, \label{eq:sphericaldelta12choice}\\
\eps_1 &\le \min\Big\{ \frac13 \Big( \bar f\big(\alpha, \GAUSSIANMEASURE(\alpha, x) - \VARPROBLEMSPHERICAL_1, x\big) - \big( \VARPROBLEMSPHERICAL_1 + \bar f\big(\alpha, \GAUSSIANMEASURE(\alpha, x), x\big) \big) \Big), \eps \Big\} \, ,\label{eq:sphericaleps1choice}\\
\eps_2 &\le\min\Big\{ \frac{1}{6} \underH \VARPROBLEMSPHERICAL_2^2, \eps \Big\} \, ,\label{eq:sphericaleps2choice} \\
\GAUSSIANMEASURE &:= \GAUSSIANMEASURE(\alpha, x) - \VARPROBLEMSPHERICAL_1 > \beta_0 = \frac1{8}\, .\label{eq:sphericalkappachoice}
\end{align}
Therefore, $\VARPROBLEMSPHERICAL, \VARPROBLEMSPHERICAL_1, \VARPROBLEMSPHERICAL_2, \eps_1, \eps_2$ above depend on $x, \eps, \alpha, \beta_0, \beta_1, \MARGIN, \lambda$. Later on we will show $\eps_1, \eps_2$ can be taken independently of $x$. First we will justify the statements implicit in the above, including positivity, and derive some useful estimates. Recall by Lemma \upref{lem:rs_formula_properties}, $\GAUSSIANMEASURE(\alpha, x) \ge \frac1{4}$. Now note:
\begin{itemize}
\item By definition of $\underH$, for small enough $\eps \le \min\big\{1, \frac{1}{12}\underH \big\}$, we have
\begin{align*}
\rho_0\big(\alpha, 1/8, x \big) &\ge \rho_0\big(\alpha, \GAUSSIANMEASURE(\alpha, x), x \big) + \underH \big( \GAUSSIANMEASURE(\alpha, x) - 1/8 \big) \\
&\ge 1-x^2+3\eps \\
&> w^2 - x^2 \\
&> 1-x^2 = \rho_0\big(\alpha, \GAUSSIANMEASURE(\alpha, x), x \big)\, .
\end{align*}
As $\rho_0(\alpha, \GAUSSIANMEASURE, x)$ is strictly decreasing on $[\beta_0, \beta_1]$ by Lemma \upref{lem:rs_formula_properties}, it follows that such a $\VARPROBLEMSPHERICAL$ exists, is unique, and that $\VARPROBLEMSPHERICAL>0$. This justifies \equpref{eq:sphericaldeltachoice}.

\item We have 
\begin{align*}
1+\eps \ge w &= \sqrt{\rho_0\big(\alpha, \GAUSSIANMEASURE(\alpha, x) - \VARPROBLEMSPHERICAL, x\big) +x^2 } \ge \sqrt{\rho_0\big(\alpha, \GAUSSIANMEASURE(\alpha, x), x \big) + x^2 + \underH \VARPROBLEMSPHERICAL} = \sqrt{1+\underH \VARPROBLEMSPHERICAL}\, ,
\end{align*}
thus as $\eps \le \frac{1}{12}\underH$, we have $\VARPROBLEMSPHERICAL < \frac1{8}$. Hence, $\GAUSSIANMEASURE(\alpha, x) - \VARPROBLEMSPHERICAL_1 \ge \frac1{4} - \VARPROBLEMSPHERICAL_1 > \frac1{8}$, justifying \equpref{eq:sphericalkappachoice}.

\item By Lemma \upref{lem:rs_formula_properties}, we have for all $\GAUSSIANMEASURE \in [\GAUSSIANMEASURE_0, \GAUSSIANMEASURE_1]$ that 
\[ \GAUSSIANMEASURE + \bar f\big(\alpha, \GAUSSIANMEASURE, x\big) > \GAUSSIANMEASURE(\alpha, x) + \bar f\big(\alpha, \GAUSSIANMEASURE(\alpha, x), x \big)\, . \]
Taking $\GAUSSIANMEASURE = \GAUSSIANMEASURE(\alpha, x) - \VARPROBLEMSPHERICAL_1 > \frac1{8} = \beta_0$ in the above justifies $\eps_1>0$ and hence \equpref{eq:sphericaleps1choice}. Moreover by the definition of $\eps_1$, the above display gives
\begin{align}
-\VARPROBLEMSPHERICAL_1 + \bar f\big(\alpha, \GAUSSIANMEASURE(\alpha, x) - \VARPROBLEMSPHERICAL_1, x\big) > \bar f\big(\alpha, \GAUSSIANMEASURE(\alpha, x), x \big) + 3\eps_1\, .\label{eq:sphericalvarproblem1}
\end{align}

\item Let $F(x, q, \rho) := \bar\Phi(x, q, \rho) - \beta (\rho + x^2)$ and let $\rho_0 = \rho_0(\alpha, \beta, x), q_0 = q_0(\alpha, \beta, x)$.
Note $\frac{\partial F}{\partial \rho}=\frac{\partial F_I}{\partial \rho}$, $\frac{\partial F}{\partial q}=\frac{\partial F_I}{\partial q}$ for any interval $I$, thus $\frac{\partial F}{\partial \rho}(x, q_0, \rho_0) = \frac{\partial F}{\partial q}(x, q_0, \rho_0) = 0$.
Now for all $\beta \in [\beta_0, \beta_1]$,
\begin{align}
\frac{\partial}{\partial \beta} \bar f(\alpha, \beta, x) &= - \big( \rho_0(\alpha, \beta, x) + x^2 \big) + \frac{\partial F}{\partial \rho}(x, q_0, \rho_0) \cdot \frac{\partial \rho_0}{\partial \beta} + \frac{\partial F}{\partial q}(x, q_0, \rho_0) \cdot \frac{\partial q_0}{\partial \beta} \notag \\
&= - \big( \rho_0(\alpha, \beta, x) + x^2 \big)\,.\label{eq:fbarfirstderiv}
\end{align}
Thus by definition of $\underH$, we have for all $\beta \in [\beta_0, \beta_1]$,
\begin{align}
\frac{\partial^2}{\partial \beta^2}\bar f(\alpha, \beta, x) = -\frac{\partial}{\partial \beta} \rho_0(\alpha, \beta, x) \ge \underH > 0\,. \label{eq:fbarsecondderiv}
\end{align}
Now in what follows, let $\beta=\beta(\alpha, x)-\VARPROBLEMSPHERICAL_1$ as per \equpref{eq:sphericalkappachoice}, thus $\beta-\VARPROBLEMSPHERICAL_2 = \beta(\alpha, x)-\VARPROBLEMSPHERICAL$ as $\VARPROBLEMSPHERICAL_1+\VARPROBLEMSPHERICAL_2=\VARPROBLEMSPHERICAL$. 
By our choice of $w$, $\VARPROBLEMSPHERICAL$ from \equpref{eq:sphericaldeltachoice}, we obtain that for this choice of $\beta=\beta(\alpha, x)-\VARPROBLEMSPHERICAL_1$,
\begin{align}
\bar f(\alpha, \GAUSSIANMEASURE, x) &\ge \bar f(\alpha, \GAUSSIANMEASURE-\VARPROBLEMSPHERICAL_2, x) + \VARPROBLEMSPHERICAL_2 \frac{\partial}{\partial \GAUSSIANMEASURE} \bar f(\alpha,\GAUSSIANMEASURE-\VARPROBLEMSPHERICAL_2, x)+\frac{1}2 \underH \VARPROBLEMSPHERICAL_2^2 \notag \\
&= \bar f(\alpha, \GAUSSIANMEASURE-\VARPROBLEMSPHERICAL_2, x) - \VARPROBLEMSPHERICAL_2 \big( \rho_0(\alpha,\GAUSSIANMEASURE-\VARPROBLEMSPHERICAL_2, x ) + x^2 \big) + \frac{1}2\underH \VARPROBLEMSPHERICAL_2^2 \notag \\
&= \bar f(\alpha, \GAUSSIANMEASURE-\VARPROBLEMSPHERICAL_2, x) - \VARPROBLEMSPHERICAL_2 w^2 + \frac{1}2 \underH \VARPROBLEMSPHERICAL_2^2 \notag \\
&\ge \bar f(\alpha, \GAUSSIANMEASURE-\VARPROBLEMSPHERICAL_2, x) - \VARPROBLEMSPHERICAL_2 w^2 + 3\eps_2\, .\label{eq:sphericalvarproblem2}
\end{align}
\end{itemize}
As our last preliminary, we claim we can take $\eps_1, \eps_2$, and
\begin{align}
\eps_{\RS} := \min\big\{ \eps_{\RSG}(\eps_1, \DENSITYBOUND, \beta_0, \beta_1, \MARGIN, \lambda), \eps_{\RSG}(\eps_2, \DENSITYBOUND, \beta_0, \beta_1, \MARGIN, \lambda) \big\} > 0 \label{eq:sphericalconversion_epsRS_def}
\end{align}
independently of $x$, where $\eps_{\RSG}$ is from Theorem \upref{thm:gaussianmeasurersformula}. To this end first observe
\begin{align*}
w^2 - x^2 &= \rho_0\big( \alpha, \beta(\alpha, x) - \VARPROBLEMSPHERICAL, x \big) \le \rho_0\big( \alpha, \beta(\alpha, x), x \big) + \overH \VARPROBLEMSPHERICAL = 1-x^2 + \overH \VARPROBLEMSPHERICAL\,.
\end{align*}
Since $w=1+\eps$, it follows that $\VARPROBLEMSPHERICAL > \frac{2}{\overH} \eps$, and hence $\VARPROBLEMSPHERICAL, \VARPROBLEMSPHERICAL_1, \VARPROBLEMSPHERICAL_2$ are lower bounded independently of $x$.
Next note 
\begin{align*}
\frac{\partial}{\partial \beta} \big( \beta + \bar f(\alpha, \beta, x) \big) = 1 - \big( \rho_0(\alpha, \beta, x) + x^2 \big)\,,\,\frac{\partial^2}{\partial \beta^2} \big( \beta + \bar f(\alpha, \beta, x) \big) = -\frac{\partial}{\partial \beta} \rho_0(\alpha, \beta, x) \ge \underH > 0
\end{align*}
by \equpref{eq:fbarfirstderiv}, \equpref{eq:fbarsecondderiv}. Thus as $\rho_0(\alpha, \beta(\alpha, x), x) = 1-x^2$ by Lemma \upref{lem:rs_formula_properties}, we have 
\begin{align*}
\beta(\alpha, x) - \VARPROBLEMSPHERICAL_1 + \bar f\big(\alpha, \beta(\alpha, x) - \VARPROBLEMSPHERICAL_1, x \big) &\ge \beta(\alpha, x) + \bar f\big(\alpha, \beta(\alpha, x), x \big) - \VARPROBLEMSPHERICAL_1 \frac{\partial}{\partial \beta} \big( \beta + \bar f(\alpha, \beta, x) \big)+ \frac{1}2 \underH \VARPROBLEMSPHERICAL_1^2 \\
&\ge \beta(\alpha, x) + \bar f\big(\alpha, \beta(\alpha, x), x \big) + \frac{1}2 \underH \VARPROBLEMSPHERICAL_1^2\,.
\end{align*}
Therefore
\begin{align*}
\bar f\big(\alpha, \beta(\alpha, x) - \VARPROBLEMSPHERICAL_1, x \big) \ge \VARPROBLEMSPHERICAL_1 + \bar f\big(\alpha, \beta(\alpha, x), x \big)  + \frac{\underH}2 \VARPROBLEMSPHERICAL_1^2\,.
\end{align*}
Since $\VARPROBLEMSPHERICAL > \frac{2}{\overH} \eps$, it follows from the above display and \equpref{eq:sphericaleps1choice} that $\eps_1$ can be taken independently of $x$. Recalling $\VARPROBLEMSPHERICAL_2$ is lower bounded independently of $x$ and $\underH$ does not depend on $x$, it follows that $\eps_2$ can be taken independently of $x$. 
Hence $\eps_{\RS}$ can also be taken independently of $x$.

We now turn to establishing \equpref{eq:spherical_conversion_lowerbound}. 
We consider any $\eps''$, $0 < \eps'' \le \eps_{\RS}$, and define $I$ as in \equpref{eq:sphericalconversion_I_lowerbound}. 
Note $0 \not\in \bar I$ as $|x|>\eps''$.
In the following proofs, we take $\beta = \beta(\alpha, x) - \VARPROBLEMSPHERICAL_1$ as per \equpref{eq:sphericalkappachoice}.

\paragraph{Initial Bound.}
Suppose $x \ge \eps''>0$. Recall $I = [\frac{x-\eps''}w, x+\eps'']$ where $w=1+\eps$.
By \equpref{eq:sphericalvlarger1containment} and \equpref{eq:sphericalconversioncontainment_larger1},
\begin{align}
\LHS &:= \int \one\big\{ \{\theta_1/\sqrt{\DIMENSION} \in \bar I\} \cap C \cap \{\DIMENSION \le \|\btheta\|^2 \le w^2 \DIMENSION\} \big\} \exp\big(-\GAUSSIANMEASURE \| \btheta \|^2\big)\, \rmd \btheta \notag \\
&\le \exp(-\GAUSSIANMEASURE \DIMENSION) \VOL\Big( \one\big\{\theta_1/\sqrt{\DIMENSION} \in [x-\eps'', x+\eps''] \big\} \cap C \cap \big\{ \DIMENSION \le \| \btheta \| \le w^2 \DIMENSION \big\} \Big)\notag  \\
&= \exp(-\GAUSSIANMEASURE \DIMENSION) \sqrt{\DIMENSION}\int_1^w \AREA \Big(\one\Big\{\theta_1/\sqrt{\DIMENSION} \in \Big[w \cdot \frac{x-\eps''}w, x+\eps'' \Big] \Big\} \cap C \cap v\cS^{d-1}(\sqrt{d}) \Big)\, \rmd v\notag \\
&\le \exp(-\GAUSSIANMEASURE \DIMENSION) \sqrt{\DIMENSION} \int_1^w \AREA \Big( \one\Big\{\theta_1/\sqrt{\DIMENSION} \in v\Big[\frac{x-\eps''}w, x+\eps'' \Big] \Big\} \cap vC \cap v\cS^{d-1}(\sqrt{d}) \Big)\, \rmd v \notag \\
&= \exp(-\GAUSSIANMEASURE \DIMENSION) \sqrt{\DIMENSION} \int_1^w \AREA \Big( \one\Big\{\theta_1/\sqrt{\DIMENSION} \in \Big[\frac{x-\eps''}w, x+\eps'' \Big] \Big\} \cap C \cap \cS^{d-1}(\sqrt{d}) \Big) v^{\DIMENSION-1}\, \rmd v\notag  \\
&\le \exp(-\GAUSSIANMEASURE \DIMENSION)\frac1{\sqrt{\DIMENSION}} (w^\DIMENSION-1) \AREA\Big( \one\{\theta_1/\sqrt{\DIMENSION} \in I \} \cap C \cap \cS^{d-1}(\sqrt{d}) \Big)\, .\label{eq:sphericalconversionlowerboundstep1}
\end{align}
Analogously, when $x \le -\eps''<0$, we have by \equpref{eq:sphericalvlarger1containment} and \equpref{eq:sphericalconversioncontainment_larger1} that
\begin{align*}
\LHS &\le \exp(-\GAUSSIANMEASURE \DIMENSION) \sqrt{\DIMENSION}\int_1^w \AREA \Big(\one\Big\{\theta_1/\sqrt{\DIMENSION} \in \Big[x-\eps'', w \cdot \frac{x+\eps''}w \Big] \Big\} \cap C \cap v\cS^{d-1}(\sqrt{d}) \Big)\, \rmd v \\
&\le \exp(-\GAUSSIANMEASURE \DIMENSION) \sqrt{\DIMENSION} \int_1^w \AREA \Big( \one\Big\{\theta_1/\sqrt{\DIMENSION} \in v\Big[x-\eps'', \frac{x+\eps''}{w} \Big] \Big\} \cap vC \cap v\cS^{d-1}(\sqrt{d}) \Big)\, \rmd v \\
&\le \exp(-\GAUSSIANMEASURE \DIMENSION)\frac1{\sqrt{\DIMENSION}} (w^\DIMENSION-1) \AREA\Big( \one\{\theta_1/\sqrt{\DIMENSION} \in I \} \cap C \cap \cS^{d-1}(\sqrt{d}) \Big)\,.
\end{align*}
This yields the same inequality as \equpref{eq:sphericalconversionlowerboundstep1} in the $x<0$ case.
(Note here that $0 \ge \frac{x+\eps''}{w} > x+\eps''$ as $x\le -\eps''<0$.)
Next, note
\begin{align*}
\LHS = \expressionI - \expressionII - \expressionIII \, ,
\end{align*}
where 
\begin{align*}
\expressionI &:= \int \one\big\{\{\theta_1/\sqrt{\DIMENSION} \in \bar I \} \cap C \big\} \exp(-\GAUSSIANMEASURE \| \btheta \|^2)\, \rmd \btheta\, ,\\
\expressionII &:= \int \one\big\{\{\theta_1/\sqrt{\DIMENSION} \in \bar I \} \cap C \cap \{ \| \btheta \| \le \sqrt{\DIMENSION}\} \big\} \exp(-\GAUSSIANMEASURE \| \btheta \|^2)\, \rmd \btheta\, , \\
\expressionIII &:= \int \one\big\{\{\theta_1/\sqrt{\DIMENSION} \in \bar I \} \cap C \cap \{ \| \btheta \|  > w\sqrt{\DIMENSION}\} \big\}\exp(-\GAUSSIANMEASURE \| \btheta \|^2)\, \rmd \btheta\, .
\end{align*}
We now lower bound $\expressionI$ and upper bound $\expressionII, \expressionIII$ via Theorem \upref{thm:gaussianmeasurersformula}.

\paragraph{Lower bounding $\expressionI$:} As $\eps'' \le \eps_{\RS}$, by definition of $\eps_{\RS}$ in \equpref{eq:sphericalconversion_epsRS_def}, and as $\beta = \beta(\alpha, x) \in [\beta_0, \beta_1]$ by \equpref{eq:sphericalkappachoice}, Theorem \upref{thm:gaussianmeasurersformula} applied with $\eps=\frac12 \eps_1 \wedge \eps_2$ and the interval $\bar I$ implies that for $d \ge d(\eps_1, \eps_2, \eps'', \DENSITYBOUND, \GAUSSIANMEASURE_0, \GAUSSIANMEASURE_1, \MARGIN, \lambda)$ and probability at least $1-\exp\big(-d / K(\eps_1, \eps_2, \DENSITYBOUND, \GAUSSIANMEASURE_0, \GAUSSIANMEASURE_1, \MARGIN, \lambda) \big)$,
\begin{align}
\expressionI &= \int \one\big\{\{\theta_1/\sqrt{\DIMENSION} \in \bar I \} \cap C \big\} \exp(-\GAUSSIANMEASURE \| \btheta \|^2)\, \rmd \btheta \notag \\
&\ge \max\Big\{ \exp \Big( \big( \RSG_{\bar I}(n/d, \GAUSSIANMEASURE, x) -  \eps_1/2 \big)d \Big), \exp \Big( \big( \RSG_{\bar I}(n/d, \GAUSSIANMEASURE, x) -  \eps_2/2 \big) \Big) d\Big\} \notag \\
&\ge \max\Big\{ \exp \Big( \big( \RSG_{\bar I}(\alpha, \GAUSSIANMEASURE, x) -  \eps_1 \big) d\Big), \exp \Big( \big( \RSG_{\bar I}(\alpha, \GAUSSIANMEASURE, x) -  \eps_2 \big) d\Big) \Big\}\, .\label{eq:sphericalconversion_lowerbd_lowerbdI}
\end{align}
Here, we used the second part of Proposition \upref{prop:nicersuniquesol} in the last inequality, following identical reasoning as the proof for the upper bound to convert from $\RSG_{\bar I}(n/d, \beta, x)$ to $\RSG_{\bar I}(\alpha, \beta, x)$.

\paragraph{Upper bounding $\expressionII$:}
As $\eps'' \le \eps_{\RS} \le \eps_{\RSG}(\eps_1, \DENSITYBOUND, \beta_0, \beta_1, \MARGIN, \lambda)$ and $\beta \in [\beta_0, \beta_1]$, by Theorem \upref{thm:gaussianmeasurersformula}, for $d \ge d(\eps_1, \eps'', \DENSITYBOUND, \GAUSSIANMEASURE_0, \GAUSSIANMEASURE_1, \MARGIN, \lambda)$ and probability at least $1-\exp\big(-d / K(\eps_1, \DENSITYBOUND, \GAUSSIANMEASURE_0, \GAUSSIANMEASURE_1, \MARGIN, \lambda) \big)$,
\begin{align}
\expressionII &= \int \one\big\{\{\theta_1/\sqrt{\DIMENSION} \in \bar I \} \cap C \cap \{ \| \btheta \| \le \sqrt{\DIMENSION}\} \big\} \exp(-\GAUSSIANMEASURE \| \btheta \|^2)\, \rmd \btheta \notag \\
&= \int \one\big\{\{\theta_1/\sqrt{\DIMENSION} \in \bar I \} \cap C \cap \{ \| \btheta \| \le \sqrt{\DIMENSION}\} \big\} \exp\big( -(\GAUSSIANMEASURE + \VARPROBLEMSPHERICAL_1) \| \btheta \|^2 \big) \exp\big( \VARPROBLEMSPHERICAL_1  \| \btheta \|^2 \big)\, \rmd \btheta \notag \\
&\le \exp(\VARPROBLEMSPHERICAL_1 d)\int \one\big\{\{\theta_1/\sqrt{\DIMENSION} \in \bar I \} \cap C \big\} \exp\big( -(\GAUSSIANMEASURE + \VARPROBLEMSPHERICAL_1) \| \btheta \|^2 \big)\, \rmd \btheta \notag \\
&\le \exp(\VARPROBLEMSPHERICAL_1 d) \cdot \exp \Big( \big( \RSG_{\bar I}(n/d, \GAUSSIANMEASURE+\VARPROBLEMSPHERICAL_1, x) + \eps_1/2 \big) d \Big) \notag \\
&\le \exp(\VARPROBLEMSPHERICAL_1 d) \cdot \exp \Big( \big( \RSG_{\bar I}(\alpha, \GAUSSIANMEASURE+\VARPROBLEMSPHERICAL_1, x) + \eps_1 \big) d \Big) \notag \\
&= \exp\Big( \big( \RSG_{\bar I}(\alpha, \GAUSSIANMEASURE+\VARPROBLEMSPHERICAL_1, x) + \VARPROBLEMSPHERICAL_1 + \eps_1 \big)d \Big)\,.\label{eq:sphericalconversion_lowerbd_upperbdIIinit}
\end{align}
Here we again use Proposition \upref{prop:nicersuniquesol} to convert from $\RSG_{\bar I}(n/d, \beta, x)$ to $\RSG_{\bar I}(\alpha, \beta, x)$. 
Now by \equpref{eq:sphericalvarproblem1} and as $\bar f = \RSG_{\bar I} - \frac12 \log(2\pi e)$ as $0 \not\in \bar I$, combining \equpref{eq:sphericalconversion_lowerbd_lowerbdI}, \equpref{eq:sphericalconversion_lowerbd_upperbdIIinit} yields that for $d \ge d(\eps_1, \eps_2, \eps'', \alpha, \GAUSSIANMEASURE_0, \GAUSSIANMEASURE_1, \MARGIN, \lambda)$, with probability at least $1-\exp\big(d / K(\eps_1, \eps_2, \DENSITYBOUND, \GAUSSIANMEASURE_0, \GAUSSIANMEASURE_1, \MARGIN, \lambda) \big)$, 
\begin{align}\label{eq:sphericalconversion_lowerbd_upperbdII}
\expressionI &\ge \exp \Big( \big( \RSG_{\bar I}(\alpha, \GAUSSIANMEASURE, x) -  \eps_1 \big) d\Big) \notag \\
&\ge \exp(\eps_1 d)\exp\Big( \big( \RSG_{\bar I}(\alpha, \GAUSSIANMEASURE+\VARPROBLEMSPHERICAL_1, x) + \VARPROBLEMSPHERICAL_1 + \eps_1 \big) d\Big) \ge 4\expressionII\,. 
\end{align}

\paragraph{Upper bounding $\expressionIII$:} 
As $\eps'' \le \eps_{\RS} \le \eps_{\RSG}(\eps_2, \DENSITYBOUND, \beta_0, \beta_1, \MARGIN, \lambda)$ and $\beta \in [\beta_0, \beta_1]$, by Theorem \upref{thm:gaussianmeasurersformula}, for $d \ge d(\eps_2, \eps'', \DENSITYBOUND, \GAUSSIANMEASURE_0, \GAUSSIANMEASURE_1, \MARGIN, \lambda)$ and probability at least $1-\exp\big(-d / K(\eps_2, \DENSITYBOUND, \GAUSSIANMEASURE_0, \GAUSSIANMEASURE_1, \MARGIN, \lambda)\big)$,
\begin{align}
\expressionIII &= \int \one\big\{\{\theta_1/\sqrt{\DIMENSION} \in \bar I \} \cap C \cap \{ \| \btheta \|  > w\sqrt{\DIMENSION}\} \big\}\exp(-\GAUSSIANMEASURE \| \btheta \|^2)\, \rmd \btheta \notag \\
&= \int \one\big\{\{\theta_1/\sqrt{\DIMENSION} \in \bar I \} \cap C \cap \{ \| \btheta \|  > w\sqrt{\DIMENSION}\} \big\} \exp(- (\GAUSSIANMEASURE - \VARPROBLEMSPHERICAL_2) \| \btheta \|^2) \exp(-\VARPROBLEMSPHERICAL_2 \| \btheta\|^2)\, \rmd \btheta \notag \\
&\le \exp(-\VARPROBLEMSPHERICAL_2 w^2 d) \int \one\big\{\{\theta_1/\sqrt{\DIMENSION} \in \bar I\} \cap C \big\} \exp(- (\GAUSSIANMEASURE - \VARPROBLEMSPHERICAL_2) \| \btheta \|^2)\, \rmd \btheta \notag \\
&\le \exp(-\VARPROBLEMSPHERICAL_2 w^2 d) \cdot \exp \Big( \big( \RSG_{\bar I}(n/d, \GAUSSIANMEASURE-\VARPROBLEMSPHERICAL_2, x) + \eps_2/2 \big) d \Big) \notag \\
&\le \exp(-\VARPROBLEMSPHERICAL_2 w^2 d) \cdot \exp \Big( \big( \RSG_{\bar I}(\alpha, \GAUSSIANMEASURE-\VARPROBLEMSPHERICAL_2, x) + \eps_2 \big) d \Big) \notag \\
&= \exp\Big( \big(\RSG_{\bar I}(\alpha, \GAUSSIANMEASURE-\VARPROBLEMSPHERICAL_2, x) - \VARPROBLEMSPHERICAL_2 w^2 + \eps_2 \big) d\Big)\, .\label{eq:sphericalconversion_lowerbd_upperbdIIIinit}
\end{align}
Here, we again used Proposition \upref{prop:nicersuniquesol}.
By \equpref{eq:sphericalvarproblem2} and as $\bar f = \RSG_{\bar I} - \frac12 \log(2\pi e)$, we obtain that for $d \ge d(\eps_1, \eps_2, \eps'', \DENSITYBOUND, \GAUSSIANMEASURE_0, \GAUSSIANMEASURE_1, \MARGIN, \lambda)$ and probability at least $1-\exp\big(-d / K(\eps_1, \eps_2, \DENSITYBOUND, \GAUSSIANMEASURE_0, \GAUSSIANMEASURE_1, \MARGIN, \lambda) \big)$, 
\begin{align}
\expressionI &\ge \exp \Big( \big( \RSG_{\bar I}(\alpha, \GAUSSIANMEASURE, x) -  \eps_2 \big) d \Big) \notag \\
&\ge \exp(\eps_2 d) \exp\Big( \big(\RSG_{\bar I}(\alpha, \GAUSSIANMEASURE-\VARPROBLEMSPHERICAL_2, x) -\VARPROBLEMSPHERICAL_2 w^2 + \eps_2 \big) d\Big) \ge 4\expressionIII\,.\label{eq:sphericalconversion_lowerbd_upperbdIII}
\end{align}

\paragraph{Finish.} Combining \equpref{eq:sphericalconversion_lowerbd_upperbdII}, \equpref{eq:sphericalconversion_lowerbd_upperbdIII}, we obtain for $d \ge d(\eps_1, \eps_2, \eps'', \DENSITYBOUND, \GAUSSIANMEASURE_0, \GAUSSIANMEASURE_1, \MARGIN, \lambda)$, with probability at least $1-\exp\big(-d / K(\eps_1, \eps_2, \DENSITYBOUND, \GAUSSIANMEASURE_0, \GAUSSIANMEASURE_1, \MARGIN, \lambda)\big)$, we have
\begin{align*}
\LHS = \expressionI - \expressionII - \expressionIII \ge \expressionI - \frac14 \expressionI - \frac14 \expressionI = \frac12 \expressionI\,.
\end{align*}
Combining the above display with \equpref{eq:sphericalconversionlowerboundstep1} and \equpref{eq:sphericalconversion_lowerbd_lowerbdI} gives for $d \ge d(\eps_1, \eps_2, \eps'', \DENSITYBOUND, \GAUSSIANMEASURE_0, \GAUSSIANMEASURE_1, \MARGIN, \lambda)$, with probability at least $1-\exp\big(-d / K(\eps_1, \eps_2, \DENSITYBOUND, \GAUSSIANMEASURE_0, \GAUSSIANMEASURE_1, \MARGIN, \lambda)\big)$,
\begin{align}
\AREA\Big( \one\{\theta_1/\sqrt{\DIMENSION} \in I \} \cap C \cap \cS^{d-1}(\sqrt{d}) \Big) &\ge \frac{\sqrt{\DIMENSION}}{w^d-1} \cdot \exp(\GAUSSIANMEASURE d) \cdot \frac12 \expressionI \notag \\
&\ge \frac{\sqrt{\DIMENSION}}2 \exp\Big( \big(- \log w + \GAUSSIANMEASURE + \RSG_{\bar I}(\alpha, \GAUSSIANMEASURE, x) - \eps \big) d \Big)\, ,\label{eq:sphericalconversionlowerboundstep2}
\end{align}
where we recall that $\eps_1, \eps_2 \le \eps$. Recall $\GAUSSIANMEASURE \in [\beta_0, \beta_1]$, so by Lemma \upref{lem:rs_formula_properties} and as $0 \not\in \bar I$,
\[\GAUSSIANMEASURE + \RSG_{\bar I}(\alpha, \GAUSSIANMEASURE, x) \ge \inf_{q \in [0,1)} \Phi(x, q) + \frac12 \log(2e\pi)\, . \]
Since $\log w \le w-1 \le \eps$, we obtain for the same $d$ and the same probability as \equpref{eq:sphericalconversionlowerboundstep2},
\begin{align}\label{eq:sphericalconversionlowerbound_step3}
\mu_d\Big( \one\{\theta_1/\sqrt{\DIMENSION} \in I \} \cap C \Big) \AREA(\cS^{d-1}(\sqrt{d})) \ge \frac{\sqrt{\DIMENSION}}2 \exp\Big( \Big(\inf_{q \in [0,1)} \Phi(x, q) + \frac{1}2\log (2 e \pi) - 2\eps \Big) d \Big)\,. 
\end{align}
Again recall $\lim_{d\rightarrow\infty} \frac1d \log \AREA(\cS^{d-1}(\sqrt{d})) = \frac12 \log(2e\pi)$, and that $\eps_1, \eps_2$ depend only on $\eps$, $\alpha$, $\GAUSSIANMEASURE_0, \GAUSSIANMEASURE_1$, $\MARGIN$, $\lambda$ and do not depend on $x$ as justified earlier. Rearranging \equpref{eq:sphericalconversionlowerbound_step3} yields for $d \ge d(\eps, \eps'', \alpha, \GAUSSIANMEASURE_0, \GAUSSIANMEASURE_1, \MARGIN, \lambda)$ and probability at least $1-\exp\big(-d / K(\eps, \alpha, \GAUSSIANMEASURE_0, \GAUSSIANMEASURE_1, \MARGIN, \lambda) \big)$,
\begin{align*}
\frac1d \log \mu_d\Big( \one\{\theta_1/\sqrt{\DIMENSION} \in I \}\cap C \Big) \ge \inf_{q \in [0,1)} \Phi(x, q) - 4\eps\, .
\end{align*}
Taking $\eps \leftarrow \eps/4$ completes the proof of \equpref{eq:spherical_conversion_lowerbound}.

\subsection{Proof of Theorem \upref{thm:sphericalmeasurersformula} for interpolators: $\MARGIN \ge 0$}\label{subsec:sphericalpf_positivemargin}
Again, we let $\bar I := [x-\eps'', x+\eps'']$, and we still have $0 \not\in \bar I$ as $|x| > \eps''$. 
In this proof many details are analogous to Appendix \upref{subsec:sphericalpf_negativemargin} above and we only highlight the differences. 
The notation in the following is also identical to Appendix \upref{subsec:sphericalpf_negativemargin}.
Since $\MARGIN \ge 0$, \equpref{eq:v_plane_logic} now yields the comparison inequalities
\begin{align}\label{eq:intersection_containment_positivemargin}
vC \subseteq C \text{ for } v \ge 1\quad , \quad vC \supseteq C \text{ for } v \le 1\,.
\end{align}
Now for any interval $[a,b]$, we have
\begin{align}
&\text{For $1 \le v \le w$}:\quad v[a, b] \subseteq [a, bw]\text{ if $b \ge a \ge 0$}\quad , \quad  v[a,b] \subseteq [aw, b]\text{ if $a \le b \le 0$}\,,\label{eq:sphericalconversioncontainment_larger1_positivemargin} \\
&\text{For $w \le v \le 1$}: \quad v[a, b] \supseteq [a, bw]\text{ if $b \ge a \ge 0$}\quad , \quad v[a,b] \supseteq [aw, b]\text{ if $a \le b \le 0$}\,.\label{eq:sphericalconversioncontainment_smaller1_positivemargin}
\end{align} 
\paragraph{Proof of \equpref{eq:spherical_conversion_upperbound}:}
Let $w=1+\eps''$ for the proof of this upper bound.
For any $\beta \in [\beta_0, \beta_1]$, we obtain the following from \equpref{eq:intersection_containment_positivemargin}, \equpref{eq:sphericalconversioncontainment_larger1_positivemargin}. 
When $x > \eps''$ we have for $d \ge d(\eps'')$,
\begin{align*}
\LHS &:= \int \one\big\{\{\theta_1/\sqrt{\DIMENSION} \in \bar I\}  \cap C \big\} \exp(-\GAUSSIANMEASURE \| \btheta \|^2)\, \rmd \btheta \\
&\ge \exp(-\GAUSSIANMEASURE w^2 \DIMENSION)\sqrt{\DIMENSION}\int_1^w \AREA \Big(\one\Big\{\theta_1/\sqrt{\DIMENSION} \in [x-\eps'', x+\eps'']\Big\} \cap C \cap v\cS^{d-1}(\sqrt{d}) \Big)\, \rmd v \\
&\ge \exp(-\GAUSSIANMEASURE w^2 \DIMENSION)\sqrt{\DIMENSION}\int_1^w \AREA \Big(\one\Big\{\theta_1/\sqrt{\DIMENSION} \in v \Big[ x-\eps'', \frac{x+\eps''}{w} \Big]\Big\}  \cap vC \cap v\cS^{d-1}(\sqrt{d}) \Big)\, \rmd v \\
&\ge \exp(-\GAUSSIANMEASURE w^2 \DIMENSION)\, \AREA \Big(\one\{\theta_1/\sqrt{\DIMENSION} \in I\} \cap C \cap \cS^{d-1}(\sqrt{d}) \Big)\,   \frac1{2\sqrt{\DIMENSION}}\,.
\end{align*}
Similarly when $x<-\eps''$, for $d \ge d(\eps'')$ we have
\begin{align*}
\LHS &\ge \exp(-\GAUSSIANMEASURE w^2 \DIMENSION)\sqrt{\DIMENSION}\int_1^w \AREA \Big(\one\Big\{\theta_1/\sqrt{\DIMENSION} \in [x-\eps'', x+\eps'']\Big\} \cap C \cap v\cS^{d-1}(\sqrt{d}) \Big)\, \rmd v \\
&\ge \exp(-\GAUSSIANMEASURE w^2 \DIMENSION)\sqrt{\DIMENSION}\int_1^w \AREA \Big(\one\Big\{\theta_1/\sqrt{\DIMENSION} \in v \Big[ \frac{x-\eps''}w, x+\eps'' \Big]\Big\} \cap vC \cap v\cS^{d-1}(\sqrt{d}) \Big)\, \rmd v \\
&\ge \exp(-\GAUSSIANMEASURE w^2 \DIMENSION)\, \AREA \Big(\one\{\theta_1/\sqrt{\DIMENSION} \in I\} \cap C \cap \cS^{d-1}(\sqrt{d}) \Big)\,   \frac1{2\sqrt{\DIMENSION}}\,.
\end{align*}
In either case, applying Theorem \upref{thm:gaussianmeasurersformula} with the interval $\bar I$, we have for $d \ge d(\eps, \eps'', \DENSITYBOUND, \beta_0, \beta_1, \MARGIN, \lambda)$ and with probability at least $1-\exp\big(- d / K(\eps, \DENSITYBOUND, \beta_0, \beta_1, \MARGIN, \lambda) \big)$,
\begin{align}\label{eq:sphericalconversion_upperbound_whpbound}
&\frac1d \log\bigl( \LHS \big) \in \Big[ \RSG_{\bar I}(n/d, \GAUSSIANMEASURE, x) - \eps, \RSG_{\bar I}(n/d, \GAUSSIANMEASURE, x) + \eps \Big]\, .
\end{align}
We apply the second part of Proposition \upref{prop:nicersuniquesol} identically as the proof in the $\MARGIN < 0$ case to convert the first argument of $\RSG_{\bar I}(n/d, \GAUSSIANMEASURE, x)$ to $\RSG_{\bar I}(\alpha, \GAUSSIANMEASURE, x)$ for $d \ge d(\eps, \DENSITYBOUND, \beta_0, \beta_1, \MARGIN, \lambda)$.
We obtain for the same $d$ and the same probability as \equpref{eq:sphericalconversion_upperbound_whpbound},
\begin{align*}
&\frac1d \log \mu_d\Big(\one\{\theta_1/\sqrt{\DIMENSION} \in I\} \cap C \Big) \le \GAUSSIANMEASURE w^2 + \RSG_{\bar I}(\alpha, \GAUSSIANMEASURE, x) + K\eps - \frac1d \log \AREA(\cS^{d-1}(\sqrt{d})) + \frac1d \log(2\sqrt{d})\,.
\end{align*}
We now take $\beta=\beta(\alpha, x)$ in the above, use Lemma \upref{lem:rs_formula_properties} and the fact that $0 \not\in \bar I$, and finish identically to the proof of the upper bound in the $\MARGIN < 0$ case.
Here we account for the presence of $\GAUSSIANMEASURE w^2$ rather than $\GAUSSIANMEASURE$ by noting $\beta (w^2-1) \le 3\eps'' \beta_1 \le 30 \eps$.
This proves \equpref{eq:spherical_conversion_upperbound} for the $\MARGIN \ge 0$ case.

\paragraph{Proof of \equpref{eq:spherical_conversion_lowerbound}:}
We let $w = 1 - \eps$. 
We define $\underH = \underH(\alpha, \MARGIN, \lambda)$, $\overH = \overH(\alpha, \MARGIN, \lambda)$ identically as \equpref{eq:underLdef}, \equpref{eq:overLdef} and define $\bar f$ identically as \equpref{eq:sphericalconversion_barf_def}.
We now let
\begin{align*}
\VARPROBLEMSPHERICAL = \VARPROBLEMSPHERICAL(\eps, \alpha, \beta_0, \beta_1, \MARGIN, \lambda) := \min\Big\{\frac14, \frac{1}{2\overH} \eps \Big\}\,.
\end{align*}
By definition of $\overH$ and as $\frac{\partial}{\partial \beta}\rho_0(\alpha, \beta, x)<0$ by Lemma \upref{lem:rs_formula_properties}, it follows that uniformly in $x \in [-1, 1]$,
\begin{align*}
\big| \rho_0(\alpha, \beta(\alpha, x), x) - \rho_0(\alpha, \beta(\alpha, x) + 2\VARPROBLEMSPHERICAL, x) \big| \le \eps\,,
\end{align*}
and hence by Lemma \upref{lem:rs_formula_properties},
\begin{align}
\rho_0(\alpha, \beta(\alpha, x) + 2\VARPROBLEMSPHERICAL, x) \ge 1-x^2-\eps\,.\label{eq:spherical_conversion_marginpositive_rhoineq}
\end{align}
Also note $2\VARPROBLEMSPHERICAL+\beta(\alpha, x) < \beta_1=10$, where we use that $\beta(\alpha, x) \le 9$. 
By \equpref{eq:fbarfirstderiv}, \equpref{eq:fbarsecondderiv}, and Lemma \upref{lem:rs_formula_properties}, we obtain
\begin{align}
\bar f(\alpha, \beta(\alpha, x) + \VARPROBLEMSPHERICAL, x) &\ge \bar f(\alpha, \beta(\alpha, x), x) + \VARPROBLEMSPHERICAL \frac{\partial}{\partial \beta} \bar f(\alpha, \beta(\alpha, x), x) + \frac{1}2 \underH \VARPROBLEMSPHERICAL^2 \notag \\
&= \bar f(\alpha, \beta(\alpha, x), x) - \VARPROBLEMSPHERICAL \big( \rho_0(\alpha, \beta(\alpha, x), x) + x^2 \big) + \frac{1}2 \underH\VARPROBLEMSPHERICAL^2 \notag \\
&= \bar f(\alpha, \beta(\alpha, x), x) - \VARPROBLEMSPHERICAL + \frac{1}2 \underH \VARPROBLEMSPHERICAL^2\,. \label{eq:spherical_conversion_marginpositive_ineq_1}
\end{align}
Similarly we obtain
\begin{align}
\bar f(\alpha, \beta(\alpha, x) + \VARPROBLEMSPHERICAL, x) &\ge \bar f(\alpha, \beta(\alpha, x) + 2\VARPROBLEMSPHERICAL, x) - (-\VARPROBLEMSPHERICAL) \big( \rho_0(\alpha, \beta(\alpha, x) + 2\VARPROBLEMSPHERICAL, x) + x^2 \big) + \frac{1}2 \underH \VARPROBLEMSPHERICAL^2 \notag \\
&\ge \bar f(\alpha, \beta(\alpha, x) + 2\VARPROBLEMSPHERICAL, x) + \VARPROBLEMSPHERICAL (1-\eps) + \frac{1}2 \underH \VARPROBLEMSPHERICAL^2 \notag \\
&\ge \bar f(\alpha, \beta(\alpha, x) + 2\VARPROBLEMSPHERICAL, x) + \VARPROBLEMSPHERICAL w^2 + \frac{1}2 \underH \VARPROBLEMSPHERICAL^2\,,\label{eq:spherical_conversion_marginpositive_ineq_w}
\end{align}
where we used \equpref{eq:spherical_conversion_marginpositive_rhoineq} and $w=1-\eps$ combined with $\eps<1$. 
Finally, we let 
\begin{align}
\eps_1, \eps_2 &:= \min\Big\{ \frac16 \underH \VARPROBLEMSPHERICAL^2, \eps \Big\}\,, \label{eq:spherical_conversion_marginpositive_eps12choice} \\
\eps_{\RS} &:= \min\big\{ \eps_{\RSG}(\eps_1, \DENSITYBOUND, \beta_0, \beta_1, \MARGIN, \lambda), \eps_{\RSG}(\eps_2, \DENSITYBOUND, \beta_0, \beta_1, \MARGIN, \lambda) \big\} > 0\,.\label{eq:sphericalconversion_epsRS_def_marginpositive}
\end{align}
We now turn to establishing \equpref{eq:spherical_conversion_lowerbound}. 
We consider any $\eps''$, $0 < \eps'' \le \eps_{\RS}$, and define $I$ as in \equpref{eq:sphericalconversion_I_lowerbound}. 
Thus $0 \not \in \bar I$ as $|x| > \eps''$.
We now let $\beta=\beta(\alpha, x)+\VARPROBLEMSPHERICAL$ for the rest of this proof.
First suppose $x \ge \eps'' > 0$. 
Here we have by \equpref{eq:intersection_containment_positivemargin}, \equpref{eq:sphericalconversioncontainment_smaller1_positivemargin} that
\begin{align}
\LHS &:= \int \one\big\{ \{\theta_1/\sqrt{\DIMENSION} \in \bar I\} \cap C \cap \{w^2 \DIMENSION \le \|\btheta\|^2 \le \DIMENSION\} \big\} \exp\big(-\GAUSSIANMEASURE \| \btheta \|^2\big)\, \rmd \btheta \notag \\
&\le \exp(-\GAUSSIANMEASURE w^2 \DIMENSION) \VOL\Big( \one\Big\{\theta_1/\sqrt{\DIMENSION} \in \Big[x-\eps'', w \cdot \frac{x+\eps''}w \Big] \Big\} \cap C \cap \{ w^2 \DIMENSION \le \| \btheta \| \le \DIMENSION \} \Big)\notag  \\
&\le \exp(-\GAUSSIANMEASURE w^2 \DIMENSION) \sqrt{\DIMENSION} \int_w^1 \AREA \Big( \one\Big\{\theta_1/\sqrt{\DIMENSION} \in v\Big[x-\eps'', \frac{x+\eps''}w \Big] \Big\} \cap vC \cap v\cS^{d-1}(\sqrt{d}) \Big)\, \rmd v \notag \\
&\le \exp(-\GAUSSIANMEASURE w^2 \DIMENSION)\frac1{\sqrt{\DIMENSION}} \AREA\Big( \one\{\theta_1/\sqrt{\DIMENSION} \in I \} \cap C \cap \cS^{d-1}(\sqrt{d}) \Big)\, .\label{eq:sphericalconversionlowerboundstep1_marginpositive}
\end{align}
Analogously, when $x \le -\eps''<0$, we have 
\begin{align*}
\LHS &\le \exp(-\GAUSSIANMEASURE w^2 \DIMENSION) \sqrt{\DIMENSION}\int_w^1 \AREA \Big(\one\Big\{\theta_1/\sqrt{\DIMENSION} \in \Big[w \cdot \frac{x-\eps''}w, x+\eps'' \Big] \Big\} \cap C \cap v\cS^{d-1}(\sqrt{d}) \Big)\, \rmd v \\
&\le \exp(-\GAUSSIANMEASURE w^2 \DIMENSION) \sqrt{\DIMENSION} \int_w^1 \AREA \Big( \one\Big\{\theta_1/\sqrt{\DIMENSION} \in v\Big[ \frac{x-\eps''}w, x+\eps'' \Big] \Big\} \cap vC \cap v\cS^{d-1}(\sqrt{d}) \Big)\, \rmd v \\
&\le \exp(-\GAUSSIANMEASURE w^2 \DIMENSION)\frac1{\sqrt{\DIMENSION}} \AREA\Big( \one\{\theta_1/\sqrt{\DIMENSION} \in I \} \cap C \cap \cS^{d-1}(\sqrt{d}) \Big)\,.
\end{align*}
(Note here that $0 \ge x-\eps'' > \frac{x-\eps''}w$.)
Now note
\begin{align*}
\LHS = \expressionI - \expressionII - \expressionIII \, ,
\end{align*}
where 
\begin{align*}
\expressionI &:= \int \one\big\{\{\theta_1/\sqrt{\DIMENSION} \in \bar I \} \cap C \big\} \exp(-\GAUSSIANMEASURE \| \btheta \|^2)\, \rmd \btheta\, ,\\
\expressionII &:= \int \one\big\{\{\theta_1/\sqrt{\DIMENSION} \in \bar I \} \cap C \cap \{ \| \btheta \| \ge \sqrt{\DIMENSION}\} \big\} \exp(-\GAUSSIANMEASURE \| \btheta \|^2)\, \rmd \btheta\, , \\
\expressionIII &:= \int \one\big\{\{\theta_1/\sqrt{\DIMENSION} \in \bar I \} \cap C \cap \{ \| \btheta \| \le w\sqrt{\DIMENSION}\} \big\}\exp(-\GAUSSIANMEASURE \| \btheta \|^2)\, \rmd \btheta\, .
\end{align*}
We again lower bound $\expressionI$ and upper bound $\expressionII, \expressionIII$ via Theorem \upref{thm:gaussianmeasurersformula}. First, by Theorem \upref{thm:gaussianmeasurersformula} with interval $\bar I$ and the second part of Proposition \upref{prop:nicersuniquesol}, as $\eps'' \le \eps_{\RS}$ and $\beta \in [\beta_0, \beta_1]$, we have for $d \ge d(\eps_1, \eps_2, \DENSITYBOUND, \beta_0, \beta_1, \MARGIN, \lambda)$ and probability at least $1 - \exp\big( - d / K(\eps_1, \eps_2, \DENSITYBOUND, \beta_0, \beta_1, \MARGIN, \lambda)\big)$,
\begin{align*}
\expressionI \ge \max\Big\{ \exp \Big( \big( \RSG_{\bar I}(\alpha, \GAUSSIANMEASURE, x) -  \eps_1 \big) d \Big), \exp \Big( \big( \RSG_{\bar I}(\alpha, \GAUSSIANMEASURE, x) -  \eps_2 \big) d \Big) \Big\}\, .
\end{align*}
By Theorem \upref{thm:gaussianmeasurersformula} and Proposition \upref{prop:nicersuniquesol}, for the same $d$ and with the same probability,
\begin{align*}
\expressionII &= \int \one\big\{\{\theta_1/\sqrt{\DIMENSION} \in \bar I \} \cap C \cap \{ \| \btheta \| \ge \sqrt{\DIMENSION}\} \big\} \exp(-(\GAUSSIANMEASURE-\VARPROBLEMSPHERICAL) \| \btheta \|^2) \exp(-\VARPROBLEMSPHERICAL \| \btheta \|^2) \, \rmd \btheta \\
&\le \exp(-\VARPROBLEMSPHERICAL d) \int \one\big\{\{\theta_1/\sqrt{\DIMENSION} \in \bar I \} \cap C \big\} \exp(-(\GAUSSIANMEASURE-\VARPROBLEMSPHERICAL) \| \btheta \|^2) \, \rmd \btheta \\
&\le \exp\Big( \big( \RSG_{\bar I}(\alpha, \beta - \VARPROBLEMSPHERICAL, x) - \VARPROBLEMSPHERICAL + \eps_1 \big) d \Big)\,.
\end{align*}
Analogously, for the same $d$ and with the same probability,
\begin{align*}
\expressionIII &= \int \one\big\{\{\theta_1/\sqrt{\DIMENSION} \in \bar I \} \cap C \cap \{ \| \btheta \| \le w\sqrt{\DIMENSION}\} \big\}\exp(- (\GAUSSIANMEASURE + \VARPROBLEMSPHERICAL) \| \btheta \|^2) \exp(\VARPROBLEMSPHERICAL \| \btheta \|^2)\, \rmd \btheta \\
&\le \exp\Big( \big( \RSG_{\bar I}(\alpha, \beta+\VARPROBLEMSPHERICAL, x) + \VARPROBLEMSPHERICAL w^2 + \eps_2 \big) d \Big)\,.
\end{align*}
Recall $\beta-\VARPROBLEMSPHERICAL = \beta(\alpha, x)$ and $\bar f = \RSG_{\bar I} - \frac12 \log(2\pi e)$ as $0 \not \in \bar I$. Thus applying \equpref{eq:spherical_conversion_marginpositive_ineq_1}, \equpref{eq:spherical_conversion_marginpositive_ineq_w} to upper bound $\expressionII$, $\expressionIII$ respectively gives for $d \ge d(\eps_1, \eps_2, \DENSITYBOUND, \beta_0, \beta_1, \MARGIN, \lambda)$ and probability at least $1 - \exp\big( - d / K(\eps_1, \eps_2, \DENSITYBOUND, \beta_0, \beta_1, \MARGIN, \lambda)\big)$, we have $\expressionII, \expressionIII \le \frac14 \expressionI$. Therefore as $\eps_1, \eps_2 \le \eps$, for the same $d$ and with the same probability,
\begin{align}
\LHS \ge \frac12 \expressionI \ge \frac12 \exp\Big( \big(\RSG_{\bar I}(\alpha, \beta, x) - \eps\big) d\Big)\,. \notag
\end{align}
Via the same steps as the $\MARGIN < 0$ case, combining the above display with \equpref{eq:sphericalconversionlowerboundstep1_marginpositive} yields that for the same $d$ and with the same probability,
\begin{align*}
\frac1d \log \mu_d\Big(\one\{\theta_1/\sqrt{\DIMENSION} \in I\} \cap C \Big) + \frac1d \log \AREA(\cS^{d-1}(\sqrt{d})) \ge \beta w^2 + \RSG_{\bar I}(\alpha, \beta, x) - \eps + \frac1d \log \Big(\frac{\sqrt{d}}2\Big)\,.
\end{align*}
Note that as $\beta \le \beta_1 = 10$, we have $\beta w^2 \ge \beta - 3\beta \eps \ge \beta - 30 \eps$.
Finally, as $0 \not\in \bar I$, we recall $\beta + \RSG_{\bar I}(\alpha, \beta, x) \ge \inf_{q \in [0,1)} \Phi(x, q)$ by Lemma \upref{lem:rs_formula_properties}. Combining these steps with the above display and taking $\eps \leftarrow \eps/K$ for universal constant $K$, \equpref{eq:spherical_conversion_lowerbound} follows.

\subsection{Proof of Theorem \upref{thm:posteriorrsformula} for posterior}\label{subsec:sphericalpf_posterior}
Much of this proof is identical to the above and we will only highlight the differences.
We again define $\bar I := [x-\eps'', x+\eps'']$ and follow the same notation as Appendices \upref{subsec:sphericalpf_negativemargin}, \upref{subsec:sphericalpf_positivemargin}.
Similarly as before, we parametrize $\R^d$ by the scalar $v = \|\theta\|/\sqrt{d}$ and $\theta \in \cS^{d-1}(\sqrt{d})$. 

By Lipschitzness of $u(x)$, instead of the comparison inequalities relating $vC$ and $C$, we directly study how $u\big( (s_i \theta_1 + \langle g_i , \bar\theta\rangle) / \sqrt{d} \big)$ changes in $v$ via the following Lemma.
\begin{lemma}\label{lem:posterior_spherical_continuity}
With probability at least $1 - \exp(-d/K_{\upref{lem:posterior_spherical_continuity}})$, we have for all $v \ge 0$,
\begin{align*}
\sup_{\theta \in \cS^{d-1}(\sqrt{d})} \Big| \sum_{i=1}^n u\big( S_i(v \theta) \big) - \sum_{i=1}^n u\big( S_i(\theta) \big) \Big| \le K_{\upref{lem:posterior_spherical_continuity}} |v-1| d\,,
\end{align*}
where $K_{\upref{lem:posterior_spherical_continuity}} \ge 1$ depends on $\lambda$. Here probability is over the $s_i, g_i$.
\end{lemma}
\begin{proof}[Proof of Lemma \upref{lem:posterior_spherical_continuity}]
Note $S_i(v\theta) = vS_i(\theta)$. As $u(\cdot)$ is $D$-Lipschitz, we have
\begin{align}\label{eq:posterior_spherical_lem_eq1}
&\sup_{\theta \in \cS^{d-1}(\sqrt{d})} \Big| \sum_{i=1}^n u\big( vS_i(\theta) \big) - \sum_{i=1}^n u\big( S_i(\theta) \big) \Big| \notag \\
&\quad \le D|v-1| \sup_{\theta \in \cS^{d-1}(\sqrt{d})} \sum_{i=1}^n  \big| S_i(\theta) \big| \notag \\
&\quad \le D|v-1| \Big(  \sup_{\theta \in \cS^{d-1}(\sqrt{d})} \sum_{i=1}^n |s_i \theta_1/\sqrt{d}| + \sup_{\theta \in \cS^{d-1}(\sqrt{d})} \sum_{i=1}^n \big| \langle g_i, \bar \theta \rangle / \sqrt{d} \big|\Big) \notag \\
&\quad := D|v-1| \big(\expressionI + \expressionII \big)\,.
\end{align}
We now establish high-probability upper bounds on $\expressionI, \expressionII$.
Starting with $\expressionI$, observe for any $\theta \in \cS^{d-1}(\sqrt{d})$ that $\theta_1/\sqrt{d} \le 1$. Hence $\sup_{\theta \in \cS^{d-1}(\sqrt{d})} \sum_{i=1}^n |s_i \theta_1/\sqrt{d}|  \le \sum_{i=1}^n |s_i|$. Recall the $s_i^2$ are sub-Exponential with parameter $\signalsubexp = K(\lambda)$, and hence the $|s_i|$ are $K(\lambda)$ sub-Gaussian. As $n \le 2d$, Hoeffding's Inequality implies with probability at least $1 - \exp \big(-d/K(\lambda) \big)$, we have
\begin{align}\label{eq:posterior_spherical_lem_eq2}
\expressionI \le K(\lambda) d\,.
\end{align}

We now upper bound $\expressionII$.
Since $\|\bar \theta\| \le \sqrt{d}$, $\expressionII$ is maximized for $\|\bar \theta\|=\sqrt{d}$. Shifting $d$ by 1, it suffices to consider
\begin{align*}
F(G) := \sup_{\theta \in \cS^{d-1}(\sqrt{d})} \sum_{i=1}^n \big| \langle g_i, \theta \rangle / \sqrt{d} \big|
\end{align*}
where now $g_i \sim N(0, I_d)$ and where $G$ refers to the $n \times d$ Gaussian matrix with rows $g_i$.

Consider two Gaussian matrices $G_0 = \{g_i^0\}_{i=1}^n, G_1 = \{g_i^1\}_{i=1}^n$ whose rows $g_i^0, g_i^1$ are i.i.d. from $N(0, I_d)$. If $F(G_1) \ge F(G_0)$, by continuity of $\sum_{i=1}^n \big| \langle g_i^1, \theta \rangle / \sqrt{d} \big|$ in $\theta$, we have for some $\theta' \in \cS^{d-1}(\sqrt{d})$,
\begin{align*}
F(G_1) - F(G_0) = \sum_{i=1}^n \big| \langle g_i^1, \theta' \rangle / \sqrt{d} \big| - F(G_0) &\le \sum_{i=1}^n \big| \langle g_i^1, \theta' \rangle / \sqrt{d} \big| - \sum_{i=1}^n \big| \langle g_i^0, \theta' \rangle / \sqrt{d} \big| \\
&\le \sum_{i=1}^n \big|\langle g_i^1 - g_i^0, \theta' \rangle / \sqrt{d} \big| \\
&\le \sum_{i=1}^n \| g_i^1 - g_i^0 \| \\
&\le \sqrt{n} \Big( \sum_{i=1}^n \sum_{j=1}^d (g_{ij}^1 - g_{ij}^0)^2 \Big)^{1/2}\,.
\end{align*}
An analogous bound applies if $F(G_0) \ge F(G_1)$. Therefore $F(G)$ is $\sqrt{n}$-Lipschitz as a function of $G$ measured w.r.t. Euclidean norm. Clearly the density of $G$ is 1 strongly log-concave. It follows by concentration of Lipschiz functions of strongly log-concave measures, Theorem \upref{thm:logconcaveconcentrationbasics}, that
\begin{align}\label{eq:posterior_spherical_lem_logconcaveconcentration}
\P\Big( F(G) \ge \E\big[ F(G) \big] + t) \le \exp\Big( - \frac{t^2}{K n} \Big)\,,
\end{align}
where probability is over $G$.
Next we upper bound $\E\big[ F(G) \big]$. Writing $u = (u_1, \ldots, u_n)$ where each $u_i \in \{\pm 1\}$, observe that 
\begin{align*}
F(G) = \frac1{\sqrt{d}} \sup_{\theta \in \cS^{d-1}(\sqrt{d}), u \in \{\pm 1\}^n } \sum_{i=1}^n u_i \langle g_i, \theta\rangle = \frac1{\sqrt{d}} \sup_{\theta \in \cS^{d-1}(\sqrt{d}), u \in \{\pm 1\}^n } \sum_{i=1}^n \sum_{j=1}^d u_i \theta_j g_{ij}\,.
\end{align*}
Define the centered Gaussian processes $X_{\theta, u} := \sum_{i=1}^n \sum_{j=1}^d u_i \theta_j g_{ij}$, $Y_{\theta, u} := \| u\|_2 \langle g, \theta \rangle + \|\theta\|_2 \langle h, u \rangle$, where $u = (u_1, \ldots, u_n)^T$, and where $g \sim N(0, I_d)$, $h \sim N(0, I_n)$ are independent.  
For any $\theta, \theta' \in \cS^{d-1}(\sqrt{d})$, $u, u' \in \{\pm 1\}^n$, observe
\begin{align*}
\E\big[ (X_{\theta, u} - X_{\theta', u'})^2 \big] &= \|\theta\|^2 \|u\|^2 + \|\theta'\|^2 \|u'\|^2 - 2 \langle u, u' \rangle \langle \theta, \theta' \rangle\,, \\
\E\big[ (Y_{\theta, u} - Y_{\theta', u'})^2 \big] &= 2\|\theta\|^2 \|u \|^2 + 2\|\theta'\|^2 \|u'\|^2 - 2\big( \|u \| \|u'\| \langle \theta, \theta' \rangle + \|\theta\|\|\theta'\| \langle u, u' \rangle \big)\,.
\end{align*}
As $\|\theta\|=\|\theta'\|=\sqrt{d}, \|u\|=\|u'\|=\sqrt{n}$ and $\big(n-\langle u,u'\rangle \big) \big( d - \langle \theta, \theta' \rangle \big) \ge 0$, it follows that 
\begin{align*}
\E\big[ (X_{\theta, u} - X_{\theta', u'})^2 \big] \le \E\big[ (Y_{\theta, u} - Y_{\theta', u'})^2 \big]\,.
\end{align*}
By the Sudakov-Fernique inequality (see e.g. Theorem 7.2.8, \cite{vershynin2018high}), noting $\|u\|=\sqrt{n}, \|\theta\|=\sqrt{d}$ always holds for $u \in \{\pm 1\}^n, \theta \in \cS^{d-1}(\sqrt{d})$, we obtain
\begin{align*}
\sqrt{d} \E\big[F(G) \big] &= \E\Big[ \sup_{\theta \in \cS^{d-1}(\sqrt{d}), u \in \{\pm 1\}^n } X_{\theta, u} \Big] \\
&\le \E\Big[ \sup_{\theta \in \cS^{d-1}(\sqrt{d}), u \in \{\pm 1\}^n } Y_{\theta, u} \Big] \\
&\le \sqrt{n} \E\big[  \sup_{\theta \in \cS^{d-1}(\sqrt{d})} \langle g, \theta \rangle\big] + \sqrt{d} \E\big[ \sup_{u \in \{\pm 1\}^n} \langle h, u\rangle \big] \\
&= \sqrt{nd} \E\big[ \| g\|_2 \big] + \sqrt{d} \E\big[ \|h\|_1 \big] \\
&= O(n^{1/2} d)\,.
\end{align*}
Using $n \le 2d$, it follows that $\E\big[ F(G) \big] \le O(d)$. Combining with \equpref{eq:posterior_spherical_lem_logconcaveconcentration} gives that with probability at least $1 - \exp(-d/K)$, we have 
\begin{align}\label{eq:posterior_spherical_lem_eq3}
\expressionII \le K d\,.
\end{align}
Now combining \equpref{eq:posterior_spherical_lem_eq1}, \equpref{eq:posterior_spherical_lem_eq2}, \equpref{eq:posterior_spherical_lem_eq3} with a Union Bound, the Lemma follows.
\end{proof}
We now have the technical ingredients in place to prove Theorem \upref{thm:posteriorrsformula}.

\paragraph{Proof of \equpref{eq:posterior_spherical_conversion_upperbound}:}
Let $K'(\lambda) := K_{\upref{lem:posterior_spherical_continuity}}(\lambda)$, where $K'$ is as in the statement of Theorem. We thus consider $\eps'' \le \min\{ \eps / K_{\upref{lem:posterior_spherical_continuity}}(\lambda), \eps_{\RSG}\}$ and let $I$ be as per \equpref{eq:posterior_sphericalconversion_I_upperbound}. Let $w=1+\eps''$. 

For any $\beta \in [\beta_0, \beta_1]$, by \equpref{eq:sphericalconversioncontainment_larger1_positivemargin}, we obtain the following, very similarly to the proof of \equpref{eq:spherical_conversion_upperbound}.
When $x > \eps''$, we have for $d \ge d(\eps'')$,
\begin{align*}
\LHS &:= \int \one\{\theta_1/\sqrt{\DIMENSION} \in \bar I\} \exp\Big( \sum_{i=1}^n u\big(S_i(\theta)\big)-\GAUSSIANMEASURE \| \btheta \|^2 \Big)\, \rmd \btheta \\
&\ge \int \one \Big\{ \big\{\theta_1/\sqrt{\DIMENSION} \in [x-\eps'', x+\eps''] \big\} \cap \big\{d \le \|\btheta \|^2 \le w^2 d \big\} \Big\} \exp\Big( \sum_{i=1}^n u\big(S_i(\theta)\big)-\GAUSSIANMEASURE \| \btheta \|^2 \Big)\, \rmd \btheta \\
&\ge \exp(-\GAUSSIANMEASURE w^2 \DIMENSION)\sqrt{\DIMENSION}\int_1^w \int_{\cS^{d-1}(\sqrt{d})} \one\big\{v\theta_1/\sqrt{\DIMENSION} \in [x-\eps'', x+\eps'']\big\} v^{d-1} \exp\Big( \sum_{i=1}^n u\big(S_i(v\theta)\big) \Big) \, \rmd \theta \rmd v \\
&\ge \exp(-\GAUSSIANMEASURE w^2 \DIMENSION)\sqrt{\DIMENSION}\int_1^w \int_{\cS^{d-1}(\sqrt{d})} \one\Big\{\theta_1/\sqrt{\DIMENSION} \in \Big[ x-\eps'', \frac{x+\eps''}{w} \Big] \Big\} v^{d-1} \exp\Big( \sum_{i=1}^n u\big(S_i(v\theta)\big) \Big) \, \rmd \theta \rmd v \,,
\end{align*}
where the last inequality uses \equpref{eq:sphericalconversioncontainment_larger1_positivemargin}.
When $x < -\eps''$, we have for $d \ge d(\eps'')$, now using \equpref{eq:sphericalconversioncontainment_larger1_positivemargin} in the $a \le b \le 0$ case,
\begin{align*}
\LHS \ge \exp(-\GAUSSIANMEASURE w^2 \DIMENSION)\sqrt{\DIMENSION}\int_1^w \int_{\cS^{d-1}(\sqrt{d})} \one\Big\{\theta_1/\sqrt{\DIMENSION} \in \Big[ \frac{x-\eps''}{w}, x+\eps''\Big] \Big\} v^{d-1} \exp\Big( \sum_{i=1}^n u\big(S_i(v\theta)\big) \Big) \, \rmd \theta \rmd v \,.
\end{align*}
By Lemma \upref{lem:posterior_spherical_continuity}, with probability at least $1-\exp\big( -d / K(\lambda) \big)$, we have for all $\theta \in \cS^{d-1}(\sqrt{d})$ and all $v, 1 \le v \le w$ that 
\begin{align*}
\Big| \sum_{i=1}^n u\big( S_i(v \theta) \big) - \sum_{i=1}^n u\big( S_i(\theta) \big) \Big| \le K_{\upref{lem:posterior_spherical_continuity}} |v-1| d \le \eps d\,,
\end{align*}
where we use that $1 \le v \le w \le 1+\eps/K_{\upref{lem:posterior_spherical_continuity}}$. Combining the above displays and recalling the definition of $I$ in \equpref{eq:posterior_sphericalconversion_I_upperbound} yields that for $d \ge d(\eps'')$, with probability at least $1-\exp\big( d / K(\lambda) \big)$,
\begin{align}
\LHS &\ge \exp \big(-\beta w^2 d - \eps d \big) \sqrt{d} \int_1^w \int_{\cS^{d-1}(\sqrt{d})} \one\big\{  \theta_1/\sqrt{\DIMENSION} \in I \big\} v^{d-1} \exp\Big(\sum_{i=1}^n u\big( S_i(\theta) \big) \Big)\, \rmd \theta \rmd v \notag \\
&\ge \frac{\exp \big(-\beta w^2 d - \eps d \big)}{2\sqrt{d}} \int_{\cS^{d-1}(\sqrt{d})} \one\big\{  \theta_1/\sqrt{\DIMENSION} \in I \big\} \exp\Big(\sum_{i=1}^n u\big( S_i(\theta) \big) \Big)\, \rmd \theta\,.\notag 
\end{align}
That is, 
\begin{align}\label{eq:posterior_spherical_upperbd_main}
\frac1d \log Z(I) \le \frac1d \log(\LHS) + \frac1d \log(2\sqrt{d}) + \beta w^2 + \eps\,.
\end{align}

With \equpref{eq:posterior_spherical_upperbd_main}, the rest of the proof of \equpref{eq:posterior_spherical_conversion_upperbound} now follows identically as the proof of \equpref{eq:spherical_conversion_upperbound}, the one difference being that we apply Theorem \upref{thm:posteriorgaussianmeasurersformula} rather than Theorem \upref{thm:posteriorrsformula} to control $\frac1d \log(\LHS)$.

\paragraph{Proof of \equpref{eq:posterior_spherical_conversion_lowerbound}:} Let $K_{\upref{lem:posterior_spherical_continuity}} = K_{\upref{lem:posterior_spherical_continuity}}(\lambda) \ge 1$, and $w=1- \frac{\eps}{K_{\upref{lem:posterior_spherical_continuity}}}$. 
Recall that Lemma \upref{lem:rs_formula_properties}, Proposition \upref{prop:nicersuniquesol} still hold here, where the RS equations are given in terms of $u(x)$ rather than the hard indicator $\one\{\cdot \ge \MARGIN\}$.
Thus we may perform the following steps.

In particular, define $\underH, \overH$ as in \equpref{eq:underLdef}, \equpref{eq:overLdef}. We now define 
\begin{align}\label{eq:posterior_sphericalconversion_barf_def}
\bar f(\alpha, \beta, x) := \bar\Phi\big(x, q_0(\alpha, \beta, x), \rho_0(\alpha, \beta, x)\big) - \beta\big(\rho_0(\alpha, \beta, x) + x^2)\,,
\end{align}
where consider $\bar\Phi(x,q,\rho)$ given with $u(x)$ (rather than the hard indicator function $\one\{\cdot \ge \MARGIN\}$ for the interpolators).
We now let 
\begin{align*}
\VARPROBLEMSPHERICAL = \VARPROBLEMSPHERICAL(\eps, \alpha, \beta_0, \beta_1, \lambda) := \min\Big\{\frac14, \frac{1}{2\overH K_{\upref{lem:posterior_spherical_continuity}}} \eps \Big\}\,.
\end{align*}
Analogously as the proof of \equpref{eq:spherical_conversion_marginpositive_rhoineq}, we have uniformly in $x \in [-1, 1]$:
\begin{align}
\rho_0(\alpha, \beta(\alpha, x) + 2\VARPROBLEMSPHERICAL, x) \ge 1-x^2-\frac{\eps}{K_{\upref{lem:posterior_spherical_continuity}}}\,.\label{eq:posterior_spherical_conversion_marginpositive_rhoineq}
\end{align}
Again analogously the proof of \equpref{eq:spherical_conversion_marginpositive_ineq_1}, \equpref{eq:spherical_conversion_marginpositive_ineq_w}, now using \equpref{eq:posterior_spherical_conversion_marginpositive_rhoineq} and $\eps<1, K_{\upref{lem:posterior_spherical_continuity}} \ge 1$, we obtain
\begin{align}
\bar f(\alpha, \beta(\alpha, x) + \VARPROBLEMSPHERICAL, x) &\ge \bar f(\alpha, \beta(\alpha, x), x) - \VARPROBLEMSPHERICAL + \frac{1}2 \underH \VARPROBLEMSPHERICAL^2\,, \label{eq:posterior_spherical_conversion_marginpositive_ineq_1} \\
\bar f(\alpha, \beta(\alpha, x) + \VARPROBLEMSPHERICAL, x) &\ge \bar f(\alpha, \beta(\alpha, x) + 2\VARPROBLEMSPHERICAL, x) + \VARPROBLEMSPHERICAL w^2 + \frac{1}2 \underH \VARPROBLEMSPHERICAL^2\,.\label{eq:posterior_spherical_conversion_marginpositive_ineq_w}
\end{align}
Finally, we again let
\begin{align}
\eps_1 = \eps_2 &:= \min\Big\{ \frac16 \underH \VARPROBLEMSPHERICAL^2, \eps \Big\}\,, \label{eq:posterior_spherical_conversion_marginpositive_eps12choice} \\
\eps_{\RS} &:= \min\big\{ \eps_{\RSG}(\eps_1, \DENSITYBOUND, \beta_0, \beta_1, \lambda), \eps_{\RSG}(\eps_2, \DENSITYBOUND, \beta_0, \beta_1, \lambda) \big\} > 0\,,\label{eq:posterior_sphericalconversion_epsRS_def}
\end{align}
where now $\eps_{\RSG}$ comes from Theorem \upref{thm:posteriorgaussianmeasurersformula}.

We now turn to establishing \equpref{eq:posterior_spherical_conversion_lowerbound}. 
Again, consider any $\eps''$, $0 < \eps'' \le \eps_{\RS}$, and define $I$ as in \equpref{eq:posterior_sphericalconversion_I_lowerbound}. 
Thus $0 \not \in \bar I$ as $|x| > \eps''$.
We now let $\beta=\beta(\alpha, x)+\VARPROBLEMSPHERICAL$ for the rest of this proof.
First suppose $x \ge \eps'' > 0$. 
Here we have by \equpref{eq:sphericalconversioncontainment_smaller1_positivemargin} that
\begin{align}
\LHS &:= \int \one\big\{ \{\theta_1/\sqrt{\DIMENSION} \in \bar I\} \cap \{w^2 \DIMENSION \le \|\btheta\|^2 \le \DIMENSION\} \big\} \exp\Big( \sum_{i=1}^n u\big( S_i(\theta) \big) -\GAUSSIANMEASURE \| \btheta \|^2\Big)\, \rmd \btheta \notag \\
&\le \exp(-\GAUSSIANMEASURE w^2 \DIMENSION) \sqrt{d} \int_w^1 \int_{\cS^{d-1}(\sqrt{d})} \one\big\{v \theta_1/\sqrt{\DIMENSION} \in [x-\eps'', x+\eps''] \big\} v^{d-1} \exp\Big( \sum_{i=1}^n u\big( S_i(v\theta) \big) \Big)\, \rmd \btheta \rmd v \notag \\
&\le \exp(-\GAUSSIANMEASURE w^2 \DIMENSION) \sqrt{\DIMENSION} \int_w^1 \int_{\cS^{d-1}(\sqrt{d})} \one\Big\{ \theta_1/\sqrt{d} \in \Big[ x-\eps'', \frac{x+\eps''}w \Big] \Big\} v^{d-1} \exp\Big( \sum_{i=1}^n u\big( S_i(v\theta) \big) \Big)\, \rmd \btheta \rmd v \, .\notag
\end{align}
Analogously, when $x \le -\eps''<0$, we have 
\begin{align*}
\LHS &\le \exp(-\GAUSSIANMEASURE w^2 \DIMENSION) \sqrt{\DIMENSION} \int_w^1 \one\Big\{\theta_1/\sqrt{\DIMENSION} \in \Big[ \frac{x-\eps''}w, x+\eps'' \Big] \Big\} v^{d-1} \exp\Big( \sum_{i=1}^n u\big( S_i(v\theta) \big) \Big)\, \rmd \btheta \rmd v \, .
\end{align*}
By Lemma \upref{lem:posterior_spherical_continuity}, with probability at least $1-\exp\big( -d / K(\lambda) \big)$, we have for all $\theta \in \cS^{d-1}(\sqrt{d})$ and all $v, w \le v \le 1$ that 
\begin{align*}
\Big| \sum_{i=1}^n u\big( S_i(v \theta) \big) - \sum_{i=1}^n u\big( S_i(\theta) \big) \Big| \le K_{\upref{lem:posterior_spherical_continuity}} |v-1| d \le \eps d\,,
\end{align*}
where we use that $1-\eps/K_{\upref{lem:posterior_spherical_continuity}}=w \le v \le 1$. 
Note as $K_{\upref{lem:posterior_spherical_continuity}} \ge 1$, $w \ge 1-\eps$. Therefore by the definition of $I$ in \equpref{eq:posterior_sphericalconversion_I_lowerbound},
\begin{align*}
\Big[ x-\eps'', \frac{x+\eps''}w \Big] \subseteq I \text{ when } x \ge \eps''\,,\, \Big[ \frac{x-\eps''}w, x+\eps'' \Big] \subseteq I \text{ when }x \le -\eps''\,.
\end{align*}
Combining the above displays now yields that with probability at least $1-\exp\big( d / K(\lambda) \big)$,
\begin{align}
\LHS &\le \exp \big(-\beta w^2 d + \eps d \big) \sqrt{d} \int_w^1 \int_{\cS^{d-1}(\sqrt{d})} \one\big\{  \theta_1/\sqrt{\DIMENSION} \in I \big\} v^{d-1} \exp\Big(\sum_{i=1}^n u\big( S_i(\theta) \big) \Big)\, \rmd \theta \rmd v \notag \\
&\le \frac{\exp \big(-\beta w^2 d + \eps d \big)}{\sqrt{d}} \int_{\cS^{d-1}(\sqrt{d})} \one\big\{  \theta_1/\sqrt{\DIMENSION} \in I \big\} \exp\Big(\sum_{i=1}^n u\big( S_i(\theta) \big) \Big)\, \rmd \theta\,.\notag 
\end{align}
That is, 
\begin{align}\label{eq:posterior_spherical_lowerbd_main}
\frac1d \log Z(I) \ge \frac1d \log(\LHS) + \frac1d \log \sqrt{d} + \beta w^2 - \eps\,.
\end{align}

With \equpref{eq:posterior_spherical_lowerbd_main}, the rest of the proof of \equpref{eq:posterior_spherical_conversion_lowerbound} now follows identically as the proof of \equpref{eq:spherical_conversion_lowerbound} in the $\MARGIN \ge 0$ case. We again write 
\begin{align*}
\LHS = \expressionI - \expressionII - \expressionIII\,,
\end{align*}
where now $\one\{C\}$ in the definition of $\expressionI, \expressionII, \expressionIII$ is replaced by $\exp\big( \sum_{i=1}^n u(S_i(\theta))\big)$. (Observe that in the proof of \equpref{eq:spherical_conversion_lowerbound}, $C$ takes the same form, for $u(x) = \one\{x \ge \MARGIN\}$.)
We now apply Theorem \upref{thm:posteriorgaussianmeasurersformula} to establish high-probability formulas for $\expressionI, \expressionII, \expressionIII$ up to accuracy $\eps_1, \eps_2$, where $\eps_1, \eps_2$ are as defined in \equpref{eq:posterior_spherical_conversion_marginpositive_eps12choice}. Combined with \equpref{eq:posterior_spherical_conversion_marginpositive_ineq_1}, \equpref{eq:posterior_spherical_conversion_marginpositive_ineq_w}, this proves that
\begin{align*}
\frac14 \expressionI \ge \expressionII, \expressionIII\,.
\end{align*}
Combining with \equpref{eq:posterior_spherical_lowerbd_main} and finishing as in the proof of \equpref{eq:spherical_conversion_lowerbound} now establishes \equpref{eq:posterior_spherical_conversion_lowerbound}.

\section{Technical results for proof of Theorems \ref{thm:mollifiedhamiltonianfreeenergy}, \ref{thm:posterior_mollifiedhamiltonianfreeenergy}}\label{sec:interpolationconcentration}
Here we prove the intermediate technical results Lemma \upref{lem:interpolationkeyconcentration}, Lemma \upref{lem:derivativeintbounded}, and Proposition \upref{prop:logconcaveconcentration} on the interpolation argument and concentration of measure, used in the proof of Theorems \upref{thm:mollifiedhamiltonianfreeenergy} and \upref{thm:posterior_mollifiedhamiltonianfreeenergy}. Throughout Appendix \upref{sec:interpolationconcentration}, for the interpolators we let $D$ be as stated in Theorem \upref{thm:mollifiedhamiltonianfreeenergy}, and for the posterior we recall that $|u^{l}| \le D=D(\lambda)$ for $l=1,2,3,4$.
Also here in Appendix \upref{sec:interpolationconcentration}, as per the proof for the interpolators, we will always write dependence on $\MARGIN$ in $K$ when it is present in the following proofs; note there is no such dependence for the posterior. 

We first prove Lemma \upref{lem:interpolationkeyconcentration}, which lets us `transfer' expectations of the form $\nu_{t,v}(f)$ for $0 \le v \le 1$ to expectations $\nu_{t,1}(f)$ when $v=1$, using some of the ideas from cavity in $\SAMPLES$. 
We then finish the cavity in $\SAMPLES$ argument and execute the cavity in $\DIMENSION$ argument to prove Lemma \upref{lem:derivativeintbounded}, which gives us control over the derivatives $\big| \frac{\rmd}{\rmd t} \nu_t\big( \theta_{\DIMENSION}^1 \theta_{\DIMENSION}^2 \big) \big|$, $\big| \frac{\rmd}{\rmd t} \nu_t\big( (\theta_{\DIMENSION}^1)^2 \big) \big|$ arising in the interpolation argument.
Both these arguments will crucially use Proposition \upref{prop:logconcaveconcentration} on concentration of the interpolating Gibbs measure $\nu_{t,v}(\cdot)$, which we prove later in Appendix \upref{subsec:logconcaveconcentration}. 
The proof of Proposition \upref{prop:logconcaveconcentration} will crucially use classical results on concentration of Lipschitz functions of strongly log-concave measures; note the quadratic component $-\beta \| \theta \|^2$ we added when relaxing to the Gaussian measure is strongly concave and $u(\cdot)$ in Theorem \upref{thm:mollifiedhamiltonianfreeenergy} is concave.

\subsection{Proof of Lemma \upref{lem:interpolationkeyconcentration}}\label{subsec:pf_interpolationkeyconcentration}
Our first step will be to establish the following on the derivative of $\nu_{t,v}$ w.r.t. $v$. 
\begin{lemma}\label{lem:cavityinMintegrationbyparts}
Consider any $t \in [0,1]$ and a function $B_v$ that is either 1, $u'(S_{\SAMPLES, t, v}^1) u'(S_{\SAMPLES, t, v}^2)$, or $u''(S_{\SAMPLES, t, v}^1) + u'(S_{\SAMPLES, t, v}^1)^2$. Then for any $f:(\R^{\DIMENSION-1})^\REPLICAS \rightarrow \R$, $\Big| \frac{\partial}{\partial v} \nu_{t,v}(fB_v) \Big|$ is upper bounded by
\begin{align*}
&K\Bigg( \sum_{\REPLICAIDX \le \REPLICAS+1} \nu_{t,v}\big( \big|f \big| \big|R_{\REPLICAIDX,\REPLICAIDX} - \rho_{\SAMPLES, \DIMENSION}\big| \big) + \sum_{1 \le \REPLICAIDX_1 \neq \REPLICAIDX_2 \le \REPLICAS+2} \nu_{t,v}\big( \big|f \big| \big|R_{\REPLICAIDX_1,\REPLICAIDX_2} - q_{\SAMPLES, \DIMENSION}\big| \big) + \Big(\eps''+\frac1{\DIMENSION} \Big) \nu_{t,v}(f^2)^{1/2} \Bigg)\,,
\end{align*}
where $K$ depends on $D, \GAUSSIANMEASURE_0, \GAUSSIANMEASURE_1, \SIGNALSHORTDEPENDENCE, \REPLICAS$, and also on $\MARGIN$ for the interpolators.
\end{lemma}
For proving Lemma \upref{lem:interpolationkeyconcentration} we will only need Lemma \upref{lem:cavityinMintegrationbyparts} with $B_v = 1$. 
However the full generality where $B_v$ equals 1, $u'(S_{\SAMPLES, t, v}^1) u'(S_{\SAMPLES, t, v}^2) = u'(S_v^1) u'(S_v^2)$, or $u''(S_v^1) + u'(S_v^1)^2$ will be needed in proving Lemma \upref{lem:derivativeintbounded} later, specifically in the proof of Lemma \upref{lem:cavityinMrbarr}.
\begin{proof}[Proof of Lemma \upref{lem:cavityinMintegrationbyparts}]
Throughout this proof, we use the shorthand $S_v^{\ell} = S_{n,t,v}^{\ell}$. Recall
\begin{align*}
\nu_{t,v}(f B_v) &= \E \frac{\big\langle f B_v \exp\big( \sum_{\REPLICAIDX \le \REPLICAS} u(S^{\REPLICAIDX}_v) \big) \big\rangle_{t,\sim} }{\big\langle \exp\big( u(S^{\REPLICAS+1}_v) \big) \big\rangle_{t,\sim}^\REPLICAS } = \E \Big[ \frac{g(v)}{h(v)} \Big]\,,
\end{align*}
where we let 
\begin{align*}
g(v) := \Big\langle f B_v \exp\big( \sum_{\REPLICAIDX \le \REPLICAS} u(S^{\REPLICAIDX}_v) \big) \Big\rangle_{t,\sim}\,,\, h(v) := \Big\langle \exp\big( u(S^{\REPLICAS+1}_v) \big) \Big\rangle_{t,\sim}^\REPLICAS\,.
\end{align*}
Note the denominator of the Gibbs average defining $\langle \cdot \rangle_{t, \sim}$ does not depend on $v$. We thus compute
\begin{align*}
\frac{\partial}{\partial v} h(v) = \REPLICAS \Big\langle \exp\big( u(S^{\REPLICAS+1}_v) \big) \Big\rangle_{t,\sim}^{\REPLICAS-1}\, \Big\langle \exp\big( u(S^{\REPLICAS+1}_v) \big) u'(S^{\REPLICAS+1}_v) \cdot \frac{\partial}{\partial v} S^{\REPLICAS+1}_v \Big\rangle_{t,\sim}\,,
\end{align*}
and
\begin{align*}
\frac{\partial}{\partial v} g(v) = \Big\langle fB_v \exp\Big( \sum_{\REPLICAIDX \le \REPLICAS} u(S^{\REPLICAIDX}_v) \Big) \Big( \sum_{\REPLICAIDX \le \REPLICAS} u'(S^{\REPLICAIDX}_v) \cdot \frac{\partial}{\partial v} S^{\REPLICAIDX}_v + r(v) \Big) \Big\rangle_{t,\sim} \,,
\end{align*}
where 
\begin{align*}
r(v) := \begin{cases}  0 &: B_v = 1\,, \\
u''(S^1_v) u'(S^2_v) \cdot \frac{\partial}{\partial v} S^1_v + u'(S^1_v) u''(S^2_v) \cdot \frac{\partial}{\partial v} S^2_v &: B_v = u'(S_v^1) u'(S_v^2)\,, \\
u'''(S^1_v) \cdot \frac{\partial}{\partial v} S^1_v + 2u'(S^1_v) u''(S^1_v) \cdot \frac{\partial}{\partial v} S^1_v &: B_v = u''(S_v^1) + u'(S_v^1)^2\,.
\end{cases}
\end{align*}
Note
\begin{align*}
\frac{\partial}{\partial v}S^{\REPLICAIDX}_v &= \Big( \frac{\theta_1^{\REPLICAIDX}}{\sqrt{\DIMENSION}} - x\Big)s_{\SAMPLES} + \frac1{2\sqrt{v\DIMENSION}} \Big( \sum_{2 \le \DIMENSIONIDX \le \DIMENSION-1} g_{\SAMPLES,\DIMENSIONIDX} \theta_\DIMENSIONIDX^\REPLICAIDX + \sqrt{t} g_{\SAMPLES,\DIMENSION} \theta_{\DIMENSION}^\REPLICAIDX \Big)  - \frac1{2\sqrt{1-v}} \big( \sqrt{q_{\SAMPLES, \DIMENSION}}\QUENCHEDDISORDER + \sqrt{\rho_{\SAMPLES, \DIMENSION} - q_{\SAMPLES, \DIMENSION}}\ANNEALEDDISORDER^\REPLICAIDX \big)\,.
\end{align*}
Collecting common terms into the `signal part' $\expressionI$ corresponding to $s_{\SAMPLES}$ and the rest of the expression $\expressionII$ for $\frac{\partial}{\partial v}S^{\REPLICAIDX}_v$, we obtain
\begin{align*}
\frac{\partial}{\partial v}\nu_{t,v}(f B_v) &= \expressionI + \expressionII\,,
\end{align*}
where $\expressionI, \expressionII$ are defined by
\begin{align*}
\expressionI &:= \nu_{t,v}\Bigg( fB_v \Bigg(\sum_{\REPLICAIDX \le \REPLICAS} u'(S^{\REPLICAIDX}_v) \cdot \Big(\frac{\theta_1^\REPLICAIDX}{\sqrt{\DIMENSION}} - x \Big) s_{\SAMPLES}  + r_1(v) \Bigg)\Bigg) - \REPLICAS \nu_{t,v}\Bigg( fB_v \cdot u'(S_v^{\REPLICAS+1}) \cdot \Big( \frac{\theta_1^{\REPLICAS+1}}{\sqrt{\DIMENSION}}-x \Big)s_{\SAMPLES} \Bigg)\,,
\end{align*}
where 
\begin{align*}
r_1(v) &:= \begin{cases} 0 &: B_v = 1\,, \\
u''(S^1_v) u'(S^2_v) \cdot \Big(\frac{\theta_1^1}{\sqrt{\DIMENSION}} - x \Big)s_{\SAMPLES}  + u'(S^1_v) u''(S^2_v) \cdot \Big(\frac{\theta_1^2}{\sqrt{\DIMENSION}} - x \Big)s_{\SAMPLES}  &: B_v = u'(S^1_v) u'(S^2_v)\,, \\ 
\big( u'''(S^1_v) + 2u'(S^1_v) u''(S^1_v) \big) \cdot \Big(\frac{\theta_1^1}{\sqrt{\DIMENSION}} - x \Big)s_{\SAMPLES} &: B_v = u''(S^1_v) + u'(S^1_v)^2\,,
\end{cases}
\end{align*}
and
\begin{align*}
\expressionII &:= \nu_{t,v}\Bigg( fB_v \Big( \sum_{\REPLICAIDX \le \REPLICAS} u'(S^{\REPLICAIDX}_v) \bar{S}_v^{\REPLICAIDX} + r_2(v) \Big) \Bigg) - \REPLICAS \nu_{t,v}\Big( fB_v u'(S^{\REPLICAS+1}_v) \bar{S}_v^{\REPLICAS+1}\Big)\,,
\end{align*}
where
\begin{align*}
\bar{S}_v^\REPLICAIDX &:= \frac1{2\sqrt{v\DIMENSION }} \Big( \sum_{2 \le \DIMENSIONIDX \le \DIMENSION-1} g_{\SAMPLES, \DIMENSIONIDX} \theta_j^{\REPLICAIDX} + \sqrt{t} g_{\SAMPLES, \DIMENSION} \theta_{\DIMENSION}^\REPLICAIDX \Big) - \frac1{2\sqrt{1-v}} \big( \sqrt{q_{\SAMPLES, \DIMENSION}}\QUENCHEDDISORDER + \sqrt{\rho_{\SAMPLES, \DIMENSION} - q_{\SAMPLES, \DIMENSION}}\ANNEALEDDISORDER^\REPLICAIDX \big)\,, \\
r_2(v) &:= \begin{cases} 0 &: B_v = 1\,, \\
u''(S^1_v) u'(S^2_v) \bar{S}^1_v + u'(S^1_v) u''(S^2_v) \bar{S}^2_v &: B_v = u'(S^1_v) u'(S^2_v)\,, \\ 
\big( u'''(S^1_v) + 2u'(S^1_v) u''(S^1_v) \big) \bar{S}^1_v &: B_v = u''(S^1_v) + u'(S^1_v)^2\,.
\end{cases}
\end{align*}
First we upper bound $\big|\expressionI\big|$, via an analogous derivation as in the proof of Proposition \upref{prop:interpolationperconstraint} (specifically, the derivation of the upper bound on $\big|\expressionI\big|$ therein).
Since $\exp(\cdot) \ge 0$, we have 
\begin{align*}
\big|\expressionI\big| \le K(L) \sum_{1 \le \REPLICAIDX \le \REPLICAS+1} \E \frac{ \Big\langle \big|f s_n \big( \frac{\theta_1^{\REPLICAIDX}}{\sqrt{d}} - x \big) h_{\REPLICAIDX}\big| \cdot \exp \big( \sum_{1 \le \REPLICAIDX \le \REPLICAS+1} u(S_v^\REPLICAIDX) \big) \Big\rangle_{t, \sim}  }{ \big( \big\langle \exp u(S_{\SAMPLES, t, v}^{\REPLICAS+2})\big\rangle_{t, \sim} \big)^{\REPLICAS+1} }\,,
\end{align*}
where $h_{\REPLICAIDX}$ is a sum of products of $u^{(l')}(s_v^{1}), \ldots, u^{(l')}(s_v^{\REPLICAS+1})$ for $l'=1,2,3$, given by the above definition for $\expressionI$. Thus we have $|h_{\REPLICAIDX}| \le K(D, L)$ by our bounds on the first 4 derivatives of $u(\cdot)$. Due to the indicator $\one\{\theta^\REPLICAIDX_1/\sqrt{d} \in I\}$ in the definition of the Gibbs measure $\langle \cdot \rangle_{t, \sim}$ from \equpref{eq:interpolationtildemeasure}, which is applied for all $1 \le \REPLICAIDX \le \REPLICAS+1$ since each $f s_n \big( \frac{\theta_1^{\REPLICAIDX}}{\sqrt{d}} - x \big) h_{\REPLICAIDX}$ is a function of all the $\REPLICAS+1$ replicas, we can upper bound the numerator in the above expression as follows:
\begin{align*}
\Big\langle \Big|f s_n \Big( \frac{\theta_1^{\REPLICAIDX}}{\sqrt{d}} - x \Big) h_{\REPLICAIDX}\Big| \cdot \exp \Big( \sum_{1 \le \REPLICAIDX \le \REPLICAS+1} u(S_v^\REPLICAIDX) \Big) \Big\rangle_{t, \sim} &\le \eps'' K(D, L) \cdot \Big\langle \big| f \big| \big| s_n \big| \exp \Big( \sum_{1 \le \REPLICAIDX \le \REPLICAS+1} u(S_v^\REPLICAIDX) \Big) \Big\rangle_{t, \sim}\,.
\end{align*}
Thus by Cauchy-Schwarz,
\begin{align}
\big|\expressionI\big| &\le \eps'' K(D, L) \cdot  \nu_{t,v}\big( |f| |s_{\SAMPLES}| \big) \le \eps'' K(D, \REPLICAS) \cdot \nu_{t,v}(s_{\SAMPLES}^2)^{1/2} \nu_{t,v}(f^2)^{1/2}\,.\notag
\end{align}
Since $s_n$ does not depend on the replicas $\theta^1, \ldots, \theta^L$ or the annealed disorder $W^1, \ldots, W^L$, we have $\nu_{t,v}(s_n^2)^{1/2} = \E[s_n^2]^{1/2} \le K(\lambda)$. Combining with the above yields
\begin{align}\label{eq:int_by_parts_boundI}
\big|\expressionI\big| \le \eps'' K(D, \SIGNALSHORTDEPENDENCE, \REPLICAS) \cdot \nu_{t,v}(f^2)^{1/2}\,.
\end{align}

Next we upper bound $\big|\expressionII\big|$. We simplify $\expressionII$ via Guassian Integration by Parts w.r.t. the $g_{\SAMPLES,\DIMENSIONIDX}, g_{\SAMPLES,\DIMENSION}, z, \ANNEALEDDISORDER^\REPLICAIDX$. 
We give the computation in the case where $B_v = u''(S^1_v) + u'(S^1_v)^2$, the computation when $B_v = 1$ or $B_v = u'(S_v^1) u'(S_v^2)$ is analogous. This computation is analogous to the proof of Lemma 3.2.8, \cite{talagrand2010mean}; this is because the derivatives of $S_v^{\REPLICAIDX}$ w.r.t. the $g_{n,j}, g_{n,d}, z, \ANNEALEDDISORDER^{\REPLICAIDX}$ are the same as when the $S_v^{\REPLICAIDX}$ are defined without the signal part. A similar computation is detailed on the proof of Lemma 2.3.2, p.\ 172 of \cite{talagrand2010mean}. Spelling out the computation in full detail when $B_v = u''(S^1_v) + u'(S^1_v)^2$, we obtain
\begin{equation}\label{eq:biginterpolation_intbyparts1}
\begin{aligned}
&\nu_{t,v}\Bigg( fB_v \Big( \sum_{\REPLICAIDX \le \REPLICAS} u'(S^{\REPLICAIDX}_v) \bar{S}_v^{\REPLICAIDX} + r_2(v) \Big) \Bigg) \\
&\quad = \frac12 \nu_{t,v}\Big( f (R_{1,1} - \rho_{n,d}) u'(S_v^1) \big(u'''(S_v^1) + 2u'(S_v^1)u''(S_v^1)\big) \Big)  \\
&\quad\qquad + \frac12 \sum_{1\le \REPLICAIDX \le \REPLICAS}\nu_{t,v}\Big( f(R_{1,1}-\rho_{n,d}) \big( u''(S_v^{\REPLICAIDX})+u'(S_v^{\REPLICAIDX})^2\big) B_v \Big) \\
&\quad\qquad + \frac12 \sum_{2 \le \REPLICAIDX \le \REPLICAS} \nu_{t,v}\Big( f(R_{1,l}-q_{n,d}) u'(S_v^{\REPLICAIDX}) \big( u'''(S_v^1) + 2u'(S_v^1) u''(S_v^1)\big) \Big) \\
&\quad\qquad + \frac12 \sum_{1 \le \REPLICAIDX \neq \REPLICAIDX' \le L} \nu_{t,v}\Big( f (R_{l, l'} - q_{n,d}) u'(S_v^{\REPLICAIDX}) u'(S_v^{\REPLICAIDX'}) B_v \Big) \\
&\quad\qquad - \frac{L}2 \sum_{1 \le \REPLICAIDX \le L} \nu_{t,v}\Big( f (R_{l,L+1}-q_{n,d}) u'(S_v^{\REPLICAIDX}) u'(S_v^{L+1}) B_v \Big) + \err_1\,,
\end{aligned}
\end{equation}
and
\begin{equation}\label{eq:biginterpolation_intbyparts2}
\begin{aligned}
\nu_{t,v}\Big( fB_v r_2(v) \Big) &= \frac12 \nu_{t,v}\Bigg( f (R_{1,1} - \rho_{n,d}) \Big( \big(u'''(S_v^1) + 2u'(S_v^1)u''(S_v^1) \big)^2 \\
&\qquad\qquad \qquad \qquad\qquad\qquad + \big(u''''(S_v^1) + 2 u''(S_v^1)^2 + 2u'(S_v^1) u'''(S_v^1) \big) B_v \\
&\qquad\qquad \qquad \qquad\qquad\qquad + \big( u'''(S_v^1) + 2u'(S_v^1) u''(S_v^1) \big) u'(S_v^1) B_v \Big) \Bigg) \\
&\qquad + \frac12 \sum_{2 \le \REPLICAIDX \le \REPLICAS} \nu_{t,v}\Big( f (R_{1,l} - q_{n,d}) \big( u'''(S_v^1)+2u'(S_v^1) u''(S_v^1) \big) u'(S_v^1) B_v \Big) \\
&\qquad - \frac{L}2 \nu_{t,v}\Big( f (R_{1,L+1} - q_{n,d}) \big( u'''(S_v^1)+2u'(S_v^1) u''(S_v^1) \big) u'(S_v^{L+1})\Big) + \err_2\,,
\end{aligned}
\end{equation}
and
\begin{equation}\label{eq:biginterpolation_intbyparts3}
\begin{aligned}
\nu_{t,v}\Big( fB_v u'(S^{\REPLICAS+1}_v) \bar{S}_v^{\REPLICAS+1}\Big) &= \frac12 \nu_{t,v}\Big( f (R_{L+1,L+1}-\rho_{n,d}) \big(u''(S_v^{L+1} + u'(S_v^L)^2) \big) B_v \Big) \\
&\qquad + \frac12 \nu_{t,v}\Big( f (R_{1,L+1}-q_{n,d}) u'(S_v^{L+1}) \big( u'''(S_v^1)+2u'(S_v^1)u''(S_v^1) \big) \Big) \\
&\qquad + \frac12 \sum_{1 \le \REPLICAIDX \le L} \nu_{t,v}\Big( f (R_{l,L+1}-q_{n,d}) u'(S_v^{\REPLICAIDX}) u'(S_v^{L+1}) B_v \Big) \\
&\qquad - \frac{L+1}2 \nu_{t,v}\Big( f (R_{L+1,L+2}-q_{n,d}) u'(S_v^{L+1}) u'(S_v^{L+2}) B_v \Big) + \err_3\,,
\end{aligned}
\end{equation}
where $\err_1, \err_2, \err_3$ are such that
\begin{align*}
|\err_1|, |\err_2|, |\err_3| \le \frac{1-t}{2d} \cdot K\nu_{t,v}\big( |f| \max_{1 \le \REPLICAIDX \le L+2} (\theta_d^\REPLICAIDX)^2 \big)\,.
\end{align*}
In particular, $\err_1, \err_2, \err_3$ arise because the last coordinate $\theta_d^{\REPLICAIDX}$ is interpolated via the parameter $t \in [0,1]$.
By the second part of Proposition \upref{prop:logconcaveconcentration} and Cauchy-Schwarz we have 
\begin{align}\label{eq:bigintbyparts_errbound}
|\err_1|, |\err_2|, |\err_3| \le \frac{K(D, \beta_0, \beta_1, \MARGIN, \SIGNALSHORTDEPENDENCE, \REPLICAS)}d \cdot \nu_{t,v}(f^2)^{1/2}\,.
\end{align}
Recall that 
\begin{align*}
\expressionII = \nu_{t,v}\Bigg( fB_v \Big( \sum_{\REPLICAIDX \le \REPLICAS} u'(S^{\REPLICAIDX}_v) \bar{S}_v^{\REPLICAIDX} + r_2(v) \Big) \Bigg) - \REPLICAS \nu_{t,v}\Big( fB_v u'(S^{\REPLICAS+1}_v) \bar{S}_v^{\REPLICAS+1}\Big)\,.
\end{align*} 
Thus combining \equpref{eq:biginterpolation_intbyparts1}, \equpref{eq:biginterpolation_intbyparts2}, \equpref{eq:biginterpolation_intbyparts3} and using the bound \equpref{eq:bigintbyparts_errbound} on $|\err_1|, |\err_2|, |\err_3|$, we obtain
\begin{align*}
\Big| \expressionII \Big| &\le K(D, \GAUSSIANMEASURE_0, \GAUSSIANMEASURE_1, \MARGIN, \SIGNALSHORTDEPENDENCE, \REPLICAS)\Big( \sum_{\REPLICAIDX \le \REPLICAS+1} \nu_{t,v}\big( \big|f \big| \big|R_{\ell, \ell} - \rho_{\SAMPLES, \DIMENSION}\big| \big) \\
&\qquad \qquad \qquad \qquad \qquad \qquad \qquad  + \sum_{1 \le \REPLICAIDX_1 \neq \REPLICAIDX_2 \le \REPLICAS+2} \nu_{t,v}\big( \big|f \big| \big|R_{\REPLICAIDX_1,\REPLICAIDX_2} - q_{\SAMPLES, \DIMENSION}\big| \big) + \frac1{\DIMENSION} \cdot \nu_{t,v}(f^2)^{1/2} \Big)\,.
\end{align*}
Again, we obtain an analogous bound as the above when $B_v = 1$ or $B_v = u'(S_v^1) u'(S_v^2)$, with the same arguments. Combining this bound with the bound \equpref{eq:int_by_parts_boundI} on $|\expressionI|$ proves Lemma \upref{lem:cavityinMintegrationbyparts}.
\end{proof}

We will also use the following variant of Gr\"{o}nwall's Inequality:
\begin{lemma}[Lemma A.13.1 in \cite{talagrand2010mean}]\label{lem:gronwallvariant}
Suppose a differentiable function $g:[0,1] \rightarrow \R$, $g \ge 0$ is such that $\Big|\frac{\rmd}{\rmd t} g(t) \Big| \le c_1 g(t) + c_2$ for some $c_1,c_2 \ge 0$. Then
\begin{align*}
g(t) \le \exp\big( c_1(1-t) \big) \Big( g(1) + \frac{c_2}{c_1} \Big)\,.
\end{align*}
\end{lemma}

Now we return to the proof of Lemma \upref{lem:interpolationkeyconcentration}. 
By the second part of Proposition \upref{prop:logconcaveconcentration}, 
\begin{align*}
\nu_{t,v}\Big( \Big\{ \big|R_{\ell, \ell} - \rho_{\SAMPLES, \DIMENSION}\big| \ge K(D, \GAUSSIANMEASURE_0, \GAUSSIANMEASURE_1, \MARGIN, \SIGNALSHORTDEPENDENCE) \Big\} \Big) &\le K\exp(-N)\,, \\
\nu_{t,v}\Big( \Big\{ \big|R_{\REPLICAIDX_1,\REPLICAIDX_2} - q_{\SAMPLES, \DIMENSION}\big| \ge K(D, \GAUSSIANMEASURE_0, \GAUSSIANMEASURE_1,  \MARGIN, \SIGNALSHORTDEPENDENCE) \Big\} \Big) &\le K\exp(-N)\text{ for }\REPLICAIDX_1 \neq \REPLICAIDX_2\,.
\end{align*}
Here to bound $\rho_{n,d}, q_{n,d}$ we apply Proposition \upref{prop:logconcaveconcentration} with $t=v=1$. Additionally, the second pat of Proposition \upref{prop:logconcaveconcentration} combined with Lemma \upref{lem:mgftomomentslemma} implies that
\begin{align*}
\nu_{t,v}\big( (R_{\ell, \ell} - \rho_{\SAMPLES, \DIMENSION})^4 \big)^{1/4}\,,\,\nu_{t,v}\big( (R_{\REPLICAIDX_1,\REPLICAIDX_2} - q_{\SAMPLES, \DIMENSION})^4 \big)^{1/4} \le K(D, \GAUSSIANMEASURE_0, \GAUSSIANMEASURE_1, \MARGIN, \SIGNALSHORTDEPENDENCE)\,.
\end{align*}
Now H\"{o}lder's Inequality implies that
\begin{align*}
&\nu_{t,v}\big( |f| \big|R_{\ell, \ell} - \rho_{\SAMPLES, \DIMENSION}\big| \big) \\
&\quad \le K(D,\GAUSSIANMEASURE_0,\GAUSSIANMEASURE_1,  \MARGIN, \SIGNALSHORTDEPENDENCE)\nu_{t,v}\big( |f| \one\{ |R_{\ell, \ell} - \rho_{\SAMPLES, \DIMENSION}| \le K(D, \GAUSSIANMEASURE_0, \GAUSSIANMEASURE_1,  \MARGIN, \SIGNALSHORTDEPENDENCE) \} \big) \\
&\qquad \qquad + K(D,\GAUSSIANMEASURE_0,\GAUSSIANMEASURE_1,  \MARGIN, \SIGNALSHORTDEPENDENCE) \nu_{t,v}\big( f^2 \big)^{1/2}  \nu_{t,v}\big( \big|R_{\ell, \ell} - \rho_{\SAMPLES, \DIMENSION}\big|^4 \big)^{1/4}  \nu_{t,v}\big( \one\{ |R_{\ell, \ell} - \rho_{\SAMPLES, \DIMENSION}| \le K(D, \GAUSSIANMEASURE_0, \GAUSSIANMEASURE_1,  \MARGIN, \SIGNALSHORTDEPENDENCE)\} \big)^{1/4} \\
&\quad \le K(D,\GAUSSIANMEASURE_0,\GAUSSIANMEASURE_1,  \MARGIN, \SIGNALSHORTDEPENDENCE) \nu_{t,v}\big( |f| \big) + K(D,\GAUSSIANMEASURE_0,\GAUSSIANMEASURE_1,  \MARGIN, \SIGNALSHORTDEPENDENCE) \exp(-N) \nu_{t,v}\big( f^2 \big)^{1/2}\, .
\end{align*}
The same argument as above establishes the same upper bound on $\nu_{t,v}\big( |f| \big|R_{\REPLICAIDX_1,\REPLICAIDX_2} - q_{\SAMPLES, \DIMENSION}\big| \big)$.

Now applying Lemma \upref{lem:cavityinMintegrationbyparts} with $B_v = 1$, and combining with the above, we obtain that
\begin{align*}
\Big| \frac{\partial}{\partial v} \nu_{t,v}(f) \Big| &\le K(D, \GAUSSIANMEASURE_0, \GAUSSIANMEASURE_1, \MARGIN, \SIGNALSHORTDEPENDENCE, \REPLICAS)\Bigg( \sum_{\REPLICAIDX \le \REPLICAS+1} \nu_{t,v}\big( |f| \big|R_{\ell, \ell} - \rho_{\SAMPLES, \DIMENSION}\big| \big) \\
&\qquad \qquad \qquad \qquad \qquad \qquad \qquad + \sum_{1 \le \REPLICAIDX_1 \neq \REPLICAIDX_2 \le \REPLICAS+2} \nu_{t,v}\big( |f| \big|R_{\REPLICAIDX_1,\REPLICAIDX_2} - \rho_{\SAMPLES, \DIMENSION}\big| \big) \\
&\qquad \qquad \qquad \qquad \qquad \qquad \qquad + \Big(\eps''+\frac1{\DIMENSION} \Big) \nu_{t,v}(f^2)^{1/2} \Bigg) \\
&\le K(D, \GAUSSIANMEASURE_0, \GAUSSIANMEASURE_1, \MARGIN, \SIGNALSHORTDEPENDENCE, \REPLICAS) \Bigg( \nu_{t,v}\big( |f| \big) + \Big(\eps'' + \frac1{\DIMENSION} + \exp(-\DIMENSION) \Big) \nu_{t,v}\big( f^2 \big)^{1/2} \Bigg)\, .
\end{align*}
For $f \ge 0$ we thus have
\begin{align*}
\Big| \frac{\partial}{\partial v} \nu_{t,v}(f) \Big| &\le K(D, \GAUSSIANMEASURE_0, \GAUSSIANMEASURE_1, \MARGIN, \SIGNALSHORTDEPENDENCE, \REPLICAS) \Bigg( \nu_{t,v}\big( f \big) + \Big(\eps'' + \frac1{\DIMENSION}  \Big) \sup_{t,v \in [0,1]^2}  \nu_{t,v}\big( f^2 \big)^{1/2} \Bigg)\, .
\end{align*}
Recalling $\nu_t \equiv \nu_{t,1}$, we conclude Lemma \upref{lem:interpolationkeyconcentration} via Lemma \upref{lem:gronwallvariant}.

\subsection{Proof of Lemma \upref{lem:derivativeintbounded}}\label{subsec:pf_derivativeintbounded}
First we explicitly compute the derivative $\frac{\rmd}{\rmd t}\nu_{t}(f)$ using Gaussian Integration by Parts:
\begin{lemma}\label{lem:derivativewrttexpression}
For a function $f:(\R^{\DIMENSION-1})^\REPLICAS \rightarrow \R$, we have 
\begin{align*}
\frac{\rmd}{\rmd t}\nu_{t}(f) = \expressionI_f + \expressionII_f + \expressionIII_f + \expressionIV_f + \expressionV_f\,,
\end{align*}
where we define the quantities $\expressionI_f, \expressionII_f, \expressionIII_f, \expressionIV_f, \expressionV_f$ as follows:
\begin{align*}
\expressionI_f &:= \frac{\alpha}2 \Big( \sum_{\REPLICAIDX \le \REPLICAS} \nu_t\big( (\theta_{\DIMENSION}^\REPLICAIDX)^2 (u''(S^{\REPLICAIDX}_{\SAMPLES, t}) + u'(S^{\REPLICAIDX}_{\SAMPLES, t})^2) f \big) - \REPLICAS \nu_t\big( (\theta_{\DIMENSION}^\REPLICAS)^2 (u''(S^{\REPLICAS+1}_{\SAMPLES, t}) + u'(S^{\REPLICAS+1}_{\SAMPLES, t})^2) f \big) \Big)\,,\\
\expressionII_f &:= \alpha \Big( \sum_{1 \le \REPLICAIDX_1 < \REPLICAIDX_2 \le \REPLICAS} \nu_t\big( \theta_{\DIMENSION}^{\REPLICAIDX_1} \theta_{\DIMENSION}^{\REPLICAIDX_2} u'(S^{\REPLICAIDX_1}_{\SAMPLES, t}) u'(S^{\REPLICAIDX_2}_{\SAMPLES,t}) f \big) - \REPLICAS \sum_{\REPLICAIDX \le \REPLICAS} \nu_t\big( \theta_{\DIMENSION}^{\REPLICAIDX} \theta_{\DIMENSION}^{\REPLICAS+1} u'(S^{\REPLICAIDX}_{\SAMPLES, t}) u'(S^{\REPLICAS+1}_{\SAMPLES,t}) f \big) \\
&\qquad\qquad + \frac{\REPLICAS(\REPLICAS+1)}2 \nu_t\big( \theta_{\DIMENSION}^{\REPLICAS+1} \theta_{\DIMENSION}^{\REPLICAS+2} u'(S^{\REPLICAS+1}_{\SAMPLES, t}) u'(S^{\REPLICAS+2}_{\SAMPLES,t}) f \big) \Big)\,, \\
\expressionIII_f &:= -r \Big( \sum_{1 \le \REPLICAIDX_1 < \REPLICAIDX_2 \le n} \nu_t\big( \theta_{\DIMENSION}^{\REPLICAIDX_1} \theta_{\DIMENSION}^{\REPLICAIDX_2} f \big) - \REPLICAS \sum_{\REPLICAIDX \le \REPLICAS} \nu_t\big( \theta_{\DIMENSION}^{\REPLICAIDX} \theta_{\DIMENSION}^{\REPLICAS+1} f \big)  + \frac{\REPLICAS(\REPLICAS+1)}2 \nu_t\big( \theta_{\DIMENSION}^{\REPLICAS+1} \theta_{\DIMENSION}^{\REPLICAS+2} f \big) \Big)\,,\\
\expressionIV_f &:= -\frac{r}2 \Big( \sum_{\REPLICAIDX \le \REPLICAS} \nu_t\big( (\theta_{\DIMENSION}^{\REPLICAIDX})^2 f \big) - \REPLICAS \nu_t\big( (\theta_{\DIMENSION}^{\REPLICAS+1})^2 f \big) \Big)\,,\\
\expressionV_f &:= \frac{r-\bar{r}}2 \Big( \sum_{\REPLICAIDX \le \REPLICAS} \nu_t\big( (\theta_{\DIMENSION}^{\REPLICAIDX})^2 f \big) - \REPLICAS \nu_t\big( (\theta_{\DIMENSION}^{\REPLICAS+1})^2 f \big) \Big)\,.
\end{align*}
\end{lemma}
\begin{proof}[Proof of Lemma \upref{lem:derivativewrttexpression}]
Note 
\begin{align*}
\frac{\rmd}{\rmd t}\nu_{t}(f) &= \sum_{\REPLICAIDX \le \REPLICAS} \nu_{t}\Big( f \frac{\partial}{\partial t} H^\REPLICAIDX_{t,\SAMPLES} \Big) - \REPLICAS \nu_{t}\Big( f \frac{\partial}{\partial s}H^{\REPLICAS+1}_{t, \SAMPLES} \Big)\,, \\
\frac{\partial}{\partial t} H^{\REPLICAIDX}_{t,\SAMPLES} &= \sum_{\SAMPLEIDX \le \SAMPLES} u'(S_{\SAMPLEIDX,t}^\REPLICAIDX) \cdot \frac1{2\sqrt{t\DIMENSION}} g_{\SAMPLEIDX,\DIMENSION}\theta^\REPLICAIDX_{\DIMENSION} - \sqrt{\frac{r}{1-t}} \theta^\REPLICAIDX_{\DIMENSION} Y + \frac{r-\bar{r}}2 \big(\theta^\REPLICAIDX_{\DIMENSION}\big)^2\,.
\end{align*}
The result now follows by symmetry w.r.t. $\SAMPLEIDX, 1 \le \SAMPLEIDX \le \SAMPLES$, and then applying Gaussian Integration by Parts w.r.t. $g_{\SAMPLES,\DIMENSION}, z''$. (The argument is analogous to the proof of Proposition 3.2.1, \cite{talagrand2010mean}.)
\end{proof}

Next collecting like terms from Lemma \upref{lem:derivativewrttexpression}, we obtain
\begin{equation}\label{eq:approxrealRSequationsderivexpression}
\begin{aligned}
\Big| \frac{\rmd}{\rmd t}\nu_{t}(f) \Big| &\le \frac12 \sum_{\REPLICAIDX \le \REPLICAS} \Big( \alpha \nu_t\big( (\theta_{\DIMENSION}^{\REPLICAIDX})^2 (u''(S^{\REPLICAIDX}_{\SAMPLES, t}) + u'(S^{\REPLICAIDX}_{\SAMPLES, t})^2) f \big) - \bar{r} \nu_t\big( (\theta_{\DIMENSION}^{\REPLICAIDX})^2 f \big) \Big) \\
&\qquad + \sum_{1 \le \REPLICAIDX_1 < \REPLICAIDX_2 \le \REPLICAS} \Big( \alpha \nu_t\big( \theta_{\DIMENSION}^{\REPLICAIDX_1} \theta_{\DIMENSION}^{\REPLICAIDX_2} u'(S^{\REPLICAIDX_1}_{\SAMPLES, t}) u'(S^{\REPLICAIDX_2}_{\SAMPLES,t}) f \big) - r \nu_t\big( \theta_{\DIMENSION}^{\REPLICAIDX_1} \theta_{\DIMENSION}^{\REPLICAIDX_2} f \big) \Big) \\
&\qquad - \REPLICAS \sum_{\REPLICAIDX \le \REPLICAS} \Big( \alpha \nu_t\big( \theta_{\DIMENSION}^{\REPLICAIDX} \theta_{\DIMENSION}^{\REPLICAS+1} u'(S^{\REPLICAIDX}_{\SAMPLES, t}) u'(S^{\REPLICAS+1}_{\SAMPLES, t}) f \big)  - r \nu_t\big( \theta_{\DIMENSION}^{\REPLICAIDX} \theta_{\DIMENSION}^{\REPLICAS+1} f \big) \Big) \\
&\qquad + \frac{\REPLICAS(\REPLICAS+1)}2 \Big( \alpha \nu_t\big( \theta_{\DIMENSION}^{\REPLICAS+1} \theta_{\DIMENSION}^{\REPLICAS+2} u'(S^{\REPLICAS+1}_{\SAMPLES, t}) u'(S^{\REPLICAS+2}_{\SAMPLES, t}) f \big) - r \nu_t\big( \theta_{\DIMENSION}^{\REPLICAS+1} \theta_{\DIMENSION}^{\REPLICAS+2} f \big) \Big)\,.
\end{aligned}
\end{equation}
We now upper bound the right hand side in the above.
This is where we use the particular choice of $r,\bar{r}$ from \equpref{eq:interpolationrrbardef1}, \equpref{eq:interpolationrrbardef2}.
Specifically, we will employ cavity in $\SAMPLES$ over the parameter $v$ to prove:
\begin{lemma}\label{lem:cavityinMrbarr}
Let $f = \big(\theta_{\DIMENSION}^1\big)^2$ or $\theta_{\DIMENSION}^1 \theta_{\DIMENSION}^2$. For any $t \in [0,1]$ and any $\REPLICAIDX \le \REPLICAS+2$, $1 \le \REPLICAIDX_1 \neq \REPLICAIDX_2 \le \REPLICAS+2$, 
\begin{align*}
&\Big| \alpha \nu_t\big( (\theta_{\DIMENSION}^{\REPLICAIDX})^2 (u''(S^{\REPLICAIDX}_{\SAMPLES, t}) + u'(S^{\REPLICAIDX}_{\SAMPLES, t})^2) f \big) - \bar{r} \nu_t\big( (\theta_{\DIMENSION}^{\REPLICAIDX})^2 f \big) \Big| \\
&\qquad \qquad \qquad \vee \Big| \alpha \nu_t\big( \theta_{\DIMENSION}^{\REPLICAIDX_1} \theta_{\DIMENSION}^{\REPLICAIDX_2} u'(S^{\REPLICAIDX_1}_{\SAMPLES, t}) u'(S^{\REPLICAIDX_2}_{\SAMPLES, t}) f \big) - r\nu_t\big( \theta_{\DIMENSION}^{\REPLICAIDX_1} \theta_{\DIMENSION}^{\REPLICAIDX_2} f \big) \Big|  \\
&\quad \le  K \Bigg(  \nu_t\big( (R_{1,2} - q_{\SAMPLES, \DIMENSION})^2 \big)^{1/2} + \nu_t\big( (R_{1,1} - \rho_{\SAMPLES, \DIMENSION})^2 \big)^{1/2} + \Big(\eps''+\frac1{\DIMENSION}\Big)\Bigg)\,,
\end{align*}
where $K$ depends on $D, \GAUSSIANMEASURE_0, \GAUSSIANMEASURE_1, \SIGNALSHORTDEPENDENCE, \REPLICAS$, and also on $\MARGIN$ for the interpolators.
\end{lemma}
\begin{proof}[Proof of Lemma \upref{lem:cavityinMrbarr}]
Recall $S_{\SAMPLES, t}^\REPLICAIDX = S_{\SAMPLES, t, 1}^\REPLICAIDX$. First, we show 
\begin{align}
\alpha \nu_{t,0}\big( (\theta_{\DIMENSION}^{\REPLICAIDX})^2 (u''(S^{\REPLICAIDX}_{\SAMPLES, t, 0}) + u'(S^{\REPLICAIDX}_{\SAMPLES, t, 0})^2) f \big) &= \bar{r} \nu_{t,0}\big( (\theta_{\DIMENSION}^{\REPLICAIDX})^2 f \big) \, \label{eq:cavityinMrbarrstep0pt1} \\
\alpha \nu_{t,0}\big( \theta_{\DIMENSION}^{\REPLICAIDX_1} \theta_{\DIMENSION}^{\REPLICAIDX_2} u'(S^{\REPLICAIDX_1}_{\SAMPLES, t, 0}) u'(S^{\REPLICAIDX_2}_{\SAMPLES, t, 0}) f \big) &= r\nu_{t,0}\big( \theta_{\DIMENSION}^{\REPLICAIDX_1} \theta_{\DIMENSION}^{\REPLICAIDX_2} f \big) \,. \label{eq:cavityinMrbarrstep0pt2}
\end{align}
We prove \equpref{eq:cavityinMrbarrstep0pt1}, the proof of \equpref{eq:cavityinMrbarrstep0pt2} being analogous. Note $S_{\SAMPLES, t, 0}^\REPLICAIDX = x s_n + \sqrt{q_{n,d}} \QUENCHEDDISORDER + \sqrt{\rho_{n,d} - q_{n,d}} \ANNEALEDDISORDER^{\REPLICAIDX}$ does not depend on $\btheta$, $Y$, $g_{i,j}$, $s_i$ for $i \le n-1$ and that $H_{t, \SAMPLES-1, \DIMENSION}$ does not depend on $\ANNEALEDDISORDER$, $z$, or $s_n$. Therefore 
\begin{align*}
&\Big\langle (\theta_{\DIMENSION}^{\REPLICAIDX})^2 (u''(S^{\REPLICAIDX}_{\SAMPLES, t, 0}) + u'(S^{\REPLICAIDX}_{\SAMPLES, t, 0})^2)f\,  \exp\Big( \sum_{\REPLICAIDX \le \REPLICAS} u(S^{\REPLICAIDX}_{\SAMPLES, t, 0} ) \Big) \Big\rangle_{t,\sim} \\
&\quad = \E_{\ANNEALEDDISORDER}\Big[ \big(u''(S^{\REPLICAIDX}_{\SAMPLES, t, 0}) + u'(S^{\REPLICAIDX}_{\SAMPLES, t, 0})^2 \big)\, \exp\Big( \sum_{\REPLICAIDX \le \REPLICAS} u(S^{\REPLICAIDX'}_{\SAMPLES, t, 0} ) \Big) \Big]\, \big\langle (\theta_d^{\REPLICAIDX})^2 f \big\rangle_{t, \sim}\, .
\end{align*}
By the same rationale, $\nu_{t,0}\big( (\theta_d^{\REPLICAIDX})^2 f) = \E \big\langle (\theta_d^{\REPLICAIDX})^2 f \big\rangle_{t, \sim}$. Consequently, we have
\begin{align*}
&\nu_{t,0}\Big( (\theta_{\DIMENSION}^{\REPLICAIDX})^2 \big(u''(S^{\REPLICAIDX}_{\SAMPLES, t, 0}) + u'(S^{\REPLICAIDX}_{\SAMPLES, t, 0})^2\big) f \Big) \\
&\quad = \E\Bigg[ \frac{ \E_{\ANNEALEDDISORDER}[( u''(S_{n,t,0}^{\REPLICAIDX}) + u'(S_{n,t,0}^{\REPLICAIDX})^2 ) \exp u(S_{n,t,0}^{\REPLICAIDX})] }{\E_{\ANNEALEDDISORDER}[\exp u(S_{n,t,0}^{L+1})]} \cdot \big\langle (\theta_d^{\REPLICAIDX})^2 f \big\rangle_{t, \sim} \Bigg] \\
&\quad = \E_{z, s_n}\Bigg[ \frac{ \E_{\ANNEALEDDISORDER}[( u''(S_{n,t,0}^{\REPLICAIDX}) + u'(S_{n,t,0}^{\REPLICAIDX})^2 ) \exp u(S_{n,t,0}^{\REPLICAIDX})] }{\E_{\ANNEALEDDISORDER}[\exp u(S_{n,t,0}^{\REPLICAIDX})]} \Bigg] \cdot \E_{Y, g_{i, j}, s_i:i \le n-1}\Big[ \big\langle (\theta_d^{\REPLICAIDX})^2 f \big\rangle_{t, \sim} \Big] \\
&\quad = \frac{\bar{r}}{\alpha} \nu_{t, 0}\Big((\theta_d^{\REPLICAIDX})^2 f \Big)\, .
\end{align*}
This proves \equpref{eq:cavityinMrbarrstep0pt1}. Again, \equpref{eq:cavityinMrbarrstep0pt2} follows with an analogous argument.

In light of \equpref{eq:cavityinMrbarrstep0pt1}, \equpref{eq:cavityinMrbarrstep0pt2}, recalling $\nu_t \equiv \nu_{t,1}$, we next decompose
\begin{equation}\label{eq:cavityinMrbarrstep1pt1}
\begin{aligned}
&\Big| \alpha \nu_{t}\big( (\theta_{\DIMENSION}^{\REPLICAIDX})^2 (u''(S^{\REPLICAIDX}_{\SAMPLES, t}) + u'(S^{\REPLICAIDX}_{\SAMPLES, t})^2) f \big) - \bar{r} \nu_{t, 1}\big( (\theta_{\DIMENSION}^{\REPLICAIDX})^2 f \big) \Big| \\
&\quad\le \Big| \alpha \nu_{t, 1}\big( (\theta_{\DIMENSION}^{\REPLICAIDX})^2 (u''(S^{\REPLICAIDX}_{\SAMPLES, t, 1}) + u'(S^{\REPLICAIDX}_{\SAMPLES, t, 1})^2) f \big) - \alpha \nu_{t,0}\big( (\theta_{\DIMENSION}^{\REPLICAIDX})^2 (u''(S^{\REPLICAIDX}_{\SAMPLES, t, 0}) + u'(S^{\REPLICAIDX}_{\SAMPLES, t, 0})^2) f \big) \Big| \\
&\qquad \qquad \qquad + \Big|  \bar{r} \nu_{t,0}\big( (\theta_{\DIMENSION}^{\REPLICAIDX})^2 f \big) - \bar{r} \nu_{t,1}\big( (\theta_{\DIMENSION}^{\REPLICAIDX})^2 f \big) \Big|\,, 
\end{aligned}
\end{equation}
and similarly we decompose
\begin{equation}\label{eq:cavityinMrbarrstep1pt2}
\begin{aligned}
&\Big| \alpha \nu_{t}\big( \theta_{\DIMENSION}^{\REPLICAIDX_1} \theta_{\DIMENSION}^{\REPLICAIDX_2} u'(S^{\REPLICAIDX_1}_{\SAMPLES, t}) u'(S^{\REPLICAIDX_2}_{\SAMPLES, t}) f \big) - r\nu_{t}\big( \theta_{\DIMENSION}^{\REPLICAIDX_1} \theta_{\DIMENSION}^{\REPLICAIDX_2} f \big) \Big|  \\
&\quad\le \Big| \alpha \nu_{t,1}\big( \theta_{\DIMENSION}^{\REPLICAIDX_1} \theta_{\DIMENSION}^{\REPLICAIDX_2} u'(S^{\REPLICAIDX_1}_{\SAMPLES, t}) u'(S^{\REPLICAIDX_2}_{\SAMPLES, t}) f \big) - \alpha \nu_{t,0}\big( \theta_{\DIMENSION}^{\REPLICAIDX_1} \theta_{\DIMENSION}^{\REPLICAIDX_2} u'(S^{\REPLICAIDX_1}_{\SAMPLES, t}) u'(S^{\REPLICAIDX_2}_{\SAMPLES, t}) f \big) \Big|  \\
&\qquad \qquad \qquad + \Big| r\nu_{t,0}\big( \theta_{\DIMENSION}^{\REPLICAIDX_1} \theta_{\DIMENSION}^{\REPLICAIDX_2} f \big) - r\nu_{t,1}\big( \theta_{\DIMENSION}^{\REPLICAIDX_1} \theta_{\DIMENSION}^{\REPLICAIDX_2} f \big) \Big|\,. 
\end{aligned}
\end{equation}
Note $\alpha \le \DENSITYBOUND$ and $\big|r\big|, \big|\bar{r}\big| \le \DENSITYBOUND \cdot K(D)$ as $|u'|, |u''| \le D$. Now consider any $v \in [0,1]$. We apply Lemma \upref{lem:cavityinMintegrationbyparts} with the function $(\theta_{\DIMENSION}^{\REPLICAIDX})^2 f$, recalling here that $f = \big(\theta_{\DIMENSION}^1\big)^2$ or $\theta_{\DIMENSION}^1 \theta_{\DIMENSION}^2$. 
Applying Cauchy-Schwarz on each of the terms arising from Lemma \upref{lem:cavityinMintegrationbyparts}, and using Proposition \upref{prop:logconcaveconcentration} and Lemma \upref{lem:mgftomomentslemma}, we obtain that $\nu_{t,v}(f^2)^{1/2} \le K(D, \beta_0, \beta_1, \MARGIN, \SIGNALSHORTDEPENDENCE)$. Therefore we get
\begin{align*}
&\Big| \frac{\partial}{\partial v} \nu_{t, v} \big( (\theta_{\DIMENSION}^{\REPLICAIDX})^2 (u''(S^{\REPLICAIDX}_{\SAMPLES, t, 1}) + u'(S^{\REPLICAIDX}_{\SAMPLES, t, 1})^2) f \big) \Big|\, ,\, \Big| \frac{\partial}{\partial v} \nu_{t,v} \big( (\theta_{\DIMENSION}^{\REPLICAIDX})^2 f \big) \Big| \\
&\quad \le K(D, \beta_0, \beta_1, \MARGIN, \SIGNALSHORTDEPENDENCE, \REPLICAS) \Bigg(  \nu_{t,v}\big( (R_{1,2} - q_{\SAMPLES, \DIMENSION})^2 \big)^{1/2} + \nu_{t,v}\big( (R_{1,1} - \rho_{\SAMPLES, \DIMENSION})^2 \big)^{1/2} + \Big(\eps''+\frac1{\DIMENSION}\Big)\Bigg)\, ,
\end{align*}
where we use symmetry between replicas. 
Next applying Lemma \upref{lem:interpolationkeyconcentration} with $f=R_{1,1}-\rho_{n,d}$ and $f=R_{1,2}-q_{n,d}$, and then applying the second part of Proposition \upref{prop:logconcaveconcentration} with Lemma \upref{lem:mgftomomentslemma}, we obtain
\begin{align*}
\nu_{t,v}\big( (R_{1,2} - q_{\SAMPLES, \DIMENSION})^2 \big)^{1/2} \le K(D, \beta_0, \beta_1, \MARGIN, \SIGNALSHORTDEPENDENCE)\Bigg( \nu_t\big( (R_{1,2} - q_{\SAMPLES, \DIMENSION})^2 \big)^{1/2} + \Big(\eps''+\frac1d\Big)\Bigg)\,,
\end{align*}
and analogously for $\nu_{t,v}\big( (R_{1,1} - \rho_{\SAMPLES, \DIMENSION})^2 \big)^{1/2}$. Combining the above two displays, we obtain
\begin{align*}
&\Big| \frac{\partial}{\partial v} \nu_{t, v} \big( (\theta_{\DIMENSION}^{\REPLICAIDX})^2 (u''(S^{\REPLICAIDX}_{\SAMPLES, t, 1}) + u'(S^{\REPLICAIDX}_{\SAMPLES, t, 1})^2) f \big) \Big|\, ,\, \Big| \frac{\partial}{\partial v} \nu_{t,v} \big( (\theta_{\DIMENSION}^{\REPLICAIDX})^2 f \big) \Big| \\
&\quad \le K(D, \beta_0, \beta_1, \MARGIN, \SIGNALSHORTDEPENDENCE, \REPLICAS)\Bigg(  \nu_t\big( (R_{1,2} - q_{\SAMPLES, \DIMENSION})^2 \big)^{1/2} + \nu_t\big( (R_{1,1} - \rho_{\SAMPLES, \DIMENSION})^2 \big)^{1/2} + \Big(\eps''+\frac1{\DIMENSION}\Big)\Bigg)\, .
\end{align*}
Similarly, we have
\begin{align*}
&\Big| \frac{\partial}{\partial v} \nu_{t,v} \big( \theta_{\DIMENSION}^{\REPLICAIDX_1} \theta_{\DIMENSION}^{\REPLICAIDX_2} u'(S^{\REPLICAIDX_1}_{\SAMPLES, t}) u'(S^{\REPLICAIDX_2}_{\SAMPLES, t}) f \big) \Big|\, ,\, \Big| \frac{\partial}{\partial v} \nu_{t,v} \big( \theta_{\DIMENSION}^{\REPLICAIDX_1} \theta_{\DIMENSION}^{\REPLICAIDX_2}  f \big) \Big|\\
&\quad \le K(D, \beta_0, \beta_1, \MARGIN, \SIGNALSHORTDEPENDENCE, \REPLICAS) \Bigg(  \nu_t\big( (R_{1,2} - q_{\SAMPLES, \DIMENSION})^2 \big)^{1/2} + \nu_t\big( (R_{1,1} - \rho_{\SAMPLES, \DIMENSION})^2 \big)^{1/2} + \Big(\eps''+\frac1{\DIMENSION}\Big)\Bigg)\, .
\end{align*}
Combining the above with \equpref{eq:cavityinMrbarrstep1pt1}, \equpref{eq:cavityinMrbarrstep1pt2}, and noting $\alpha \le \DENSITYBOUND \le 2$, the Lemma follows.
\end{proof}

Combining Lemma \upref{lem:cavityinMrbarr} and \equpref{eq:approxrealRSequationsderivexpression}, it follows that for $f = \big(\theta_{\DIMENSION}^1\big)^2$ or $\theta_{\DIMENSION}^1 \theta_{\DIMENSION}^2$, we have 
\begin{align*}
\Big| \frac{\rmd}{\rmd t}\nu_{t}(f) \Big| \le K(D, \GAUSSIANMEASURE_0, \GAUSSIANMEASURE_1, \MARGIN, \SIGNALSHORTDEPENDENCE, \REPLICAS) \Bigg(  \nu_t\big( (R_{1,2} - q_{\SAMPLES, \DIMENSION})^2 \big)^{1/2} + \nu_t\big( (R_{1,1} - \rho_{\SAMPLES, \DIMENSION})^2 \big)^{1/2} + \Big(\eps''+\frac1{\DIMENSION}\Big)\Bigg)\,.
\end{align*}
By Proposition \upref{prop:logconcaveconcentration}, we have 
\begin{align*}
\nu_t\big( (R_{1,1} - \rho_{\SAMPLES, \DIMENSION})^2 \big)^{1/2} &\le K \nu_t\big( (R_{1,1} - \nu_t(R_{1,1}) )^2 \big)^{1/2} + 2 \big| \nu_t(R_{1,1}) - \nu_1(R_{1,1}) \big| \\
&\le \frac{K(D, \GAUSSIANMEASURE_0, \GAUSSIANMEASURE_1, \MARGIN, \SIGNALMEDIUMDEPENDENCE)}{\sqrt{\DIMENSION}} + 2 \big| \nu_t(R_{1,1}) - \nu_1(R_{1,1}) \big|\,,
\end{align*}
and similarly for $R_{1,2}$. Here we use that $\rho_{\SAMPLES, \DIMENSION} = \nu_1(R_{1,1})$, $q_{\SAMPLES, \DIMENSION} = \nu_1(R_{1,2})$. Thus combining the above two displays, to complete the proof of Lemma \upref{lem:derivativeintbounded}, it suffices to upper bound
\begin{align*}
\Big| \frac{\rmd}{\rmd t} \nu_t(R_{1,1}) \Big|\, ,\, \Big| \frac{\rmd}{\rmd t} \nu_t(R_{1,2}) \Big| \le \frac{K(D, \GAUSSIANMEASURE_0, \GAUSSIANMEASURE_1, \MARGIN, \SIGNALMEDIUMDEPENDENCE)}{\sqrt{d}} \,.
\end{align*}
To this end, considering $f = R_{1,2}$ and following the notation of Lemma \upref{lem:derivativewrttexpression}, we note that
\begin{align*}
\expressionI_f = \expressionI_f - \nu_t(R_{1,2}) = \expressionI_{f-\nu_t(R_{1,2})}\,,
\end{align*}
because $\nu_t(R_{1,2})$ is just a constant (note that for any constant function $c$, $\expressionI_c=0$). Similarly $\expressionII_f = \expressionII_{f-\nu_t(R_{1,2})}$, $\expressionIII_f = \expressionIII_{f-\nu_t(R_{1,2})}$, $\expressionIV_f = \expressionIV_{f-\nu_t(R_{1,2})}$, $\expressionV_f = \expressionV_{f-\nu_t(R_{1,2})}$. However note that $\alpha \le \DENSITYBOUND$ and $r,\bar{r} \le K(D, \DENSITYBOUND)$. Thus by Proposition \upref{prop:logconcaveconcentration}, since $|u'|, |u''| \le D$, $\alpha \le \DENSITYBOUND \le 2$, and by H\"{o}lder's Inequality used on each individual summand comprising the following expressions, we obtain
\begin{align*}
&\expressionI_{f-\nu_t(R_{1,2})}\,,\,\expressionII_{f-\nu_t(R_{1,2})}\,,\,\expressionIII_{f-\nu_t(R_{1,2})}\,,\,\expressionIV_{f-\nu_t(R_{1,2})}\,,\,\expressionV_{f-\nu_t(R_{1,2})} \le \frac{K(D, \GAUSSIANMEASURE_0, \GAUSSIANMEASURE_1, \MARGIN, \SIGNALMEDIUMDEPENDENCE)}{\sqrt{d}}\,.
\end{align*}
The same derivation applies for $f = R_{1,1}$, thus completing the proof of Lemma \upref{lem:derivativeintbounded}.

\subsection{Proof of Proposition \ref{prop:logconcaveconcentration}}\label{subsec:logconcaveconcentration}
Here, we prove Proposition \upref{prop:logconcaveconcentration}. 
For only Appendix \upref{subsec:logconcaveconcentration}, we define the following Gibbs average $\langle \cdot \rangle$:
\begin{align*}
\langle f \rangle := \frac{\big\langle f \exp \big( \sum_{1 \le \REPLICAIDX \le \REPLICAS} u(S_{\SAMPLES, t, v}^\REPLICAIDX) \big) \big\rangle_{t, \sim} }{\big\langle \exp u(S_{\SAMPLES, t, v}^{L+1})\big\rangle_{t, \sim}^\REPLICAS }\,.
\end{align*}
Consequently by \equpref{eq:nu_t_v_def}, we have $\nu_{t,v}(f) = \E \langle f \rangle$.
Note we can rewrite $\langle f \rangle$ as
\begin{align}
\frac1{Z'^L} \int \E_{\ANNEALEDDISORDER}\Bigg[ f \prod_{1 \le l \le L} \one\{\theta_1^l/\sqrt{d} \in I\} \exp\Big( u\big(S^l_{n,t,v}(\btheta^l)\big) + H_{t,n-1,d}(\btheta^l) - \beta \| \btheta^l \|^2 \Big) \Bigg] \, \rmd\btheta^1 \cdots \rmd\btheta^L\,,\label{eq:gibbsmeasurerewriteconcentration}
\end{align}
where 
\begin{align}
Z' &= \int_{\R^d} \E_{\ANNEALEDDISORDER}\Big[ \one\{\theta_1/\sqrt{d} \in I\}  \exp\Big( u(S_{n,t,v}) + H_{t,n-1,d}(\btheta) - \beta \| \btheta \|^2 \Big) \Big]\, \rmd\btheta\,.\label{eq:Z'def}
\end{align}
Note the constraint $\prod_{1 \le l \le L} \one\{\theta_1^l/\sqrt{d} \in I\}$ is convex, $u$ is concave, the $S_{n,t,v}$ and $S_{i,t}$ are affine in $\btheta$, and the rest of $H_{t,n-1,d}$ is log-concave or affine in $\btheta$. It follows that the measure $\langle \cdot \rangle$ is $\beta$-strongly log-concave in $\btheta$, a fact that will be exploited repeatedly via concentration of Lipschitz functions of strongly log-concave measures (Theorem \upref{thm:logconcaveconcentrationbasics}) in what follows.

\subsubsection{First part of Proposition \upref{prop:logconcaveconcentration}}
We prove the first part of Proposition \upref{prop:logconcaveconcentration} by directly combining the following Propositions \upref{prop:concentrationoverlapfixedgibbs}, \upref{prop:concentrationmeangibbswrtdisorder}. First, we have the following concentration property of $f = R_{1, 1}$ or $R_{1, 2}$ about their average $\langle f \rangle$ with the quenched disorder fixed:
\begin{proposition}\label{prop:concentrationoverlapfixedgibbs}
Consider $u(\cdot)$ as in Theorem \upref{thm:mollifiedhamiltonianfreeenergy} for the interpolators or Theorem \upref{thm:posterior_mollifiedhamiltonianfreeenergy} for the posterior.
Then for $d \ge d_{\upref{prop:concentrationoverlapfixedgibbs}}(\eps'')$ and all $k \le d$, we have
\begin{align*}
\Big\langle \big(R_{1, 1} - \langle R_{1, 1}\rangle\big)^{2k} \Big\rangle\, ,\, \Big\langle  \big(R_{1, 2} - \langle R_{1, 2}\rangle \big)^{2k} \Big\rangle \le \Big( \frac{ k K B^{\star}}{d^2} \Big)^k\, ,
\end{align*}
where $K$ depends on $D, \beta_0, \beta_1, \SIGNALSHORTDEPENDENCE$ and also on $\MARGIN$ for the interpolators, and where
\begin{equation}\label{eq:Bstardef}
\begin{aligned}
B^{\star} &:= \Big( \sum_{i=1}^n \| \bar{g}_i \|^2 \Big)^{1/2} + \frac1{16d} \sum_{j=2}^d \Big\{ \Big(\sum_{i=1}^{n-1} g_{ij}\Big)^2 + g_{nj}^2 \Big\} \\
&\qquad \qquad + \Big(n \sum_{i=1}^n s_i^2 \Big)^{1/2}  + \frac1{8d\signalsubexp} \Big\{ \Big( \sum_{i=1}^{n-1} s_i \Big)^2 + s_n^2 \Big\} + Y^2 + \QUENCHEDDISORDER^2 + d\,.
\end{aligned}
\end{equation}
\end{proposition}
Second, we have concentration of $\langle f \rangle$ about $\E \langle f \rangle$:
\begin{proposition}\label{prop:concentrationmeangibbswrtdisorder} 
Consider $u(\cdot)$ as in Theorem \upref{thm:mollifiedhamiltonianfreeenergy} for the interpolators or Theorem \upref{thm:posterior_mollifiedhamiltonianfreeenergy} for the posterior. Then for $k \le \frac{d}4$, we have
\begin{align*}
\E\Big[ \big( \big\langle R_{1,1} \big\rangle - \E \big\langle R_{1,1} \big\rangle \big)^{2k} \Big]~,~ \E\Big[ \big( \big\langle R_{1,2} \big\rangle - \E \big\langle R_{1,2} \big\rangle \big)^{2k} \Big] &\le \Big( \frac{K k}d \Big)^k\,,
\end{align*}
where $K$ depends on $D, \beta_0, \beta_1, \SIGNALMEDIUMDEPENDENCE$, and also on $\MARGIN$ for the posterior.
\end{proposition}
Combining Proposition \upref{prop:concentrationoverlapfixedgibbs} and Proposition \upref{prop:concentrationmeangibbswrtdisorder} yields the first part of Proposition \upref{prop:logconcaveconcentration}. For the rest of the proof of the first part of Proposition \upref{prop:logconcaveconcentration}, we assume $u(\cdot)$ satisfies the conditions of Theorem \upref{thm:mollifiedhamiltonianfreeenergy}.
Central in the proofs of both Proposition \upref{prop:concentrationoverlapfixedgibbs} and Proposition \upref{prop:concentrationmeangibbswrtdisorder} is concentration of Lipschitz functions w.r.t strongly log-concave measures, as written in Theorem \upref{thm:logconcaveconcentrationbasics}.
We will also use the following fact which follows by a standard symmetrization argument, see (3.22) of \cite{talagrand2010mean}.
\begin{lemma}\label{lem:symmetrization}
For any $m \in \R$, any function $f$ such that the following is defined, and any probability measure $\mu$, we have 
\begin{align*}
\int \Big(f - \int f \rmd \mu \Big)^{2k}\, \rmd \mu \le 2^{2k} \int (f-m)^{2k}\, \rmd \mu\,.
\end{align*}
\end{lemma}

\paragraph{Proof of Proposition \upref{prop:concentrationoverlapfixedgibbs}:}
We now prove Proposition \upref{prop:concentrationoverlapfixedgibbs}.
First, we show that $B^{\star}$ can be controlled on the exponential scale:
\begin{lemma}\label{lem:lem310}
For all $k \le d$, we have
\begin{align*}
\E\big[ B^{\star k} \big] \le \big( K d \big)^k\,,
\end{align*}
where $K$ depends on $\SIGNALMEDIUMDEPENDENCE$.
\end{lemma}
\begin{proof}
Our strategy will be to upper bound $\E\big[ \exp B^\star \big]$ and then convert this to a bound on the $k$-th moment of $B^\star$ by Lemma \upref{lem:mgftomomentslemma}. Since $n \le 2d$, and by applying Cauchy-Schwarz on the $(s_i)$,
\begin{align*}
B^{\star} &\le \frac1{8n} \Big( \sum_{i=1}^n \| \bar g_i \|^2 \Big) + 2n + \frac1{16d} \sum_{j=2}^d \Big\{ \Big(\sum_{i=1}^{n-1} g_{ij}\Big)^2 + g_{nj}^2 \Big\} + \frac{1}{2\signalsubexp} \Big( \sum_{i=1}^n s_i^2 \Big) + \signalsubexp n \\ 
&\qquad \qquad + \frac{Y^2}4 + \frac{\QUENCHEDDISORDER^2}4 + Kd \\
&= \expressionI + \expressionII + \expressionIII + \frac{Y^2}4 + \frac{Z^2}4 + (\signalsubexp + 2) n + Kd\,,
\end{align*}
where we define
\begin{align*}
\expressionI := \frac1{8n} \Big( \sum_{i=1}^n \| \bar g_i \|^2 \Big)\,,\, \expressionII := \frac1{16d} \sum_{j=2}^d \Big\{ \Big(\sum_{i=1}^{n-1} g_{ij}\Big)^2 + g_{nj}^2 \Big\} \,,\, \expressionIII := \frac{1}{2\signalsubexp} \Big( \sum_{i=1}^n s_i^2 \Big)\,.
\end{align*}
By independence of $g_{ij}, s_i, Y, Z$, and as $n \le 2d$, we have 
\begin{align}
&\E\big[ \exp B^\star \big] \notag \\
&\quad \le \E\big[ \exp \big( \expressionI + \expressionII \big) \big] \E\big[ \exp \expressionIII \big] \E\big[ \exp Y^2/4] \E\big[ \exp Z^2/4] \exp\big( K(\lambda) d \big) \notag \\
&\quad \le \E\big[ \exp 2 \expressionI ] \E\big[ \exp 2 \expressionII ] \E\big[ \exp \expressionIII \big] \E\big[ \exp Y^2/4] \E\big[ \exp Z^2/4] \exp\big( K(\lambda) d \big)\,,\label{eq:boundingmgfBstar_0}
\end{align}
where the last step used Cauchy-Schwarz. Note $\E\big[ \big( \exp \frac{Y^2}{4} \big) \big]$, $\E\big[ \big( \exp \frac{\QUENCHEDDISORDER^2}{4} \big) \big] \le K$ as $Y, Z \sim N(0, 1)$.
It remains to upper bound $\E\big[ \exp 2\expressionI \big]$, $\E\big[ \exp 2\expressionII \big]$, and $\E\big[ \exp \expressionIII \big]$. 
\begin{itemize}
\item Upper bounding $\E\big[ \exp 2\expressionI \big]$: By a direct calculation, we have the following, as stated as (A.11) in \cite{talagrand2010mean}:
\begin{align*}
\E\Big[ \Big(\exp \frac{g_{ij}^2}{4n} \Big)\Big] \le \Big( 1-\frac1{2n} \Big)^{-1/2}\,.
\end{align*}
Thus the independence of the $(g_{ij})$ implies
\begin{align*}
\E\big[ \exp 2\expressionI \big] \le \Big( 1-\frac1{2n} \Big)^{-nd/2} \le K^d \le \exp (Kd)\,.
\end{align*}

\item Upper bounding $\E\big[ \exp 2\expressionII \big]$: By independence of the $(g_{ij})$, we have $\sum_{i=1}^{n-1} g_{ij} \sim N(0, n-1)$ for all $2 \le j \le d$, and thus 
\begin{align*}
\frac1{8d} \Big( \sum_{i=1}^{n-1} g_{ij} \Big)^2 + \frac1{8d} g_{nj}^2 \sim \frac{1}{8d} N(0, n-1)^2 + \frac1{8d} N(0, 1)^2 \sim \frac{n}{8d} N(0, 1)^2 \implies 2 \expressionII \sim \frac{n}{8d} \chi^2_{d-1}\,,
\end{align*}
where $\chi^2_{d-1}$ denotes a $\chi^2$ distribution with $d-1$ degrees of freedom. Since $n \le 2d$, standard bounds on the MGF of a $\chi^2$ random variable now gives
\begin{align*}
\E\big[ \exp 2\expressionII \big] = \Big(1-2 \cdot \frac{n}{8d} \Big)^{-\frac{d-1}2} \le 2^{\frac{d-1}2} \le \exp(Kd)\,.
\end{align*}

\item Upper bounding $\E\big[ \exp \expressionIII \big]$: Since the $s_i^2$ are $\signalsubexp$ sub-Exponential, by e.g. Proposition 2.7.1 of \cite{vershynin2018high}), we have $\E\big[ \exp s_i^2 / 2\signalsubexp \big] \le K$. As the $(s_i)$ are independent, we obtain 
\begin{align*}
\E\big[ \exp \expressionIII \big] \le K^n \le \exp(K d)\,.    
\end{align*}
\end{itemize}
Combining the bounds in the above three cases together with \equpref{eq:boundingmgfBstar_0} and recalling $\signalsubexp$ depends only on $\lambda$ implies
\begin{align}
\E \big[ \exp B^{\star} \big] &\le \exp\big( K(\lambda) d \big)\,,\label{eq:boundmgfBstar}
\end{align}
and the conclusion now follows from Lemma \upref{lem:mgftomomentslemma}.
\end{proof}

We next lower bound the partition function $Z'$ in order to apply Lemma \upref{lem:logconcaveboundmomentslemma} to control the relevant overlaps.
\begin{lemma}\label{lem:lem316}
For $d \ge d_{\upref{lem:lem316}}(\eps'')$, letting $Z'$ be as per \equpref{eq:Z'def}, we have 
\begin{align*}
Z' \ge \exp\big( -K_{\upref{lem:lem316}}  B^{\star} \big)\,,
\end{align*}
where $K_{\upref{lem:lem316}}$ depends on $D, \beta_0, \beta_1, \SIGNALSHORTDEPENDENCE$, and also on $\MARGIN$ for the posterior.
\end{lemma}
\begin{proof}
We will prove this Lemma in two stages. First we will prove it when $v=t=1$, which lets us upper bound $|q_{n,d}|$. Then we use this upper bound on $|q_{n,d}|$ to complete the proof for general $t,v\in [0,1]$. Note the conditions on $u$ from Theorem \upref{thm:mollifiedhamiltonianfreeenergy} give $u(x) \ge -D(\kappa)(|x|+1)$ for the interpolators; similarly, we have $u(x) \ge -D(\lambda)(|x|+1)$ for the posterior. In an abuse of notation, we will write this bound as $u(x) \ge -D(\kappa, \lambda)(|x|+1)$, where we emphasize that there is no dependence on $\kappa$ for the posterior.

We begin with steps of the proof that are common to both stages.
Rewrite
\begin{align*}
Z' &= \int_{\frac1{\sqrt{d}} \theta_1 \in I} \int_{\R^{d-1}} \E_{\ANNEALEDDISORDER}\Big[ \exp\Big( u(S_{n,t,v}) + \sum_{1 \le i \le n-1} u(S_{i,t}) \\
&\qquad \qquad \qquad \qquad \qquad \qquad + \theta_d\sqrt{(1-t)r} Y - \frac{(1-t)(r-\bar{r})}2 \theta_d^2 - \beta \| \bar{\btheta} \|^2 - \beta \theta_1^2 \Big)\Big] \, \rmd \overline{\btheta} \rmd \theta_1\,.
\end{align*}
Notice $\bar{r} \le r \le K(D)$ as per their definitions in \equpref{eq:psidef}, \equpref{eq:psibardef}, thus $\frac{(1-t)(r-\bar{r})}2 \le \frac{r-\bar{r}}2 \le K_{\upref{lem:lem316}}(D)$. Letting $\beta' := K_{\upref{lem:lem316}}(D)+\beta$, we consider the density $\gamma$ on $\R^{d-1}$ with density $(\GAUSSIANMEASURE'/\pi)^{\frac{d-1}2} \exp(-\GAUSSIANMEASURE' \|\overline{\btheta}\|^2)$. 
We observe that
\begin{align*}
&\int_{\R^{d-1}} \E_{\ANNEALEDDISORDER}\Big[ \exp\Big( u(S_{n,t,v}) + \sum_{1 \le i \le n-1} u(S_{i,t}) + \theta_d\sqrt{(1-t)r} Y - \frac{(1-t)(r-\bar{r})}2 \theta_d^2 - \beta \| \bar{\btheta} \|^2 - \beta \theta_1^2 \Big)\Big] \, \rmd \overline{\btheta} \\
&\quad\ge \exp(-\GAUSSIANMEASURE \theta_1^2) \Big( \frac{\pi}{\GAUSSIANMEASURE'}\Big)^{\frac{d-1}2 } \int_{\R^{d-1}} \E_{\ANNEALEDDISORDER}\Big[ \exp\Big( u(S_{n,t,v}) + \sum_{1 \le i \le n-1} u(S_{i,t}) + \theta_d\sqrt{(1-t)r}Y \Big) \Big]\, \rmd \gamma\big(\bar{\btheta}) \\
&\quad \ge \exp(-\GAUSSIANMEASURE \theta_1^2) \Big( \frac{\pi}{\GAUSSIANMEASURE'}\Big)^{\frac{d-1}2 } \exp\Bigg(\int_{\R^{d-1}} \E_{\ANNEALEDDISORDER}\Big[ u(S_{n,t,v}) + \sum_{1 \le i \le n-1} u(S_{i,t}) + \theta_d\sqrt{(1-t)r}Y \Big] \, \rmd \gamma\big(\bar{\btheta})\Bigg)\,.
\end{align*}
Next by rotational invariance of $\gamma$, we have the following bound, as stated in (3.24) in \cite{talagrand2010mean}:
\begin{align}
\int \big| \langle x, \overline{\btheta}\rangle \big| \rmd \gamma\big(\bar{\btheta}) = \sqrt{\frac1{\pi \GAUSSIANMEASURE'}} \| x \| \le \frac1{\sqrt{\GAUSSIANMEASURE'}} \| x \| \quad \text{for all} \quad x \in \R^{d-1}\, .  \label{eq:gaussianmeasureorthogonalitybound}  
\end{align} 
Also define $\bar{g_i}' := \Big( g_{i,2}, \ldots, g_{i,d-1}, \sqrt{t} g_{i,d} + \frac{\sqrt{d}}{n} \cdot \sqrt{(1-t)r}Y \Big)$. Using \equpref{eq:gaussianmeasureorthogonalitybound}, it follows that
\begin{align*}
&\int_{\R^{d-1}} \E_{\ANNEALEDDISORDER}\Big[ u(S_{n,t,v}) + \sum_{1 \le i \le n-1} u(S_{i,t}) + \theta_d\sqrt{(1-t)r}Y \Big] \, \rmd \gamma\big(\bar{\btheta}) \\
&\quad \ge -D(\kappa, \lambda) \Bigg( \int_{\R^{d-1}} \E_{\ANNEALEDDISORDER}\Bigg[ \Big| \sqrt{v} \Big(\frac1{\sqrt{\DIMENSION}} \sum_{2 \le \DIMENSIONIDX \le \DIMENSION-1} g_{\SAMPLES, \DIMENSIONIDX} \theta_{\DIMENSIONIDX} + \sqrt{\frac{t}{\DIMENSION}} g_{\SAMPLES, \DIMENSION} \theta_{\DIMENSION} \Big) + \sqrt{1-v} \big( \sqrt{q_{\SAMPLES, \DIMENSION}} \QUENCHEDDISORDER + \sqrt{\rho_{\SAMPLES, \DIMENSION} -q_{\SAMPLES, \DIMENSION}}\ANNEALEDDISORDER  \big) \Big| \\
&\qquad \qquad \qquad \qquad \qquad \qquad \qquad \qquad + \sum_{1 \le i \le n-1} \Big| \frac{s_i \theta_1}{\sqrt{d}} + \frac1{\sqrt{d}} \sum_{2 \le j \le d-1} g_{i,j} \theta_j + \sqrt{\frac{t}{d}} g_{i,d} \theta_d \Big| \\
&\qquad \qquad \qquad \qquad \qquad \qquad \qquad \qquad + \Big| v \cdot \frac{s_n \theta_1}{\sqrt{d}} + (1-v) s_n x \Big| + n + \theta_d\sqrt{(1-t)r}Y \Bigg] \, \rmd \gamma\big(\bar{\btheta}) \Bigg) \\
&\quad\ge -D(\kappa, \lambda) \Bigg(n + \frac1{\sqrt{d}} \sum_{1 \le i \le n-1} \int \big| \langle \bar{g}_i', \bar{\btheta} \big| \rangle \rmd \gamma(\bar{\btheta}) + \sqrt{\frac{v}{d}} \int \big| \langle \bar{g}_n', \bar{\btheta} \rangle \big| \rmd \gamma(\bar{\btheta}) \\
&\qquad \qquad \qquad \qquad \qquad + \frac{|\theta_1|}{\sqrt{d}} \sum_{1 \le i \le n} |s_i| + |\QUENCHEDDISORDER| \sqrt{(1-v)q_{n,d}} + (1-v) |s_n| + K\sqrt{(1-v)(\rho_{n,d}-q_{n,d})} \Bigg) \\
&\quad\ge -D(\kappa, \lambda) \Bigg( n + \frac1{\sqrt{d \beta'}} \sum_{1 \le i \le n} \|\bar{g}_i'\| + \frac{|\theta_1|}{\sqrt{d}} \sum_{1 \le i \le n} |s_i| + |\QUENCHEDDISORDER| \sqrt{(1-v) q_{n,d}} + |s_n| + K\sqrt{\rho_{n,d}-q_{n,d}}  \Bigg)\,.
\end{align*}
Here, in the second inequality we took expectation with respect to $\ANNEALEDDISORDER$ and used the bound $\E_{\ANNEALEDDISORDER}\big[|\ANNEALEDDISORDER|\big] \le K$, and in the third inequality we used \equpref{eq:gaussianmeasureorthogonalitybound}.

Next, remark that $\sum_{1 \le i \le n} \big|s_i\big| \le \big( n \sum_{1 \le i \le n} s_i^2 \big)^{1/2}$. Also by Triangle Inequality and as $|r| \le K(D)$, we have
\begin{align*}
\frac1{\sqrt{d}} \sum_{1 \le i \le n} \| \bar{g_i}' \| &\le \frac1{\sqrt{d}}  \sum_{1 \le i \le n} \Big( \| \bar{g}_i \| + \frac{\sqrt{d(1-t)r}}{n} |Y| \Big) \le \Big( \frac{n}{d} \sum_{1 \le i \le n} \| \bar{g}_i \|^2 \Big)^{1/2} + K(D) |Y|\,.
\end{align*}

Now specializing to when $t=v=1$, combining the above displays, recalling the definition $\beta'
= K_{\upref{lem:lem316}}(D)+\beta$, and using that $\frac{\theta_1}{\sqrt{d}} \in I \subseteq [-1,1]$, $n \le \DENSITYBOUND d \le 2d$ and the definition of $B^{\star}$ from \equpref{eq:Bstardef}, we obtain
\begin{align*}
&\int_{\R^{d-1}} \E_{\ANNEALEDDISORDER}\Big[ \exp\Big( u(S_{n,1,1}) + \sum_{1 \le i \le n-1} u(S_{i,1}) \beta \| \bar{\btheta} \|^2 - \beta \theta_1^2 \Big)\Big] \, \rmd \overline{\btheta} \\
&\quad \ge \exp(-\GAUSSIANMEASURE \theta_1^2) \Big( \frac{\pi}{\GAUSSIANMEASURE'}\Big)^{\frac{d-1}2 } \exp\Bigg( -D(\kappa, \lambda) \Big( n + \frac1{\sqrt{d \beta'}} \sum_{1 \le i \le n} \|\bar{g}_i'\| + \frac{|\theta_1|}{\sqrt{d}} \sum_{1 \le i \le n} |s_i| + |s_n| \Big) \Bigg) \\
&\quad \ge \exp\Bigg( -K(D, \beta_0, \beta_1, \MARGIN) \Big( d + \Big( \sum_{1 \le i \le n} \|\bar{g}_i\|^2 \Big)^{1/2} + \Big( d \sum_{1 \le i \le n} s_i^2 \Big)^{1/2} + |Y| \Big) \Bigg) \\
&\quad \ge \exp\Big( -K(D, \beta_0, \beta_1, \MARGIN) B^{\star} \Big)\,.
\end{align*}
Using that $|\sqrt{d} I| \ge 2\eps''$, that $\eps'' \ge \exp(-d)$ if $d \ge d_{\upref{lem:lem316}}(\eps'')$, and that $B^{\star} \ge d$, we thus obtain when $t=v=1$,
\begin{align*}
Z' \ge \int_{\frac1{\sqrt{d}} \theta_1 \in I } \exp\Big( -K(D, \beta_0, \beta_1, \MARGIN) B^{\star} \Big)\,\rmd\theta_1 \ge \exp\Big( -K(D, \beta_0, \beta_1, \MARGIN) B^{\star} \Big)\,.
\end{align*}
The condition of Lemma \upref{lem:logconcaveboundmomentslemma} thus applies, and combining it with Lemma \upref{lem:lem310} yields that 
\begin{align}
0 \le |q_{n,d}| \le \rho_{n,d} &= \frac1d \E\Big[ \big\langle \| \bar{\btheta} \|^2 \big\rangle \Big] \notag  \\
&\le \frac1d \E\Big[ K(D, \beta_0, \beta_1, \MARGIN, \SIGNALSHORTDEPENDENCE) B^{\star} \Big] \notag \\
&\le \frac1d \cdot K(D, \beta_0, \beta_1, \MARGIN, \SIGNALSHORTDEPENDENCE) d = K(D, \beta_0, \beta_1, \MARGIN, \SIGNALSHORTDEPENDENCE)\,.\label{eq:upperboundrhowithinpf}
\end{align}

We now return to the proof for general $t,v\in[0,1]^2$. Again by the aforementioned upper bounds on $\sum_{1 \le i \le n}|s_i|$ and $\frac1{\sqrt{d}} \sum_{1 \le i \le n} \| \bar{g}_i'\|$, recalling the definition $\beta'
= K_{\upref{lem:lem316}}(D)+\beta$, and using that $\frac{\theta_1}{\sqrt{d}} \in I \subseteq [-1,1]$, $n \le \DENSITYBOUND d \le 2d$ and the definition of $B^{\star}$, we have 
\begin{align*}
&\int_{\R^{d-1}} \E_{\ANNEALEDDISORDER}\Big[ \exp\Big( u(S_{n,t,v}) + \sum_{1 \le i \le n-1} u(S_{i,t}) \beta \| \bar{\btheta} \|^2 - \beta \theta_1^2 \Big)\Big] \, \rmd \overline{\btheta} \\
&\quad \ge \exp(-\GAUSSIANMEASURE \theta_1^2) \Big( \frac{\pi}{\GAUSSIANMEASURE'}\Big)^{\frac{d-1}2 } \exp\Bigg( -D(\kappa, \lambda) \Big( n + \frac1{\sqrt{d \beta'}} \sum_{1 \le i \le n} \|\bar{g}_i'\| + \frac{|\theta_1|}{\sqrt{d}} \sum_{1 \le i \le n} |s_i| \\
&\qquad \qquad \qquad \qquad \qquad \qquad \qquad \qquad \qquad \qquad + |\QUENCHEDDISORDER|\sqrt{(1-v)q_{n,d}} + K\sqrt{\rho_{n,d}-q_{n,d}} + |s_n| \Big) \Bigg) \\
&\quad \ge \exp\Bigg( -K(D, \beta_0, \beta_1, \MARGIN, \SIGNALSHORTDEPENDENCE) \Big( d + \Big( \sum_{1 \le i \le n} \|\bar{g}_i\|^2 \Big)^{1/2} + \Big( d \sum_{1 \le i \le n} s_i^2 \Big)^{1/2} + |Y| + |\QUENCHEDDISORDER| \Big) \Bigg) \\
&\quad \ge \exp\Big( -K(D, \beta_0, \beta_1, \MARGIN, \SIGNALSHORTDEPENDENCE) B^{\star} \Big)\,.
\end{align*}
Here, we used \equpref{eq:upperboundrhowithinpf} to upper bound $|\rho_{n,d}|, |q_{n,d}|$ and hence $K\sqrt{\rho_{n,d}-q_{n,d}}$. 

Using that $|\sqrt{d} I| \ge 2\eps''$, that $\eps'' \ge \exp(-d)$ if $d \ge d_{\upref{lem:lem316}}(\eps'')$, and that $B^{\star} \ge d$, we thus obtain for general $t,v\in[0,1]^2$ that
\begin{align*}
Z' \ge \int_{\frac1{\sqrt{d}} \theta_1 \in I } \exp\Big( -K(D, \beta_0, \beta_1, \MARGIN, \SIGNALSHORTDEPENDENCE) B^{\star} \Big)\,\rmd\theta_1 \ge \exp\Big( -K(D, \beta_0, \beta_1, \MARGIN, \SIGNALSHORTDEPENDENCE) B^{\star} \Big)\,.
\end{align*}
This proves Lemma \upref{lem:lem316}.
\end{proof}

We now have the ingredients needed to prove Proposition \upref{prop:concentrationoverlapfixedgibbs}.
\begin{proof}[Proof of Proposition \upref{prop:concentrationoverlapfixedgibbs}]
First, we claim that 
\begin{align}
\big\langle \| \bar\btheta \|^{2k} \big\rangle \le \big\langle \| \bar\btheta \|^{4k} \big\rangle^{1/2} &\le \big( K(D, \beta_0, \beta_1, \MARGIN, \SIGNALSHORTDEPENDENCE)\,B^{\star}\big)^k \,, \label{eq:concentrationoverlapfixedgibbs_initbounds1} \\
\Big\langle \exp\Big( \frac{\GAUSSIANMEASURE}2 \|\bar\btheta\|^2 \Big) \Big\rangle &\le \exp\big( K_{\upref{lem:logconcaveboundmomentslemma}}(K_{\upref{lem:lem316}}, D, \GAUSSIANMEASURE_0, \beta_1, \MARGIN, \SIGNALSHORTDEPENDENCE)\, B^{\star} \big)\,. \label{eq:concentrationoverlapfixedgibbs_initbounds2}
\end{align}
Note that the functions we take averages with respect to in the Gibbs measure $\langle \cdot \rangle$ above are a function of $L=1$ replica. Next, note we have $Z' \ge \exp(-K_{\upref{lem:lem316}} B^{\star})$ by Lemma \upref{lem:lem316}. For the interpolators or for the posterior in the logistic case where $u \le 0$, \equpref{eq:concentrationoverlapfixedgibbs_initbounds1}, \equpref{eq:concentrationoverlapfixedgibbs_initbounds2} follow immediately by Lemma \upref{lem:logconcaveboundmomentslemma}, as the corresponding $U(\cdot)$ from \equpref{eq:lem315_measureform} is such that $U \le 0$ and all the $a_j=0$. For the posterior in the GMM case, since $u(x) = \sqrt{\lambda} x$, we note Lemma \upref{lem:lem315} applies where
\begin{align*}
a_0 &= \sqrt{1-v}\big( \sqrt{q_{n,d}} Z + \sqrt{\rho_{n,d} - q_{n,d}} W \big) + s_n \cdot (1-v) x\,, \\
a_1 &= \sqrt{\frac{\lambda}{d}} \Big( \sum_{i=1}^{n-1} s_i + s_n v\Big)\,,\\
a_j &= \sqrt{\frac{\lambda}{d}} \Big( \sum_{i=1}^{n-1} g_{ij} + \sqrt{v} g_{nj} \Big)\text{ for all } 2 \le j \le d-1\,, \\
a_d &= \sqrt{\frac{\lambda t}{d}} \Big( \sum_{i=1}^{n-1} g_{ij} + \sqrt{v} g_{nj} \Big) + \sqrt{(1-t)r} Y\,, \\
U(\theta) &= -\frac{(1-t)(r-\bar r)}2 \theta_d^2 \le 0\,.
\end{align*}
By \equpref{eq:upperboundrhowithinpf}, using that $\E_W[\exp K W] \le K$ we have
\begin{align*}
\log \E_W\big[ \exp a_0 \big] \le K(D, \beta_0, \beta_1, \MARGIN, \lambda) \big(|Z|+|s_n| \big) + K(D, \beta_0, \beta_1, \MARGIN, \lambda)\,.
\end{align*}
It follows that 
\begin{align*}
\frac1{2\beta_0} \sum_{1 \le j \le d} a_j^2 + \log \E_W\big[ \exp a_0 \big] \le K(D, \beta_0, \beta_1, \MARGIN, \lambda) B^{\star}\,.
\end{align*}
Combining the above display with Lemma \upref{lem:logconcaveboundmomentslemma} and Lemma \upref{lem:lem316} now establishes \equpref{eq:concentrationoverlapfixedgibbs_initbounds1}, \equpref{eq:concentrationoverlapfixedgibbs_initbounds2}.

We first prove the desired result for $R_{1,2}$. Let $m = \big\langle \bar\btheta \big\rangle$, thus $\langle R_{1,2} \rangle = \frac{\| m \|^2}d$, and
\begin{equation*}
\big| R_{1,2} - \langle R_{1,2} \rangle \big| \le \Big| \frac{\langle \bar{\btheta}^1, \bar{\btheta}^2 \rangle}{d} - \frac{\langle \bar{\btheta^1}, m \rangle }{d} \Big| + \Big| \frac{\langle \bar{\btheta}^1, m \rangle}{N} - \frac{\langle m, m \rangle}{N} \Big|\,.
\end{equation*}
For fixed $\btheta^1$, the map $f(x) = \frac{\langle \bar\btheta^1, x \rangle}d$ is Lipschitz with constant $\frac{\| \bar\btheta^1\|}d$. Letting $\mu_I(\btheta^l)$ denote the Gibbs measure w.r.t. only the $l$-th replica, which is $\beta$ strongly log-concave, Theorem \upref{thm:logconcaveconcentrationbasics} gives
\[ \int \Big(  \frac{\langle \bar\btheta^1, \bar\btheta^2\rangle }{d} - \frac{\langle \bar\btheta^1, m \rangle }{d} \Big)^{2k} \rmd \mu_I(\btheta^2) \le \Big( \frac{k\, K(\GAUSSIANMEASURE)\, \|\bar\btheta^1\|^2}{d^2} \Big)^k\,.\]
We now integrate this inequality for $\btheta^1$ with respect to $\mu_I(\btheta^1)$. By \equpref{eq:concentrationoverlapfixedgibbs_initbounds1},
\begin{align}
\Big\langle \Big(  \frac{\langle \bar\btheta^1, \bar\btheta^2\rangle }{d} - \frac{\langle \bar\btheta^1, m \rangle }{d} \Big)^{2k} \Big\rangle &= \int \Big(  \frac{\langle \bar\btheta^1, \bar\btheta^2\rangle }{d} - \frac{\langle \bar\btheta^1, m \rangle }{d} \Big)^{2k} \rmd \mu_I(\btheta^1) \rmd \mu_I(\btheta^2) \notag \\
&\le \Big( \frac{k\, K(\GAUSSIANMEASURE)}{d^2} \Big)^k \big\langle \|\bar\btheta^1\|^{2k} \big\rangle \notag \\
&\le \Big( \frac{k\, K(D, \beta_0, \beta_1, \MARGIN, \SIGNALSHORTDEPENDENCE)\, B^{\star}}{d^2} \Big)^k\,.\label{eq:logconcaveR12pt1}
\end{align}
Likewise, the map $f(x) = \frac{\langle x, m \rangle}d$ is Lipschitz with constant $\frac{\| m \|}d$. Thus Theorem \upref{thm:logconcaveconcentrationbasics} gives 
\[ \Big\langle \Big( \frac{\langle \bar\btheta^1, m \rangle}{d}  -\frac{\langle m, m \rangle}{d} \Big)^{2k} \Big\rangle = \int \Big( \frac{\langle \bar\btheta^1, m \rangle}{d}  -\frac{\langle m, m \rangle}{d} \Big)^{2k} \rmd \mu_I(\btheta^1) \le \Big( \frac{k\,K(\GAUSSIANMEASURE)\,\| m \|^2}{d^2} \Big)^k~. \]
Again by Lemma \upref{lem:logconcaveboundmomentslemma} and Jensen's Inequality repeatedly, this gives
\[ \|m\|^{2k} = \big\| \langle \bar\btheta \rangle \big\|^{2k} \le \big\langle \| \bar\btheta \| \big\rangle^{2k} \le \big\langle \| \bar\btheta \|^{2k} \big\rangle \le \big( K(D, \beta_0, \beta_1, \MARGIN, \SIGNALSHORTDEPENDENCE)\,B^{\star}\big)^k\,.\]
Combining with the above display yields
\begin{align}
\Big\langle \Big( \frac{\langle \bar\btheta^1, m \rangle}{d}  -\frac{\langle m, m \rangle}{d} \Big)^{2k} \Big\rangle &\le \Big( \frac{k\,K(D, \beta_0, \beta_1, \MARGIN, \SIGNALSHORTDEPENDENCE)\,B^{\star} }{d^2} \Big)^k\,.\label{eq:logconcaveR12pt2}
\end{align}
As $(a+b)^{2k} \le 2^{2k-1}(a^{2k}+b^{2k})$, combining \equpref{eq:logconcaveR12pt1} and \equpref{eq:logconcaveR12pt2}, we obtain the desired upper bound for $\big\langle (  R_{1,2} - \langle R_{1,2} \rangle )^{2k} \big\rangle$.

We now prove the result for $R_{1,1}$. The idea is to truncate $R_{1,1}$ into a part $f$ where we have the requisite bounds, and another part $\phi$ which is 0 with high probability. Specifically, define the parameter 
\begin{align*}
a := \frac{K_{\upref{lem:logconcaveboundmomentslemma}}(K_{\upref{lem:lem316}}, D, \GAUSSIANMEASURE_0, \beta_1, \MARGIN, \SIGNALSHORTDEPENDENCE) \sqrt{6B^{\star}}}{\sqrt{\beta_0}} > 0
\end{align*}
controlling the truncation, and in terms of $a$ define
\[ f(\btheta) = \min\Big\{ \frac{\| \bar\btheta \|^2}{d}, \frac{a^2}{d} \Big\} = \frac{\big( \min\{ \| \bar\btheta \|, a \} \big)^2}{d}\,.\]
As the cylinder $\big\{\btheta:\| \bar\btheta \| \le a\big\}$ is convex, it follows that $f$ is Lipschitz with constant $\frac{2a}d$. Thus Theorem \upref{thm:logconcaveconcentrationbasics} gives
\begin{align}
\Big\langle \big(f - \langle f \rangle \big)^{2k} \Big\rangle \le \Big( \frac{k\,K(\GAUSSIANMEASURE)\,a^{2}}{d^2} \Big)^k\,.\label{eq:concentrationR11eq0}
\end{align}
Define
\begin{align*}
\phi(\btheta) := \frac{\| \bar\btheta \|^2}{d} - f(\btheta) \ge 0\,.   
\end{align*}
Thus $\phi(\btheta)  \le \frac{\| \bar\btheta \|^2}{d} \, \one\{\|\bar\btheta\| \ge a\}$.
By Lemma \upref{lem:symmetrization} and Cauchy–Schwarz,
\begin{align}
\Big\langle \big(\phi - \langle \phi \rangle \big)^{2k} \Big\rangle &\le 2^{2k} \big\langle \phi^{2k} \big\rangle \notag \\
&\le 2^{2k} \Big\langle \Big(\frac{\| \bar\btheta \|^2}{d}\Big)^{2k} \, \one\{\|\bar\btheta\| \ge a\}\Big\rangle \notag \\
&\le 2^{2k} \Big\langle \Big(\frac{\| \bar\btheta \|^2}{d}\Big)^{4k} \Big\rangle^{1/2} \Big\langle \one\{\|\bar\btheta\| \ge a\}^{2k}\Big\rangle^{1/2} \notag \\
&= 2^{2k} \Big\langle \Big(\frac{\| \bar\btheta \|^2}{d}\Big)^{4k} \Big\rangle^{1/2} \Big\langle \one\{\|\bar\btheta\| \ge a\}\Big\rangle^{1/2} \,.\label{eq:concentrationR11eq1}
\end{align}
By \equpref{eq:concentrationoverlapfixedgibbs_initbounds1}, \equpref{eq:concentrationoverlapfixedgibbs_initbounds2}, Markov's Inequality and the choice of $a$ gives 
\begin{align}
\Big\langle \one\{\|\bar\btheta\| \ge a\}\Big\rangle &= \P\Big( \exp\Big( \frac{\GAUSSIANMEASURE}2 \|\bar\btheta\|^2 \Big)  \ge \exp\Big( \frac{\beta a^2}2 \Big) \Big) \notag \\
&\leq 
\exp\Big(  K_{\upref{lem:logconcaveboundmomentslemma}}(K_{\upref{lem:lem316}}, D, \GAUSSIANMEASURE_0, \beta_1, \MARGIN, \SIGNALSHORTDEPENDENCE)\, B^{\star} - \frac{\GAUSSIANMEASURE a^2}{2} \Big) \le \exp\big(-2B^{\star}\big)\,. \label{eq:concentrationR11eq3} 
\end{align}
Plugging in \equpref{eq:concentrationoverlapfixedgibbs_initbounds1}, \equpref{eq:concentrationR11eq3} into \equpref{eq:concentrationR11eq1}, we obtain
\[ \Big\langle \big(\phi - \langle \phi \rangle \big)^{2k} \Big\rangle \leq \exp\big(-B^{\star}\big) \Big( \frac{K(D, \GAUSSIANMEASURE_0, \beta_1, \MARGIN, \SIGNALSHORTDEPENDENCE)\, B^{\star}}{d} \Big)^{2k}\,. \]
Since $R_{1,1} = \| \bar\btheta \|^2 / d = f + \phi$, using that $(a + b)^{2k} \leq 2^{2k}(a^{2k} + b^{2k})$, combining with \equpref{eq:concentrationR11eq0} and recalling our choice of $a$, we obtain
\begin{align*}
\Big\langle \big(R_{1, 1} - \big\langle R_{1, 1} \big\rangle \big)^{2k} \Big\rangle &\leq 2^{2k} \Bigg( \Big\langle \big(f - \langle f \rangle \big)^{2k} \Big\rangle + \Big\langle \big(\phi - \langle \phi \rangle \big)^{2k} \Big\rangle  \Bigg) \\
&\leq \Big( \frac{k\,K(D, \GAUSSIANMEASURE_0, \beta_1, \MARGIN, \SIGNALSHORTDEPENDENCE)\,B^{\star}}{d^2} \Big)^k + \exp\big(-B^{\star}\big) \Big( \frac{K(D, \GAUSSIANMEASURE_0, \beta_1, \MARGIN, \SIGNALSHORTDEPENDENCE)\, B^{\star}}{d} \Big)^{2k}\,.
\end{align*}
Note for $x \ge 0$ that $x^k \le k^k \exp(x)$, thus $\exp(-x) \le \left( k / x \right)^k$ for $x>0$. 
Clearly $B^{\star} > 0$. 
Hence
\begin{align*}
\exp\big(-B^{\star}\big) \Big( \frac{K(D, \GAUSSIANMEASURE_0, \beta_1, \MARGIN, \SIGNALSHORTDEPENDENCE)\, B^{\star}}{d} \Big)^{2k} &\le \Big( \frac{k}{B^{\star}} \Big)^k \Big( \frac{K(D, \GAUSSIANMEASURE_0, \beta_1, \MARGIN, \SIGNALSHORTDEPENDENCE)\, B^{\star}}{d} \Big)^{2k} \\
&\le \Big( \frac{k\, K(D,\GAUSSIANMEASURE_0, \beta_1, \MARGIN, \SIGNALSHORTDEPENDENCE)\, B^{\star}}{d^2} \Big)^k\,.  
\end{align*}
Combining the above two displays gives the desired upper bound on $\big\langle (R_{1, 1} - \big\langle R_{1, 1} \big\rangle )^{2k} \big\rangle$.
\end{proof}

\paragraph{Proof of Proposition \upref{prop:concentrationmeangibbswrtdisorder}}
We will follow the strategy of \cite{talagrand2010mean} of considering $\langle R_{1, 1} \rangle$ and $\langle R_{1, 2}\rangle$ as functions of the disorder $(s_j)_{j \le n} \times (g_{i, j})_{2 \le i \le n, 1 \le j \le d} \times (\ANNEALEDDISORDER^l)_{1 \le l \le L} \times \QUENCHEDDISORDER \times Y$. The density of the disorder is strongly log-concave by Assumption \upref{ass:disorder_technical_assumption}, and to apply Theorem \upref{thm:logconcaveconcentrationbasics} on concentration of Lipschitz functions of strongly log-concave measures, we need to establish control on the Lipschitz constant of these functions.
\begin{lemma}\label{lem:concentrationofmeanslipschitzexpression}
Let $f=R_{1,1}$ or $R_{1,2}$ and let $\grad$ denote the derivative of the function $\langle f \rangle$ w.r.t. the disorder. Let 
\begin{align*}
u(\btheta) := \Big( u'\big(S_{1,t}(\btheta)\big), \ldots, u'\big(S_{n-1,t}(\btheta)\big), u'\big(S_{n,t,v}(\btheta)\big) \Big)^T \in \R^n\,,
\end{align*}
and analogously define $u^l(\btheta^l)$ for each replica $l$. Then $\| \grad \|^2$ is upper bounded by
\begin{align*}
&K\Bigg( \big\langle \dot{f} \dot{\theta^1_d} \big\rangle^2 + \frac1d \Big\| \Big\langle \big( \btheta^1 - \langle \btheta^1 \rangle\big) \dot{f} u'\big(S_{n,t,v}^1(\btheta^1)\big) \Big\rangle \Big\|^2 + \frac1d \sum_{i \le n-1} \Big\| \Big\langle \big( \btheta^1 - \langle \btheta^1 \rangle\big) \dot{f} u'\big(S_{i,t}^1(\btheta^1)\big) \Big\rangle \Big\|^2 \\
&\qquad\qquad \qquad + \Big( 1 + \frac1d \big\langle \| \btheta^1 \|^2 \big\rangle \Big) \Big\| \Big\langle \big( u^1(\btheta^1) - \big\langle u^1(\btheta^1) \big\rangle \big) \dot{f} \Big\rangle \Big\|^2 \Bigg)\,,
\end{align*}
where $K$ depends on $D, \beta_0, \beta_1, \SIGNALSHORTDEPENDENCE$ and also on $\MARGIN$ for the interpolators, and where we use the notation
\begin{align*}
\dot{f} := f - \langle f \rangle\,,\,\dot{\btheta^l} := \btheta^l - \langle \btheta^l\rangle\,,\,\dot{\theta_j^l} := \theta_j^l - \langle \theta_j^l\rangle\,,
\end{align*}
and
\begin{align*}
\dot{u}'(S^l) := u'(S^l) - \langle u'(S^l)\rangle\quad\text{ for }\quad S^l=S^l_{i,t}(\btheta^l)\text{ or }S^l_{n,t,v}(\btheta^l)\,.
\end{align*}
\end{lemma}
\begin{proof}
First, consider $x$ an element of the quenched disorder, i.e. not one of the $\ANNEALEDDISORDER^l$. For just this proof, we let $H_{t, v, n, d}(\btheta^l) := u\big(S^l_{n,t,v}(\btheta^{l})\big) + H_{t,n-1,d}(\btheta^{l}) - \beta \| \btheta^{l} \|^2$. By \equpref{eq:gibbsmeasurerewriteconcentration}, we remark that 
\begin{align*}
\frac{\partial}{\partial x} \langle f \rangle
&= \frac1{Z'^L} \int_{\R^d} \E_{\ANNEALEDDISORDER}\Bigg[ f \sum_{1 \le l \le L}\Big\{ \one\{\theta_1^{l}/\sqrt{d} \in I\} \exp\big( H_{t, v, n, d}(\btheta^l) \big) \cdot \frac{\partial}{\partial x} H_{t, v, n, d}(\btheta^l) \\
&\qquad \qquad \qquad \qquad \cdot \prod_{1 \le l' \neq l \le L} \one\{\theta_1^{l'}/\sqrt{d} \in I\} \exp\big( H_{t, v, n, d}(\btheta^{l'}) \big) \Big\} \Bigg] \, \rmd\btheta^1 \cdots \rmd\btheta^L \\
&\qquad \qquad - \frac{L}{Z'^{L+1}} \int_{\R^d} \E_{\ANNEALEDDISORDER}\Bigg[ f \prod_{1 \le l \le L} \one\{\theta_1^l/\sqrt{d} \in I\} \exp\big( H_{t,v,n,d}(\btheta^l) \big) \Bigg]\, \rmd\btheta^1 \cdots \rmd\btheta^L \\
&\qquad \qquad \qquad \qquad \cdot \int_{\R^d} \E_{\ANNEALEDDISORDER}\Big[ \one\{\theta^{L+1}_1/\sqrt{d} \in I\} \exp\big( H_{t,v,n,d}(\btheta^{L+1}) \big) \cdot \frac{\partial}{\partial x} H_{t,v,n,d}(\btheta^{L+1}) \Big]\, \rmd\btheta^{L+1} \\
&= \sum_{1 \le l \le L} \Big\langle f \frac{\partial}{\partial x} H_{t,v,n,d}(\btheta^l) \Big\rangle - \langle f \rangle \Big\langle \frac{\partial}{\partial x} H_{t,v,n,d}(\btheta) \Big\rangle\,.
\end{align*}
Note that the $\theta^1, \ldots, \theta^L$ such that there exists some $l'$ with $\frac1{\sqrt{d}} \theta^{l'} \not\in I$ make no contribution to $\langle f \rangle$. Thus $\langle f \rangle = \langle f' \rangle$ where $f'$ is the restriction of $f$ to $\prod_{1 \le l \le L} \one\{\theta_1^l/\sqrt{d} \in I\}$. Consequently we have 
\begin{align*}
\Big\langle \frac{\partial}{\partial x} H_{t,v,n,d}(\btheta) \Big\rangle = \begin{cases}
\Big\langle u'\big(S_{i,t}(\btheta)\big) \cdot \frac{\theta_1}{\sqrt{d}} \Big\rangle &: x = s_i, i \le n-1 \\
\Big\langle u'\big( S_{n,t,v}(\btheta) \big) \cdot \Big( \frac{v\theta_1}{\sqrt{d}} + (1-v)x \Big) \Big\rangle &: x = s_n \\
\Big\langle u'\big(S_{i,t}(\btheta)\big) \cdot \frac{\theta_j}{\sqrt{d}} \Big\rangle &: x = g_{i,j}, i \le n-1, j \le d-1 \\
\Big\langle u'\big(S_{i,t}(\btheta)\big) \cdot \theta_d\sqrt{\frac{t}{d}} \Big\rangle &: x = g_{i,d}, i \le n-1 \\
\Big\langle u'\big(S_{n,t,v}(\btheta)\big) \cdot \theta_j \sqrt{\frac{v}{d}} \Big\rangle &: x = g_{n, j}, j \le d-1 \\
\Big\langle u'\big(S_{n,t,v}(\btheta)\big) \cdot \theta_d \sqrt{\frac{vt}{d}} \Big\rangle &: x = g_{n, d} \\
\Big\langle u'\big(S_{n,t,v}(\btheta)\big) \cdot \sqrt{(1-v) q_{n,d}} \Big\rangle &: x = \QUENCHEDDISORDER \\
\Big\langle \theta_d\sqrt{(1-t) r} \Big\rangle  &: x = Y\,.
\end{cases}
\end{align*}
An analogous derivation applies also for $\Big\langle f \frac{\partial}{\partial x} H_{t,v,n,d}(\btheta^l) \Big\rangle$.

Now consider $x = \ANNEALEDDISORDER^l$ for some $1 \le l \le L$. We remark that 
\begin{align*}
&\frac{\partial}{\partial \ANNEALEDDISORDER^l} \E_{\ANNEALEDDISORDER^l}\Big[ \one\{\theta_1^l/\sqrt{d}\in I\} \exp\big(H_{t,v,n,d}(\btheta)\big) \Big] \\
&\quad = -\E_{\ANNEALEDDISORDER^l}\Big[ \one\{\theta_1^l/\sqrt{d}\in I\} \exp\big(H_{t,v,n,d}(\btheta)\big) \Big] \\
&\qquad \qquad \qquad + \E_{\ANNEALEDDISORDER^l}\Big[ \one\{\theta_1^l/\sqrt{d}\in I\} \exp\big(H_{t,v,n,d}(\btheta)\big) \cdot u'\big(S_{n,t,v}^l(\btheta^l)\big) \frac{\partial}{\partial \ANNEALEDDISORDER^l} u\big(S_{n,t,v}^l(\btheta^l)\big) \Big] \\
&\quad = -\E_{\ANNEALEDDISORDER^l}\Big[ \one\{\theta_1^l/\sqrt{d}\in I\} \exp\big(H_{t,v,n,d}(\btheta)\big) \Big] \\
&\qquad \qquad \qquad + \E_{\ANNEALEDDISORDER^l}\Big[ \one\{\theta_1^l/\sqrt{d}\in I\} \exp\big(H_{t,v,n,d}(\btheta)\big) \cdot u'\big(S_{n,t,v}^l(\btheta^l)\big) \sqrt{(1-v)(\rho_{n,d}-q_{n,d})} \Big]\,,
\end{align*}
where the last step follows from the precense of the constraint $\one\{\theta_1^l/\sqrt{d}\in I\}$ and analogous logic as above. Since $f=R_{1,1}$ or $R_{1,2}$ does not depend on the $\ANNEALEDDISORDER^l$, it follows from an identical derivation as above that
\begin{align*}
\frac{\partial}{\partial \ANNEALEDDISORDER^l} \langle f \rangle &= \sum_{1 \le l \le L} -\langle f \rangle + \langle f \rangle + \sqrt{(1-v)(\rho_{n,d}-q_{n,d})} \Big( \Big\langle f u'\big(S_{n,t,v}^l(\btheta^l)\big) \Big\rangle - \langle f \rangle \Big\langle u'\big(S_{n,t,v}^l(\btheta^l)\big) \Big\rangle \Big) \\
&= \sum_{1 \le l \le L} \sqrt{(1-v)(\rho_{n,d}-q_{n,d})} \Big( \Big\langle f u'\big(S_{n,t,v}^l(\btheta^l)\big) \Big\rangle - \langle f \rangle \Big\langle u'\big(S_{n,t,v}^l(\btheta^l)\big) \Big\rangle \Big) \,.
\end{align*}
Note $\langle \dot{f} \rangle=0$. Thus for any $1 \le l \le L$,
\begin{align*}
\big\langle \dot{f} u'(S^l) \big\rangle = \big\langle \dot{f} \dot{u}'(S^l) \big\rangle\,,\,\big\langle \dot{f} \theta^l_d \big\rangle = \big\langle \dot{f} \dot{\theta^l_d} \big\rangle\,,
\end{align*}
and for any $1 \le l \le L$ and any $1 \le j \le d$,
\begin{align*}
\big\langle \dot{f} \theta^l_j u'(S^l) \big\rangle &= \big\langle \dot{f} \dot{\theta^l_j} u'(S^l) \big\rangle + \big\langle \theta^l_j \big\rangle \big\langle \dot{f} u'(S^l) \big\rangle = \big\langle \dot{f} \dot{\theta^l_j} u'(S^l) \big\rangle + \big\langle \theta^l_j \big\rangle \big\langle \dot{f} \dot{u}'(S^l) \big\rangle \,.
\end{align*}
Since $t, v \in [0,1]^2$, $L \le 2$ for $f=R_{1,1}, R_{1,2}$, and as $f=R_{1,1}, R_{1,2}$ are symmetric in $\btheta^1, \btheta^2$, combining all of the displays above implies that 
\begin{align*}
\| \grad \|^2 \le K \Big(\expressionI + \expressionII + \expressionIII\Big)\,
\end{align*}
where
\begin{align*}
\expressionI &:= \frac1d \sum_{1 \le j \le d} \big\langle \dot{f}\dot{\theta^1_j} u'\big(S^1_{n,t,v}(\btheta^1) \big)\big\rangle^2 + \frac1d \sum_{i \le n-1, 1 \le j \le d} \big\langle \dot{f}\dot{\theta^1_j} u'\big(S_{i,t} (\btheta^1) \big) \big\rangle^2\,, \\
\expressionII &:= \Big( \frac1d \sum_{1 \le j \le d}\langle \theta^1_j \rangle^2 \Big) \Big( \big\langle \dot{f} \dot{u}'\big(S^1_{n,t,v} (\btheta^1) \big) \big\rangle^2 + \sum_{i \le n-1} \big\langle \dot{f} \dot{u}'\big(S_{i,t} (\btheta^1) \big) \big\rangle^2 \Big)\,, \\
\expressionIII &:= \max\big\{ x^2, q_{n,d}, r, \rho_{n,d}-q_{n,d} \big\}\Big( \big\langle \dot{f} \dot{u}'\big(S^1_{n,t,v} (\btheta^1) \big) \big\rangle^2 + \big\langle \dot{f} \dot{\theta^1_d} \big\rangle^2 \Big)\,.
\end{align*}
Note $|x| \le 1$ and that $r \le K(D)$. Moreover by \equpref{eq:upperboundrhowithinpf}, we have $|q_{n,d}|, |\rho_{n,d}| \le K(D, \beta_0, \beta_1, \kappa, \SIGNALSHORTDEPENDENCE)$. Also remark that
\begin{align*}
\frac1d \sum_{1 \le j \le d}\langle \theta^1_j \rangle^2 \le \frac1d \sum_{1 \le j \le d} \langle (\theta^1_j)^2 \rangle = \frac1d \big\langle \| \btheta^1 \|^2 \big\rangle\,.    
\end{align*}
Next, for $S^1$ either $S^1_{n,t,v}(\btheta^1)$ or $S_{i,t}(\btheta^1)$, we observe that 
\begin{align*}
&\sum_{1 \le j \le d} \Big\langle \dot{f}\dot{\theta^1_j} u'\big(S^1(\btheta^1) \big)\Big\rangle^2 \\
&\quad = \sum_{1 \le j \le d} \Big\langle \dot{f}(\btheta^1, \ldots, \btheta^L) \dot{f}(\btheta^{L+1}, \ldots, \btheta^{2L}) u'\big(S^1(\btheta^1) \big) u'\big(S^1(\btheta^{L+1}) \big) \dot{\theta^1_j} \dot{\theta^{L+1}_j} \Big\rangle \\
&\quad = \Big\langle \dot{f}(\btheta^1, \ldots, \btheta^L) \dot{f}(\btheta^{L+1}, \ldots, \btheta^{2L}) u'\big(S^1(\btheta^1) \big) u'\big(S^1(\btheta^{L+1}) \big) \cdot \big(\btheta^1 - \langle \btheta^1 \rangle \big)^T \big(\btheta^{L+1} - \langle \btheta^{L+1} \rangle \big) \Big\rangle \\
&\quad = \Big\| \Big\langle \big( \btheta^1 - \langle \btheta^1 \rangle\big) \dot{f} u'(S^1) \Big\rangle \Big\|^2 \,.
\end{align*}
Analogously, we have
\begin{align*}
&\big\langle \dot{f} \dot{u}'\big(S^1_{n,t,v} (\btheta^1) \big) \big\rangle^2 + \sum_{i \le n-1} \big\langle \dot{f} \dot{u}'\big(S_{i,t} (\btheta^1) \big) \big\rangle^2 \\
&\quad = \Big\langle \dot{f}(\btheta^1, \ldots, \btheta^L) \dot{f}(\btheta^{L+1}, \ldots, \btheta^{2L}) \cdot \Big( u\big( S^1(\btheta^1) \big) - \big\langle u\big( S^1(\btheta^1) \big) \big\rangle \Big)^T  \Big( u\big( S^1(\btheta^{L+1}) \big) - \big\langle u\big( S^1(\btheta^{L+1}) \big) \big\rangle \Big) \Big\rangle \\
&\quad = \Big\| \Big\langle \big( u\big( S^1(\btheta^1) \big) - \big\langle u\big( S^1(\btheta^1) \big) \big\rangle \big) \dot{f} \Big\rangle \Big\|^2\,. 
\end{align*}
Using these observations, the Lemma follows from combining terms in the above expressions $\expressionI, \expressionII, \expressionIII$. 
\end{proof}

We now aim to upper bound the expressions appearing in the bound from Lemma \upref{lem:concentrationofmeanslipschitzexpression}. 
\begin{lemma}\label{lem:lem3113}
For any function $f$ of $\btheta^1, \ldots, \btheta^L$ mapping to $\R$ (possibly depending on the disorder), we have
\[ \Big\| \Big\langle \big( \btheta^1 - \langle \btheta^1 \rangle\big) f \Big\rangle \Big\|^2 \le K \langle f^2 \rangle\,,\]
and consequently we also have $\big\langle f \dot{\theta^1_d} \big\rangle^2 \le K \langle f^2 \rangle$, where $K$ depends on $\beta_0, \beta_1$.
\end{lemma}
\begin{proof}
Consider any $y\in \R^d$. By Cauchy-Schwarz,
\begin{align*}
\Big| \big\langle \Big( \btheta^1 - \langle \btheta^1 \rangle\Big) f \big\rangle^T y \Big| &= \Big\langle f  \big( \btheta^1 - \langle \btheta^1 \rangle\big)^T y \Big\rangle \le \big\langle f^2 \big\rangle^{1/2} \Big\langle \Big(\big( \btheta^1 - \langle \btheta^1 \rangle\big)^T y\Big)^2\Big\rangle^{1/2}\,.
\end{align*}
By Theorem \upref{thm:logconcaveconcentrationbasics} applied with the function $(\btheta^1)^T y$ of $\btheta^1, \ldots, \btheta^L$ which is $\| y \|$-Lipschitz,
\begin{align*}
\Big\langle \Big(\big( \btheta^1 - \langle \btheta^1 \rangle\big)^T y\Big)^2\Big\rangle \le K(\beta_0, \beta_1) \|y\|^2\,.
\end{align*}
Hence, we have 
\begin{align*}
\Big| \Big\langle \big( \btheta^1 - \langle \btheta^1 \rangle\big) f \Big\rangle^T y \Big| &\le K(\beta_0, \beta_1) \big\langle f^2 \big\rangle^{1/2} \|y\|\,,
\end{align*}
and since $y$ is arbitrary the first part of the Lemma follows. The second part of this Lemma follows upon observing that analogously as in the proof of Lemma \upref{lem:concentrationofmeanslipschitzexpression}, we have 
\begin{align*}
\big\langle f \dot{\theta^1_d} \big\rangle^2 \le \sum_{1 \le j \le d} \big\langle f \dot{\theta^1_j} \big\rangle^2 = \Big\| \Big\langle \big( \btheta^1 - \langle \btheta^1 \rangle \big) f\Big\rangle \Big\|^2\,.
\end{align*}
Applying the first part of this Lemma then finishes the proof.
\end{proof}

To bound the other term in Lemma \upref{lem:concentrationofmeanslipschitzexpression}, define the matrix $M$ as follows: for all $1 \le i \le n$ and $1 \le j \le d$, we let
\begin{align}
(M)_{i, j} &:= \begin{cases} s_i &: 1 \le i \le n-1, j=1 \\ g_{i, j} &: 1 \le i \le n-1, 2 \le j \le d-1 \\ \sqrt{t} g_{i, d} &: 1 \le i \le n-1, j=d \\ vs_n &: i=n, j=1 \\ \sqrt{v} g_{n,j} &: i=n, 2 \le j \le d-1 \\ \sqrt{vt} g_{n,d} &: i=n, j=d\,.\end{cases} \label{eq:overlapconcentrationmatMdef} \\
B' &:= \| M \|_{\OPNORM}\,.\label{eq:overlapconcentrationB'def}
\end{align}
We now control $B'$ as follows:
\begin{lemma}\label{lem:overlapconcentrationopnormcontrol}
We have $B'^2 \le K d$ with probability at least $1-K \exp(-d)$, where $K$ depends on $\lambda$.
\end{lemma}
\begin{proof}
Note $B'^2 = \| M M^T \|_{\OPNORM}$. Consider any unit vector $y \in \R^n$. Note for all $1 \le j \le d$ that
\begin{align}\label{eq:concentration_controlopnorm_0}
\big(M^{T} y\big)_j &= \begin{cases}
\sum_{i=1}^{n-1} s_i y_i + v s_n y_n &: j=1 \\
\sum_{i=1}^{n-1} g_{i, j} y_i + \sqrt{v} g_{n, j} y_n &: 2 \le j \le d-1\\
\sqrt{t} \sum_{i=1}^{n-1} g_{i, d} y_i + \sqrt{tv} g_{n, d} y_n &: j=d\,.
\end{cases}
\end{align}
Let $\bar{M}$ be the $n \times d-1$ matrix comprising of the last $d-1$ columns of $M$; note this matrix does not involve the $s_i$, which we handle in a separate step. Observe that 
\begin{align}\label{eq:concentration_controlopnorm_1}
\sum_{j=2}^{d-1} \Big( \sum_{i=1}^{n-1} g_{i, j} y_i + \sqrt{v} g_{n, j} y_n \Big)^2 + \Big( \sqrt{t} \sum_{i=1}^{n-1} g_{i, d} y_i + \sqrt{tv} g_{n, d} y_n \Big)^2 = \| \bar{M}^T y \|^2 \le \|\bar{M}\|_{\OPNORM}^2\,.
\end{align}
By concentration of operator norm for matrices with i.i.d. sub-Gaussian entries, e.g. Theorem 4.4.5 of \cite{vershynin2018high}, as $t,v\in[0,1]^2$ we have $\|\bar{M}\|_{\OPNORM} \le K\sqrt{d}$ with probability at least $1-K \exp(-d)$.
Also, as $v \in [0,1]$,
\begin{align}\label{eq:concentration_controlopnorm_2}
\Big( \sum_{i=1}^{n-1} s_i y_i + v s_n y_n \Big)^2 \le \|y\|^2 \| s \|^2 = \sum_{i=1}^n s_i^2\,.
\end{align}
As each $s_i^2$ is $\signalsubexp$ sub-Exponential and $n \le 2d$, the latter is at most $K(\signalsubexp) d$ with probability at least $1-K \exp(-d)$ by Bernstein's Inequality. Combining \equpref{eq:concentration_controlopnorm_0} with the aforementioned high-probability bounds for \equpref{eq:concentration_controlopnorm_1} and \equpref{eq:concentration_controlopnorm_2}, the Lemma follows.
\end{proof}
Next note for any sequences $(x_j)_{j \le d}$, $(y_i)_{i \le n}$, we have 
\[ \sum_{i \le n, j \le d} M_{i, j} x_j y_i \le B' \Big( \sum_{j \le d} x_j^2 \Big)^{1/2} \Big( \sum_{i \le n} y_i^2 \Big)^{1/2}\,. \]
Now applying the above inequality for $x_j = \sum_{i \le n} M_{i, j} y_i$, we obtain that
\begin{align}
\Big( \sum_{j \le d} \Big( \sum_{i \le n} M_{i, j} y_i \Big)^2 \Big)^{1/2} \le B' \Big( \sum_{i \le n}  y_i^2 \Big)^{1/2}\,.\label{eq:opnormboundtal347}  
\end{align}
\begin{lemma}\label{lem:tallem3114}
Consider any sequence $y = (y_i)_{i \le n} \in \R^n$. Then the function $f(\btheta^1, \ldots, \btheta^L) = \big\langle u^1(\btheta^1), y \big\rangle$ is Lipschitz with constant $B' D \| y\| / \sqrt{d}$.
\end{lemma}
\begin{proof}
Clearly $f$ only depends on $\btheta^1$ among $\theta^1, \ldots, \theta^L$. Note for all $1 \le j \le d$ that by construction of $M_{i,j}$,
\[ \frac{\partial}{\partial \theta^1_j} \big\langle u^1(\btheta^1), y \big\rangle = \frac1{\sqrt{d}} \Big( \sum_{1 \le i \le n-1} y_i u''\big(S_{i,t} (\btheta^1) \big) M_{i,j} + y_n u''\big(S_{n,t,v} (\btheta^1) \big) M_{n,j} \Big)\,.\]
Applying \equpref{eq:opnormboundtal347} for the sequence $\Big( y_1 u''\big(S_{1,t}(\btheta^1)_, \ldots, y_{n-1} u''\big(S_{n-1,t}(\btheta^1)\big), y_n u''\big(S_{n,t,v}(\btheta^1)\big) \Big)$ and recalling $|u''| \le D$ now gives that 
\[ \Big( \sum_{j \le d} \Big( \frac{\partial f}{\partial \theta^1_j} \Big)^2 \Big)^{1/2} \le \frac{B'}{\sqrt{d}} \Big( \sum_{i \le n-1} y_i^2 u''\big(S_{i,t} (\btheta^1)\big)^2 + y_n^2 u''\big(S_{n,t,v}(\btheta^1)\big)^2 \Big)^{1/2} \le \frac{B' D \| y\|}{\sqrt{d}}\,. \]
Since $f$ only depends on $\btheta^1$, the conclusion follows.
\end{proof}
\begin{lemma}\label{lem:tallem3115}
For any function $f$ of $\btheta^1, \ldots, \btheta^L$ mapping to $\R$ (possibly depending on the disorder), we have 
\[ \Big\| \Big\langle \big( u^1(\btheta^1) - \big\langle u^1(\btheta^1) \big\rangle \big) f \Big\rangle \Big\|^2 \le \frac{K B'^2 }{d} \langle f^2 \rangle\,, \]
where $K$ depends on $D, \GAUSSIANMEASURE$.
\end{lemma}
\begin{proof}
Consider any vector $y \in \R^n$. By Cauchy-Schwarz,
\begin{align*}
\Big\langle \big( u^1(\btheta^1) - \big\langle u^1(\btheta^1) \big\rangle \big) f \Big\rangle^T y &= \Big\langle f \big( u^1(\btheta^1) - \big\langle u^1(\btheta^1) \big\rangle \big)^T y \Big\rangle \le \langle f^2 \rangle^{1/2} \Big\langle \Big( \big( u^1(\btheta^1) - \big\langle u^1(\btheta^1) \big\rangle \big)^T y \Big)^2 \Big\rangle^{1/2}\,.
\end{align*}
Consider the function $\btheta \rightarrow \big\langle u(\btheta), y \big\rangle$. By Lemma \upref{lem:tallem3114}, this function is Lipschitz with constant $B' D \| y \| / \sqrt{d}$.
Applying Theorem \upref{thm:logconcaveconcentrationbasics}, we obtain
\[ \Big\langle \Big( \big( u^1(\btheta^1) - \big\langle u^1(\btheta^1) \big\rangle \big)^T y \Big)^2 \Big\rangle \le \frac{K(\beta) B'^2 D^2 \|y\|^2}{d}\,.\]
Since $y \in \R^n$ was arbitrary, the desired conclusion follows.
\end{proof}

\begin{corollary}\label{corr:overlapslipschitzonC}
For $f=R_{1,1}, R_{1,2}$, defining $\grad$ as in Lemma \upref{lem:concentrationofmeanslipschitzexpression}, we have 
\begin{align*}
\| \grad \|^2 \le K \Bigg( \frac{B^{\star}}{d^2} + \Big(1+\frac{B^{\star}}d\Big) \frac{B^{\star} B'^2 }{d^3} \Bigg)\,,
\end{align*}
where $K$ depends on $D, \beta_0, \beta_1, \SIGNALSHORTDEPENDENCE$ and also on $\MARGIN$ for the interpolators.
\end{corollary}
\begin{proof}
Applying the bound from Lemma \upref{lem:concentrationofmeanslipschitzexpression}, it remains to upper bound 
\begin{align*}
\expressionI &:= \big\langle \dot{f} \dot{\theta_d^1} \big\rangle^2 + \frac1{d} \big\| \langle \big( \btheta^1 - \langle \btheta^1 \rangle\big) \dot{f} u'(S^1_{n,t,v}) \rangle \big\|^2 + \frac1{d} \sum_{i \le n-1} \big\| \langle \big( \btheta^1 - \langle \btheta^1 \rangle\big) \dot{f} u'(S^1_{i,t}) \rangle \big\|^2 \,, \\
\expressionII &:= \Big( 1 + \frac1d \big\langle \| \btheta^1 \|^2 \big\rangle \Big) \Big\| \Big\langle \big( u^1(\btheta^1) - \big\langle u^1(\btheta^1) \big\rangle \big) \dot{f} \Big\rangle \Big\|^2\,.
\end{align*}
First, we upper bound $\expressionI$. To this end, applying Lemma \upref{lem:lem3113} with the function $\dot{f}$ or $\dot{f} u'(S^1)$ where $S^1 = S^1_{i,t}$ or $S^1_{n,t,v}$, we have 
\begin{align*}
\expressionI \le K(\beta_0, \beta_1) \Big( \langle \dot{f}^2 \rangle + \frac1d \langle \dot{f}^2 u'(S^1_{n,t,v})^2 \rangle + \frac1d \sum_{i \le n-1} \langle \dot{f}^2 u'(S^1_{i,t})^2 \rangle \Big) \le \frac{ K(D, \beta_0, \beta_1, \MARGIN, \lambda) B^{\star}}{d^2}\,,
\end{align*}
where we use $|u'| \le D$, Proposition \upref{prop:concentrationoverlapfixedgibbs} to control $\langle \dot{f}^2 \rangle$, and that $n/d \le \DENSITYBOUND = 2$. 

Next, we upper bound $\expressionII$. First, by \equpref{eq:concentrationoverlapfixedgibbs_initbounds1},
\begin{align*}
\Big( 1 + \frac1d \big\langle \| \btheta^1 \|^2 \big\rangle \Big) \le 1 + \frac{K(D, \beta_0, \beta_1, \kappa, \lambda) B^{\star}}d\,.
\end{align*}
By Lemma \upref{lem:tallem3115} and Proposition \upref{prop:concentrationoverlapfixedgibbs},
\begin{align*}
\Big\| \Big\langle \big( u^1(\btheta^1) - \big\langle u^1(\btheta^1) \big\rangle \big) \dot{f} \Big\rangle \Big\|^2 \le \frac{K(D, \beta) B'^2}d \langle \dot{f}^2 \rangle \le \frac{K(D, \beta_0, \beta_1, \MARGIN, \lambda) B^{\star} B'^2}{d^3} \,.
\end{align*}
This yields an upper bound on $\expressionII$, and combining the upper bounds on $\expressionI$ and $\expressionII$ establishes Corollary \upref{corr:overlapslipschitzonC}.
\end{proof}

We now have done the necessary preparation to prove Proposition \upref{prop:concentrationmeangibbswrtdisorder}.
\begin{proof}[Proof of Proposition \upref{prop:concentrationmeangibbswrtdisorder}]
For this proof, we let $f= \langle R_{1,1} \rangle$ or $\langle R_{1,2} \rangle$. 
Similarly as the proof of Theorem \upref{thm:posteriorgaussianmeasurersformula}, we consider the space $\mathcal{S} = \R^{nd + L + 2}$ where the first $n$ coordinates correspond to the $s_j$, the next $(n-1)d$ coordinates correspond to $g_{i,j}$, the next $L$ coordinates correspond to $\ANNEALEDDISORDER^L$, and the last two coordinates correspond to $\QUENCHEDDISORDER$ and $Y$. 
We endow $\mathcal{S}$ with the product measure $\gamma$ where the measure on each coordinate is given by law of each corresponding part of the disorder, that is, $s_j, g_{i,j}, \ANNEALEDDISORDER^L, \QUENCHEDDISORDER$ or $Y$. Thus integration w.r.t. $\gamma$ corresponds to taking expectation w.r.t. the disorder.

Consider the set 
\begin{align*}
C := \{ x \in \mathcal{S}\,:\, B^{\star} \le K(\SIGNALSHORTDEPENDENCE)d\,,\, B'^2 \le K(\SIGNALSHORTDEPENDENCE)d\} \,,
\end{align*}
where by \equpref{eq:boundmgfBstar} and Lemma \upref{lem:overlapconcentrationopnormcontrol}, $K(\SIGNALSHORTDEPENDENCE)$ has been chosen large enough so that
\begin{align*}
\P(C^c) \le K(\SIGNALSHORTDEPENDENCE) \exp(-d)\,,    
\end{align*}
where probability here is w.r.t. $\gamma$.
Moreover, $C$ is a convex subset of $\mathcal{S}$, because $B^{\star}$ is a convex function of the disorder, the operator norm of a matrix is a convex function of its entries, and the entries of $M$ are an affine transformation of the disorder.

Next, let $\gamma'$ be a probability measure defined on $C$ with density proportional to that of $\gamma$. By Assumption \upref{ass:disorder_technical_assumption}, $\gamma'$ is $\min\{1, \SIGNALGAUSSIAN\}$ strongly log-concave. By Corollary \upref{corr:overlapslipschitzonC}, $f$ is Lipschitz with constant $\frac{K(D, \beta_0, \beta_1, \MARGIN, \SIGNALSHORTDEPENDENCE)}{\sqrt{d}}$ on $C$. Let $m = \int f \rmd \gamma'$. Since $C$ is convex and since $\gamma'$ is supported on $C$, Theorem \upref{thm:logconcaveconcentrationbasics} gives that for all $k \ge 1$,
\begin{align*}
\int \big(f - m)^{2k} \rmd \gamma' \le \Big(\frac{K(D, \beta_0, \beta_1, \MARGIN, \SIGNALMEDIUMDEPENDENCE) k}{d}\Big)^k\,.
\end{align*}
Thus by definition of $\gamma'$,
\begin{align*}
\E\big[ (f-m)^{2k} \one\{C\} \big] &= \int_C (f-m)^{2k} \rmd \gamma = \gamma(C) \int_C (f-m)^{2k} \rmd\gamma' \le \Big(\frac{K(D, \beta_0, \beta_1, \MARGIN, \SIGNALMEDIUMDEPENDENCE) k}{d}\Big)^k\,.
\end{align*}
We now combine Lemma \upref{lem:lem316} and Lemma \upref{lem:logconcaveboundmomentslemma} together as in the proof of \equpref{eq:concentrationoverlapfixedgibbs_initbounds1}, and then take expectations and applying Lemma \upref{lem:lem310}. Using that $f=R_{1,1}$ or $R_{1,2}$, we obtain
\begin{align*}
\E\big[ f^{4k} \big] \le K(D, \beta_0, \beta_1, \MARGIN, \SIGNALSHORTDEPENDENCE)^k\,.
\end{align*}
Next we upper bound $m$. By Lemma \upref{lem:lem316} and Lemma \upref{lem:logconcaveboundmomentslemma}, we have $|f| \le \frac{K(D, \beta_0, \beta_1, \MARGIN, \SIGNALSHORTDEPENDENCE) B^{\star}}d$. Hence by definition of $C^{\star}$,
\begin{align*}
|m| = \Big| \int f \rmd \gamma' \Big| \le K(D, \beta_0, \beta_1, \MARGIN, \SIGNALSHORTDEPENDENCE)\,.
\end{align*}
It follows from the above two displays that $\E\big[ (f-m)^{4k} \big] \le K(D, \beta_0, \beta_1, \MARGIN, \SIGNALSHORTDEPENDENCE)^k$. Thus
\begin{align*}
\E\big[ (f-m)^{2k} \big] &= \E\big[ (f-m)^{2k} \one\{C\} \big] + \E\big[ (f-m)^{2k} \one\{C^c\} \big] \\
&\le \Big(\frac{K(D, \beta_0, \beta_1, \MARGIN, \SIGNALMEDIUMDEPENDENCE) k}{d}\Big)^k + \E\big[ (f-m)^{4k} \big] \P(C^c) \\
&\le \Big(\frac{K(D, \beta_0, \beta_1, \MARGIN, \SIGNALMEDIUMDEPENDENCE) k}{d}\Big)^k + K(D, \beta_0, \beta_1, \MARGIN, \SIGNALSHORTDEPENDENCE)^k \cdot K(\SIGNALSHORTDEPENDENCE) \exp(-d) \\
&\le \Big(\frac{K(D, \beta_0, \beta_1, \MARGIN, \SIGNALMEDIUMDEPENDENCE) k}{d}\Big)^k\,,
\end{align*}
where the last inequality uses that $\exp(-d) \le (k/d)^k$, valid for all $k, d > 0$. The result now follows from symmetrization, Lemma \upref{lem:symmetrization}.
\end{proof}
As noted before, combining Proposition \upref{prop:concentrationoverlapfixedgibbs} and Proposition \upref{prop:concentrationmeangibbswrtdisorder} proves the first part of Proposition \upref{prop:logconcaveconcentration}.

\subsubsection{Second part of Proposition \upref{prop:logconcaveconcentration}}
Recall that we have \equpref{eq:concentrationoverlapfixedgibbs_initbounds2} by Lemma \upref{lem:lem316} and Lemma \upref{lem:logconcaveboundmomentslemma}. 
Now by H\"{o}lder's Inequality and \equpref{eq:boundmgfBstar}, taking $K(D, \GAUSSIANMEASURE_0, \GAUSSIANMEASURE_1, \MARGIN, \SIGNALSHORTDEPENDENCE)$ appropriately and combining with \equpref{eq:concentrationoverlapfixedgibbs_initbounds2}, we obtain
\begin{align*}
\log \nu_{t,v}\Big( \exp \Big( \frac{\| \bar{\btheta} \|^2}{K(D, \GAUSSIANMEASURE_0, \GAUSSIANMEASURE_1, \MARGIN, \SIGNALSHORTDEPENDENCE} \Big) \Big) &= \log \E\Big\langle \exp \Big( \frac{\| \bar{\btheta} \|^2}{K(D, \GAUSSIANMEASURE_0, \GAUSSIANMEASURE_1, \MARGIN, \SIGNALSHORTDEPENDENCE)} \Big) \Big\rangle \\
&\le \log \E\Big\langle \exp \Big( \frac{ \beta \| \bar{\btheta} \|^2}{2} \Big) \Big\rangle^{\frac1{1+K_{\upref{lem:logconcaveboundmomentslemma}}(K_{\upref{lem:lem316}}, D, \beta_0, \beta_1, \MARGIN, \lambda)}} \\
&\le \log \E \exp(B^{\star}) \\
&\le K(\SIGNALSHORTDEPENDENCE)\, d\,.    
\end{align*}
It remains to prove that 
\begin{align}\label{eq:second_part_concentration_step1}
\log \nu_{t,v}\Big( \exp \Big( \frac{\theta_{\DIMENSION}^2}{K(D, \GAUSSIANMEASURE_0, \GAUSSIANMEASURE_1, \MARGIN, \SIGNALSHORTDEPENDENCE)} \Big) \Big) \le K(D, \GAUSSIANMEASURE_0, \GAUSSIANMEASURE_1, \MARGIN, \SIGNALSHORTDEPENDENCE)\,.    
\end{align}
To this end, we cite the following Lemma from \cite{talagrand2010mean}. This function is stated for concave $T \le 0$ and $\beta'>0$ in \cite{talagrand2010mean}, however the proof therein applies verbatim under the conditions stated below.
\begin{lemma}[Lemma 3.2.5, \cite{talagrand2010mean}]\label{lem:lastcoordinateconvexityprinciple}
Consider a concave function $T(\btheta)$ defined on $\R^d$ such that the following integrals are well-defined. Consider any $(a_j)_{j \le d}$, $\beta > 0, \beta' \ge 0$ and any convex set $C \subseteq \R^d$. Define the measure $\mu$ on $\R^d$ by defining for all $B \subseteq \R^d$,
\begin{align*}
\mu(B) \propto \int_{B \cap C} \exp\Big( T(\btheta) - \beta \| \btheta \|^2 - \beta' \btheta_d^2 + \sum_{j \le d} a_j \theta_j \Big)\, \rmd\btheta\,.
\end{align*}
Let 
\begin{align*}
C' := \big\{ \theta_d \in \R\,:\,\exists \btheta' \in \R^{d-1}\text{ s.t. }(\btheta', \theta_d) \in C\}\,.
\end{align*}
Then the following function $f$ on $C'$ is concave, where $f$ is defined by
\begin{align*}
f(\theta_d) := \log \int_{\theta = (\btheta', \theta_d) \in C} \exp\Big( T(\btheta) - \beta \sum_{j \le d-1} \theta_j^2 + \sum_{j \le d-1} a_j \theta_j \Big)\, \rmd\btheta'\,.
\end{align*}
Moreover, letting
\begin{align*}
w(\theta) := f(\theta) - (\beta + \beta') \theta^2 + a_d \theta\,,
\end{align*}
the law of $\theta_d$ under $\mu$ is the measure on $C'$ with density proportional to $\exp w(\theta)$.
\end{lemma}
We will also need the following corollary of the Pr\'ekopa-Leindler Inequality:
\begin{lemma}\label{lem:logconcavemarginallemma}
Consider a log-concave distribution $\btheta$ and and a random variable $\ANNEALEDDISORDER$ with a log-concave density. Then for any function $T(\btheta, \ANNEALEDDISORDER)$ that is jointly concave in $\btheta$ and $\ANNEALEDDISORDER$, the function $\btheta \rightarrow \log \E_{\ANNEALEDDISORDER}\big[ \exp T(\btheta, \ANNEALEDDISORDER) \big]$ is concave in $\btheta$.
\end{lemma}
\begin{proof}
Letting $\varphi(\ANNEALEDDISORDER)$ denote the density of $\ANNEALEDDISORDER$, it follows that $S(\btheta, \ANNEALEDDISORDER) = \exp T(\btheta, \ANNEALEDDISORDER) \varphi(\ANNEALEDDISORDER)$ is log-concave. Thus it remains to prove that $\btheta \rightarrow \int S(\btheta, \ANNEALEDDISORDER) \rmd \ANNEALEDDISORDER$ is log-concave. Considering any $\lambda \in [0,1]$ and any $\btheta_1, \btheta_2$, we have by log-concavity of $S(\btheta, \ANNEALEDDISORDER)$ that for any $\ANNEALEDDISORDER$,
\begin{align*}
S(\lambda \btheta_1 + (1-\lambda) \btheta_2, \ANNEALEDDISORDER) \ge S(\btheta_1, \ANNEALEDDISORDER)^\lambda S(\btheta_2, \ANNEALEDDISORDER)^{1-\lambda}\,.
\end{align*}
It follows from the Pr\'ekopa-Leindler Inequality  (see \cite{prekopa1973logarithmic}) that 
\begin{align*}
\int S(\lambda \btheta_1 + (1-\lambda) \btheta_2, \ANNEALEDDISORDER) \rmd\ANNEALEDDISORDER \ge \Big(\int S(\btheta_1, \ANNEALEDDISORDER) \rmd\ANNEALEDDISORDER\Big)^{\lambda} \Big(\int S(\btheta_2, \ANNEALEDDISORDER) \rmd\ANNEALEDDISORDER\Big)^{1-\lambda}\,.
\end{align*}
This proves that $\int S(\btheta, \ANNEALEDDISORDER) \rmd \ANNEALEDDISORDER$ is log-concave, hence the conclusion.
\end{proof}
We now complete the proof of the second part of Proposition \upref{prop:logconcaveconcentration}. By Theorem \upref{thm:logconcaveconcentrationbasics} applied with the 1-Lipschitz function $\btheta \rightarrow \theta_d$, we have
\begin{align}\label{eq:second_part_concentration_step2}
\Big\langle \exp \Big( \frac{(\theta_d - \langle \theta_d \rangle)^2}{K(\beta_0)} \Big) \Big\rangle \le K\,.
\end{align}
Next, note by convexity of the function $\exp x^2$ that for $f_1, f_2$ satisfying $\nu_{t,v}\big( \exp(f_1^2 / K) \big) \le K$, $\nu_{t,v}\big( \exp(f_2^2 / K) \big) \le K$, we have the upper bound $\nu_{t,v}\big( \exp\big( (f_1+f_2)^2 / K \big) \big) \le K$ (here $K$ can differ). Thus taking expectations of \equpref{eq:second_part_concentration_step2} and applying this observation, to establish \equpref{eq:second_part_concentration_step1}, it thus remains to show that
\begin{align}
\E\Big[ \exp \Big( \frac{\langle \theta_d \rangle^2}{K(D, \beta_0, \beta_1, \MARGIN, \SIGNALSHORTDEPENDENCE)} \Big) \Big] \le K(D, \beta_0, \beta_1, \MARGIN, \SIGNALSHORTDEPENDENCE)\,.\label{eq:logconcavelastcoordgoal}
\end{align} 
To this end, we apply Lemma \upref{lem:lastcoordinateconvexityprinciple}. Define the convex set $C = \R^d \cap \one\{\theta_1/\sqrt{d}\in I\}$, and let
\begin{align*}
T_0(\btheta, \ANNEALEDDISORDER) := \sum_{1 \le i \le n} u(S_{i,t}) + u(S_{n,t,v})\,.    
\end{align*}
Since $S_{i,t}$ is affine in $\btheta$ and as $S_{n,t,v}$ is jointly affine in $\btheta$ and $\ANNEALEDDISORDER$, since $u$ is concave, it follows that $T_0(\btheta, \ANNEALEDDISORDER)$ is concave jointly in $\btheta, \ANNEALEDDISORDER$. 
Thus letting $T(\btheta) := \log \E_{\ANNEALEDDISORDER}\big[ \exp T_0(\btheta) \big]$, $T$ is concave by Lemma \upref{lem:logconcavemarginallemma}. 

Next, let $\beta' = \frac{(1-t)(r-\bar{r})}2 \ge 0$, where the inequality follows as $r \ge \bar{r}$. Define $a_j$ for $1 \le j \le d$ by letting $a_j = 0$ for all $j \le d-1$ and $a_d = \sqrt{(1-t) r}Y$.
Defining $\mu$ from Lemma \upref{lem:lastcoordinateconvexityprinciple} with these choices, as dependence on $\ANNEALEDDISORDER$ in $\nu_{t,v}$ is only via $S_{n,t,v}(\btheta)$, as $T(\btheta)$ is defined with expectation w.r.t. $W$, and as the number of replicas $L=1$ here, it follows that $\mu$ is exactly the Gibbs average $\langle \cdot \rangle$.

Note $C'=\big\{ \theta \in \R\,:\,\exists \btheta' \in \R^{d-1}\text{ s.t. }(\btheta', \theta) \in C\} = \R$. Hence by Lemma \upref{lem:lastcoordinateconvexityprinciple}, the following function is concave on all of $\R$:
\begin{align*}
f(\theta_d) &:= \log \int_{\theta = (\btheta', \theta_d) \in \R^d \cap \one\{\theta_1/\sqrt{d}\in I\}} \E_{\ANNEALEDDISORDER}\Big[ \exp\Big( T_0(\btheta, \ANNEALEDDISORDER) - \beta \sum_{j \le d-1} \theta_j^2 \Big) \Big]\, \rmd\btheta'\,.
\end{align*}
Moreover, by Lemma \upref{lem:lastcoordinateconvexityprinciple}, we know that letting 
\begin{align*}
w(\theta) := f(\theta) - \Big(\beta + \frac{(1-t)(r-\bar{r})}2\Big) \theta^2 + \sqrt{(1-t)r}Y \theta\,,
\end{align*}
the law of $\theta_d$ under $\mu$ is the measure on $\R$ with density proportional to $\exp w(\theta)$. 

As $r \ge \bar{r}$, we know $\beta + \frac{(1-t)(r-\bar{r})}2 > 0$, and thus $w(\theta)$ is $\beta$ strongly concave and attains a unique maximum $\theta^{\star}$ on $\R$. By Lemma \upref{lem:realRSeqsolsboundedtechnicalhelper}, it follows that $\big\langle (\theta_d - \theta^{\star})^2 \big\rangle \le K(\beta_0)$. By Jensen and Triangle Inequality, we thus have $\langle \theta_d \rangle \le K(\beta_0) + |\theta^{\star}|$ (for a different $K$). 

Note that $w'(\theta^{\star})=0$. Since $w$ is $\beta$ strongly concave, $|w'(0)| = |w'(\theta^{\star}) - w'(0)| \ge 2\beta |\theta^{\star}|$. Note $w'(0) = f'(0) + \sqrt{(1-t)r} Y$. Combining the above steps yields
\begin{align}\label{eq:concentration_logconcave_end1}
\big\langle \theta_d \big\rangle \le K(\beta_0) \Big(1+|f'(0)|+\sqrt{(1-t)r}|Y|\Big)\,.
\end{align}
Note since $Y$ is a standard Gaussian, we have $\exp (Y^2) \le K$.
By \equpref{eq:concentration_logconcave_end1}, since $r \le K(D)$, it suffices to prove the following to establish \equpref{eq:logconcavelastcoordgoal} and hence complete the proof:
\begin{align}\label{eq:concentration_logconcave_end2}
\E \Big[ \exp \Big( \frac{f'(0)^2}{K(D, \beta_0, \beta_1, \MARGIN, \SIGNALSHORTDEPENDENCE)} \Big) \Big] \le K(D, \beta_0, \beta_1, \MARGIN, \SIGNALSHORTDEPENDENCE)\,.
\end{align}
Let $\langle \cdot \rangle_f$ denote the Gibbs average defined as follows: for any test function $h$,
\begin{align*}
\langle h \rangle_f &= \frac1{\bar Z} \int_{\theta=(\btheta', \theta_d) \in \R^d \cap \one\{\theta_1/\sqrt{d}\in I\}} \E_{\ANNEALEDDISORDER}\Big[ h \exp\Big( T_0(\theta, W) - \beta \sum_{j \le d-1} \theta_j^2 + \sqrt{(1-t)r}Y \theta_d \Big) \Big]\, \rmd\btheta\,,
\end{align*}
where $\bar Z$ denotes the corresponding normalizing constant.
Explicitly calculating yields
\begin{align}
f'(\theta_d) &= \Bigg\langle \sum_{1 \le i \le n-1} u'(S_{i,t}) \sqrt{\frac{t}d} g_{i,d} + u'(S_{n,t,v}) \sqrt{\frac{tv}d } g_{n,d} \Bigg\rangle_f\,.\label{eq:concentration_logconcave_writeasmixture}
\end{align}
Notice that when $\theta_d=0$, $S_{i,t}$, $S_{n,t,v}$, $\sqrt{(1-t)r}Y \theta_d=0$, and therefore $T_0(\theta, W)$ do not depend on $(g_{i,d})_{1 \le i \le n}, Y$. Hence when $\theta_d=0$, $\bar Z$ also does not depend on $(g_{i,d})_{1 \le i \le n}, Y$.
\equpref{eq:concentration_logconcave_writeasmixture} thus implies that $f'(0)$ can be written as a linear combination of $(g_{i,d})_{1 \le i \le n}$ with coefficients independent of $(g_{i,d})_{1 \le i \le n}$. Specifically, these coefficients are Gibbs averages of $u'(S_{i,t}) \sqrt{\frac{t}{d}}$, $u'(S_{n,t,v}) \sqrt{\frac{tv}{d}}$ (where we consider $S_{i,t}, S_{n,t,v}$ with $\theta_d=0$) w.r.t. $\langle \cdot \rangle_f$.

Since the $(g_{i,d})_{1 \le i \le n} \sim N(0,1)$, it follows that $f'(0)$ is a Gaussian. Let $\bar{\E}$ denote expectation in the $g_{i,d}$. As $|u'| \le D$, the corresponding coefficients on the $g_{i,d}$ are most $D / \sqrt{d}$ in magnitude. By the independence of the $g_{i,d}$, it follows that 
\begin{align*}
\bar\E\big[ f'(0)^2 \big] \le \frac{n}d \cdot D^2 \le K(D)\,,
\end{align*}
where we use that $n/d \le \DENSITYBOUND \le 2$. Since $f'(0)$ is Gaussian, it follows from standard upper bounds for the MGF of a Gaussian with a given variance (see (A.11) in \cite{talagrand2010mean}) that for suitable $K(D, \beta_0, \beta_1, \MARGIN, \SIGNALSHORTDEPENDENCE)$,
\begin{align*}
\E \Big[ \exp \Big( \frac{f'(0)^2}{K(D, \beta_0, \beta_1, \MARGIN, \SIGNALSHORTDEPENDENCE)} \Big) \Big] = \E \bar\E \Big[ \exp \Big( \frac{f'(0)^2}{K(D, \beta_0, \beta_1, \MARGIN, \SIGNALSHORTDEPENDENCE)} \Big)\Big] \le 2\,.
\end{align*}
This establishes \equpref{eq:concentration_logconcave_end2} and hence establishes Proposition \upref{prop:logconcaveconcentration}.

\section{Properties of the RS equations}\label{sec:RSequationsunique}
Here, we establish several important properties of the RS equations. 
Specifically, we prove Proposition \upref{prop:realRSeqsolsbounded} in Appendix \upref{subsec:real_rs_bounded}, Proposition \upref{prop:nicersuniquesol} in Appendix \upref{subsec:nice_rs_uniquesol_pf}, and Lemma \upref{lem:rs_formula_properties} in Appendix \upref{subsec:spherical_conversion_RS_pfs}.
Throughout this section, recall that $\pdfnormal$ denotes the PDF of a standard univariate normal $N(0,1)$ and $\pdfsignal$ denotes the PDF of $S \sim \cD$.

All three of these results are critical, and are proved through similar ideas. 
Discussing the proof of Proposition \upref{prop:nicersuniquesol} for the interpolators for concreteness, we must establish that the RS equations exhibit a particular convex-concave structure.
Talagrand's proof of this convex-concave structure in Section 3.3, \cite{talagrand2010mean} proceeds by showing each part of the expression for the relevant second derivatives is of the correct sign. 
This relies crucially on the fact that $\MARGIN \ge 0$ and that there is no signal in the data.

Here we instead proceed by showing the part of the expression for the relevant second derivatives whose signs are not correct are small in magnitude, and therefore do not eliminate the desired convex-concave structure.
We must also bound similar terms in proving Proposition \upref{prop:realRSeqsolsbounded}.
The magnitude of these terms all can be upper bounded by $\alpha \E\Big[f\big(\frac{\MARGIN - xS - \sqrt{q} Z}{\sqrt{\rho-q}}\big) \Big]$ for a function $f(x) \le O(|x|^2+1)$ that depends on the Inverse Mills' Ratio from \equpref{eq:mills}, where expectation is over $Z, S$.

The magnitude of the expression for each such part can be directly shown to be of order $O\Big(\alpha \cdot \frac{\MARGIN^2+q+1}{\rho-q}\Big)$ for all $\MARGIN \in \R$.
However, to establish our result for all $\alpha \le \alpha_0$ for all $\MARGIN < 0$, such a bound is not sufficient.
The idea is instead as follows. 
When $\frac{\MARGIN - xS - \sqrt{q} Z}{\sqrt{\rho-q}} \le 0$, the relevant quantities are simple to control by properties of the Inverse Mills' Ratio.
Else, we let $\MARGIN' := \MARGIN - xS$ be the `effective margin' and analyze the sign of $\MARGIN'$. 
\begin{itemize}
    \item If $|S| < |\MARGIN|$, as $|x| \le 1$, we must have $\MARGIN' < 0$. In this case, observe that the probability over $Z$ that $\frac{\MARGIN - xS - \sqrt{q} Z}{\sqrt{\rho-q}} = \frac{\MARGIN' - \sqrt{q} Z}{\sqrt{\rho-q}} > 0$ is exponentially small in $|\MARGIN'|^2/q$.
    
    \item Else, the probability over $S$ that $|S| \ge |\MARGIN|$ is exponentially small in $|\MARGIN|^2$ as the law of $S^2$ is sub-Exponential.
\end{itemize}
Combining both these bounds enables us to upper bound $\E\Big[f\big(\frac{\MARGIN - xS - \sqrt{q} Z}{\sqrt{\rho-q}}\big) \Big] $ independently of $\MARGIN$ for $\MARGIN < 0$, which is sufficient for our purposes. 

\subsection{Proof of Proposition \upref{prop:realRSeqsolsbounded}}\label{subsec:real_rs_bounded}
We consider an arbitrary solution $(q, \rho, r, \bar{r})$ of \equpref{eq:realRSequations} and show that $(q, \rho) \in [0, \Csol) \times [0, \Csol)$ and  $\frac{1}{\rho-q} < \Csol$.
This is sufficient to establish Proposition \upref{prop:realRSeqsolsbounded}, up to showing $q > 0$. To this end, note $q=0$ implies $r=0$. Now as $\lambda>0$, we have $u' > 0$ on a set of positive Lebesgue measure. It follows that for every $Z, S$, there exists a set of $W$ with positive Lebesgue measure such that $u'(\APPROXCONSTRAINT) \exp u(\APPROXCONSTRAINT) > 0$, where $\APPROXCONSTRAINT$ is given in terms of $W, Z, S$ by \equpref{eq:approxconstraint}. Since $u' \exp u \ge 0$ pointwise, it follows that $r > 0$ by definition of \equpref{eq:realRSequations}. Hence $q=0$ yields contradiction, so as $q \ge 0$, we must have $q>0$. 

\paragraph{Proof for posterior.} The proof of Proposition \upref{prop:realRSeqsolsbounded} for the posterior is direct. Recall that $|u^l| \le D=D(\lambda)$ for $l=1,2,3,4$. Thus by definition of $\psi_{\alpha}$, $\bar{\psi}_{\alpha}$ in \equpref{eq:interpolationrrbardef1}, \equpref{eq:interpolationrrbardef2}, we have $|r|, |\bar r| \le K(D)$. Also note $r \ge 0$ by definition of $\psi_{\alpha}$. Since $r \ge \bar r$, the definition of the system \equpref{eq:realRSequations} implies $0 \le q, \rho \le K(D)$. 
Next, we note $\rho, \rho-q \ge \frac1{2\beta+r - \bar r} \ge \frac1{K(D)}$.
This proves Proposition \upref{prop:realRSeqsolsbounded} for the posterior.

\paragraph{Proof for interpolators.}
The rest of Appendix \upref{subsec:real_rs_bounded} is now devoted to the proof of Proposition \upref{prop:realRSeqsolsbounded} for the interpolators, which poses significantly more difficulties.
We first state and prove the following Lemmas which we will need in the following proof.
In the following, we recall the definition of $\APPROXCONSTRAINT = \APPROXCONSTRAINT(q, \rho)$ from \equpref{eq:approxconstraint}.
\begin{lemma}\label{lem:calculaterealrsequatios}
Letting $s \sim \cD$, define
\[ v(y) = \log \E_{\ANNEALEDDISORDER}\big[ \exp u\big(x\SIGNALDISORDER + y +\sqrt{\rho-q}\ANNEALEDDISORDER\big)\big]\, .\]
Then we have $v' \ge 0$ and 
\[ r=\alpha \E\big[ v'\big(\sqrt{q} \QUENCHEDDISORDER\big)^2\big]\, ,\, \bar{r}-r=\alpha\E\big[ v''\big(\sqrt{q} \QUENCHEDDISORDER\big) \big]\,.\]
\end{lemma}
\begin{proof}
The fact that $v' \ge 0$ follows directly as $u' \ge 0$. The equalities written above follow from direct calculation, using that $r = \psi_{\alpha}(q, \rho), \bar{r} = \bar{\psi}_{\alpha}(q, \rho)$.
\end{proof}

\begin{lemma}\label{lem:RSboundedupperboundr}
Letting $Y := \frac{\MARGIN - x\SIGNALDISORDER - \sqrt{q}\QUENCHEDDISORDER}{\sqrt{\rho-q}}$ and $K_{\upref{lem:RSboundedupperboundr}} = 20$, we have 
\[ \Bigg( \frac{\E_{\ANNEALEDDISORDER}\big[ \ANNEALEDDISORDER \exp u(\APPROXCONSTRAINT) \big]}{\E_{\ANNEALEDDISORDER}\big[ \exp u(\APPROXCONSTRAINT) \big]} \Bigg)^2 \le K_{\upref{lem:RSboundedupperboundr}} + Y^2 \one\{Y \ge 0\}\,.\]
\end{lemma}
\begin{proof}
If $Y \le 1$, then as $u \le 0$ and $u(\APPROXCONSTRAINT)=0$ for $\APPROXCONSTRAINT \ge \MARGIN$, we have 
\begin{align*}
\E_{\ANNEALEDDISORDER}\big[ \ANNEALEDDISORDER \exp u(\APPROXCONSTRAINT) \big] &\le \E_{\ANNEALEDDISORDER}\big[ |\ANNEALEDDISORDER| \big] = \sqrt{\frac{2}{\pi}}\,, \\
\E_{\ANNEALEDDISORDER}\big[ \exp u(\APPROXCONSTRAINT) \big] &\ge \P_{\ANNEALEDDISORDER}(\APPROXCONSTRAINT \ge \MARGIN) = \P_{\ANNEALEDDISORDER}(\ANNEALEDDISORDER \ge Y) \ge \P_{\ANNEALEDDISORDER}(\ANNEALEDDISORDER \ge 1)\,.
\end{align*} 
As $K_{\upref{lem:RSboundedupperboundr}} \ge 20$, we obtain an upper bound of $K_{\upref{lem:RSboundedupperboundr}}$ in this case.

Else suppose $Y \ge 1$. First as $u' \ge 0$, observe that 
\begin{align*}
\E_{\ANNEALEDDISORDER}\big[ \ANNEALEDDISORDER \exp u(\APPROXCONSTRAINT) \big] = \E_{\ANNEALEDDISORDER}\big[ \exp u(\APPROXCONSTRAINT) \cdot u'(\APPROXCONSTRAINT) \big] \cdot \sqrt{\rho-q} \ge 0\,.    
\end{align*}
Thus we may upper bound
\begin{align*}
\E_{\ANNEALEDDISORDER}\big[ \ANNEALEDDISORDER \exp u(\APPROXCONSTRAINT) \big] &= \E_{\ANNEALEDDISORDER}\big[ \ANNEALEDDISORDER \one\{ \ANNEALEDDISORDER \le Y \} \exp u(\APPROXCONSTRAINT) \big] + \E_{\ANNEALEDDISORDER}\big[ \ANNEALEDDISORDER \one\{ \ANNEALEDDISORDER \ge Y \} \exp u(\APPROXCONSTRAINT) \big]  \\
&\le Y \E_{\ANNEALEDDISORDER}\big[ \one\{ \ANNEALEDDISORDER \le Y \} \exp u(\APPROXCONSTRAINT) \big] + \E_{\ANNEALEDDISORDER}\big[ \ANNEALEDDISORDER \one\{ \ANNEALEDDISORDER \ge Y \} \big] \\
&\le Y \E_{\ANNEALEDDISORDER}\big[ \one\{ \ANNEALEDDISORDER \le Y \} \exp u(\APPROXCONSTRAINT) \big] + \frac1{\sqrt{2\pi}} \exp(-Y^2/2)\,.
\end{align*}
We next lower bound, using that $u(\APPROXCONSTRAINT) = 0$ for $\ANNEALEDDISORDER \ge Y$,
\begin{align*}
\E_{\ANNEALEDDISORDER}\big[ \exp u(\APPROXCONSTRAINT) \big] &= \E_{\ANNEALEDDISORDER}\big[ \one\{ \ANNEALEDDISORDER \le Y \} \exp u(\APPROXCONSTRAINT) \big] + \E_{\ANNEALEDDISORDER}\big[ \one\{ \ANNEALEDDISORDER \ge Y \} \exp u(\APPROXCONSTRAINT) \big] \\
&= \E_{\ANNEALEDDISORDER}\big[ \one\{ \ANNEALEDDISORDER \le Y \} \exp u(\APPROXCONSTRAINT) \big] + \P(\ANNEALEDDISORDER \ge Y) \\
&\ge \E_{\ANNEALEDDISORDER}\big[ \one\{ \ANNEALEDDISORDER \le Y \} \exp u(\APPROXCONSTRAINT) \big] + \frac1{\sqrt{2\pi}} \cdot \frac{Y}{Y^2+1} \exp(-Y^2/2)\,.
\end{align*}
Here the lower bound on $\P(\ANNEALEDDISORDER \ge Y)$ is standard, see e.g. \cite{vershynin2018high}.
We thus obtain 
\begin{align*}
0 \le \frac{\E_{\ANNEALEDDISORDER}\big[ \ANNEALEDDISORDER \exp u(\APPROXCONSTRAINT) \big]}{\E_{\ANNEALEDDISORDER}\big[ \exp u(\APPROXCONSTRAINT) \big]} &\le \frac{Y \E_{\ANNEALEDDISORDER}\big[ \one\{ \ANNEALEDDISORDER \le Y \} \exp u(\APPROXCONSTRAINT) \big] + \frac1{\sqrt{2\pi}} \exp(-Y^2/2)}{\E_{\ANNEALEDDISORDER}\big[ \one\{ \ANNEALEDDISORDER \le Y \} \exp u(\APPROXCONSTRAINT) \big] + \frac1{\sqrt{2\pi}} \cdot \frac{Y}{Y^2+1} \exp(-Y^2/2)} \\
&\le Y+\frac1{Y}\,,
\end{align*}
where we use the following inequality stated on p. 238 of \cite{talagrand2010mean} that for $a, b > 0$, $
\frac{aY+b}{a+\frac{Y}{1+Y^2} b} \le Y+\frac1{Y}$.
Hence we have, as $Y \ge 1$ in this case,
\begin{align*}
\Big( \frac{\E_{\ANNEALEDDISORDER}\big[ \ANNEALEDDISORDER \exp u(\APPROXCONSTRAINT) \big]}{\E_{\ANNEALEDDISORDER}\big[ \exp u(\APPROXCONSTRAINT) \big]} \Big)^2 \le Y^2 + 2 + \frac1{Y^2} \le Y^2+3\,. 
\end{align*}
This completes the proof of the Lemma.
\end{proof}

\begin{lemma}\label{lem:realRSeqsolsboundedtechnicalcheck}
For any $q \le \rho$, we have $\psi_{\alpha}(q, \rho) \ge \bar{\psi}_{\alpha}(q, \rho)$. In particular, we have $r \ge \bar{r}$.
\end{lemma}
\begin{proof}
When $q=\rho$ the desired inequality follows as $u'' \le 0$. Else when $q<\rho$, it suffices to show
\begin{align*}
&\frac{\E_{\ANNEALEDDISORDER}\big[ (\ANNEALEDDISORDER^2 - 1) \exp u(\APPROXCONSTRAINT)\big] }{\E_{\ANNEALEDDISORDER}\big[ \exp u(\APPROXCONSTRAINT)\big]} \le \Big( \frac{\E_{\ANNEALEDDISORDER}\big[ \ANNEALEDDISORDER \exp u(\APPROXCONSTRAINT) \big]}{\E_{\ANNEALEDDISORDER}\big[ \exp u(\APPROXCONSTRAINT) \big]} \Big)^2 \\
&\iff \E_{\ANNEALEDDISORDER}\big[ \ANNEALEDDISORDER^2 \exp u(\APPROXCONSTRAINT) \big] - \frac{\E_{\ANNEALEDDISORDER}\big[ \ANNEALEDDISORDER \exp u(\APPROXCONSTRAINT) \big]^2}{\E_{\ANNEALEDDISORDER}\big[ \exp u(\APPROXCONSTRAINT) \big]} \le \E_{\ANNEALEDDISORDER}\big[ \exp u(\APPROXCONSTRAINT) \big]\,.
\end{align*}
Letting $w(y) = u(x\SIGNALDISORDER + \sqrt{q}\QUENCHEDDISORDER + \sqrt{\rho-q} \cdot y) - \frac{y^2}2$ which is $2$ strongly-concave as $u'' \le 0$, Lemma \upref{lem:realRSeqsolsboundedtechnicalhelper} yields for some real $y'$,
\begin{align*}
\E_{\ANNEALEDDISORDER}\big[ \exp u(\APPROXCONSTRAINT) \big] &\ge \E_{\ANNEALEDDISORDER}\big[ (\ANNEALEDDISORDER-y')^2 \exp u(\APPROXCONSTRAINT) \big] \ge \E_{\ANNEALEDDISORDER}\big[ \ANNEALEDDISORDER^2 \exp u(\APPROXCONSTRAINT) \big] - \frac{\E_{\ANNEALEDDISORDER}\big[ \ANNEALEDDISORDER \exp u(\APPROXCONSTRAINT) \big]^2}{\E_{\ANNEALEDDISORDER}\big[ \exp u(\APPROXCONSTRAINT) \big]}\,,
\end{align*}
completing the proof of the Lemma.
\end{proof}

Now we return to the proof of Proposition \upref{prop:realRSeqsolsbounded}. 
\paragraph{Bounding $q, \rho$:} First, clearly $0 \le q \le \rho$ by definition of the system \equpref{eq:realRSequations}. It remains to upper bound $q, \rho$. Now, note as $r \ge \bar{r}$ by Lemma \upref{lem:realRSeqsolsboundedtechnicalcheck},
\begin{align}\label{eq:boundedsols_q_init}
\rho - q = \frac1{2\GAUSSIANMEASURE + r - \bar{r}} \le \frac1{2\beta_0} \,.
\end{align}
It thus remains to upper bound $q$. If $0 \le q \le 1$ we obtain
\begin{align*}
0 \le q, \rho \le \frac{1}{2\beta_0} + 1 < \Csol\,,
\end{align*}
yielding the required bound.
We now suppose that $q > 1$.
Now by definition of \equpref{eq:realRSequations}, we have 
\begin{align}
r = q (2\beta + r - \bar{r})^2 = \frac{q}{(\rho-q)^2}\,.\label{eq:realRSboundedeq1}
\end{align}
Next by Lemma \upref{lem:RSboundedupperboundr}, letting $Y := \frac{\MARGIN - x\SIGNALDISORDER - \sqrt{q}\QUENCHEDDISORDER}{\sqrt{\rho-q}}$, we have
\begin{align}
r &= \psi_{\alpha}\big(q,\rho\big) = \frac{\alpha}{\rho - q} \E \Big( \frac{\E_{\ANNEALEDDISORDER}\big[ \ANNEALEDDISORDER \exp u(\APPROXCONSTRAINT) \big]}{\E_{\ANNEALEDDISORDER}\big[ \exp u(\APPROXCONSTRAINT) \big]} \Big)^2 \le \frac{\alpha}{\rho - q}\Big( K_{\upref{lem:RSboundedupperboundr}} + \E\big[ Y^2 \one\{Y \ge 0\}\big] \Big)\, .\label{eq:realRSbounded_rbound}
\end{align}
Now we upper bound $\E\big[ Y^2 \one\{Y \ge 0\}\big]$.

\paragraph{Upper bound for $\MARGIN \ge 0$.}
Using $\E[S^2] \le \signalsubexp$, $\E[Z^2]=1$, $x^2 \le 1$, we have
\begin{align}\label{eq:bound_sols_positivemargin}
\E\big[ Y^2 \one\{Y \ge 0\}\big] \le \E[Y^2] &\le \frac{3}{\rho-q} \E\big[ \MARGIN^2 + x^2 S^2 + qZ^2 \big] \le \frac{3(\MARGIN^2 + \signalsubexp + q)}{\rho-q}\,.
\end{align}

\paragraph{Upper bound for $\MARGIN < 0$.}
Here as we would like to establish this result for all $\alpha \le \alpha_0$ at most a universal constant, the proof idea is more complicated.
Define the `effective margin'
\begin{align}\label{eq:effective_margin_def}
\MARGIN' := \MARGIN - x\SIGNALDISORDER\,.
\end{align}
We now analyze the sign of $\MARGIN'$. Note:
\begin{itemize}
    \item If $|\SIGNALDISORDER| < |\MARGIN|$, as $|x| \le 1$, we must have $\MARGIN'=\MARGIN - x\SIGNALDISORDER < 0$. Then for $Y < 0$ to occur, we must have $|Z| \ge \frac{|\MARGIN'|}{\sqrt{q}}$, the probability of which over $Z$ is exponentially unlikely in $|\MARGIN'|/\sqrt{q}$.
    \item Else, the probability over $S$ that $|\SIGNALDISORDER| \ge |\MARGIN|$ is exponentially unlikely in $|\MARGIN|^2$ as the law of $\SIGNALDISORDER^2$ is sub-Exponential. 
\end{itemize}
Thus in both cases we can upper bound $\E\big[ Y^2 \one\{Y \ge 0\}\big]$. We now execute this idea. Write
\begin{align*}
\E\big[ Y^2 \one\{Y \ge 0\}\big] &= \int_{s, z} \frac{(\MARGIN' - z\sqrt{q})^2}{\rho - q} \one\{ \MARGIN' \ge z\sqrt{q} \} \pdfnormal(z) \pdfsignal(s)\, \rmd z \rmd s \\
&= \int_{s:|s| < |\MARGIN|} \int_z \frac{(\MARGIN' - z\sqrt{q})^2}{\rho - q} \one\{\MARGIN' \ge z\sqrt{q} \}  \pdfnormal(z) \pdfsignal(s)\, \rmd z \rmd s \\
&\qquad \qquad + \int_{s:|s| \ge |\MARGIN|} \int_z \frac{(\MARGIN' - z\sqrt{q})^2}{\rho-q} \one\{\MARGIN' \ge z\sqrt{q} \}  \pdfnormal(z) \pdfsignal(s)\, \rmd z \rmd s\, .
\end{align*}
If $|s| < |\MARGIN|$, then $|xs| < |\MARGIN|$, so $\MARGIN' = \MARGIN-xs < 0$ as $\MARGIN < 0$. Recall $q \ge 1$, otherwise the proof is already completed. Thus $\MARGIN' \ge z\sqrt{q}$ implies that $|z| \ge \frac{|\MARGIN'|}{\sqrt{q}}$. Hence 
\begin{align*}
& \int_{s:|s| < |\MARGIN|} \int_z \frac{(\MARGIN' - z\sqrt{q})^2}{\rho - q} \one\{\MARGIN' \ge z\sqrt{q} \}  \pdfnormal(z) \pdfsignal(s)\, \rmd z \rmd s \\
&\quad\le \frac{q}{\rho-q} \int_{s:|s| < |\MARGIN|} \int_z 2\Big( \frac{|\MARGIN'|^2}{q} + z^2 \Big) \cdot \exp(-|\MARGIN'|^2/4q) \cdot \frac1{\sqrt{2\pi}} \exp(-z^2/4) \pdfsignal(s)\, \rmd z \rmd s  \\
&\quad \le \frac{2q}{\rho-q} \Big(\int_s \pdfsignal(s) \, \rmd s \Big) \Bigg( \int_z \Big( \frac{|\MARGIN'|^2}{q} + z^2 \Big) \cdot \exp(-|\MARGIN'|^2/4q) \exp(-z^2/4)\, \rmd z \Bigg) \\
&\quad \le \frac{10q}{\rho-q}\, .
\end{align*}
When $|s| \ge |\MARGIN|$, we employ a similar strategy, but now we analyze the probability in $S$ rather than $Z$. Since $q \ge 1$, we may bound
\begin{align*}
&\int_{s:|s| \ge |\MARGIN|} \int_z \frac{(\MARGIN' - z\sqrt{q})^2}{\rho-q} \one\{\MARGIN' \ge z\sqrt{q} \}  \pdfnormal(z) \pdfsignal(s)\, \rmd z \rmd s \\
&\quad \le \frac{2q}{\rho-q} \int_{s:|s| \ge |\MARGIN|} \int_z \Big( z^2 + \frac{2(\MARGIN^2+s^2 )}{q} \Big) \pdfnormal(z) \pdfsignal(s)\, \rmd z \rmd s \\
&\quad \le \frac{2q}{\rho-q} \Big( \int_z z^2 \pdfnormal(z)\, \rmd z + 2  \int_s s^2 \pdfsignal(s)\, \rmd s + 2\MARGIN^2 \int_{s:|s| \ge |\MARGIN|} \pdfsignal(s)\, \rmd s  \Big) \\
&\quad \le \frac{2q}{\rho-q} \Big(1 + 2 \signalsubexp + 2\MARGIN^2 \P\Big( S^2 > \MARGIN^2 \Big) \Big) \\
&\quad \le \frac{2q}{\rho-q} \Big(1 + 2 \signalsubexp + 4\MARGIN^2 \exp\Big( -\frac{\MARGIN^2}{ \signalsubexp} \Big) \Big) \\
&\quad \le \frac{2q}{\rho-q} \Big( 1 + 3.5  \signalsubexp \Big) \,.
\end{align*}
Here we used that $S^2$ is $\signalsubexp$ sub-Exponential to upper bound $\P( S^2 > \MARGIN^2 )$, and used Assumption \upref{ass:disorder_technical_assumption} to upper bound $\E[S^2] \le \signalsubexp$.
Combining the above displays yields
\begin{align}\label{eq:bound_sols_negativemargin}
\E\big[ Y^2 \one\{Y \ge 0\}\big] \le \frac{q}{\rho-q} \Big( 12 + 7  \signalsubexp  \Big)\,.
\end{align}

Having upper bounded $\E\big[ Y^2 \one\{Y \ge 0\}\big]$ in both cases, in the $\MARGIN>0$ case we obtain via \equpref{eq:realRSbounded_rbound} and \equpref{eq:bound_sols_positivemargin}, using that $q \ge 1$ by assumption here (else we have already finished the proof) and that $\signalsubexp \ge 1$,
\begin{align*}
q = r(\rho-q)^2 &\le (\rho-q)^2 \cdot \frac{\alpha}{\rho-q} \Big( K_{\upref{lem:RSboundedupperboundr}} + \E\big[ Y^2 \one\{Y \ge 0\}\big] \Big) \\
&\le \alpha K_{\upref{lem:RSboundedupperboundr}} (\rho-q) + 3 \alpha q \big(\MARGIN^2 + \signalsubexp + 1 \big)\,.
\end{align*}
Similarly in the $\MARGIN < 0$ case we have by \equpref{eq:realRSbounded_rbound} and \equpref{eq:bound_sols_negativemargin},
\begin{align*}
q = r(\rho-q)^2 &\le \alpha K_{\upref{lem:RSboundedupperboundr}} (\rho-q) + \alpha q \big( 12 + 7  \signalsubexp  \big)\, .
\end{align*}
Since $\DENSITYBOUND \le \min\big\{ \frac{1}{2 K_{\upref{lem:RSboundedupperboundr}}}, \frac1{24 + 14  \signalsubexp  } \big\}$ when $\MARGIN < 0$ and $\DENSITYBOUND \le \frac1{3(\MARGIN^2 + \signalsubexp + 1)}$ when $\MARGIN \ge 0$, combining with \equpref{eq:boundedsols_q_init} yield that in either case, we have
\[ q \le \frac{q}2 + \frac{\rho-q}2 = \frac{\rho}2 \le \frac{q}2 + \frac1{4\beta_0}\,,\quad\text{thus}\quad q \le \frac{1}{2\beta_0}\, .\]
Consequently we obtain $\rho \le \frac{1}{\beta_0}$ from \equpref{eq:boundedsols_q_init}, and thus $0 \le q, \rho < \Csol$.

\paragraph{Upper bounding $\frac1{\rho-q}$:} From the definition of the system \equpref{eq:realRSequations}, we have
\[ \frac1{\rho - q} = 2\GAUSSIANMEASURE + r - \bar{r} = 2\beta - \alpha \E\big[ v''\big(\sqrt{q}\QUENCHEDDISORDER\big) \big]\, .\]
Next by Gaussian Integration by Parts and as $v' \ge 0$, we have 
\begin{align*}
\E\big[ v''\big(\sqrt{q}\QUENCHEDDISORDER\big) \big] &= \frac1{\sqrt{q}} \E\big[ \QUENCHEDDISORDER v'\big(\sqrt{q}\QUENCHEDDISORDER\big) \big] \\
&\ge \frac1{\sqrt{q}} \E\big[ \one\{\QUENCHEDDISORDER \le 0\} \cdot \QUENCHEDDISORDER v'\big(\sqrt{q}\QUENCHEDDISORDER\big) \big] \\
&\ge -\frac1{\sqrt{q}} \E\big[ \QUENCHEDDISORDER^2 \one\{\QUENCHEDDISORDER \le 0\} \big]^{1/2} \cdot \E\big[ v'\big(\sqrt{q}\QUENCHEDDISORDER\big)^2 \big]^{1/2} \\
&= -\frac1{\sqrt{2q}} \E\big[ v'\big(\sqrt{q}\QUENCHEDDISORDER\big)^2 \big]^{1/2}\,.
\end{align*}
Combining the above with Lemma \upref{lem:calculaterealrsequatios} gives
\[ \frac1{\rho - q} \le 2\GAUSSIANMEASURE + \frac{\alpha}{\sqrt{2q}}\E\big[ v'\big(\sqrt{q}\QUENCHEDDISORDER\big)^2 \big]^{1/2} = 2\GAUSSIANMEASURE + \sqrt{\frac{\alpha r}{2q} }\, .\]
Combining the above with \equpref{eq:realRSboundedeq1} and using $\alpha \le \DENSITYBOUND \le 1/2$ yields
\[ \frac1{\rho-q} \le 2\GAUSSIANMEASURE + \sqrt{\frac{\alpha}{2}} \cdot \frac1{\rho-q}\,,\quad\text{thus}\quad\frac1{\rho - q} \le 4\GAUSSIANMEASURE \le 4\beta_1 < \Csol\, .\]

\subsection{Proof of Proposition \upref{prop:nicersuniquesol}}\label{subsec:nice_rs_uniquesol_pf}
Throughout Appendix \upref{subsec:nice_rs_uniquesol_pf}, $I$ is held fixed, and so we denote $F_I$ by $F$ for simplicity of notation. As done in Chapter 3 of \cite{talagrand2010mean}, consider the transformation
\[ y = \frac{q}{\rho-q}\,.\]
Hence given $y$, we set $q=\frac{y \rho}{1+y}$, and have $\rho > q$.
Define
\begin{align}
G(y,\rho) &:= F\Big(\frac{y \rho}{1+y},\rho \Big) \notag\,.
\end{align}
Thus we have
\begin{equation}\label{eq:G_transformed_def}
\begin{aligned}
G(y, \rho) &= \alpha \E \Big[ \log \E_{\ANNEALEDDISORDER} \exp u\Big( x \SIGNALDISORDER + \sqrt{\frac{y \rho}{1+y}} \QUENCHEDDISORDER + \sqrt{\frac{\rho}{1+y}} \ANNEALEDDISORDER \Big)\Big]  + \frac{y}2 + \frac12 \log \rho - \frac12 \log(1+y) \\
&\qquad \qquad - \beta \big(\rho + r_I(x)\big)\,,
\end{aligned}
\end{equation}
where $u(\cdot)$ is given by $\exp u(x) = \one\{x \ge \MARGIN\}$ for the interpolators, and by \equpref{eq:defuSgmm} in the GMM case or \equpref{eq:defuSlogistic} in the logistic case for the posterior. Consider the system 
\begin{align}
\frac{\partial G}{\partial \rho}=\frac{\partial G}{\partial y}=0\,.\label{eq:RStransformed}
\end{align}
Note this system does not depend on the $-\beta r(x)$ term in the definition of $G(y, \rho)$.

\begin{lemma}\label{lem:RSuniquesoltransform}
The system \equpref{eq:niceRSeqs} has a unique solution $(q_0, \rho_0) \in [0, \Csol] \times [\frac1{\Csol}, \Csol]$ with $q_0 < \rho_0$ iff \equpref{eq:RStransformed} has a unique solution $(y_0, \rho_0) \in [0, \infty) \times [\frac1{\Csol}, \Csol]$.
\end{lemma}
\begin{proof}
First, notice the transformation $y=\frac{q}{\rho-q}$ maps any $(q,\rho) \in [0, \Csol] \times [\frac1{\Csol}, \Csol]$ with $q<\rho$ to some $(y,\rho) \in [0, \infty) \times [\frac1{\Csol}, \Csol]$. Also for any $(y,\rho) \in [0, \infty) \times [\frac1{\Csol}, \Csol]$, notice the transformation $q = \frac{y\rho}{1+y}$ is such that $0 \le q < \rho \le \Csol$. Thus  the transformation $q = \frac{y\rho}{1+y}$ maps any $(y,\rho) \in [0, \infty) \times [\frac1{\Csol}, \Csol]$ to some $(q,\rho) \in [0, \Csol] \times [\frac1{\Csol}, \Csol]$ with $q<\rho$.

Thus, it suffices to show that $(q_0, \rho_0) \in [0, \Csol] \times [\frac1{\Csol}, \Csol]$ with $q_0<\rho_0$ is a solution to \equpref{eq:niceRSeqs} iff $(y_0, \rho_0)\in [0, \infty) \times [\frac1{\Csol}, \Csol]$ is a solution to \equpref{eq:RStransformed}, where $y_0 = \frac{q_0}{\rho_0-q_0}$. Establishing that $(q_0, \rho_0)$, $(y_0, \rho_0)$ belong in the appropriate domains has already been established above. To show that this transformation preserves solutions, note
\begin{align*}
\frac{\partial G}{\partial \rho}(y, \rho) &= \frac{\partial F}{\partial \rho}\Big( \frac{y\rho}{1+y}, \rho \Big) + \frac{\partial F}{\partial q}\Big( \frac{y\rho}{1+y}, \rho \Big) \cdot \frac{\partial q}{\partial \rho}(y, \rho) \quad,\quad\frac{\partial G}{\partial y}(y, \rho) = \frac{\partial F}{\partial q}\Big( \frac{y\rho}{1+y}, \rho \Big) \cdot \frac{\partial q}{\partial y}(y, \rho) \,.
\end{align*}
Note $\frac{\partial q}{\partial y}(y_0, \rho_0) = \frac{\rho_0}{(1+y_0)^2} > 0$ as $\rho_0 > \frac1{\Csol} > 0$. Thus if $(y_0, \rho_0)$ solves \equpref{eq:RStransformed}, then we have $\frac{\partial F}{\partial q}(q_0, \rho_0)=0$, and consequently $\frac{\partial F}{\partial \rho}(q_0, \rho_0)=0$. Moreover if $(q_0, \rho_0)$ solves \equpref{eq:niceRSeqs}, plugging into the equations above directly implies that $(y_0, \rho_0)$ solves \equpref{eq:RStransformed}. This proves the Lemma.
\end{proof}
Thus in what follows, we now work with \equpref{eq:RStransformed}.

\paragraph{Preliminary calculations for interpolators.}
We first explicitly write the system \equpref{eq:RStransformed} and establish some of its useful properties.
Since we consider $\rho_0 \ge \frac1{\Csol} > 0$, we may divide by $\rho$ in the following. 
Now let 
\begin{align}
\INVERSEMILLSEXPRESSION := \sqrt{\frac{1+y}{\rho}}(\MARGIN - x\SIGNALDISORDER) - \sqrt{y}\QUENCHEDDISORDER = \sqrt{\frac{1+y}{\rho}} \MARGIN' - \sqrt{y}\QUENCHEDDISORDER \,,\label{eq:RS_Vexpression_def}
\end{align}
where here and in the following, as in the proof of Proposition \upref{prop:realRSeqsolsbounded}, we define the `effective margin' $\MARGIN' := \MARGIN - x \SIGNALDISORDER$ as per \equpref{eq:effective_margin_def}.
Hence letting $\cN(\cdot)$ denote the complementary c.d.f. of the standard normal as per \equpref{eq:complementary_cdf_normal}, we have
\begin{align}\label{eq:convert_to_complementarycdf}
\log \P_{\ANNEALEDDISORDER} \Big( xS + \sqrt{\frac{y \rho}{1+y}}Z + \sqrt{\frac{\rho}{1+y}}W \ge \MARGIN \Big) \equiv \log \cN\Big( \sqrt{\frac{1+y}{\rho}} (\MARGIN - x\SIGNALDISORDER) - \sqrt{y}\QUENCHEDDISORDER \Big)\,.
\end{align}
Recalling the definition of the Inverse Mills' Ratio $\cA(x)$ in \equpref{eq:mills}, we record a few of its properties.
\begin{lemma}[Lemma 3.3.7 of \cite{talagrand2010mean} and Lemma 17, \cite{el2022algorithmic}]\label{lem:millsratioproperties}
We have:
\begin{enumerate}
\item $\cA(x) \ge x$ for all $x \in \R$.
\item $\cA'(x) = \cA(x)^2 - x\cA(x) \ge 0$ for all $x \in \R$.
\item $x \cA(x) \cA'(x) \le \cA(x)^2$ for all $x \in \R$.
\item $x\cA(x) \le 1+x^2$ for all $x \in \R$.
\item $\lim_{x\rightarrow\infty}\frac{\cA(x)}{x} = 1$.
\item $\cA'(x), \cA''(x) \le \MILLSBOUND$ for a universal constant $\MILLSBOUND \ge 1$, $0 \le \cA(x) \le 1$ for all $x<0$, and $0 \le \cA(x) \le 2 x + 1$ for all $x \ge 0$.
\end{enumerate}
\end{lemma}
\begin{proof}
All the above are proved in Lemma 3.3.7 of \cite{talagrand2010mean} and Lemma 17, \cite{el2022algorithmic} except the bound $\cA''(x) \le \MILLSBOUND$. 
This follows from noting as in the proof of Lemma 17, \cite{el2022algorithmic} the bound $\cA''(x) = \cA(x) P(x, \cA(x))$ for a polynomial $P(\cdot,\cdot)$. 
It thus suffices to control $\cA''(x)$ for $x$ sufficiently large. 
The $x \rightarrow +\infty$ case follows from Corollary 1.6 of \cite{pinelis2019exact}, while the $x \rightarrow -\infty$ case follows from noting that $\cA(x) \rightarrow 0$ exponentially fast as $x \rightarrow -\infty$.
\end{proof}

We now explicitly write the system \equpref{eq:RStransformed}. By explicit calculation and \equpref{eq:convert_to_complementarycdf}, we have
\begin{align}
\frac{\partial G}{\partial \rho} &= \frac{\alpha}2 \E\Big[ (\MARGIN - x\SIGNALDISORDER) \frac{(1+y)^{1/2}}{\rho^{3/2}} \cA(\INVERSEMILLSEXPRESSION)\Big] + \frac1{2\rho} - \GAUSSIANMEASURE = \frac{\alpha}2 \E\Big[ \MARGIN' \frac{(1+y)^{1/2}}{\rho^{3/2}} \cA(\INVERSEMILLSEXPRESSION)\Big] + \frac1{2\rho} - \GAUSSIANMEASURE \,,\label{eq:partialGrho} \\
\frac{\partial G}{\partial y} &= \alpha \E\Big[ \Big( -\frac{\MARGIN - x\SIGNALDISORDER}{2\rho^{1/2} (1+y)^{1/2}} + \frac{Z}{2y^{1/2}} \Big) \cA(\INVERSEMILLSEXPRESSION) \Big] + \frac{y}{2(1+y)} \,.
\end{align}
Applying Gaussian Integration by Parts and Lemma \upref{lem:millsratioproperties},
\begin{align*}
\E\Big[ \frac{Z}{\sqrt{y}} \cA(\INVERSEMILLSEXPRESSION) \Big] = -\E\Big[ \cA'(\INVERSEMILLSEXPRESSION) \Big] &= \E\Big[ -\cA(\INVERSEMILLSEXPRESSION)^2 + \INVERSEMILLSEXPRESSION \cA(\INVERSEMILLSEXPRESSION) \Big] \\
&= \E \Big[ -\cA(\INVERSEMILLSEXPRESSION)^2 \Big] + \E\Big[ \Big( (\MARGIN - x\SIGNALDISORDER)\sqrt{\frac{1+y}{\rho}} - \sqrt{y}Z \Big) \cA(\INVERSEMILLSEXPRESSION) \Big]\,.
\end{align*}
Rearranging yields
\begin{align*}
\frac{Z}{\sqrt{y}} \E[\cA(\INVERSEMILLSEXPRESSION)] = \frac1{1+y}\E\Big[-\cA(\INVERSEMILLSEXPRESSION)^2 + (\MARGIN - x\SIGNALDISORDER)\sqrt{\frac{1+y}{\rho}} \cA(\INVERSEMILLSEXPRESSION) \Big]\,. 
\end{align*}
We thus obtain
\begin{align}
\frac{\partial G}{\partial y} = -\frac{\alpha}{2(1+y)} \E\Big[ \cA(\INVERSEMILLSEXPRESSION)^2 \Big] + \frac{y}{2(1+y)} \,.\label{eq:partialGy}
\end{align}
Following the argument of \cite{shcherbina2003rigorous}, \cite{talagrand2010mean}, we establish that \equpref{eq:partialGrho}, \equpref{eq:partialGy} have a unique solution by considering $\frac{\partial^2 G}{\partial^2 \rho}$ and the sign of $\frac{\partial G}{\partial y}$, and arguing that the signs work out to yield a saddle point and thus a unique solution to \equpref{eq:RStransformed}. 
We first compute from \equpref{eq:partialGrho}, 
\begin{align}
\frac{\partial^2 G}{\partial \rho^2} &= -\frac{3\alpha(1+y)^{1/2} } {4 \rho^{5/2}} \E\Big[ \MARGIN' \cA(\INVERSEMILLSEXPRESSION) \Big] + \frac{\alpha (1+y)^{1/2}}{2 \rho^{3/2}} \E\Big[ \MARGIN' \cA'(\INVERSEMILLSEXPRESSION) \cdot -\frac12 (1+y)^{1/2} \rho^{-3/2} \MARGIN' \Big] - \frac1{2\rho^2} \notag \\
&= -\frac{3\alpha (1+y)^{1/2}}{4 \rho^{5/2}} \E\Big[ \MARGIN' \cA(\INVERSEMILLSEXPRESSION) \Big] - \frac{\alpha (1+y)}{4 \rho^3} \E\Big[ \MARGIN'^2 \cA'(\INVERSEMILLSEXPRESSION) \Big] - \frac1{2\rho^2}\,.\label{eq:partialGrhosecondderiv}
\end{align}
Now in preparation on our work on the sign of $\frac{\partial G}{\partial y}$, we define
\begin{align}
g(y, \rho) := y-\alpha \E \big[ \cA(\INVERSEMILLSEXPRESSION)^2 \big]\,,\quad\text{thus}\quad\frac{\partial G}{\partial y}(y, \rho) = \frac{g(y, \rho)}{2(1+y)}\,.\label{eq:partialGpartialyhelperfunc}
\end{align}
We will need to show $\frac{\partial g}{\partial y} > 0$ in an appropriate domain, which we do in Lemma \upref{lem:ycontrolled}. In preparation, we compute
\begin{align*}
\frac{\partial g}{\partial y} &= 1-\alpha \E \Big[ 2\cA(\INVERSEMILLSEXPRESSION) \cA'(\INVERSEMILLSEXPRESSION) \Big( \frac{\MARGIN'}{2\sqrt{\rho(1+y)}} - \frac12 \QUENCHEDDISORDER y^{-1/2} \Big)\Big] \\
&= 1 - \frac{\alpha}{\sqrt{\rho(1+y)}} \E\Big[ \MARGIN' \cA(\INVERSEMILLSEXPRESSION) \cA'(\INVERSEMILLSEXPRESSION) \Big] + \frac{\alpha}{\sqrt{y}} \E\Big[ \QUENCHEDDISORDER \cA(\INVERSEMILLSEXPRESSION) \cA'(\INVERSEMILLSEXPRESSION) \Big]\,.
\end{align*}
Notice $\frac{\partial V}{\partial \QUENCHEDDISORDER} = -\sqrt{y}$, so Gaussian Integration by Parts gives
\[\E\Big[ \QUENCHEDDISORDER \cA(\INVERSEMILLSEXPRESSION) \cA'(\INVERSEMILLSEXPRESSION) \Big] = -\sqrt{y} \E\Big[ \cA'^2(\INVERSEMILLSEXPRESSION)+\cA(\INVERSEMILLSEXPRESSION)\cA''(\INVERSEMILLSEXPRESSION) \Big]\,. \]
Thus
\begin{align}
\frac{\partial g}{\partial y} = 1 - \frac{\alpha}{\sqrt{\rho(1+y)}} \E\Big[ \MARGIN' \cA(\INVERSEMILLSEXPRESSION) \cA'(\INVERSEMILLSEXPRESSION) \Big] - \alpha \E\Big[ \cA'^2(\INVERSEMILLSEXPRESSION)+\cA(\INVERSEMILLSEXPRESSION)\cA''(\INVERSEMILLSEXPRESSION) \Big]\,.\label{eq:partialGpartialysecond}
\end{align}
To control the signs of the above, since we do not have non-negative margin -- $\MARGIN \ge 0$ or $\MARGIN' = \MARGIN - x\SIGNALDISORDER \ge 0$ is not uniformly true -- the signs of each component of the derivative do not work out directly in the same way as in \cite{shcherbina2003rigorous}. Rather we show that for $\alpha \le \DENSITYBOUND$, the part of the derivatives that do not have the correct sign are small in magnitude. We will carry this argument out in the proofs to follow.

\paragraph{Preliminary calculations for posterior.}
We now define
\begin{align}\label{eq:RSeqs_barV_def_posterior}
\bar V := xS + \sqrt{\frac{y \rho}{1+y}}Z + \sqrt{\frac{\rho}{1+y}}W\,.
\end{align}
Here in Appendix \upref{subsec:nice_rs_uniquesol_pf}, we let $\langle \cdot \rangle$ denote a Gibbs average w.r.t. $\E_W \exp u(\bar V)$. In what follows, when the derivatives $u^{l}$ are written, they are taken with argument $\bar V$ unless otherwise stated. We obtain by Gaussian Integration by Parts,
\begin{align}
\frac{\partial G}{\partial \rho} &= \frac{\alpha}2 \E \big\langle u'' + u'^2 \big\rangle - \frac{\alpha y}{2(1+y)} \E \big\langle u' \big\rangle^2 + \frac1{2\rho} - \beta\,, \label{eq:posterior_partialGrho}\\
\frac{\partial G}{\partial y} &= -\frac{\alpha \rho}{2(1+y)^2} \E \big\langle \exp u' \big\rangle^2 + \frac{y}{2(1+y)}\,. \label{eq:posterior_partialGy}
\end{align}
We will follow the same strategy as in the interpolators of establishing the saddle point structure of the RS equations. As such we compute
\begin{align}\label{eq:posterior_partialGrho_twice}
\frac{\partial^2 G}{\partial \rho^2} &= \frac{\alpha}4\Big( \expressionI + \expressionII + \expressionIII + \expressionIV \Big) -\frac1{2\rho^2} \,,
\end{align}
where 
\begin{equation}\label{eq:posterior_secondderiv_expressions}
\begin{aligned}
\expressionI &:= \frac1{1+y} \E\Big[ \big\langle u'''' + 4u''' u' + 3u''^2 + 6 u'' u'^2 + u'^4 \big\rangle - \big\langle u'' + u'^2 \big\rangle^2 \Big]\,, \\
\expressionII &:= \frac{y}{1+y} \E \Big[ \big\langle u'''' + 4u''' u' + 3u''^2 + 6 u'' u'^2 + u'^4 \big\rangle + 2 \big\langle u''+ u'^2 \big\rangle \big\langle u' \big\rangle^2 \\
&\qquad \qquad \qquad \qquad \qquad - \big\langle u''' + 3u'' u' + u'^3 \big\rangle \big\langle u' \big\rangle - \big\langle u''' + 3u'' u' + u'^3 \big\rangle \big\langle u'' + u'^2 \big\rangle \Big]\,, \\
\expressionIII &:= \frac{-y}{2(1+y)^2} \E\Big[ \big\langle u''' + 3u'' u' + u'^3 \big\rangle\big\langle u' \big\rangle  - \big\langle u'' + u'^2 \big\rangle \big\langle u' \big\rangle^2 \Big] \,, \\
\expressionIV &:= \frac{-y^2}{2(1+y)^2} \E\Big[ \big\langle u'' + u'^2 \big\rangle^2 + \big\langle u''' + 3u''u' + u'^3 \big\rangle \big\langle u' \big\rangle + 3 \big\langle u' \big\rangle^4 - 5\big\langle u'' + u'^2 \big\rangle \big\langle u' \big\rangle^2 \Big]\,.
\end{aligned}
\end{equation}
We also define similarly to \equpref{eq:partialGpartialyhelperfunc},
\begin{align}\label{eq:posterior_RS_g_def}
g(y) := y(1+y) - \alpha \rho \E \big\langle \exp u' \big\rangle^2\,,\quad\text{thus}\quad \frac{\partial G}{\partial y}(y, \rho) = \frac{g(y)}{2(1+y)^2}\,.
\end{align}
Again in Lemma \upref{lem:ycontrolled}, we show $\frac{\partial g}{\partial y}>0$. In preparation, we compute
\begin{align}\label{eq:posterior_secondderiv_y}
\frac{\partial g}{\partial y} = 1 + 2y - \frac{\alpha \rho}{(1+y)^2} \Big( \expressionI + \expressionII \Big)\,,
\end{align}
where
\begin{equation}\label{eq:posterior_secondderiv_y_expressions}
\begin{aligned}
\expressionI &:= \E\Big[  -\big\langle u''' + 3u'' u' + u'^3 \big\rangle \big\langle u' \big\rangle + \big\langle u'' + u'^2 \big\rangle \big\langle u' \big\rangle^2 \Big] \,, \\
\expressionII &:= \E\Big[ \big\langle u'' + u'^2 \big\rangle^2 + \big\langle u''' + 3u''u' + u'^3 \big\rangle \big\langle u' \big\rangle + 3 \big\langle u' \big\rangle^4 - 5\big\langle u'' + u'^2 \big\rangle \big\langle u' \big\rangle^2 \Big]
\end{aligned}
\end{equation}

\paragraph{Bounding relevant quantities for the interpolators.}
To study the system \equpref{eq:RStransformed} for the interpolators, we will need to bound various quantities involving the Inverse Mills' Ratio. This step is not necessary to study \equpref{eq:RStransformed} for the posterior.
To this end, we first control $\E[\cA(\INVERSEMILLSEXPRESSION)^2]$ in terms of $y$. We then use this to upper bound $y$. Finally we use these two upper bounds to upper bound the quantities of interest, namely $\E\big[\cA(\INVERSEMILLSEXPRESSION)\big], \E\big[\cA(\INVERSEMILLSEXPRESSION)^2\big], \E\big[|\MARGIN - x\SIGNALDISORDER|\cA(\INVERSEMILLSEXPRESSION)\big]$.

\begin{lemma}\label{lem:controlexpectationofAWinit}
Consider $(\alpha,\GAUSSIANMEASURE)\in [0, \DENSITYBOUND] \times [\beta_0, \beta_1]$. Then for any $y \ge 0$ and $\rho \ge \frac1{\Csol}$, 
\begin{align*}
\E\Big[ \cA(\INVERSEMILLSEXPRESSION)^2 \Big] \le K_{\upref{lem:controlexpectationofAWinit}} y + K_{\upref{lem:controlexpectationofAWinit}} \,,
\end{align*}
where 
\begin{align*}
K_{\upref{lem:controlexpectationofAWinit}} = \begin{cases} 12 \Csol  \signalsubexp +  \signalsubexp + 51 &: \MARGIN < 0\,,\\ 24 \Csol (\MARGIN^2 + \signalsubexp) + 12 &: \MARGIN \ge 0\,.\end{cases}
\end{align*}
\end{lemma}
\begin{proof}
We perform a similar argument as in the proof of Proposition \upref{prop:realRSeqsolsbounded}. 
We break into the two cases $\MARGIN \ge 0$ and $\MARGIN < 0$.

\paragraph{When $\MARGIN \ge 0$.}
By Lemma \upref{lem:millsratioproperties}, $0\le \cA(x) \le 2|x|+1$ for all $x \in \R$. 
Consequently using $x^2 \le 1$ and $\rho \ge \frac1{\Csol}$, we have
\begin{align*}
\cA(\INVERSEMILLSEXPRESSION)^2 \le (2|\INVERSEMILLSEXPRESSION|+1)^2 \le 6\INVERSEMILLSEXPRESSION^2 + 3 &\le 12 \MARGIN'^2 \cdot \frac{1+y}{\rho} + 12 y Z^2 + 3 \\
&\le 24 \Csol (\MARGIN^2 + S^2) (y+1) + 12 y Z^2 + 3\,.
\end{align*}
As $\E[S^2] \le \signalsubexp$, $\E[Z^2]=1$, we obtain
\begin{align*}
\E\Big[ \cA(\INVERSEMILLSEXPRESSION)^2 \Big] \le 24 \Csol (\MARGIN^2 + \signalsubexp) (y+1) + 12y + 3 \le K_{\upref{lem:controlexpectationofAWinit}} y + K_{\upref{lem:controlexpectationofAWinit}}\,. 
\end{align*}

\paragraph{When $\MARGIN < 0$.}
Recalling that $\MARGIN' = \MARGIN - x\SIGNALDISORDER$, write
\begin{align*}
\E\Big[ \cA(\INVERSEMILLSEXPRESSION)^2 \Big] &= \int_{s, z} \cA(\INVERSEMILLSEXPRESSION)^2 \pdfnormal(z) \pdfsignal(s)\, \rmd z \rmd s \\
&= \int_{s:|s| < |\MARGIN| } \int_z \cA(\INVERSEMILLSEXPRESSION)^2 \pdfnormal(z) \pdfsignal(s)\, \rmd z \rmd s + \int_{s:|s| \ge |\MARGIN| } \int_z \cA(\INVERSEMILLSEXPRESSION)^2 \pdfnormal(z) \pdfsignal(s)\, \rmd z \rmd s\,.
\end{align*}

We separately upper bound each of these two terms:
\begin{itemize}
\item If $|s| < |\MARGIN|$, then as $|xs| < |\MARGIN|$ and $\MARGIN < 0$, we must have $\MARGIN' = \MARGIN - xs < 0$. Now, consider the set $E := \{ (s, z)\,:\, \INVERSEMILLSEXPRESSION \le 1\}$. Hence for $(s, z) \in E$, $\cA(\INVERSEMILLSEXPRESSION) \le \frac1{\sqrt{2\pi} \cN(1)}$. 

Otherwise consider $(s, z) \not\in E$, so here we have $\INVERSEMILLSEXPRESSION>1$. First note that in this case, $y > 0$; if $y=0$, then as $\MARGIN' < 0$, we have $1 < V = \MARGIN' \sqrt{\frac{1+y}{\rho}} - z\sqrt{y} < 0$ which is a contradiction. Thus recalling $\rho>0$,
\[ 1 < \INVERSEMILLSEXPRESSION = \MARGIN' \sqrt{\frac{1+y}{\rho}} - z\sqrt{y} \implies z \le \MARGIN'\sqrt{\frac{1+y}{y\rho}}  - \frac1{\sqrt{y}}\,. \]
Hence as $\MARGIN' < 0$ and as $y, \rho > 0$, we have 
\begin{align*}
z^2 \ge |\MARGIN'|^2 \cdot \frac{1+y}{y \rho} + \frac1{y} \ge \frac{|\MARGIN'|^2}{\rho}\,.
\end{align*}
Moreover when $\INVERSEMILLSEXPRESSION > 1$, by Lemma \upref{lem:millsratioproperties}, 
\[ \cA(\INVERSEMILLSEXPRESSION)^2 \le (2\INVERSEMILLSEXPRESSION+1)^2 \le 8\INVERSEMILLSEXPRESSION^2 + 2 \le 16\Big( |\MARGIN'|^2 \cdot \frac{1+y}{\rho} + z^2 y \Big)+2\,.\]
Consequently we can upper bound 
\begin{align*}
&\int_{s:|s| < |\MARGIN|} \int_z \cA(\INVERSEMILLSEXPRESSION)^2 \pdfnormal(z) \pdfsignal(s)\, \rmd z \rmd s \\
&\quad= \int_{(s, z) \not\in E, |s| < |\MARGIN|} \cA(\INVERSEMILLSEXPRESSION)^2 \pdfnormal(z) \pdfsignal(s)\, \rmd z \rmd s + \int_{(s, z) \in E, |s| < |\MARGIN|} \cA(\INVERSEMILLSEXPRESSION)^2 \pdfnormal(z) \pdfsignal(s)\, \rmd z \rmd s \\
&\quad = \int_{(s, z) \not\in E, |s| < |\MARGIN|} \Big(\frac{16 |\MARGIN'|^2 (y+1)}{\rho} + z^2 y + 2 \Big) \pdfnormal(z) \pdfsignal(s)\, \rmd z \rmd s + \frac1{2\pi \cN(1)^2} \\
&\quad \le 16 (y+1) \int_{\R^2} \frac{|\MARGIN'|^2}{\rho} \exp\Big( -\frac{|\MARGIN'|^2}{4\rho} \Big) \cdot \frac1{\sqrt{2\pi}} \exp(-z^2/4) \pdfsignal(s)\, \rmd z \rmd s + y + 2 + \frac1{2\pi \cN(1)^2} \\
&\quad \le 35 y + 43\,.
\end{align*}
Here we used that $\E[Z^2] = 1$, and upper bounded $\frac{\MARGIN'^2}{\rho} \exp\Big( -\frac{\MARGIN'^2}{4\rho} \Big)$ by using that for all $x \ge 0$, we have $x \exp(-x/4) \le K$ where $K$ is a universal constant.

\item Else if $|s| \ge |\MARGIN|$, note by Lemma \upref{lem:millsratioproperties} that 
\begin{align*}
\cA(\INVERSEMILLSEXPRESSION)^2 \le 1 + (2\INVERSEMILLSEXPRESSION+1)^2  \le 8\INVERSEMILLSEXPRESSION^2 + 3 &\le 16\Big((\MARGIN - xs)^2 \cdot \frac{1+y}{\rho} + z^2 y\Big) + 3 \\
&\le 32 \Csol (\MARGIN^2 +  s^2) (y+1) + 16 z^2 y + 3\,.    
\end{align*}
Thus as $\SIGNALDISORDER^2$ is $\signalsubexp$ sub-Exponential and $\E[\SIGNALDISORDER^2] \le \signalsubexp$, 
\begin{align*}
&\int_{s:|s| \ge |\MARGIN|} \int_z \cA(\INVERSEMILLSEXPRESSION)^2 \pdfnormal(z) \pdfsignal(s)\, \rmd z \rmd s \\
&\quad \le 32 \Csol (y+1) \MARGIN^2 \int_{s:|s| \ge |\MARGIN|} \int_z \pdfnormal(z) \pdfsignal(s)\, \rmd z \rmd s +  (y+1) \int_{\R^2} s^2 \pdfnormal(z) \pdfsignal(s)\, \rmd z \rmd s \\
&\qquad \qquad \qquad + 16 y \int_{\R^2} z^2 \pdfnormal(z) \pdfsignal(s)\, \rmd z \rmd s + 3 \\
&\quad = 32 \Csol (y+1) \MARGIN^2 \cdot \P\Big(S^2 \ge \MARGIN^2 \Big) +  \signalsubexp (y+1) + 16y+3 \\
&\quad \le 32 \Csol (y+1) \MARGIN^2 \cdot \exp\Big( -\frac{\MARGIN^2}{ \signalsubexp}\Big) +  \signalsubexp (y+1) + 16y+3 \\
&\quad \le y(12 \Csol  \signalsubexp +  \signalsubexp + 16 ) + (12 \Csol  \signalsubexp +  \signalsubexp + 3)\,.
\end{align*}
\end{itemize}
Summing the above bounds yields the desired upper bound on $\E\big[ \cA(\INVERSEMILLSEXPRESSION)^2 \big]$.
\end{proof}

\begin{lemma}\label{lem:ybounded}
Suppose $(\alpha,\GAUSSIANMEASURE) \in [0, \DENSITYBOUND] \times [\beta_0, \beta_1]$. For any $(y,\rho)$ such that $y \ge 0$, $\rho \ge \frac1{\Csol}$ that satisfies $\frac{\partial G}{\partial y}(y, \rho)=0$ -- in particular, to any solution to \equpref{eq:RStransformed} -- we must have $y \le 1$.
\end{lemma}
\begin{proof}
By our expression \equpref{eq:partialGy} for $\frac{\partial G}{\partial y}$ and by Lemma \upref{lem:controlexpectationofAWinit}, we obtain
\[ y = \alpha \E\big[ \cA(\INVERSEMILLSEXPRESSION)^2 \big] \le \alpha (K_{\upref{lem:controlexpectationofAWinit}} y + K_{\upref{lem:controlexpectationofAWinit}}) \le \frac12 y + \frac12 \implies y \le 1\,, \]
where we used that $\alpha \le \frac1{2 K_{\upref{lem:controlexpectationofAWinit}}}$.
\end{proof}

The last preliminary step we need is the following Lemma, which lets us control various expectations of $\cA(\INVERSEMILLSEXPRESSION)$.
\begin{lemma}\label{lem:controlexpectationofAW}
Suppose $(\alpha,\GAUSSIANMEASURE)\in [0, \DENSITYBOUND] \times [\beta_0, \beta_1]$. Then we have for any $(y,\rho) \in [0,1] \times [\frac1{\Csol}, \Csol]$ that
\[ 0 \le \E\big[ \cA(\INVERSEMILLSEXPRESSION) \big]\,,\, \E\big[ \cA(\INVERSEMILLSEXPRESSION)^2 \big]\,,\, \E\big[ \big| \MARGIN' \cA(\INVERSEMILLSEXPRESSION) \big| \big] \le K_{\upref{lem:controlexpectationofAW}}\,,\]
where
\begin{align*}
K_{\upref{lem:controlexpectationofAW}} = \begin{cases} 42 \Csol  \signalsubexp + 13 \Csol^{3/2} + 7.5 &: \MARGIN < 0\,, \\ 8.3 \Csol^{1/2} (\MARGIN^2 + \signalsubexp) + 2.3 &: \MARGIN \ge 0\,. \end{cases}    
\end{align*}
\end{lemma}
\begin{proof}
The desired lower bounds are obvious as $\cA \ge 0$, and the first two upper bounds are immediate by combining Lemma \upref{lem:controlexpectationofAWinit} and Lemma \upref{lem:ybounded}. 
As $\cA \ge 0$, it remains to upper bound $\E\big[ |\MARGIN'| \cA(\INVERSEMILLSEXPRESSION) \big]$, which we do similarly as in the proof of Lemma \upref{lem:controlexpectationofAWinit}.
We again break into the cases $\MARGIN \ge 0$, $\MARGIN < 0$.

\paragraph{When $\MARGIN \ge 0$.}
By Lemma \upref{lem:millsratioproperties}, we have $\cA(x) \le 2|x|+1$ for all $x \in \R$. 
We thus may write
\begin{align*}
|\MARGIN'| \cA(\INVERSEMILLSEXPRESSION) &\le 2|\MARGIN'| |\INVERSEMILLSEXPRESSION| + |\MARGIN'| \\
&\le 2 \Big( |\MARGIN'|^2 \sqrt{\frac{1+y}{\rho}} + \sqrt{y} |Z| |\MARGIN'| \Big) + |\MARGIN'| \\
&\le 4 (\MARGIN^2 + x^2 S^2) \sqrt{2 \Csol} + \frac{2|Z|+1}2 \cdot \big( 2(\MARGIN^2 + x^2 S^2) + 1 \big)\,,
\end{align*}
where we upper bounded $|\MARGIN'|$ by $\frac{\MARGIN'^2+1}2$. Consequently as $S$, $Z$ are independent,
\begin{align*}
\E\big[ \big| \MARGIN' \cA(\INVERSEMILLSEXPRESSION) \big| \big] &\le 5.7 (\MARGIN^2 + \signalsubexp) \Csol^{1/2} + 2.6(\MARGIN^2 + \signalsubexp) + 2.3 \le K_{\upref{lem:controlexpectationofAW}}\,.
\end{align*}

\paragraph{When $\MARGIN < 0$.}
We write
\begin{align*}
\E\big[ |\MARGIN'| \cA(\INVERSEMILLSEXPRESSION) \big] &= \int_{s, z} |\MARGIN'| \cA(\INVERSEMILLSEXPRESSION) \pdfnormal(z) \pdfsignal(s) \, \rmd z \rmd s \\
&= \int_{s:|s| < |\MARGIN|} \int_z |\MARGIN'| \cA(\INVERSEMILLSEXPRESSION) \pdfnormal(z) \pdfsignal(s) \, \rmd z \rmd s  + \int_{s:|s| \ge |\MARGIN|} \int_z |\MARGIN'| \cA(\INVERSEMILLSEXPRESSION) \pdfnormal(z) \pdfsignal(s) \, \rmd z \rmd s\,.
\end{align*}
Again, we upper bound each of these two integrals separately.
\begin{itemize}

\item If $|s| < |\MARGIN|$, then as $|xs| < |\MARGIN|$ and $\MARGIN' < 0$, we have $\MARGIN' < 0$.
Define the sets
\begin{align*}
E_1 &:= \Big\{ (s, z)\,:\, z\sqrt{y} \ge \frac{\MARGIN'}2 \sqrt{\frac{1+y}{\rho}} \Big\}\,, \\
E_2 &:= \Big\{ (s, z)\,:\, \MARGIN' \sqrt{\frac{1+y}{\rho}} - 1 < z\sqrt{y} < \frac{\MARGIN'}2 \sqrt{\frac{1+y}{\rho}} \Big\}\,, \\
E_3 &:= \Big\{ (s, z)\,:\, z\sqrt{y} \le \MARGIN' \sqrt{\frac{1+y}{\rho}} - 1\Big\}\,.
\end{align*}
Note as $\MARGIN' < 0$ and $y, \rho \ge 0$, $E_1, E_2, E_3$ are all disjoint and partition $\R^2$. 
We thus can write 
\begin{align*}
\int_{s:|s| < |\MARGIN|} \int_z |\MARGIN'| \cA(\INVERSEMILLSEXPRESSION) \pdfnormal(z) \pdfsignal(s) \, \rmd z \rmd s &= \expressionI + \expressionII + \expressionIII\,,
\end{align*}
where
\begin{align*}
\expressionI &:= \int_{(s, z) \in E_1:|s| < |\MARGIN|} |\MARGIN'| \cA(\INVERSEMILLSEXPRESSION) \pdfnormal(z) \pdfsignal(s) \, \rmd z \rmd s \\
\expressionII &:= \int_{(s, z) \in E_2:|s| < |\MARGIN|} |\MARGIN'| \cA(\INVERSEMILLSEXPRESSION) \pdfnormal(z) \pdfsignal(s) \, \rmd z \rmd s \\
\expressionIII &:= \int_{(s, z) \in E_3:|s| < |\MARGIN|} |\MARGIN'| \cA(\INVERSEMILLSEXPRESSION) \pdfnormal(z) \pdfsignal(s) \, \rmd z \rmd s \,.
\end{align*}
We now upper bound $\expressionI, \expressionII, \expressionIII$:
\begin{enumerate}
    
\item For $\expressionI$: by definition of $E_1$, for $(s,z) \in E_1$ we have 
\begin{align*}
\INVERSEMILLSEXPRESSION = \MARGIN' \sqrt{\frac{1+y}{\rho}} - z\sqrt{y} \le \frac{\MARGIN'}{2}\sqrt{\frac{1+y}{\rho}} < 0\,.
\end{align*}
Hence $\cN(\INVERSEMILLSEXPRESSION) \ge \frac12$ for $(s,z) \in E_1$, and moreover as $\MARGIN' < 0$, we have for $(s,z) \in E_1$ that
\[ \INVERSEMILLSEXPRESSION^2 \ge \frac{|\MARGIN'|^2}{4} \cdot \frac{1+y}{\rho} \ge \frac{|\MARGIN'|^2}{4\rho} \ge \frac{|\MARGIN'|^2}{4 \Csol} \,.\]
Thus for $(s, z) \in E_1$, 
\begin{align*}
0 \le \cA(\INVERSEMILLSEXPRESSION) \le \frac{e^{-\INVERSEMILLSEXPRESSION^2/2}}{\sqrt{2\pi} \cN(\INVERSEMILLSEXPRESSION)} \le \frac2{\sqrt{2\pi}} \exp\Big( -\frac{|\MARGIN'|^2}{8\Csol} \Big)\,.
\end{align*}
Thus, 
\begin{align*}
\expressionI \le \int_{(s, z) \in E_1:|s| < |\MARGIN|} \frac{2|\MARGIN'|}{\sqrt{2\pi}} \exp\Big( -\frac{|\MARGIN'|^2}{8\Csol} \Big) \pdfnormal(z) \pdfsignal(s) \, \rmd z \rmd s \le \sqrt{\Csol}\,.
\end{align*}

\item For $\expressionII$: for $(s, z) \in E_2$, we have $\INVERSEMILLSEXPRESSION = \MARGIN' \sqrt{\frac{1+y}{\rho}} - z\sqrt{y} \le 1$. Therefore for $(s,z) \in E_2$, $\cA(\INVERSEMILLSEXPRESSION) \le \frac{1}{\sqrt{2\pi} \cN(1)}$. Furthermore, it is not possible in this case that $y=0$, as then we have $0 = z\sqrt{y} < \frac{\MARGIN'}2\sqrt{\frac{1+y}{\rho}} < 0$ as $\MARGIN' < 0$. Thus $y>0$ and so
\[ z \le \frac{\MARGIN'}2 \sqrt{\frac{1+y}{y \rho}} < 0 \implies z^2 \ge \frac{|\MARGIN'|^2}4 \cdot \frac{1+1/y}{\rho} \ge \frac{|\MARGIN'|^2}{4\Csol}\,. \]
Therefore, 
\begin{align*}
\expressionII &\le \int_{(s, z) \in E_2:|s| < |\MARGIN|} |\MARGIN'| \cdot \frac{1}{2\pi \cN(1)} \exp\Big(-\frac{|\MARGIN'|^2}{16 \Csol}\Big) \exp(-z^2/4) \pdfsignal(s) \, \rmd z \rmd s \\
&\le 1.72\sqrt{\Csol} \int_{\R^2 } \frac1{\sqrt{2\pi}} \exp(-z^2/4) \pdfsignal(s) \, \rmd z \rmd s \le 2.5\sqrt{\Csol} \,.
\end{align*}

\item For $\expressionIII$: by Lemma \upref{lem:millsratioproperties} and Lemma \upref{lem:ybounded}, 
\begin{align*}
0 \le \cA(\INVERSEMILLSEXPRESSION) &\le 2|\INVERSEMILLSEXPRESSION|+1 \le 2|\MARGIN'|\sqrt{\frac{1+y}{\rho}} + 2|z|\sqrt{y} + 1 \le 3|\MARGIN'|\sqrt{\Csol} + 2|z| + 1 \,.
\end{align*}
Moreover, here we have as $\MARGIN' \le 0$ and $\rho \le \Csol$ that
\begin{align*}
z \le \MARGIN'\sqrt{\frac{1+y}{y \rho}} - \frac1{\sqrt{y}} < 0 \implies z^2 \ge |\MARGIN'|^2 \cdot \frac{1+y}{y \rho} + \frac1{y} \ge \frac{|\MARGIN'|^2}{\Csol}\,.
\end{align*}
Therefore,
\begin{align*}
\expressionIII &\le 3\sqrt{\Csol} \int_{(s, z) \in E_3:|s| < |\MARGIN|} |\MARGIN'|^2 \exp\Big( -\frac{|\MARGIN'|^2}{4\Csol} \Big) \cdot \frac1{\sqrt{2\pi}} \exp(-z^2/4) \pdfsignal(s) \, \rmd z \rmd s \\
&\qquad + 2\int_{(s, z) \in E_3:|s| < |\MARGIN|} |z| |\MARGIN'| \exp\Big( -\frac{|\MARGIN'|^2}{4\Csol} \Big) \cdot \frac1{\sqrt{2\pi}} \exp(-z^2/4) \pdfsignal(s) \, \rmd z \rmd s \\
&\qquad + \int_{(s, z) \in E_3:|s| < |\MARGIN|} |\MARGIN'| \exp\Big( -\frac{|\MARGIN'|^2}{4\Csol} \Big) \cdot \frac1{\sqrt{2\pi}} \exp(-z^2/4) \pdfsignal(s) \, \rmd z \rmd s \\
&\le 6.25 \Csol^{3/2} + 3.15 \Csol^{1/2}\,.
\end{align*}
\end{enumerate}
Putting everything together and using that $\Csol \ge 1$, we obtain
\[ \int_{s:|s| < |\MARGIN|} \int_z |\MARGIN'| \cA(\INVERSEMILLSEXPRESSION) \pdfnormal(z) \pdfsignal(s) \, \rmd z \rmd s = \expressionI + \expressionII + \expressionIII \le 13 \Csol^{3/2}\,.\]

\item If $|s| \ge |\MARGIN|$, we obtain by Lemma \upref{lem:millsratioproperties} and Lemma \upref{lem:ybounded} that
\begin{align*}
0 \le \cA(\INVERSEMILLSEXPRESSION) \le 2\INVERSEMILLSEXPRESSION+1 \le 2|\INVERSEMILLSEXPRESSION|+1 \le 3 |\MARGIN'| \sqrt{\Csol} + 2|z| + 1\,,
\end{align*}
thus by AM-GM,
\begin{align*}
|\MARGIN'| \cA(\INVERSEMILLSEXPRESSION) \le \frac{|\MARGIN'|^2}2 + \frac32 \big(9 |\MARGIN'|^2 \Csol + 4|z|^2+1 \big) &\le |\MARGIN'|^2(13.5 \Csol + 0.5) + 6 z^2 + 1.5 \\
&\le (\MARGIN^2 + s^2  )(27 \Csol + 1) + 6 z^2 + 1.5\,.
\end{align*}
We obtain
\begin{align*}
&\int_{s:|s| \ge |\MARGIN|} \int_z \cA(\INVERSEMILLSEXPRESSION)^2 \pdfnormal(z) \pdfsignal(s)\, \rmd z \rmd s  \\
&\quad \le (27 \Csol + 1) \MARGIN^2 \int_{s:|s| \ge |\MARGIN|} \int_z \pdfnormal(z) \pdfsignal(s)\, \rmd z \rmd s + (27 \Csol + 1)  \int_{s:|s| \ge |\MARGIN|} \int_z s^2 \pdfnormal(z) \pdfsignal(s)\, \rmd z \rmd s \\
&\qquad\qquad  + \int_{s:|s| \ge |\MARGIN|} (6 z^2 + 1.5) \pdfnormal(z) \pdfsignal(s)\, \rmd z \rmd s \\
&\quad \le (27 \Csol + 1) \MARGIN^2 \P\Big(S^2 > \MARGIN^2 \Big) + (27 \Csol + 1)  \signalsubexp + 7.5 \\
&\quad \le (27 \Csol + 1) \MARGIN^2 \exp\Big( -\frac{\MARGIN^2}{ \signalsubexp} \Big) + (27 \Csol + 1)  \signalsubexp + 7.5 \\
&\quad \le 1.5(27 \Csol + 1)  \signalsubexp + 7.5 \,,
\end{align*}
where we used that $S^2$ is $\signalsubexp$ sub-Exponential and $\E[S^2] \le \signalsubexp$.
\end{itemize}
Summing the upper bounds from both cases and using that $\Csol \ge 1$, we obtain
\[ \E\big[ |\MARGIN'| \cA(\INVERSEMILLSEXPRESSION) \big] \le 42  \Csol  \signalsubexp + 13 \Csol^{3/2} + 7.5 \,.\]
This proves the Lemma in both cases $\MARGIN < 0$, $\MARGIN \ge 0$.
\end{proof}

\paragraph{Finishing the argument.}
By Lemma \upref{lem:ybounded}, to show \equpref{eq:RStransformed} has a unique solution in the desired domain, it suffices to show \equpref{eq:RStransformed} has a unique solution in $(y,\rho) \in [0,1] \times [\frac1{\Csol}, \Csol]$. 
Now, we consider the signs of the derivatives of $G$ in $y$ and $\rho$. This is done for the interpolators and posterior together in a unified manner.
\begin{lemma}\label{lem:rhocontrolled}
Suppose $(\alpha,\GAUSSIANMEASURE)\in [0, \DENSITYBOUND] \times [\GAUSSIANMEASURE_0, \GAUSSIANMEASURE_1]$. For $(y, \rho) \in [0,1] \times [\frac1{\Csol}, \Csol]$, we have $\frac{\partial^2 G}{\partial^2 \rho}(y, \rho) < 0$.
\end{lemma}
\begin{proof}
We first establish the desired result for the interpolators. Consider \equpref{eq:partialGrhosecondderiv}, which states
\begin{align*}
\frac{\partial^2 G}{\partial^2 \rho}(y, \rho) = -\frac{3\alpha (1+y)^{1/2}}{4 \rho^{5/2}} \E\Big[ \MARGIN' \cA(\INVERSEMILLSEXPRESSION) \Big] - \frac{\alpha (1+y)}{4 \rho^3} \E\Big[ \MARGIN'^2 \cA'(\INVERSEMILLSEXPRESSION) ] - \frac1{2\rho^2}\,.
\end{align*}
The second and third terms in this expression are non-positive as $\cA' \ge 0$ by Lemma \upref{lem:millsratioproperties}, and the third term is strictly negative. Finally, note by Lemma \upref{lem:ybounded} and Lemma \upref{lem:controlexpectationofAW} and since $\alpha \le \frac{2}{3 \sqrt{2\Csol} K_{\upref{lem:controlexpectationofAW}}}$, we have
\begin{align*}
\Big|-\frac{3\alpha (1+y)^{1/2}}{4 \rho^{5/2}} \E\Big[ \MARGIN' \cA(\INVERSEMILLSEXPRESSION) \Big] \Big| \le \frac{3\sqrt{2 \Csol} \alpha }{4 \rho^2} \cdot K_{\upref{lem:controlexpectationofAW}} \le \frac1{4\rho^2}\,.
\end{align*}
It follows that 
\[ \frac{\partial^2 G}{\partial \rho^2} \le \frac1{4\rho^2}-\frac1{2\rho^2}=-\frac1{2\rho^2}<0\,.\]

For the posterior, we follow the same strategy. Recall \equpref{eq:posterior_partialGrho_twice} states that 
\begin{align*}
\frac{\partial^2 G}{\partial \rho^2} &= \frac{\alpha}4\Big( \expressionI + \expressionII + \expressionIII + \expressionIV \Big) -\frac1{2\rho^2} \,,
\end{align*}
where $\expressionI, \expressionII, \expressionIII, \expressionIV$ are as defined in \equpref{eq:posterior_secondderiv_expressions}. As $1/(1+y), y/(1+y) \in [0,1]$ for $y \in [0,1]$, we have $|\expressionI|, |\expressionII|, |\expressionIII|, |\expressionIV| \le K(D)$. Since $\rho \le \Csol$ and $\alpha$ is small enough in terms of $D, \Csol$, it follows that $\frac{\partial^2 G}{\partial \rho^2}<0$ for the posterior as well.
\end{proof}

\begin{lemma}\label{lem:ycontrolled}
Suppose $(\alpha,\GAUSSIANMEASURE)\in [0, \DENSITYBOUND] \times [\GAUSSIANMEASURE_0, \GAUSSIANMEASURE_1]$. Then for any $\rho' \in [\frac1{\Csol}, \Csol]$, we have the following:
\begin{enumerate}
    \item There exists a unique $y' \in [0, 1]$ such that $\frac{\partial G}{\partial y}(y',\rho')=0$. 
    \item The one variable function $y \rightarrow G(y, \rho')$ attains its minimum on $[0,\infty)$ at $y=y'$. 
    \item We have $\frac{\partial^2 G}{\partial y^2}(y', \rho') > 0$.
    \item Defining $y' := y'(\rho')$ from 1), $y'(\rho')$ is continuous in $\rho'$ for $\rho' \in [\frac1{\Csol}, \Csol]$.
\end{enumerate}
\end{lemma}
\begin{proof}
We first prove 1), 2) of the Lemma. 
Recall $g(y, \rho)$ defined in \equpref{eq:partialGpartialyhelperfunc} for the interpolators or \equpref{eq:posterior_RS_g_def} for the posterior.
In either case, $g(y, \rho)$ is clearly jointly continuous in both arguments.
We claim that to prove 1), 2) of the Lemma, it suffices to show that $g(y, \rho')=0$ at a unique $y=y'$, and that $g(y, \rho') < 0$ for $y < y'$, $g(y, \rho') > 0$ for $y$ such that $y' < y \le 1$.

To justify why, note by \equpref{eq:partialGpartialyhelperfunc} or \equpref{eq:posterior_RS_g_def} that this implies $\frac{\partial G}{\partial y}(y,\rho')=0$ at $y=y'$, and that $\frac{\partial G}{\partial y}(y,\rho')<0$ for $y<y'$, $\frac{\partial G}{\partial y}(y,\rho')>0$ for $y' < y \le 1$. From here, we note that another $y'' \in (1,\infty)$ with $\frac{\partial G}{\partial y}(y'', \rho')=0$ would contradict Lemma \upref{lem:ybounded}. This would imply 1). 
Moreover by continuity of $\frac{\partial G}{\partial y}(\cdot, \rho')$ and as $\frac{\partial G}{\partial y}(y,\rho')>0$ for $y' < y \le 1$, since there is no $y'' \in (1,\infty)$ with $\frac{\partial G}{\partial y}(y'', \rho')=0$, this would mean that $\frac{\partial G}{\partial y}(y'', \rho') > 0$ for all $y'' \in (1, \infty)$. Together with the above, this would imply 2).

Now, we claim that $\frac{\partial g}{\partial y}(y, \rho') > 0$ for all $y \in [0, 1]$, for both the interpolators and the posterior:
\begin{itemize}
    \item For the interpolators: By \equpref{eq:partialGpartialysecond}, 
\begin{align*}
\frac{\partial g}{\partial y} &= 1 - \frac{\alpha}{\sqrt{\rho(1+y)}} \E\Big[ \MARGIN' \cA(\INVERSEMILLSEXPRESSION) \cA'(\INVERSEMILLSEXPRESSION) \Big] - \alpha \E\Big[ \cA'^2(\INVERSEMILLSEXPRESSION)+\cA(\INVERSEMILLSEXPRESSION)\cA''(\INVERSEMILLSEXPRESSION) \Big]\,.
\end{align*}
By Lemma \upref{lem:millsratioproperties} and Lemma \upref{lem:controlexpectationofAW} we have 
\begin{align*}
\Big| \E\Big[ \cA'^2(\INVERSEMILLSEXPRESSION)+\cA(\INVERSEMILLSEXPRESSION)\cA''(\INVERSEMILLSEXPRESSION) \Big] \Big| &\le \MILLSBOUND^2 + \MILLSBOUND \E[ \cA(\INVERSEMILLSEXPRESSION) ] \le \MILLSBOUND^2 + \MILLSBOUND K_{\upref{lem:controlexpectationofAW}} \,,\\
\Big| \frac1{\sqrt{\rho(1+y)}}\E\Big[ \MARGIN' \cA(\INVERSEMILLSEXPRESSION) \cA'(\INVERSEMILLSEXPRESSION) \Big] \Big| &\le \Csol^{1/2} \MILLSBOUND \E\Big[ \big|\MARGIN' \cA(\INVERSEMILLSEXPRESSION) \big| \Big] \le \Csol^{1/2} \MILLSBOUND K_{\upref{lem:controlexpectationofAW}}\,.
\end{align*}
Thus since $\alpha \le \frac14 \min\big\{ \frac1{\MILLSBOUND^2 + \MILLSBOUND K_{\upref{lem:controlexpectationofAW}}}, \frac1{\Csol^{1/2} \MILLSBOUND K_{\upref{lem:controlexpectationofAW}}} \big\}$, it follows that 
\begin{align}\label{eq:RS_interpolants_lowerbd_gderiv}
\frac{\partial g}{\partial y} &= 1 - \frac{\alpha}{\sqrt{\rho(1+y)}} \E\Big[ \MARGIN' \cA(\INVERSEMILLSEXPRESSION) \cA'(\INVERSEMILLSEXPRESSION) \Big] - \alpha \E\Big[ \cA'^2(\INVERSEMILLSEXPRESSION)+\cA(\INVERSEMILLSEXPRESSION)\cA''(\INVERSEMILLSEXPRESSION) \Big] \ge \frac12 > 0\,.
\end{align}

\item For the posterior: Recall \equpref{eq:posterior_secondderiv_y}, that $y \ge 0$ and $|u^l| \le D$ for $l=1,2,3,4$, and that $\rho' \le \Csol$, $\frac1{1+y} \le 0$. Thus since $\alpha \le \alpha_0$ and $\alpha_0$ is small enough in terms of $D, \Csol$,
\begin{align}\label{eq:posterior_RS_interpolants_lowerbd_gderiv}
\frac{\partial g}{\partial y} = 1 + 2y - \frac{\alpha \rho}{(1+y)^2} \Big( \expressionI + \expressionII \Big) \ge \frac12\,,
\end{align}
where $\expressionI, \expressionII$ are as in \equpref{eq:posterior_secondderiv_y_expressions}. We have $|\expressionI|, |\expressionII| \le K(D)$ as $|u^l| \le D$ for $l=1,2,3,4$.
\end{itemize}
Note that the bounds \equpref{eq:RS_interpolants_lowerbd_gderiv}, \equpref{eq:posterior_RS_interpolants_lowerbd_gderiv} hold uniformly in $(y, \rho) \in [0,1] \times [\frac1{\Csol}, \Csol]$. This implies there is at most one $y' \in [0,1]$ such that $g(y',\rho')=0$. Moreover, if this $y'$ exists, our above work clearly implies $g(y) < 0$ for $y < y'$, $g(y) > 0$ for $y > y'$. 

Now we show there exists such an $y' \in [0,1]$. When $y=0$, evidently $g(y) \le 0$. Now consider when $y=1$. For the interpolators, note as $\alpha < \frac1{K_{\upref{lem:controlexpectationofAW}}}$, we have by Proposition \upref{lem:controlexpectationofAW} that
\[ g(1)=1 -\alpha \E[ \cA(\INVERSEMILLSEXPRESSION)^2 ] > 0. \]
For the posterior, since $|u'| \le D, \rho' \le \Csol$, and $\alpha$ is small enough in terms of $D, \Csol$, we have
\begin{align*}
g(1) = 2 - \alpha \rho \E \big\langle \exp u' \big\rangle^2 > 0\,.
\end{align*}
Thus for both the interpolators and the posterior, there exists $y' \in [0,1]$ with $g(y', \rho')=0$ by the Intermediate Value Theorem. The above steps now prove that $g(y, \rho')=0$ at a unique $y=y'$, and that $g(y, \rho') < 0$ for $y < y'$, $g(y, \rho') > 0$ for $y$ such that $y' < y \le 1$, which as discussed earlier implies 1), 2) of the Lemma.
\newline

Next, we prove 3) of the Lemma for both the interpolators and the posterior:
\begin{itemize}
\item For the interpolators: Differentiating \equpref{eq:partialGpartialyhelperfunc}, 
\[ \frac{\partial^2 G}{\partial y^2}(y, \rho') = \frac{\frac{\partial g}{\partial y}(y, \rho') (1+y)-g(y, \rho')}{2(1+y)^2}\, .\]
Observe $g(y, \rho') \le y$. Also from the above steps, we have 
\[ \frac{\partial g}{\partial y}(y, \rho') \ge 1 - \alpha \max\big\{ \MILLSBOUND^2 + \MILLSBOUND K_{\upref{lem:controlexpectationofAW}}, \Csol^{1/2} \MILLSBOUND K_{\upref{lem:controlexpectationofAW}}\big\}\,.\]
It therefore remains to show 
\[ 1 - \alpha \max\big\{ \MILLSBOUND^2 + \MILLSBOUND K_{\upref{lem:controlexpectationofAW}}, \Csol^{1/2} \MILLSBOUND K_{\upref{lem:controlexpectationofAW}}\big\}  > \frac{y'}{y'+1}. \]
Since $\frac{y}{y+1}$ is increasing and as $y' \le 1$ by 1), it suffices to show 
\[ 1 - \alpha \max\big\{ \MILLSBOUND^2 + \MILLSBOUND K_{\upref{lem:controlexpectationofAW}}, \Csol^{1/2} \MILLSBOUND K_{\upref{lem:controlexpectationofAW}}\big\} > \frac23 \,, \]
which holds true by our condition on $\alpha$. 
\item For the posterior: Differentiating \equpref{eq:posterior_RS_g_def} and using \equpref{eq:posterior_secondderiv_y}, 
\begin{align*}
\frac{\partial^2 G}{\partial y^2}(y, \rho') = \frac{\frac{\partial g}{\partial y}(y, \rho') (1+y) - 2g(y)}{2(1+y)^3} = \frac{1+y + 2 \alpha \rho \E \big\langle \exp u' \big\rangle^2 - \frac{\alpha \rho}{(1+y)^2} \big( \expressionI + \expressionII \big)}{2(1+y)^3}\,,
\end{align*}
where $\expressionI, \expressionII$ are as in \equpref{eq:posterior_secondderiv_y_expressions}. Note $|\expressionI|, |\expressionII| \le K(D)$ as $|u^l| \le D$ for $l=1,2,3,4$. Thus as $y \ge 0, \rho \le \Csol$, and $\alpha \le \alpha_0$ is small enough in terms of $D, \Csol$, we have $\frac{\partial^2 G}{\partial y^2}(y', \rho')>0$.
\end{itemize}
This proves 3) of the Lemma in both cases.

Finally, we prove 4) of the Lemma. Note that $y(\rho')$ is a solution to the fixed point equation $g(y(\rho'), \rho') = 0$. As $\frac{\partial g}{\partial y} > 0$ holds uniformly on $[0,1] \times [\frac1{\Csol}, \Csol]$ as argued above, and as $y(\rho') \in [0,1]$ as shown in 1), we conclude 4) by the Implicit Function Theorem.
\end{proof}

The above Lemmas establish the desired convex-concave structure of the RS equations, and we are now ready to complete the proof of Proposition \upref{prop:nicersuniquesol}.
\begin{lemma}\label{lem:atmostonesol}
Consider $(\alpha,\GAUSSIANMEASURE) \in [0, \DENSITYBOUND] \times [\beta_0, \beta_1]$. Then there is at most one solution $(y,\rho) \in [0,1] \times [\frac1{\Csol}, \Csol]$ to \equpref{eq:RStransformed}. 
\end{lemma}
\begin{proof}
Suppose that there are two distinct solutions $(y_1,\rho_1)$ and $(y_2,\rho_2)$ both in the domain $[0,1] \times [\frac1{\Csol}, \Csol]$. 
Assume without loss of generality that $y_1<y_2$. 
Directly from Lemma \upref{lem:rhocontrolled} and Lemma \upref{lem:ycontrolled}, we have $y_1 \neq y_2$ and $\rho_1 \neq \rho_2$. 
Applying Taylor Expansion to degree 2 about $(y_1,\rho_1)$, since $\frac{\partial G}{\partial \rho}(y_1,\rho_1)=0$ by definition of the system of equations \equpref{eq:RStransformed} and since $\frac{\partial^2 G}{\partial \rho^2}(y_1,\cdot)<0$ from Lemma \upref{lem:rhocontrolled}, we see that $G(y_1, \rho_1) > G(y_1, \rho_2)$ (since we Taylor expand about a $\rho' \in [\min\{\rho_1, \rho_2\}, \max\{\rho_1, \rho_2\}] \subseteq [\frac1{\Csol}, \Csol]$, we may apply the above work). 
By Lemma \upref{lem:ycontrolled}, since $\frac{\partial G}{\partial y}(y_2, \rho_2)=0$, we have $G(y_1,\rho_2) > G(y_2,\rho_2)$. 
By analogous reasoning again using Lemma \upref{lem:rhocontrolled}, we have $G(y_2,\rho_2) > G(y_2,\rho_1)$. 
Again by Lemma \upref{lem:ycontrolled}, we have $G(y_2,\rho_1)>G(y_1,\rho_1)$. 
Combining these inequalities, we obtain
\[ G(y_1, \rho_1) >G(y_1,\rho_2) >G(y_2,\rho_2) >G(y_2,\rho_1)>G(y_1,\rho_1)\,,\]
a contradiction.
\end{proof}

\begin{lemma}\label{lem:solexists}
Consider $(\alpha,\GAUSSIANMEASURE) \in [0, \DENSITYBOUND] \times [\beta_0, \beta_1]$. Then there exists a solution $(y,\rho) \in [0,1] \times [\frac1{\Csol}, \Csol]$ to \equpref{eq:RStransformed}. 
\end{lemma}
\begin{proof}
By Lemma \upref{lem:ycontrolled} and Lemma \upref{lem:ybounded}, for any $\rho \in [\frac1{\Csol}, \Csol]$, there exists a unique $y(\rho) \in [0,1]$ where $\frac{\partial G}{\partial y}(y(\rho), \rho)=0$. Thus it is enough to show there exists a $\rho' \in [\frac1{\Csol}, \Csol]$ such that $\frac{\partial G}{\partial \rho}(y(\rho'), \rho')=0$. 

Let $h(\rho):=2\frac{\partial G}{\partial \rho}(y(\rho), \rho)$. By 4) of Lemma \upref{lem:ycontrolled} and using the explicit form \equpref{eq:partialGrho}, \equpref{eq:posterior_partialGrho} for $\frac{\partial G}{\partial \rho}$ for the interpolators and posterior respectively, $h(\rho)$ varies continuously in $\rho \in [\frac1{\Csol}, \Csol]$. 

We now prove that $h(1.5/\beta) < 0, h(1/3.5 \beta) > 0$ for both the interpolators and the posterior.
\begin{itemize}
    \item For the interpolators: Observe that because $y(\rho) \le 1$ by Lemma \upref{lem:ybounded}, Lemma \upref{lem:controlexpectationofAW} yields
\begin{align*}
\Big| \frac{(1+y(\rho))^{1/2}}{\rho^{3/2}} \E\Big[ \MARGIN' \cA(\INVERSEMILLSEXPRESSION) \Big] \Big| \le 2^{1/2} \Csol^{3/2} K_{\upref{lem:controlexpectationofAW}}\,.
\end{align*}
Now note by \equpref{eq:partialGrho} and as $\alpha < \frac{0.9 \beta_0}{\Csol^{3/2} K_{\upref{lem:controlexpectationofAW}}}$, we have
\[ h\Big(\frac{1.5}{\GAUSSIANMEASURE} \Big) = \frac{\alpha (1+y(\rho))^{1/2}}{\rho^{3/2}} \E\Big[ \MARGIN' \cA(\INVERSEMILLSEXPRESSION) \Big] + \frac23 \GAUSSIANMEASURE - 2\GAUSSIANMEASURE <  \alpha 2^{1/2} \Csol^{3/2} K_{\upref{lem:controlexpectationofAW}} + \frac23 \GAUSSIANMEASURE - 2\GAUSSIANMEASURE < 0\,.\]
Similarly we obtain
\[ h\Big( \frac1{3.5 \beta} \Big) \ge -\frac43 \beta + 3.5\beta - 2\beta > 0\,.\]
\item For the posterior: Again $y(\rho) \le 1$ by Lemma \upref{lem:ybounded}. Thus as $|u'|, |u''| \le D$,
\begin{align*}
\Big| \E \big\langle u'' + u' \big\rangle - \frac{y}{1+y} \E \big\langle u' \big\rangle^2 \Big| \le K(D)\,.
\end{align*}
Thus as $\alpha$ is small enough in terms of $\beta_0, \beta_1, D, \Csol$, by \equpref{eq:posterior_partialGrho},
\begin{align*}
h\Big(\frac{1.5}{\GAUSSIANMEASURE} \Big) &\le \alpha K(D) + \frac{2}3 \beta - 2 \beta < 0\quad,\quad h\Big( \frac1{3.5 \beta} \Big) \ge -\alpha K(D) + 3.5 \beta - 2 \beta > 0\,.
\end{align*}
\end{itemize}
Thus by Intermediate Value Theorem and continuity of $h(\rho)$, there is some $\rho' \in [\frac1{3.5 \beta}, 1.5 \beta]$ such that $h(\rho')=0$, that is, $\frac{\partial G}{\partial \rho}(y(\rho'), \rho')=0$. By our definition of $\Csol$, we have $[\frac1{3.5 \beta}, \frac{ 1.5}{ \beta } ] \subset [\frac1{\Csol}, \Csol]$ for all $\beta \in [\beta_0, \beta_1]$. Hence $(y(\rho'), \rho)$ is a solution in $[0,1] \times [\frac1{\Csol}, \Csol]$ to \equpref{eq:RStransformed}, and this proves the Lemma. 
\end{proof}

Finally, recall that for $(\alpha, \beta) \in [0, \DENSITYBOUND] \times [\beta_0, \beta_1]$, if there is a solution $(y,\rho) \in [0,\infty) \times [\frac1{\Csol}, \Csol]$, then this solution must be in $[0,1] \times [\frac1{\Csol}, \Csol]$ by Lemma \upref{lem:ybounded}. We conclude that \equpref{eq:RStransformed} has a unique solution $(y,\rho) \in [0,\infty) \times [\frac1{\Csol}, \Csol]$ (in particular, this solution must be in $[0,1] \times [\frac1{\Csol}, \Csol]$). By our initial remarks in this proof, it follows that \equpref{eq:niceRSeqs} have a unique solution $(q, \rho) \in [0, \Csol] \times [\frac1{\Csol}, \Csol]$, as desired. 
Finally, as this unique solution $(y, \rho) \in [0,1] \times [\frac1{\Csol}, \Csol]$, we have $\rho - q = \frac{\rho}{y+1} \ge \frac1{2\Csol}$.
This proves the first part of Proposition \upref{prop:nicersuniquesol}.

We now prove the second part of Proposition \upref{prop:nicersuniquesol}.
First, recall $\frac{\partial^2 G}{\partial^2 y}\big(y_0(\alpha, \beta, x), \rho_0(\alpha, \beta, x) \big)$ and $\frac{\partial^2 G}{\partial^2 \rho }\big(y_0(\alpha, \beta, x), \rho_0(\alpha, \beta, x) \big)$ are of opposite signs and nonzero by Lemma \upref{lem:rhocontrolled} and Lemma \upref{lem:ycontrolled}. Thus the determinant of the following Jacobian is strictly negative:
\[ \det \begin{bmatrix} \frac{\partial^2 G}{\partial^2 y}\big(y_0(\alpha, \beta, x), \rho_0(\alpha, \beta, x)\big) & \frac{\partial^2 G}{\partial y \partial \rho} \big(y_0(\alpha, \beta, x), \rho_0(\alpha, \beta, x)\big) \\  \frac{\partial^2 G}{\partial y \partial \rho} \big(y_0(\alpha, \beta, x), \rho_0(\alpha, \beta, x)\big) & \frac{\partial^2 G}{\partial^2 \rho }\big(y_0(\alpha, \beta, x), \rho_0(\alpha, \beta, x)\big)  \end{bmatrix} < 0~.\]
Since $G$ is infinitely differentiable in $x, y, \rho$, and as $\big(y_0(\alpha, \beta, x), \rho_0(\alpha, \beta, x)\big)$ denotes the unique solution $(y, \rho) \in [0, \infty) \times [\frac1{\Csol}, \Csol]$ to the system $\frac{\partial G}{\partial y} = \frac{\partial G}{\partial \rho}=0$, it follows by the Implicit Function Theorem that $y_0(\alpha, \beta, x), \rho_0(\alpha, \beta, x)$ are infinitely differentiable in $\alpha$, $\beta$ and $x$. Hence as $y_0(\alpha, \beta, x) \ge 0$, it follows that
\[q_0(\alpha, \beta, x) = \frac{y_0(\alpha, \beta, x) \rho_0(\alpha, \beta, x)}{1+y_0(\alpha, \beta, x)} \]
is also infinitely differentiable in $\alpha$, $\beta$ and $x$.

\subsection{Proof of Lemma \upref{lem:rs_formula_properties}}\label{subsec:spherical_conversion_RS_pfs}
For this proof, we make the same transformation to $G(y,\rho)$ we considered in Appendix \upref{subsec:nice_rs_uniquesol_pf}, where $G(y, \rho)$ is defined as in \equpref{eq:G_transformed_def} (again, here $I$ is held fixed so $F \equiv F_I$). 
We let $q_0 = q_0(\alpha, \beta, x)$, $\rho_0 = \rho_0(\alpha, \beta, x)$, and $y_0 = y_0(\alpha, \beta, x) = \frac{q_0}{\rho_0-q_0}$. 

We first prove 1) of the Lemma. 
Differentiating the relations $\frac{\partial G}{\partial \rho}(y_0,\rho_0)\equiv0$, $\frac{\partial G}{\partial y}(y_0,\rho_0)\equiv0$ with respect to $\beta$ (recall $y_0, \rho_0$ depend on $\beta$), the Chain Rule gives
\begin{align*}
\frac{\partial^2 G}{\partial \rho \partial y}(y_0,\rho_0)\cdot \frac{\partial y_0}{\partial \beta} + \frac{\partial^2 G}{\partial \rho^2}(y_0, \rho_0)\cdot \frac{\partial \rho_0}{\partial \beta} + \frac{\partial^2 G}{\partial \beta \partial \rho}(y_0, \rho_0) &= 0\,, \\
\frac{\partial^2 G}{\partial y^2}(y_0, \rho_0)\cdot \frac{\partial y_0}{\partial \beta} + \frac{\partial^2 G}{\partial \rho \partial y}(y_0, \rho_0) \cdot \frac{\partial \rho_0}{\partial \beta} + \frac{\partial^2 G}{\partial \beta \partial y}(y_0, \rho_0) &= 0\,,
\end{align*}
where differentiability is justified by Proposition \upref{prop:nicersuniquesol}.
Now identically as functions, we have $\frac{\partial^2 G}{\partial \beta \partial \rho}=-1$, $\frac{\partial^2 G}{\partial \beta \partial y}=0$ by the definition of $G$ in \equpref{eq:G_transformed_def}. 
Multiplying the first relation above by $\frac{\partial \rho_0}{\partial \beta}$ and the second relation above by $\frac{\partial y_0}{\partial \beta}$ and subtracting, the $\frac{\partial^2 G}{\partial \rho \partial y}$ terms cancel and we obtain
\begin{align*}
\frac{\partial \rho_0}{\partial \beta} = \Big( \frac{\partial \rho_0}{\partial \beta}\Big)^2 \cdot \frac{\partial^2 G}{\partial \rho^2}(y_0, \rho_0) - \Big( \frac{\partial y_0}{\partial \beta} \Big)^2 \cdot \frac{\partial^2 G}{\partial y^2}(y_0, \rho_0) < 0\,,
\end{align*}
where Lemma \upref{lem:rhocontrolled} and Lemma \upref{lem:ycontrolled} justify that $\frac{\partial^2 G}{\partial \rho^2}(y_0, \rho_0)<0$, $\frac{\partial^2 G}{\partial y^2}(y_0, \rho_0)>0$. This proves 1) of the Lemma.
\begin{remark}
Note that we only have $\frac{\partial^2 G}{\partial \rho^2}(y_0, \rho_0)<0$, $\frac{\partial^2 G}{\partial y^2}(y_0, \rho_0)>0$ \textup{at the solution} $y_0, \rho_0$.
\end{remark}

We next prove 2) of the Lemma.
To this end, we will show that $\rho_0(\alpha, 1/4, x) > 1-x^2$, $\rho_0(\alpha, 9, x) < 1-x^2$. The existence and uniqueness of such a $\beta = \beta(\alpha, x) \in [1/4, 9] \subseteq [\beta_0, \beta_1]$ then follows by part 1) of this Lemma, and moreover this establishes that $\beta(\alpha, x) \in (1/4, 9)$.

We first show that $\rho_0(\alpha, 1/4, x) > 1 \ge 1-x^2$. Define
\begin{align}
\expressionI &:= \begin{cases}
\E\Big[ \MARGIN' \frac{(1+y_0)^{1/2}}{\rho_0^{3/2}} \cA(\INVERSEMILLSEXPRESSION)\Big] &: \text{ for interpolators}\,, \\
\E\big\langle u''(\bar V) + u'(\bar V)^2 \big\rangle - \frac{y}{1+y} \E\big\langle u'(\bar V) \big\rangle^2 &:\text{ for posterior}\,,
\end{cases}
\end{align}
where $V$ is as in \equpref{eq:RS_Vexpression_def} and $\bar V$ is as in \equpref{eq:RSeqs_barV_def_posterior}, taken with arguments $\rho_0, y_0$.
By \equpref{eq:partialGrho} for the interpolators and \equpref{eq:posterior_partialGrho} for the posterior, we have 
\begin{align}\label{eq:rs_lemma_2_equality}
0 = \frac{\alpha}2 \expressionI + \frac1{2 \rho_0 } - \beta\,.
\end{align}
Now note
\begin{align}\label{eq:rs_lemma_2_bound}
|\expressionI| \le \begin{cases}
\frac{\sqrt{2}}{\rho_0} K_{\upref{lem:controlexpectationofAW}} \Csol^{1/2} &: \text{ for interpolators}\,, \\
2D^2+D &: \text{ for posterior}\,.
\end{cases}
\end{align}
For the posterior, this simply uses that $|u'|, |u''| \le D$. For the interpolators, we use the proof of Proposition \upref{prop:nicersuniquesol}, that $\rho_0 \ge \frac1{\Csol}$, $y_0 \ge 0$, and that $y_0 \le 1$ by Lemma \upref{lem:ybounded}. These facts combined with Lemma \upref{lem:controlexpectationofAW} imply the above bound on $\big| \expressionI \big|$.

Now as $\alpha \le \alpha_0$ and by our conditions on $\alpha_0$, we have by \equpref{eq:rs_lemma_2_bound} that $\alpha |\expressionI| / 2 < \frac1{4 \rho_0}$.
Combining with \equpref{eq:rs_lemma_2_equality} now establishes
\begin{align}\label{eq:rs_lemma_2_bd1}
\beta \in \Big( \frac1{4\rho_0}, \frac{3}{4\rho_0} \Big) \implies \frac1{4\beta} < \rho_0(\alpha, \beta, x) < \frac3{4\beta} \implies \rho_0(\alpha, 1/4, x) > 1 \ge 1-x^2\,.
\end{align}
Next, we claim that $\rho_0(\alpha, 9, x) < 1-x^2$ for $x \in [-1+\delta, 1-\delta]$. To this end, \equpref{eq:rs_lemma_2_equality}, \equpref{eq:rs_lemma_2_bound} imply that 
\begin{align}\label{eq:rs_lemma_2_bd2}
9 \le \frac1{2\rho_0(\alpha, 9, x)} + \frac{\alpha}2 |\expressionI| \implies \rho_0(\alpha, 9, x) < 1-x^2\,.
\end{align}
Here, the last inequality holds for $x \in [-1+\delta, 1-\delta]$, using the definition of $\delta$ in \equpref{eq:deltaalphafulldef} for the interpolators or \equpref{eq:posteriordeltadef} for the posterior.
As mentioned earlier, \equpref{eq:rs_lemma_2_bd1}, \equpref{eq:rs_lemma_2_bd2} yields 2) of the Lemma. \\

Now to prove 3) of the Lemma, for fixed $(\alpha, x)$, it follows from Proposition \upref{prop:nicersuniquesol} that $f(\beta)$ is infinitely differentiable in $\beta$.
Let $F(x, q, \rho) = \bar\Phi(x, q, \rho) - \beta (\rho + x^2)$.
Note $\frac{\partial F}{\partial \rho}=\frac{\partial F_I}{\partial \rho}$, $\frac{\partial F}{\partial q}=\frac{\partial F_I}{\partial q}$, thus $\frac{\partial F}{\partial \rho}(x, q_0, \rho_0) = \frac{\partial F}{\partial q}(x, q_0, \rho_0) = 0$.
Therefore,
\begin{align*}
\frac{\rmd}{\rmd \beta}f(\beta) &= 1 - \big( \rho_0 + x^2 \big) + \frac{\partial F}{\partial \rho}(x, q_0, \rho_0) \cdot \frac{\partial \rho_0}{\partial \beta} + \frac{\partial F}{\partial q}(x, q_0, \rho_0) \cdot \frac{\partial q_0}{\partial \beta}= 1 - \big( \rho_0 + x^2 \big)\,.
\end{align*}
Thus 
\begin{align*}
\frac{\rmd^2}{\rmd \beta^2} f(\beta) &= -\frac{\partial \rho_0(\alpha, \beta, x)}{\partial \beta}\,.
\end{align*}
By 1) of this Lemma, $\frac{\rmd^2}{\rmd \beta^2} f(\beta) > 0$ on $[\beta_0, \beta_1]$, and by 2) of this Lemma, $\frac{\rmd}{\rmd \beta} f\big(\beta(\alpha, x)\big)=0$. Hence 3) of the Lemma follows. \\

Next, we prove 4) of the Lemma. 
Since $\rho_0\big(\alpha, \beta(\alpha, x), x) = 1-x^2$, we have 
\begin{align*}
f\big(\beta(\alpha, x)\big) = \bar\Phi\big( x, q_0\big(\alpha, \beta(\alpha, x), x\big), 1-x^2 \big)\,.   
\end{align*}
Note as $x \in [-1+\delta, 1-\delta] \subseteq [-17/18, 17/18]$, $1 \ge 1-x^2 \ge \frac{1}{10}$. Thus $1-x^2 \in \big[\frac1{\Csol}, \Csol\big]$ as $\Csol \ge 10$.
Recalling the definition of $\bar\Phi$ in \equpref{eq:phi_master_rho}, define the following function of $y$:
\begin{align*}
\bar\Phi'(x, y, \rho) := \bar\Phi \Big(x, \frac{y\rho}{y+1}, \rho \Big)\,.    
\end{align*}
Note $\frac{\partial \bar\Phi'}{\partial y}\big( y_0\big( \alpha, \beta(\alpha, x), x \big), 1-x^2 \big) = \frac{\partial G_I}{\partial y}\big( y_0\big( \alpha, \beta(\alpha, x), x \big), \rho_0\big( \alpha, \beta(\alpha, x), x \big) \big) = 0$.
Now as a direct consequence of the proof of Lemma \upref{lem:ycontrolled} and by 2) of this Lemma, $\bar\Phi'(x, y, 1-x^2)$ attains its unique minimum on $[0, \infty)$ at $y=y_0\big( \alpha, \beta(\alpha, x), x \big)$.
By monotonicity of the transformation $q = \frac{y\rho}{y+1}$, it follows that $\bar\Phi(x, q, 1-x^2)$ attains a unique minimum on $[0, 1-x^2)$ at $q = q_0\big( \alpha, \beta(\alpha, x), x \big)$. 
Thus
\begin{align*}
f\big(\beta(\alpha, x)\big) = \inf_{q \in [0, 1-x^2)} \bar\Phi(x, q, 1-x^2) = \inf_{q \in [0,1)} \Phi(x, q)\,,    
\end{align*}
where the last step follows from the transformation $q \leftarrow \frac{q}{1-x^2}$. This proves 4) of the Lemma.\\

Finally, to prove 5) of the Lemma, we claim that $\beta(\alpha,x)$ is a continuous function of $(\alpha, x)$. 
To this end, note by \equpref{eq:partialGrho} and 2), we know that $\beta(\alpha, x)$ is the unique solution in $[\beta_0, \beta_1]$ to
\begin{align}
\GAUSSIANMEASURE = \frac{\alpha}2 \expressionI + \frac1{2(1-x^2)} \,,\label{eq:betacontinuous_fixedpteq}
\end{align}
where $\expressionI$ is taken with $\rho_0 = 1-x^2$.
Let $\bar f(\alpha, \beta, x)$ denote the right hand side of \equpref{eq:betacontinuous_fixedpteq}; this is continuous in $(\alpha, \beta, x)$ for $(\alpha, \beta, x) \in [0, \alpha_0] \times [\beta_0, \beta_1] \times [-1+\delta, 1-\delta]$ by Proposition \upref{prop:nicersuniquesol}. 
Consider any sequence $(\alpha_n, x_n) \rightarrow (\alpha^\star, x^\star)$ where $\alpha_n \in [0, \DENSITYBOUND]$, $x_n \in [-1+\delta, 1-\delta]$.
Consider the corresponding $\beta_n := \beta(\alpha_n, x_n)$.
Since $\beta_n \in [\beta_0, \beta_1]$, it follows that we can extract a convergent subsequence.
Consider any convergent subsequence $\beta_{n_k}$, with limit $\beta^{\star}$.
By the above it follows by compactness that $\beta^{\star} \in [\beta_0, \beta_1]$.
As $\bar f(\alpha, \beta, x)$ is continuous, we have
\begin{align*}
\beta^{\star} = \lim_{k\rightarrow\infty} \beta_{n_k} = \lim_{k\rightarrow\infty} \bar f(\alpha_{n_k}, \beta_{n_k}, x_{n_k}) = \bar f(\alpha^\star, \beta^\star, x^\star)\,.
\end{align*}
However as $(\alpha^\star, x^\star) \in [0, \DENSITYBOUND] \times [-1+\delta, 1-\delta]$, the uniqueness of $\beta(\alpha, x)$ in $[\beta_0, \beta_1]$ supplied by 2) implies that $\beta^{\star} = \beta(\alpha^\star, x^\star)$. 
Therefore, $\lim_{k \rightarrow \infty} \beta(\alpha_{n_k}, x_{n_k}) = \beta^{\star} = \beta\big(\lim_{k\rightarrow\infty}(\alpha_{n_k}, x_{n_k})\big)$.
Since the convergent subsequence $\beta_{n_k}$ is arbitrary, it follows that $\beta_n$ converges to $\beta(\alpha^\star, x^\star)$.
Continuity of $\beta(\alpha, x)$ follows.
Combining this claim with 4) and Proposition \upref{prop:nicersuniquesol} now proves 5).
This concludes the proof of Lemma \upref{lem:rs_formula_properties}.

\subsection{Proof of Propositions \upref{lem:RS_biv_eq_gmm}, \upref{lem:RS_biv_eq_logistic}, \upref{lem:RS_biv_bayes_gmm}, \upref{lem:RS_biv_bayes_logistic}, and Lemma \upref{lem:nishimori_logistic}}
\label{subsec:univariate_maximizer_pf}
Here we first prove these four Propositions and then prove Lemma \upref{lem:nishimori_logistic}. We will prove these four Propositions together in a unified manner as follows. We show the following: 
\begin{enumerate}
\item\label{it:part1} There is at most one pair $(x, q)$ satisfying $\frac{\partial}{\partial x} \Phi(x,q)=\frac{\partial}{\partial q} \Phi(x,q)=0$.
\item\label{it:part2} This implies that the saddle point problem is uniquely achieved.
\end{enumerate}

\paragraph{Proof of Part~\upref{it:part1}.} 
First, we show that any $(x, q)$ satisfying $\frac{\partial}{\partial x} \Phi(x,q)=0$, $\frac{\partial}{\partial q} \Phi(x,q)=0$ must satisfy the relevant bivariate system of equations given in Propositions \upref{lem:RS_biv_eq_gmm}, \upref{lem:RS_biv_eq_logistic}, \upref{lem:RS_biv_bayes_gmm}, \upref{lem:RS_biv_bayes_logistic} respectively. 
The proof follows from explicitly computing $\frac{\partial}{\partial x} \Phi(x,q)$ and $\frac{\partial}{\partial q} \Phi(x,q)$ and applying Gaussian Integration by Parts to simplify the resulting expression. We discuss the proof for Propositions \upref{lem:RS_biv_eq_gmm}, \upref{lem:RS_biv_eq_logistic} for the interpolators; the proof for Propositions \upref{lem:RS_biv_bayes_gmm}, \upref{lem:RS_biv_bayes_logistic} for the posterior is identical.
Recalling the definition of $S$ in \equpref{eq:defSgmm} in the GMM case and \equpref{eq:defSlogistic} in the logistic case, we let
\begin{align*}
V := \frac{\MARGIN - x S  - \sqrt{(1-x^2)q}Z}{\sqrt{(1-x^2)(1-q)}}\,.
\end{align*}
Explicitly differentiating $\Phi(x,q)$ w.r.t. $q$ and simplifying with Guassian Integration by Parts and the identity $\cA'(V) = \cA(V)^2 - V\cA(V)$ gives
\begin{align*}
\frac{q}{1-q} &= \alpha \E\Big[ \cA(V) \Big\{ \frac{\MARGIN-xS}{\sqrt{(1-q)(1-x^2)}} - \frac{Z}{\sqrt{q(1-q)}} \Big\} \Big] \\
&= \alpha \E\Big[ \cA(V) \Big\{ V - \sqrt{\frac{1-q}{q}} Z \Big\}\Big] \\
&= \alpha \E\big[ \cA(V) V + \cA'(V) \big] = \alpha \E\big[ \cA(V)^2 \big].
\end{align*}
Therefore, such $(x, q)$ satisfying the stationary point equations must satisfy \equpref{eq:RS_biv_eq_gmm_q} in the GMM case or \equpref{eq:RS_biv_eq_logistic_q} in the logistic case. 

We will now show the proof that such $(x, q)$ satisfying the stationary point equations must satisfy \equpref{eq:RS_biv_eq_logistic_x} in the logistic case. The proof for \equpref{eq:RS_biv_eq_gmm_x} in the GMM case is similar (in fact simpler). We differentiate $\Phi(x, q)$ w.r.t. $x$ to obtain
\begin{align}\label{eq:univariate_derivx_eq0}
x \sqrt{(1-q)(1-x^2)} = \alpha \E\Big[ \cA(V) \big(S - \MARGIN x \big) \Big] = \alpha \E\Big[ \cA(V) \big( YG - \MARGIN x \big) \Big]\,,
\end{align}
where we used the explicit form $S=YG$ from \equpref{eq:defSlogistic} (for the GMM case, we use $S=G+\sqrt{\lambda}$ as per \equpref{eq:defSgmm}). 
We now simplify $\E\big[ \cA(V) YG \big]$ via Gaussian Integration by Parts to obtain
\begin{align}\label{eq:univariate_derivx_eq1}
\E\Big[ \cA(V) YG \Big] = -\frac{x}{\sqrt{(1-x^2)(1-q)}} \E\big[ \cA'(V) \big] + \sqrt{\lambda} \E\Big[ \varphi\big(-\sqrt{\lambda}YG \big) \cA(V) \Big]\,.
\end{align}
Similar calculations, using the identity $\cA'(V) = \cA(V)^2 - V\cA(V)$, yield
\begin{equation}\label{eq:univariate_derivx_eq2}
\begin{aligned}
\E\big[ \cA'(V) \big] &= (1-x^2)(1-q) \E\big[ \cA(V)^2 \big] - \MARGIN \sqrt{(1-x^2)(1-q)} \E\big[ \cA(V) \big] \\
&\qquad + x\sqrt{\lambda (1-x^2)(1-q)} \E\Big[ \varphi\big(-\sqrt{\lambda}YG \big) \cA(V) \Big]\,.
\end{aligned}
\end{equation}
Combining \eqref{eq:univariate_derivx_eq0}, \eqref{eq:univariate_derivx_eq1}, \eqref{eq:univariate_derivx_eq2} and using the identity $\alpha \E\big[ \cA(V)^2 \big] = \frac{q}{1-q}$ and rearranging yields \equpref{eq:RS_biv_eq_logistic_x}.
For Proposition \upref{lem:RS_biv_bayes_gmm}, this establishes \equpref{eq:explicitxq0}, which directly implies \equpref{eq:explicitxq} and proves 1). For Proposition \upref{lem:RS_biv_bayes_logistic}, 1) now follows from Proposition \upref{lem:nishimori_logistic}. We now prove 1) for Propositions \upref{lem:RS_biv_eq_gmm}, \upref{lem:RS_biv_eq_logistic} for the interpolators. 

To this end, our next step is to bound the value of $x,q$ for solutions to $\frac{\partial}{\partial x} \Phi(x,q)=0$, $\frac{\partial}{\partial q} \Phi(x,q)=0$. Specifically, we will prove the following Lemma.
\begin{lemma}\label{lem:rs_biv_eq_stationary_upperbd}
For any solution $(x, q)$ to $\frac{\partial}{\partial x} \Phi(x,q)=0$, $\frac{\partial}{\partial q} \Phi(x,q)=0$, we have
\begin{align}\label{eq:rs_biv_eq_stationary_upperbd}
0 \le x, q &\le \alpha K(\lambda)\,.
\end{align}
\end{lemma}
To prove Lemma \upref{lem:rs_biv_eq_stationary_upperbd}, we will need the following useful Lemma. Its proof is analogous to that of Lemma \upref{lem:controlexpectationofAWinit} and we omit it for brevity.
\begin{lemma}\label{lem:rs_biv_eq_bdmills1}
We have for all $0 \le x < 1, 0 \le q < 1$ that
\begin{align*}
\E\big[ \cA(V)^p \big] \le K\Big( \frac{\MARGIN_+ + 1}{\sqrt{(1-x^2)(1-q)}} + \sqrt{\frac{q}{1-q}} + 1\Big)^p
\end{align*}
for $p \in \{1,2\}$, where $K>0$ depends on $\lambda$.
\end{lemma}
\begin{proof}[Proof of Lemma \upref{lem:rs_biv_eq_stationary_upperbd}]
We now complete the proof of Lemma \upref{lem:rs_biv_eq_stationary_upperbd}. We start with the bound on $x$. 
We apply Lemma \upref{lem:rs_biv_eq_bdmills1} with $p=1$. First, note $x \ge 0$ as when $x<0$, the LHS of \equpref{eq:RS_biv_eq_logistic_x}, \equpref{eq:RS_biv_eq_gmm_x} is strictly negative while the RHS is non-negative. Noting that $\varphi \in [0, 1]$ and $\cA \ge 0$, we obtain
\begin{align*}
\frac{x}{\sqrt{1-x^2}} &\le \alpha \sqrt{\lambda(1-q)} \E\big[ \cA(V) \big] \le \alpha \sqrt{\lambda(1-q)} \cdot K(\lambda) \Big( \frac{\MARGIN_+ + 1}{\sqrt{(1-x^2)(1-q)}} + \sqrt{\frac{q}{1-q}} + 1\Big)\,.
\end{align*}
As $\alpha \le \alpha_0(\lambda, \MARGIN_+)$, we now use $0 \le q \le 1$ and multiply both sides by $\sqrt{1-x^2}$ to obtain the bound in \equpref{eq:rs_biv_eq_stationary_upperbd} for $x$.
Note as $\alpha \le \alpha_0(\lambda, \MARGIN_+)$, we have $x \le \frac12$. 
We now prove the desired bound in \equpref{eq:rs_biv_eq_stationary_upperbd} for $q$.
By Lemma \upref{lem:rs_biv_eq_bdmills1} and \equpref{eq:RS_biv_eq_gmm_q} for the GMM case or \equpref{eq:RS_biv_eq_logistic_q} for the logistic case, we have 
\begin{align*}
\frac{q}{1-q} \le \alpha K(\lambda) \Big( \frac{\MARGIN_+^2 + 1}{(1-x^2)(1-q)} + \frac{q}{1-q} + 1\Big) \implies \frac{q}{1-q} \le \frac{q}{2(1-q)} + \alpha K(\lambda)  \Big(\frac{\MARGIN_+^2 + 1}{1-q}+1\Big)\,,
\end{align*}
where we used $\alpha \le \alpha_0(\lambda, \MARGIN_+)$ and $x \le \frac12$. 
The above display now implies \equpref{eq:rs_biv_eq_stationary_upperbd} for $q$.
\end{proof}

We now are in a position to finish the proof of 1). Let $f_1(x, q)$ denote the $\LHS$ minus the $\RHS$ of \equpref{eq:RS_biv_eq_gmm_x} or \equpref{eq:RS_biv_eq_logistic_x} in the GMM case and the logistic case respectively. We similarly let $f_2(x,q)$ denote the $\LHS$ minus the $\RHS$ of \equpref{eq:RS_biv_eq_gmm_q} or \equpref{eq:RS_biv_eq_logistic_q} in these two respective cases. We now establish the following Lemma \upref{lem:rs_biv_eq_jacobian_positive}, which shows that all principal minors of the Jacobian of $\big(f_1(x,q), f_2(x,q)\big)$ are strictly positive in the domain $(x, q) \in [0, \alpha K(\lambda)]^2$. Since any solution $(x,q)$ must lie in $[0, \alpha K(\lambda)]^2$ by Lemma \upref{lem:rs_biv_eq_stationary_upperbd}, injectivity of $\big(f_1(x,q), f_2(x,q)\big)$ on $(x, q) \in [0, \alpha K(\lambda)]^2$ and therefore 1) now follows from the Gale-Nikaido Theorem \cite{gale1965jacobian}.
\begin{lemma}\label{lem:rs_biv_eq_jacobian_positive}
We have for all $(x, q) \in [0, \alpha K(\lambda)]^2$,
\begin{align*}
\frac{\partial}{\partial x} f_1(x, q) \,,\, \frac{\partial}{\partial q} f_2(x, q) \,,\, \frac{\partial}{\partial x} f_1(x, q) \cdot \frac{\partial}{\partial q} f_2(x, q) -  \frac{\partial}{\partial q} f_1(x, q) \cdot \frac{\partial}{\partial x} f_2(x, q) > 0\,.
\end{align*}
\end{lemma}
\begin{proof}[Proof of Lemma \upref{lem:rs_biv_eq_jacobian_positive}]
Lemma \upref{lem:rs_biv_eq_stationary_upperbd} implies $x, q \le 1/2$, thus we have 
\begin{align*}
\frac{\partial}{\partial x} \frac{x}{\sqrt{1-x^2}}\,,\, \frac{\partial}{\partial x} \frac{x}{1-x^2}\,,\, \frac{\partial}{\partial q} \frac{q}{1-q}\,,\, \frac{\partial}{\partial q} \frac{q}{(1-q)^2} \ge c > 0\,,
\end{align*}
where $c>0$ is a universal constant. 
Since $x, q \le 1/2$, $\alpha \le \alpha_0(\lambda, \MARGIN_+)$, and $\varphi \in [0,1]$, it suffices to prove
\begin{align*}
\E\Big[ \Big| \cA(V) \cA'(V) \frac{\partial V}{\partial x} \Big| \Big], \E\Big[ \Big| \cA(V) \cA'(V) \frac{\partial V}{\partial q} \Big| \Big], \E\Big[ \Big|\cA'(V) \frac{\partial V}{\partial x} \Big| \Big], \E\Big[ \big| \cA(V) \big| + \Big| \cA'(V) \frac{\partial V}{\partial q} \Big| \Big] \le K(\lambda, \MARGIN_+)\,.
\end{align*}
We note that the corresponding $\rho = 1-x^2 \in \big[\frac1{\Csol}, \Csol\big]$ and $y = \frac{q}{\rho-q} = \frac{q}{1-x^2-q} \le 1$ as $(x, q) \in [0, \alpha K(\lambda)]^2$. Therefore the conditions of Lemma \upref{lem:controlexpectationofAW} apply, giving $\E\big[ \big| (\MARGIN - Sx) \cA(\INVERSEMILLSEXPRESSION) \big| \big] \le K(\lambda, \MARGIN_+)$.
Upon using Cauchy-Schwarz together with Lemma \upref{lem:rs_biv_eq_bdmills1} to upper bound $\E\big[ \big| S \cA(V) \big| \big]$, this implies that $\E\big[ \big| \MARGIN \cA(V) \big| \big] \le K(\lambda, \MARGIN_+)$.
The upper bounds on the first two quantities above now follow as $|\cA'| \le 1$, $0 \le x, q \le \frac12$, using the explicit forms of $\partial V/\partial x, \partial V/\partial q$, and applying the Cauchy-Schwarz Inequality together with Lemma \upref{lem:rs_biv_eq_bdmills1}.

To prove bound the third and fourth quantities above, as $\cA'(V) = \cA(V)^2 - V \cA(V)$, we have $|\cA'(V)| \le 3|V| |\cA(V)| + |\cA(V)|$ as $\cA(V) \le 2|V|+1$. Now, following the exact same proof as that of Lemma \upref{lem:controlexpectationofAW} and since $0 \le x, q \le \frac12$, we have $\E\big[ \big|\MARGIN'^{2} \cA(V) \big| \big] \le K(\lambda, \MARGIN_+)$. Using Cauchy-Schwarz and Lemma \upref{lem:rs_biv_eq_bdmills1} to upper bound $\E\big[ S^2 |\cA(V)|\big]$, it follows that $\E\big[ |\MARGIN|^2 |\cA(V)| \big] \le K(\lambda, \MARGIN_+)$. Again applying Cauchy-Schwarz, this implies $\E\big[ \big|\cA'(V) \partial V/\partial x\big| \big]$, $\E\big[ \big|\cA'(V) \partial V/\partial q\big| \big] \le K(\lambda, \MARGIN_+)$. Lemma \upref{lem:rs_biv_eq_jacobian_positive} now follows.
\end{proof}

\paragraph{Proof of Part~\upref{it:part2}.} We now complete the proof of Propositions \upref{lem:RS_biv_eq_gmm}, \upref{lem:RS_biv_eq_logistic}, \upref{lem:RS_biv_bayes_gmm}, \upref{lem:RS_biv_bayes_logistic}. 
Define the function $f(x) := \inf_{q \in [0,1)} \Phi(x,q)$. 
Note $f(x)$ is the infimum of an arbitrary family of continuous functions, and therefore is upper semicontinuous.
First notice for the posterior, as $u(x) \ge -D(\lambda)(|x|+1)$, we have
\begin{equation}\label{eq:rs_bivariate_lowerbd_u_posterior}
\begin{aligned}
&\E \log \E_{W} \exp u\Big( x S + \sqrt{(1-x^2)q}Z + \sqrt{(1-x^2)(1-q)} W \Big) \\
&\quad \ge -D(\lambda) \E \Big[ 1 + \big| x S + \sqrt{(1-x^2)q}Z + \sqrt{(1-x^2)(1-q)} W \big| \Big] \ge -K(\lambda)\,.
\end{aligned}
\end{equation}
For the interpolators, using that $\cA(x) \le 2x+1$ for all $x \ge 0$ by Lemma \upref{lem:millsratioproperties} and that $-\log \cN(V) \le \log 2$ when $V \le 0$, we obtain $-\log \cN(V) \le K(V^2_{+} + 1)$ for all $V \ge 0$ where $V_+ = \max\{V, 0\}$.
Since $(a+b)_+ \le b_+$, it follows that 
\begin{align}\label{eq:rs_bivariate_lowerbd_u_interpolants}
\E \log \cN(V) \ge -K \cdot \frac{ \MARGIN_+^2 + x^2 S^2 + (1-x^2) q Z^2}{(1-x^2)(1-q)}\,.
\end{align}
Thus for $\alpha \le \alpha_0(\lambda, \MARGIN_+)$, we have for both the posterior and the interpolators, 
\begin{align}\label{eq:rs_bivariate_lowerbd_q}
\Phi(x, q) \ge \frac{q}{4(1-q)} - \frac{\alpha K(\lambda) (\MARGIN_+^2+1)}{(1-x^2)(1-q)} + \frac12 \log(1-x^2) + \frac12 \log(1-q) - \alpha K(\lambda)
\end{align}
Taking $x=0$ in the above display, bounding $\alpha K(\lambda) (\MARGIN_+^2+1) \le 1/8$, and noting the infimum of the resulting lower bound is attained at $q=1/2$ now proves that $f(0) \ge -K(\lambda)$. 

Next, notice $f(x) \le \Phi(x, 0)$. Note for the interpolators or the posterior in the logistic case, we have $u \le 0$ and therefore $\Phi(x,0) \le \frac12 \log(1-x^2)$. For the posterior in the GMM case, we have $u(x) = \sqrt{\lambda} x$ and therefore
\begin{equation}\label{eq:rs_bivariate_upperbd_u}
\begin{aligned}
&\E \log \E_{W} \exp u\Big( x S + \sqrt{(1-x^2)q}Z + \sqrt{(1-x^2)(1-q)} W \Big) \\
&\quad \le \log \E \exp \sqrt{\lambda}\Big( x S + \sqrt{(1-x^2)q}Z + \sqrt{(1-x^2)(1-q)} W \Big) \le K(\lambda)\,.
\end{aligned}
\end{equation}
Thus in all cases, we have $\lim_{x\rightarrow \pm 1} \Phi(x, 0) = -\infty$, and justifying the existence of the following limit and the equality
\begin{align*}
\lim_{x \rightarrow \pm 1} f(x) = -\infty\,.    
\end{align*}
The above derivation also proves that $f(0) \le \Phi(0,0) \le K(\lambda)$. Now as $f(0) \ge -K(\lambda)$ and $\lim_{x \rightarrow \pm 1} f(x) = -\infty$, by the upper semicontinuity of $f$ and as $f(0) \le K(\lambda)$, the Weierstrass Extreme Value theorem implies that $\sup_{x \in (-1,1)} f(x)$ is attained by at least one $x^\star \in [-1,1]$ and that all such $x^\star$ must lie in $[-1+\delta(\lambda), 1-\delta(\lambda)]$. 

Consider any such $x^\star$. As $x^\star \in [-1+\delta(\lambda), 1-\delta(\lambda)]$ and $\Phi(x^\star, q)$ is continuous in $q$, \equpref{eq:rs_bivariate_lowerbd_q} implies that $\inf_{q \in [0,1)} \Phi(x^\star,q)$ is attained at $q \in [0, 1 - \delta'(\lambda)]$, since the $\RHS$ of \equpref{eq:rs_bivariate_lowerbd_q} goes to $+\infty$ as $q \rightarrow 1$ and as $\alpha \le \alpha_0(\lambda, \MARGIN_+)$. 
Thus $\sup_{x \in (-1,1)} \inf_{q \in [0,1)} \Phi(x,q)$ is attained by at least one $(x^\star, q^\star)$.

Now consider any $x \in (-1,1)$. For the interpolators, we explicitly compute using Gaussian Integration by Parts,
\begin{align*}
\frac{\partial}{\partial q}\Phi(x, 0) = -\frac{\alpha}2 \E\big[ \cA(V)^2 \big] < 0\,.
\end{align*}
For the posterior, a similar calculation gives
\begin{align*}
\frac{\partial}{\partial q}\Phi(x, 0) = -\frac{\alpha(1-x^2)}2 \E\Big[ \Big( \frac{\E_W[ u' \cdot \exp u ]}{\E_W [\exp u]} \Big)^2 \Big] < 0\,.
\end{align*}
It follows that all such $q^\star > 0$, and therefore any such $(x^\star, q^\star) \in (-1,1) \times (0,1)$. Hence any such $(x^\star, q^\star)$ solves the system $\frac{\partial}{\partial x} \Phi(x,q)=0$, $\frac{\partial}{\partial q} \Phi(x,q)=0$. By 1), it follows that $(x^\star, q^\star)$ is unique and satisfies $(x_{\star}, q_\star) \in [0, K(\lambda) \alpha]^2$ -- for the posterior these bounds follow directly (GMM) or by Lemma~\upref{lem:nishimori_logistic} (logistic). 
This completes the proof. \newline

\begin{proof}[Proof of Lemma~\upref{lem:nishimori_logistic}]
Let
\begin{equation}\label{eq:tag0}
G \sim N(0,1)\,,\qquad Y | G \sim \text{Rad}\big(\varphi(\sqrt{\lambda}G)\big)\,,\qquad \eta:=\varphi(-\sqrt{\lambda} YG)\, .
\end{equation}
Recall the system of equations
\begin{align}%
\label{eq:system_bayes1}
    \alpha\lambda (1-q) \E \big[\varphi\big(-\sqrt{\lambda}YG\big) R(x,q)\big] = \frac{x}{1-x^2}\,,
    ~~\mbox{and}~~\alpha \lambda (1-x^2) \E \big[R(x,q)^2\big] = \frac{q}{(1-q)^2}\,,
\end{align}
where
\begin{align}\label{eq:R1}
R(x,q) := \frac{\E_{W} \varphi'(\sqrt{\lambda}V)}{\E_{W}\varphi(\sqrt{\lambda}V)}\,,~~ V := xYG+\sqrt{(1-x^2)q}Z + \sqrt{(1-x^2)(1-q)}W\,,
\end{align}
and $G, Z ,W \sim N(0,1)$ are mutually independent. 
We wish to show that for $\alpha$ small enough this system has a unique solution which takes the form $(x,q)=(s/(1+s),s/(1+2s))$. 
Consider the parametrized curve
\[x_s:=\frac{s}{1+s}\,, \qquad q_s:=\frac{s}{1+2s}\,,\qquad s\ge 0\,,\]
which observes the Nishimori identity
\[\frac{q_s}{1-q_s}=x_s\,,\]
and
\begin{equation}\label{eq:tag1}
1-x_s^2=\frac{1+2s}{(1+s)^2}\,,
\qquad
\sqrt{(1-x_s^2)q_s}=\frac{\sqrt{s}}{1+s}\,,
\qquad
\sqrt{(1-x_s^2)(1-q_s)}=\frac{1}{\sqrt{1+s}}\,.
\end{equation}
We consider the auxiliary observation model
\[\bar G=\sqrt{s}\,G+W\,,
\qquad W\sim N(0,1)\,,\]
with $W$ independent of $(G,Y)$. Conditionally on $\bar G$,
the posterior law of $G$ before observing $Y$ is Gaussian with mean
$\frac{\sqrt{s}}{1+s}\bar G$ and variance $\frac{1}{1+s}$.
Since
\[\P(Y=y\,|\, G)=\varphi(\sqrt{\lambda} yG)\,,
\qquad y\in\{-1,+1\}\,,\]
Bayes' formula gives
\begin{equation}\label{eq:tag2}
\E[\eta \,|\, \bar G,Y]
=
\frac{
\E_W\varphi' \!\left(
\sqrt{\lambda} Y\left(
\frac{\sqrt{s}}{1+s}\bar G+\frac{1}{\sqrt{1+s}}W
\right)
\right)
}{
\E_W\varphi \!\left(
\sqrt{\lambda} Y\left(
\frac{\sqrt{s}}{1+s}\bar G+\frac{1}{\sqrt{1+s}}W
\right)
\right)
} = R(x_s,q_s)\,.
\end{equation}
Here we used the identity $\varphi(t)\varphi(-t)=\varphi'(t)$
and symmetry of the Gaussian distribution, together with the equalities displayed in~\eqref{eq:tag1}.
Consequently by the tower property,
\begin{equation}\label{eq:tag3}
\mathbb E[\eta R(x_s,q_s)]
=\E\left[\E[\eta\mid \bar G,Y]^2\right] 
=\E[R(x_s,q_s)^2].
\end{equation}

Now define
\[C(s):= \E[R(x_s,q_s)^2]\,.\]
Using~\eqref{eq:tag3}, the system of equations~\eqref{eq:system_bayes1} restricted to the curve
$s \mapsto (x_s,q_s)$ become
\[
\alpha\lambda(1-q_s)C(s)
=
\frac{x_s}{1-x_s^2}\,,
\qquad\mbox{and}\qquad
\alpha\lambda(1-x_s^2)C(s)
=
\frac{q_s}{(1-q_s)^2}\,.
\]
Both are equivalent to
\begin{equation}\label{eq:tag4}
s=\alpha\lambda C(s)\,.
\end{equation}

We show that for $\alpha$ small enough, there is a unique solution which must lie on this curve.
Let
\[A(x,q):=\E[\eta R(x,q)]\,, \qquad B(x,q):=\E[R(x,q)^2]\,.\]
Since $0<\varphi<1$ and $0<\varphi'<\varphi$, $A,B,C \in (0,1)$.
Therefore any solution of~\eqref{eq:system_bayes1} satisfies
\[
\frac{x}{1-x^2}\le \alpha\lambda\,,
\qquad
\frac{q}{(1-q)^2}\le \alpha\lambda\,.
\]
In particular,
\begin{equation}\label{eq:tag5}
x\le \alpha\lambda\,,
\qquad
q\le \alpha\lambda\,.
\end{equation}
Thus every solution lies in an arbitrarily small neighborhood of $(0,0)$ once $\alpha$ is sufficiently small.
We now use local uniqueness. Let
\[
H(x):=\frac{2x}{1+\sqrt{1+4x^2}}\,,
\qquad
K(q):=\frac{\sqrt{1+4q}-1}{\sqrt{1+4q}+1}\,,
\]
the inverses of $x\mapsto x/(1-x^2)$ on $[0,1)$ and $q\mapsto q/(1-q)^2$ on $[0,1)$, respectively. Hence the system is equivalent to
the fixed-point equation
\[
(x,q)=T_\alpha(x,q)\,,
\]
where
\begin{equation} \label{eq:fixed_point_map}
T_\alpha(x,q)
:=
\left(
H\!\left(\alpha\lambda(1-q)A(x,q)\right),
K\!\left(\alpha\lambda(1-x^2)B(x,q)\right)
\right).
\end{equation}

The functions $A$ and $B$ are locally Lipschitz in a neighborhood of $(0,0)$.
Indeed, after conditioning on $(G,Y)$, the remaining average (with respect to $Z$, Eq.~\eqref{eq:R1}) is a Gaussian
heat semigroup applied to the smooth function
\[
u\mapsto
\frac{\E_W\varphi'(\sqrt{\lambda}(u+\sigma W))}
{\E_W\varphi(\sqrt{\lambda}(u+\sigma W))}\,,
\]
(and its square) with $\sigma$ bounded away from $0$ near $(x,q)=(0,0)$. Differentiation under the Gaussian integral is justified by the boundedness and exponential decay properties of the derivatives of the logistic function. Hence there exist
$\delta>0$ and $L<\infty$, depending only on $\lambda$, such that $A$ and $B$
are $L$-Lipschitz on $[0,\delta]^2$.

Since $H$ and $K$ have bounded derivatives on bounded intervals, it follows
that, for $(x,q),(x',q')\in[0,\delta]^2$,
\[
\|T_\alpha(x,q)-T_\alpha(x',q')\|_\infty
\le
C_\lambda \alpha
\|(x,q)-(x',q')\|_\infty
\]
for some finite constant $C_\lambda$. Choosing $\alpha_0>0$ sufficiently small,
we may ensure that $C_\lambda\alpha_0<1$ and $\alpha_0\lambda<\delta$.
Then $T_\alpha$ is a contraction on $[0,\delta]^2$ for every
$\alpha\in(0,\alpha_0)$.
By~\eqref{eq:tag5}, every solution of the system~\eqref{eq:system_bayes1} lies in $[0,\delta]^2$ for
$\alpha<\alpha_0$. Hence the system has at most one solution.

On the other hand, the scalar equation~\eqref{eq:tag4} has a solution for small $\alpha$. Indeed, $0<C(s)<1$ and $C$ is continuous near $0$, while
\[
s-\alpha\lambda C(s)
\]
is negative at $s=0$ and positive at $s=\alpha\lambda$ for $\alpha$ small
enough. Moreover, by the same Lipschitz argument, the map
$s\mapsto \alpha\lambda C(s)$ is a contraction on a small interval, so this
solution is unique.

Thus the unique small-$\alpha$ solution of the two-dimensional system lies on
the curve
\[
(x,q)=\left(\frac{s}{1+s},\frac{s}{1+2s}\right)\,,
\]
and satisfies $s = \alpha\lambda C(s)$. This proves the claim.
\end{proof}


\section{Additional technical results}\label{sec:appendix}
\begin{lemma}[Lemma 3.2.2, \cite{talagrand2010mean}]\label{lem:realRSeqsolsboundedtechnicalhelper}
Consider a real $c>0$ and a univariate concave function $w(y)$ such that $w'' \le -c < 0$. Let $y'$ be the solution to $w'(y')=0$. Then we have
\begin{align*}
c \int (y - y')^2 \exp w(y)\, \rmd y \le \int \exp w(y)\, \rmd y\,.
\end{align*}
\end{lemma}

\subsection{Truncated Logarithm}\label{subsec:appendixtruncatedlog}
\begin{lemma}[Lemma 8.3.7, \cite{talagrand2011advanced}]\label{lem:truncatedlog1correctproof}
For $0 \le x, z \le 1$, we have 
\[ \Big| \log_A x - \log_A z \Big| \le \Big| \log_A\Big( \frac{x}{z} \Big) \Big|\, .\]
\end{lemma}
\begin{proof}
Suppose $z \le x$. The result is clear if $x \le \exp(-A)$, when the left hand side is 0, or if $z \ge \exp(-A)$, where the desired inequality follows from definition of logarithms. If $z \le \exp(-A) \le x$, as here $\frac{x}{z} \ge 1$ and $x, z \le 1$, we obtain
\[ | \log_A x - \log_A z | \le | \log x - \log z | = \Big| \log \frac{x}{z} \Big| = \Big| \log_A \frac{x}{z} \Big|\, .\]

Else, suppose $x \le z$. By similar arguments as above, the result is clear if $z \le \exp(-A)$ or if $x \ge \exp(-A)$, where we now note $x \ge \exp(-A) \ge z \exp(-A)$ in the latter case. If $x \le \exp(-A) \le z$, then noting $\log z \in [-A, 0]$, $\log x \le -A$, and $x, z \le 1$, we obtain
\[ | \log_A x - \log_A z | = | -A - \log z | \le A\,,\, | \log_A x - \log_A z | \le |\log x - \log z| = \Big|  \log \frac{x}{z}\Big|\, .\]
Therefore, we obtain
\[ | \log_A x - \log_A z| \le \min\Big\{ A, \Big| \log \frac{x}{z} \Big| \Big\} = \Big| \log_A  \frac{x}{z} \Big|\, \]
in this case as well.
\end{proof}
\begin{lemma}[Lemma 8.3.10, \cite{talagrand2011advanced}]\label{lem:truncatedlog2}
If $0 \le x,y,z \le 1$, we have 
\[ |\log_A xz - \log_A yz| \le |\log_A x - \log_A y|\one\{ z \ge \exp(-A) \}\, .\]
\end{lemma}
\begin{lemma}[Lemma 8.3.11, \cite{talagrand2011advanced}]\label{lem:truncatedlog3}
If $0 < y \le x \le 1$, then for any $c>0$, 
\[ |\log_A x - \log_A y| \le |\log_A y| \one\{ y \le c\} + \frac{|x-y|}c\, .\]
\end{lemma}

\end{document}